\documentclass[mestrado, pre-defesa]{packages/icmc}

\usepackage{rotating}           

\usepackage{amsthm, amsfonts, amssymb, amscd, rotating}
\usepackage{tikz-cd}
\usepackage{aligned-overset}
\usepackage{leftindex}

\DeclareMathOperator{\z}{\mathbb{Z}}
\newcommand{\q}{\mathbb{Q}}

\newcommand{\F}{\mathbb{F}}

\newcommand{\GL}{\text{GL}}

\newcommand{\PSL}{\text{PSL}}
\newcommand{\PGamma}{\text{P}\Gamma}
\newcommand{\SL}{\text{SL}}

\newcommand{\Ind}{\text{Ind}}
\newcommand{\coker}{\text{coker}}
\newcommand{\im}{\text{im}}
\newcommand{\id}{\text{id}}
\newcommand{\rank}{\text{rank}}

\newcommand{\Tot}{\text{Tot}}
\newcommand{\Tor}{\text{Tor}}

\newcommand{\arr}{\rightarrow}

\newcommand {\mtxx}[4]
{\left(\!
\begin{array}{cc}
\!\!#1 & \!\!#2 \\
\!\!#3 & \!\!#4
\end{array}\!\!
\right)
}

\usetikzlibrary{arrows}
\usetikzlibrary{positioning}
\usetikzlibrary{matrix}
\usepackage{venndiagram}
\usepackage[many]{tcolorbox}
\usepackage[dvipsnames]{xcolor}
\usetikzlibrary{babel}
\usetikzlibrary{fit, shapes.geometric, arrows, calc, math, positioning, 
decorations.pathreplacing, decorations.markings, patterns, intersections,
fadings, shadows.blur}
\usepackage{tkz-euclide}
\usepackage{pgf,tikz,pgfplots}
\usepackage{tikz-3dplot}
\pgfplotsset{compat=newest}
\usepgfplotslibrary{fillbetween}
\tikzset{
set arrow inside/.code={\pgfqkeys{/tikz/arrow inside}{#1}},
set arrow inside={end/.initial=>, opt/.initial=},
/pgf/decoration/Mark/.style={mark/.expanded=at position #1 with{
\noexpand\arrow[\pgfkeysvalueof{/tikz/arrow inside/opt}]
{\pgfkeysvalueof{/tikz/arrow inside/end}}}},
arrow inside/.style 2 args={set arrow inside={#1},
postaction={decorate,decoration={markings,Mark/.list={#2}}}},}

\tituloPT{A homologia integral de $\mathrm{SL}_2(\mathbb{Z}[1/n])$}
\tituloEN{The integral homology of $\mathrm{SL}_2(\mathbb{Z}[1/n])$}
\autor[Picinini, I. V.]{Isadora Vanzella Picinini}
\genero{F} 
\orientador[Orientador]{Prof. Dr.}{Behrooz Mirzaii}
\curso{MAT}
\data{25}{05}{2026} 
\idioma{EN} 

\textoresumo[brazil]{
    Seja $\mathbb{Z}[1/n] := \{a/n^r: a \in \mathbb{Z}, r \in \mathbb{Z}^{\geq 0}\}$ o domínio euclidiano obtido do anel de inteiros $\mathbb{Z}$ localizando em $n$. Ademais, seja $\text{SL}_2(\mathbb{Z}[1/n])$ o grupo das matrizes 2 por 2 invertíveis de determinante 1 com entradas em $\mathbb{Z}[1/n]$. Os grupos de homologia $H_k(\SL_2(\z[1/n]),\z)$ são de interesse em Teoria Geométrica de Grupos, Teoria Algébrica dos Números e K-teoria Algébrica. Nesta dissertação, estudamos estes grupos via uma sequência espectral obtida através da ação de $\SL_2(\z[1/n])$ no produto das árvores $B_p$ associadas com a valuação $p$-ádica em $\q$ onde $p$ é um fator primo de $n$.}{Homologia de grupos, grupo especial linear, Homologia equivariante, Ações de grupos em árvores}

\textoresumo[english]{
    Let $\mathbb{Z}[1/n] := \{a/n^r: a \in \mathbb{Z}, r \in \mathbb{Z}^{\geq 0}\}$ be the euclidean domain obtained from the ring of integers $\mathbb{Z}$ by localizing at $n$. Moreover, let $\text{SL}_2(\mathbb{Z}[1/n])$ be the group of all invertible 2-by-2 matrices of determinant one with entries in $\mathbb{Z}[1/n]$. The homology groups $H_k(\text{SL}_2(\mathbb{Z}[1/n]), \mathbb{Z})$ are of interest in Geometric Group Theory, Algebraic Number Theory and Algebraic K-theory. In this dissertation, we study these groups via a spectral sequence derived by the action of $\text{SL}_2(\mathbb{Z}[1/n])$ on the product of the trees $B_{p}$ associated with the $p$-adic valuation on $\q$ where $p$ is a prime factor of $n$.
    }{Homology of groups, special linear group, Equivariant homology, Group actions on trees}

\incluifichacatalografica{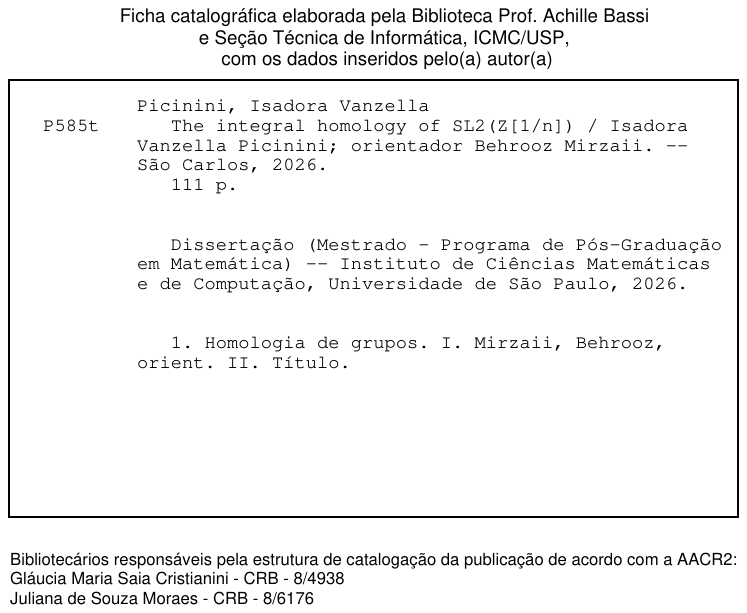}

\textoagradecimentos*{tex/pre-textual/agradecimentos}

\incluilistadefiguras

\begin{document}
\textual

\clearpage
\pagenumbering{arabic}
\setcounter{page}{1}
\chapter*{Introduction}
\markboth{Introduction}{Introduction}
\addcontentsline{toc}{chapter}{INTRODUCTION}
\label{chapter:introduction}
Let $\mathbb{Z}[1/n]$ be the following subring of the field of rational numbers $\mathbb{Q}$: $$\mathbb{Z}[1/n] := \{a/n^r: a \in \mathbb{Z}, r \in \mathbb{Z}^{\geq 0}\}.$$ This ring is an euclidean domain, obtained from the ring of integers $\mathbb{Z}$. If $n=p_1^{a_1}\cdots p_k^{a_k}$ is the prime decomposition of $n$, it is easy to check that $\mathbb{Z}[1/n]=\mathbb{Z}[1/(p_1\cdots p_k)]$.

Let $\text{SL}_2(\mathbb{Z}[1/n])$ be the group of all invertible 2-by-2 matrices of determinant one with entries in $\mathbb{Z}[1/n]$, i.e. 
$$\text{SL}_2(\mathbb{Z}[1/n]) = \left\{\begin{pmatrix}a & b \\ c & d\end{pmatrix}: a, b, c, d \in \mathbb{Z}[1/n], ad - bc = 1\right\}.$$ 

The integral homology groups of the special linear groups $\text{SL}_2(\mathbb{Z}[1/n])$, i.e. the homology groups $H_k(\text{SL}_2(\mathbb{Z}[1/n]), \mathbb{Z})$ are of interest in Geometric Group Theory \cite{brown1994}, \cite{ww1998}, Algebraic Number Theory \cite{serre1980} and Algebraic K-theory \cite{rosenberg1996}, \cite{weibel2013}. Finding the exact structure of these groups is an important and interesting problem in algebra and geometry \cite{an1998}, \cite{ww1998}. 

In this master's dissertation, we study the integral homology groups $H_k(\text{SL}_2(\mathbb{Z}[1/n]),\z)$ for a general natural number $n$. The structure of these groups have been calculated for all $k$ when $n$ is a prime number \cite{an1998}, for all $n$ if $k=1$ \cite{B-E2025}, for $n$ divisible by one of the primes 2, 3, 5, 7 or 13 if $k=2$ \cite{MRV2025}, and for most small values of $n \leq 50$ \cite{ae2014}.

Our main tool to study these homology groups will be a spectral sequence that converges to them.
Given a contractible oriented CW-complex $X$ on which a group $G$ acts, Corollary \ref{seq} shows that there exists a spectral sequence $$E^1_{p,q} = \displaystyle\bigoplus_{\sigma\in\Sigma_p}H_q(G_\sigma,\z_\sigma) \Longrightarrow H_{p+q}(G,\z)$$ where $\Sigma_p$ is a set of representatives for the $G$-orbits of $p$-cells, $G_\sigma$ is the stabilizer of a cell $\sigma$ and $\z_\sigma$ is the infinite cyclic group whose two generators correspond to the two orientations of $\sigma$ (so $g\in G_\sigma$ acts on $\z_\sigma$ as $+1$ if $g$ preserves the orientation of $\sigma$ and $-1$ otherwise). Theorem \ref{teo: maps} shows that if $X$ is finite-dimensional and $G_\sigma$ fixes $\sigma$ pointwise for all cells $\sigma$, then $\z_\sigma\simeq \z$ and the differential $d^1: E^1_{p,*}\rightarrow E^1_{p-1,*}$ is given by $$d^1\rvert_{H_*(G_\sigma,\z)}= \sum_{\tau\in\Sigma_{p-1}}[\sigma:\tau](i_{\sigma\tau})_*$$ where $(i_{\sigma\tau})_*$ is induced by the inclusion map $i_{\sigma\tau}: G_{\sigma}\hookrightarrow G_\tau$. 

To study the homology groups of $\SL_2(\z[1/n])$ we will use the product of the Bruhat-Tits Buildings $B_p$ associated with $\GL_2(\q_p)$ for the primes $p$ that divide $n$. We construct these trees using lattices. Given linearly independent vectors $v_1,v_2\in \mathbb{Q}^2$, the lattice defined by these vectors is the set $$L_{v_1v_2} = \{a_1v_1+a_2v_2 \mid a_1,a_2\in\z_{(p)}\} = \z_{(p)} v_1\oplus\z_{(p)} v_2$$ where $\z_{(p)}$ is the localization of $\z$ outside the prime ideal $(p) = p\z$. We consider an equivalence relation on the lattices where two lattices are equivalent if they are homothetic. The Structure Theorem then implies that given any two lattices $L,L'$, there exists a basis $\{u_1,u_2\}$ of $L$ such that $\{p^au_1,p^bu_2\}$ is a basis of $L'$ for some integers $a,b$. If $L,L'$ are representatives for equivalnce classes $\Lambda,\Lambda'$, we define the distance $d(\Lambda,\Lambda')$ to be $|a-b|$. We define the tree $B_p$ by letting vertices be equivalence classes of lattices.

The group $\GL_2(\q)$ acts on $B_p$ via standard representation. This action is transitive. We let $\SL_2(\z[1/p])$ act on this tree by restricting the action of $\GL_2(\q)$. Theorem \ref{teo: stab zp} shows that the stabilizer of a vertex is isomorphic to $\SL_2(\z)$ and Theorem \ref{teo: stab ep} shows that the stabilizer of an edge is isomorphic to $\Gamma_0(p)$. Finally, Theorems \ref{teo: numb vert 1} and \ref{teo: num edges 1} proves that this action has an edge as its quotient graph. 

The trees $B_p$ can be used to construct a CW-complex on which $\SL_2(\z[1/n])$ acts. If $p_1,\dots,p_k$ are the distinct prime factors of $n$, then $\text{SL}_2(\mathbb{Z}[1/n]) = \text{SL}_2(\mathbb{Z}[1/(p_1\cdots p_k)])$ and we consider the diagonal action of $\SL_2(\z[1/n])$ on the product $B_{p_1}\times \cdots\times B_{p_k}$. This complex is also a Bruhat-Tits Building. It is associated with the group $\GL_2(\q_{p_1})\times\cdots\times\GL_2(\q_{p_k})$. The cells in this complex are given by $(\sigma_1,\dots,\sigma_k)$ where each $\sigma_i$ is a cell in $B_{p_i}$. 

Theorem \ref{teo: standard cells} shows that this action can move each component cell $\sigma_i$ independently, so the quotient graph for this action is the product of the individual quotient graphs (in particular, we conclude it is a hypercube of dimension $k$). Theorem \ref{stab gen} shows that the stabilizer of an $s$-cell is isomorphic to a group $\Gamma_0(p_{i_1}\cdots p_{i_s})$ where each $p_{i_j}$ is some prime factor of $n$. Thus, there is a spectral sequence $$E_{s,t}^1 = \bigoplus H_t(\Gamma_0(p_{i_1}\cdots p_{i_s}),\z)^{2^{k-s}} \Longrightarrow H_{s+t}(\text{SL}_2(\mathbb{Z}[1/(p_1\cdots p_k)]),\z)$$ where the sum is understood to be over the $s$-element subsets of $\{p_1,\dots,p_k\}$.

We recall that the subgroup $\Gamma_0(n)$ of $\SL_2(\z)$ is given by
\[
\Gamma_0(n) :=\bigg\{ {\mtxx a b c d} \in \SL_2(\z): n\mid c \bigg\}.
\]
So to study the homology groups of $\SL_2(\z[1/n])$ via the above spectral sequence, we need to calculate the homology groups of $\Gamma_0(n)$ for any $n$.

These calculations were done for cohomology by Williams and Wisner \cite{ww1998}. Understanding the main result of this article is the topic of Chapter 3 of this dissertation. In order to calculate these homology groups, we construct a tree on which $\Gamma_0(n)$ acts. For this, we let \[
X' := \left\{\binom{w}{y} \in \z^2: \gcd(w,y)=1\right\}
\] and consider the action of $\SL_2(\z)$ on $X'$ given by multiplication on the left. For our construction, we will take $X = X'/\{\pm I\}$, that is, we identify antipodes in the lattice $\z^2$. We then consider triangles with vertices $\displaystyle\binom{w}{y},\binom{x}{z},\binom{t}{u} \in X$ such that $t = w + x, u = y + z $ and $ wz - xy = \pm 1.$ This allows us to represent triangles with $2\times 3$ matrices \[
T = \begin{pmatrix}w & x & t \\ y & z & u\end{pmatrix}
\] where any $2\times 2$ minor has determinant $\pm 1$ and one of the columns can be written as the sum of the other two.

We denote by $\Delta$ the set of triangles. We let $\SL_2(\z)$ act on $\Delta$ also by multiplication on the left and consider the restriction of this action to $\Gamma_0(n)$.

It is then possible to define a map $\varepsilon: \SL_2(\z) \rightarrow \Delta$ by adding the sum of the columns of a matrix as a third column, which is then a representative of some triangle. 
Given a triangle $T$ as before, we define three associated elements of $\SL_2(\z)$ by 
\[
T_1 = \begin{pmatrix}
w & x \\ y & z
\end{pmatrix}, \
T_2 = \begin{pmatrix}
t& -w \\ u & -y
\end{pmatrix}, \ 
T_3 = \begin{pmatrix}
x & -t \\ z & -u
\end{pmatrix}.
\]

The map $\varepsilon$ allows us to count the orbits of the $\Gamma_0(n)$-action on $\Delta$. Any triangle in $\Delta$ lifts to the three matrices above. If $T_1,T_2,T_3$ represent three different cosets in $\Gamma_0(n)$, we say $\Delta$ is of Type 1, and say it is of Type 2 otherwise. If $a(n)$ is the index of $\Gamma_0(n)$ in $\SL_2(\z)$ and $b(n)$ is the number of roots of $X^2+X+1$ mod $n$, we conclude that there are $\displaystyle\frac{a(n)+2b(n)}{3}$ distinct classes of triangles with $b(n)$ of those being of Type 2.

We construct a tree $Y$ by letting each triangle in $\Delta$ be a vertex. Here, two vertices are adjancent if the corresponding triangles are adjacent, that is, if they share an edge. The action of $\Gamma_0(n)$ on $Y$ may invert edges. To avoid this, we add a new vertex in the middle of any edge that is inverted by $\Gamma_0(n)$. These vertices are said to be Type 3. If $c(n)$ is the number of roots of $X^2+1$ mod $n$, there are $c(n)$ distinct classes of Type 3 vertices.

Proposition \ref{stab} shows that the stabilizer of an edge of $Y$ under the action of $\Gamma_0(n)$ is isomorphic to $\z/2$ and the stabilizers for vertices of Types 1, 2 and 3 are isomorphic to $\z/2,\z/6$ and $\z/4$, respectively. Theorem \ref{H_n} then uses the spectral sequence associated with this action to conclude that \[
H_k(\Gamma_0(n),\z)\simeq \begin{cases}
\z & \text{if $k=0$}\\
\z^{r(n)}\oplus G & \text{if $k=1$}\\
(\z/2)^{r(n)} & \text{if $k > 1$ is even,}\\
G & \text{if $k>1$ is odd}
\end{cases}
\]
where $\displaystyle r(n):=\frac{a(n)}{6}-\frac{2b(n)}{3}-\frac{c(n)}{2}+1$ and
\[
G := \begin{cases}
(\z/2)^{c(n)-1}\oplus (\z/3)^{b(n)}\oplus \z/4 & \text{if $c(n)>0$} \\
\z/2\oplus (\z/3)^{b(n)} & \text{if $c(n)=0$}.
\end{cases}
\]

This result allows us to conclude that the homology groups $H_k(\SL_2(\z[1/n]),\z)$ are finitely generated (Theorem \ref{teo: fint gen}) and
that the rank of $H_1(\SL_2(\z[1/n]),\z)$ is zero (Corollary \ref{cor: H1 is 0}). Moreover, for a prime $p$ we have 
$$\rank(H_k(\SL_2(\z[1/p]),\z)) = \begin{cases}
        1, & \text{ if }k=0 \\
        r(p), & \text{ if }k=2 \\
        0, & \text{ otherwise}
    \end{cases}$$ where $$r(p) = \begin{cases}
        1, & \text{ if } p=2,3\\
        \frac{p+1}{6}, &\text{ if } p\equiv 5\pmod{12}\\
        \frac{p-1}{6}, &\text{ if }p\equiv 7\pmod{12}\\
        \frac{p+7}{6}, &\text{ if }p\equiv 11\pmod{12}\\
        \frac{p-7}{6}, &\text{ if }p\equiv 1\pmod{12}
    \end{cases}$$
(Corollary \ref{cor: Hk for prime}). Furthermore, if $n=p_1\cdots p_k$, then we have the upper bound for ranks
$$\rank(H_s(\SL_2(\z[1/n]),\z)) \leq 2^{k-s+1}\sum r(p_{i_1}\cdots p_{i_{s-1}})$$
(Corollary \ref{cor: upper bound}) and the improved bound for the rank of the third homology
$$\rank(H_3(\SL_2(\z[1/n]),\z)) \leq 2^{k-2}\left(\sum_{1\leq i<j\leq k} r(p_ip_j)\right)-2^{k-1}\left(\sum_{1\leq i \leq k}r(p_i) \right)+r(p_1)$$ 
where $r(p_1)\leq r(p_i)$ for all $i=2,\dots,k$. In particular, we can conclude that if $n=pq$ with $p$ or $q$ being one of the primes $2,3,5, 7, 13$, then
 $$
 \rank(H_k(\SL_2(\z[1/(pq)]),\z))=
 \begin{cases}
1 & \text{if $k=0,2,$}\\
r(pq)-2r(p)-2r(q)+1 & \text{if $k=3$}, \\
0 & \text{otherwise.}
 \end{cases}
 $$
(Corollary \ref{r(pq)}).
These last two results on the ranks of the groups $H_3(\SL_2(\z[1/n]),\z)$ and $H_3(\SL_2(\z[1/(pq)]),\z)$ seem to be new.

Finally, using an induction argument, we show in Theorem \ref{teo: calculation H1} that for any positive integer $n$, we have 
$$H_1(\SL_2(\z[1/n]),\z) \simeq \begin{cases}
        0, &\text{ if } 2\mid n, 3\mid n \\
        \z/3, &\text{ if } 2\mid n, 3\nmid n \\
        \z/4, &\text{ if } 2\nmid n, 3\mid n \\
        \z/12, &\text{ if } 2\nmid n, 3\nmid n. 
\end{cases}$$
This last result has already been proven in \cite{B-E2025}, but the proof given in Chapter 4 is new and relies on the main results of this dissertation.

We conclude the introduction by outlining the structure of the dissertation.

In Chapter 1 we  review some of the necessary background on CW-complexes, trees and group homology necessary to understand the following chapters. 

Chapter 2 introduces spectral sequences and describes how a filtration of a double complex leads to a sequence converging to the total homology 
of that complex, which is then used to construct the spectral sequence converging to the homology of a group via its action on a contractible CW-complex. 

Considering that the homology groups of $\Gamma_0(n)$ make up the first page of the spectral sequence converging to $\SL_2(\z[1/n])$, Chapter 3 elucidates the construction of a tree due to \cite{ww1998} that allows for the calculation of such homology groups.

Finally, Chapter 4 reviews the construction of the tree $B_p$ associated with the $p$-adic valuation on $\q$ on which $\SL_2(\z[1/p])$ acts. This is a particular case of the tree of $\SL_2$ over a local field due to \cite{serre1980}. The second section of this chapter then considers the action of $\SL_2(\z[1/n])$ on the product $B_{p_1}\times \dots\times B_{p_k}$ of the trees $B_{p_i}$ where $p_1,\dots,p_k$ are the prime factors of $n$ which gives the desired spectral sequence. Section 3 then analyses this spectral sequence to obtain some results on $H_k(\SL_2(\z[1/n]),\z)$. Finally, the last section provides a description of $H_1(\SL_2(\z[1/n]),\z)$.

\chapter{Preliminaries}
\label{chapter:preliminaries}
In this chapter, we review some results on CW-complexes which are important for this dissertation. The first section considers CW-complexes more generally, whereas the second section focuses on trees, which are contractible complexes of dimension 1. The complexes appearing in the following chapters are either trees or products of trees. 

Considering that the theme of this dissertation is group homology, the last section of this chapter states the definition of the homology of a group and presents some important results on the topic. 

\section{CW-Complexes and cellular homology}

A CW-complex can be though of as a topological space that is built from $n$-disks $D^n$ starting from a set of points. For a more in-depth view of CW-complexes in general, we recommend \cite{geoghegan2008} and \cite{hatcher2002}. These are also the main references used for this section.

\begin{definicao}
    A \textit{CW-complex} $X$ is a union $\displaystyle\bigcup_{n\geq 0}X^n$ where $X^0$ is a discrete set and $X^n$ is formed from $X^{n-1}$ by attaching $n$-cells $e^n_\alpha$ using maps $\varphi^n_\alpha: S^{n-1}\rightarrow X^{n-1}$, that is, $X^n = (X^{n-1}\sqcup_\alpha D^n_\alpha)/\sim$ where $\sqcup$ denotes a disjoint sum, $D^n_\alpha$ is a copy of the \textit{$n$-disk} $D^n$ and $\sim$ is the equivalence relation defined by $x\sim \varphi^n_\alpha(x)$ whenever $x\in S^{n-1}_\alpha$.
\end{definicao}

The sets $X^n$ are called the \textit{$n$-skeletons} of $X$. The CW-complex $X$ is said to be of \textit{finite dimension} if there is some $n$ such that $X^r=X^n$ for all $r\geq n$, that is, if $X$ coincides with its $n$-skeleton for some $n$. If $n$ is the smallest integer such that $X=X^n$, we say that the \textit{dimension} of $X$ is $n$ and write $\dim X = n$. If no such $n$ exists, we say $\dim X=\infty$. 

Intuitively, the process of building a CW-complex starts with the 0-skeleton, which is a collection of points. The maps $\varphi^n_\alpha$ are maps from $S^{n-1}$ to $X^{n-1}$ which can be understood as giving instructions for each $n$-disk $D^n_\alpha$ as to which $(n-1)$-cells the border of $D^n_\alpha$ will be glued to and how this gluing will occur.

Finally, if $X$ has a finite number of cells, we say that $X$ is a \textit{finite} CW-complex.

The topology given to a CW-complex $X$ is the weak topology with respect to the $n$-skeletons $X^n$, that is, a set $A \subset X$ is open if and only if $A \cap X^n$ is
open in $X^n$ for each n.

\begin{exemplo}
    A \textit{torus} can be constructed as a CW-complex with one 0-cell, two 1-cells and a single 2-cell, as illustrated in figure \ref{fig:torus} below. This is a finite CW-complex of dimension 2.

    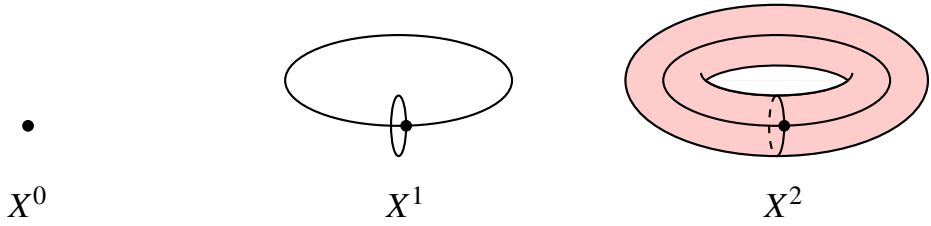
\begin{figure}[H]
        \centering
        \caption{CW-complex structure for a torus}
        \begin{tikzpicture}
            \draw[fill=black, shift={(-10,0)}] (0.1,-0.6) circle (2pt);
            
            \node[black,scale=1.1] at (-9.9,-1.6) {$X^0$};
            
            \draw[thick, shift={(-5,0)}] (0.1,-0.6) arc [x radius=0.1,y radius=0.4, start angle = 0, end angle=360];

            \draw[thick, shift={(-5,0)}] (1.5,0) arc [x radius=1.5,y radius=0.6, start angle = 0, end angle=360];

            \draw[fill=black, shift={(-5,0)}] (0.1,-0.6) circle (2pt);

            \node[black,scale=1.1] at (-4.9,-1.6) {$X^1$};
            
            \draw[black, thick, fill=red, opacity=0.2] (2,0) arc [x radius=2,y radius=1, start angle = 0, end angle=360];
            \draw[black, thick] (2,0) arc [x radius=2,y radius=1, start angle = 0, end angle=360];
            \draw[thick] (1,0.1) arc [x radius=1,y radius=0.3, start angle = 0, end angle=-180];
            \draw[thick, fill=white] (0.94281,0) arc [x radius=1,y radius=0.3, start angle = 20, end angle=160];
            \draw[thick, fill=white] (0.94281,0) arc [x radius=1,y radius=0.3, start angle = -20, end angle=-160];

            \draw[thick] (1.5,0) arc [x radius=1.5,y radius=0.6, start angle = 0, end angle=360];

            \draw[thick] (0,-1) arc [x radius=0.1,y radius=0.4, start angle = -90, end angle=90];

            \draw[thick, dashed] (0,-0.2) arc [x radius=0.1,y radius=0.4, start angle = -270, end angle=-90];

            \draw[fill=black] (0.1,-0.6) circle (2pt);

            \node[black,scale=1.1] at (0.1,-1.6) {$X^2$};
        \end{tikzpicture}
        \label{fig:torus}
        \fautor
    \end{figure}
\end{exemplo}

\begin{definicao}
    Given a CW-complex $X$ and an $n$-cell $e^n_\alpha$, the \textit{characteristic map} $h^n_\alpha: D^n \rightarrow e^n_\alpha$ is the map defined as the composition $$D^n_\alpha \hookrightarrow X^{n-1}\sqcup_{\beta}D^n_\beta\rightarrow X^n$$ where the map on the right is the quotient map. 
\end{definicao}

Note that the restriction of $h^n_\alpha$ to $S^{n-1}$ must coincide with the gluing map $\varphi^n_\alpha$. This implies that a complex can also be defined via its characteristic maps.

Given a $n$-cell $e^n_\alpha$ with characteristic map $h^n_\alpha$ in a CW-complex $X$, we denote by $\overset{\bullet}{e^n_\alpha}, \overset{\circ}{e^n_\alpha}$ the sets $h^n_\alpha(S^{n-1})$ and $h^n_\alpha(D^n\backslash S^{n-1})$. Intuitively, $\overset{\bullet}{e^n_\alpha}$ and $\overset{\circ}{e^n_\alpha}$ are the border and the interior of $e^n_\alpha$.

\begin{proposicao}\label{prp: product cells}
    Let $X,Y$ be the two CW-complexes where both $X$ and $Y$ have countably many cells. Then the product $X\times Y$ has a CW-complex structure where the $n$-cells are products of cells in the original complexes.
\end{proposicao}

\begin{demonstracao}
    Suppose $e^n_\alpha,f^m_\beta$ are the cells of $X$ and $Y$ with characteristic maps $h^n_\alpha,k^m_\beta$. Note that the union of the products $e^n_\alpha\times f^m_\beta$ is all of $X\times Y$. Since $D^n\simeq I^n$ and $I^n\times I^m\simeq I^{n+m}$, we have that $D^{n}\times D^m\simeq D^{n+m}$. Let $r_n$ be such a homeomorphism. Then the composition \begin{center}
        \begin{tikzcd}[column sep=8ex]
            D^{n+m} \arrow[r, "(r_{n+m})^{-1}"] & D^n\times D^m \arrow[r, "h^n_\alpha\times k^m_\beta"] & e^n_\alpha\times f^m_\beta
        \end{tikzcd}
    \end{center} works as a characteristic map for the cell $e^n_\alpha,f^m_\beta$. This gives $X\times Y$ a CW-complex structure.

    The condition of both $X$ and $Y$ having countably many cells comes from the fact that weak topology of the cell complex $X\times Y$ obtained from the process above is sometimes finer the product topology, and this condition suffices to guarantee the topologies will coincide (see Theorem A.6 of \cite{hatcher2002}). Another sufficient condition is for either $X$ or $Y$ to be locally compact.
\end{demonstracao}


For the remainder of this section we fix for each $n$ a homeomorphism $k_n: B^n/S^{n-1} \rightarrow S^n$. 

\begin{definicao}
    Let $e^n_\alpha$ be an $n$-cell of a CW-complex $X$. If $n=0$, an \textit{orientation} for $e^n_\alpha$ is a choice between 1 and $-1$. If $n>0$, an orientation for $e^n_\alpha$ is a choice of an equivalence class of characteristic maps under homotopy.
\end{definicao}

\begin{proposicao}
    There are only two choices of orientation for a cell $e^n_\alpha$.
\end{proposicao}

\begin{demonstracao}
    The case $n=0$ is obvious. If $h_1,h_2: (B^n,S^{n-1})\rightarrow (e^n_\alpha,\overset{\bullet}{e^n_\alpha})$ are two characteristic maps for $e^n_\alpha$, they induce maps $h_1',h_2': B^n/S^{n-1} \rightarrow e^n_\alpha/\overset{\bullet}{e^n_\alpha}$ so that the composition $$k_n\circ (h_2')^{-1}\circ h_1'\circ (k_n)^{-1}: S^n\rightarrow S^n$$ is a homeomorphism, thus, has degree $\pm 1$. So, it must be homotopic to either $\id_{S^n}$ or the canonical homeomorphism of degree $-1$. It can be shown that $h_2,h_1$ are equivalent if and only if the degree of that composition is 1 (see Lemma 2.5.5 of \cite{geoghegan2008}).

    Then, fixing a characteristic map $h$, if $h_1,h_2$ are two  characteristic maps not equivalent to $h$, we have $$k_n\circ(h_2')^{-1}\circ h_1'\circ (k_n)^{-1} = (k_n\circ(h_2')^{-1}\circ h'\circ(k_n)^{-1})\circ(k_n\circ (h')^{-1}\circ h_1'\circ (k_n)^{-1}).$$ Since the maps $k_n\circ(h_2')^{-1}\circ h'\circ(k_n)^{-1}$ and $k_n\circ (h')^{-1}\circ h_1'\circ (k_n)^{-1}$ have degree $-1$, the degree of $k_n\circ(h_2')^{-1}\circ h_1'\circ (k_n)^{-1}$ is $(-1)^2=1$. Thus, any two maps not equivalent to $h$ are equivalent amongst themselves. So the two classes are given by $h$ and any other map not equivalent to $h$.
\end{demonstracao}

\begin{definicao}
    An \textit{oriented CW-complex} is a CW-complex with a choice of orientation for each cell.
\end{definicao}

If $X$ is an oriented CW-complex and $e^n_\alpha,e^{n-1}_\beta$ are an $n$-cell and an $(n-1)$-cell of $X$ with characteristic maps $h_\alpha,h_\beta$ respectively, one can consider the following diagram, where hooked arrows represent inclusions, $s_{\alpha,\beta}$ is a composition of the inclusion $\overset{\bullet}{e^n_\alpha}\hookrightarrow X^{n-1} \rightarrow X^{n-1}/(X^{n-1}\backslash\overset{\circ}{e^{n-1}_\beta})$ and $r_\beta$ is the induced homeomorphism making the diagram below commute.  \begin{center}
    \begin{tikzcd}
        & X^{n-1} \arrow[d] & e^{n-1}_\beta \arrow[l, hook] \arrow[d] & B^{n-1} \arrow[l,"h_\beta"] \arrow[d] \\
        \overset{\bullet}{e^n_\alpha} \arrow[ur, hook] \arrow[r, "s_{\alpha,\beta}", swap] & X^{n-1}/(X^{n-1}\backslash\overset{\circ}{e^{n-1}_\beta}) & e^{n-1}_\beta/\overset{\bullet}{e^{n-1}_\beta} \arrow[l, "r_\beta"] & B^{n-1}/S^{n-2} \arrow[l, "h_\beta'"] \arrow[d, "k_{n-1}"] \\
        S^{n-1} \arrow[u, "h_\alpha"] & & & S^{n-1}
    \end{tikzcd}
\end{center}

\begin{definicao}
    Let $n\geq 2$. If $e^n_\alpha,e^{n-1}_\beta$ are an $n$-cell and an $(n-1)$-cell of an oriented CW-complex $X$ with characteristic maps $h_\alpha,h_\beta$ respectively, then the \textit{incidence number} $[e^n_\alpha: e^{n-1}_\beta]$ is defined as the degree of the map $$k_{n-1}\circ (h'_\beta)^{-1}\circ (r_\beta)^{-1}\circ s_{\alpha,\beta}\circ (h_\alpha\rvert_{S^{n-1}}): S^{n-1}\rightarrow S^{n-1}$$

    For the case $n=1$, we define the incidence number by $$[e^1_\alpha:e^0_\beta] := \begin{cases}
        \varepsilon,& \text{ if } h_\alpha(1) = e^0_\beta \text{ and } h_\alpha(-1)\neq e^0_\beta \\
        -\varepsilon,& \text{ if } h_\alpha(1) \neq e^0_\beta \text{ and } h_\alpha(-1)= e^0_\beta \\
        0,& \text{ if } h_\alpha(1) = e^0_\beta \text{ and } h_\alpha(-1)= e^0_\beta \\
        0,& \text{ if } h_\alpha(1) \neq e^0_\beta \text{ and } h_\alpha(-1)\neq e^0_\beta \\
    \end{cases}$$ where $h_\alpha$ is the characteristic map of $e^1_\alpha$ and $\varepsilon\in\{\pm1\}$ is the orientation of $e^0_\beta$
\end{definicao}

\begin{proposicao}
    The incidence number $[e^n_\alpha: e^{n-1}_\beta]$ depends on the orientation but not the choice of maps $h_\alpha,h_\beta$ in their equivalence classes.
\end{proposicao}

\begin{demonstracao}
    See Propositions 2.5.6 and 2.5.7 of \cite{geoghegan2008}.
\end{demonstracao}

Geometrically, this incidence number behaves like a "winding number" of sorts, measuring how many times the border of $e^n_\alpha$ is being wrapped around $e^{n-1}_\beta$.

Note that if $e^{n-1}_\beta$ is not a subset of $e^n_\alpha$, this map will not be surjective, so it will have degree zero.

\begin{definicao}\label{def: cell complex}
    Let $X$ be an oriented CW-complex and $R$ a commutative ring. For $n<0$, let $C_n(X,R)=0$. For $n\geq 0$, let $C_n(X,R)$ be the free $R$-module generated by the set of oriented $n$-cells of $X$ and let $\partial_n: C_n(X,R)\rightarrow C_{n-1}(X,R)$ be the map defined by $\partial_0(e^0_\alpha)=0$, and $$\partial_n e^n_\alpha = \sum_{\beta}[e^n_\alpha:e^{n-1}_\beta]e^{n-1}_\beta$$ for $n>0$.

    Then $C_\bullet(X,R) = (C_n(X,R),\partial_n)_{n\in\z}$ is a chain complex, called the cellular complex of $X$ with coefficients in $R$.
\end{definicao}

Note that since $D^n$ is compact, the image $e^n_\alpha = h_\alpha(D^n)$ under the continuous map $h_\alpha$ is also compact, so it will intersect only a finite number of $(n-1)$-cells. Thus, the sum on the border map is always finite, since $[e^n_\alpha:e^{n-1}_\beta] =0 $ if $e^{n-1}_{\beta}$ is not in $\overset{\bullet}{e^n_\alpha}$.

\begin{definicao}
    Let $X$ be an oriented CW-complex and $R$ a commutative ring. The $n$-th \textit{cellular homology} of $X$, denoted by $H_n(X,R)$, is the $n$-th homology of the cellular complex of $X$ with coefficients in $R$.
\end{definicao}

\begin{proposicao}
    The cellular homology of $X$ coincides with its singular homology.
\end{proposicao}

\begin{demonstracao}
    This proof requires a third complex for $X$, obtained by considering the $n$-th skeletons $X^n$ as a filtration of $X$ and letting $C'_n(X;R)$ be the relative singular homology with $R$-coefficients $H_n^S(X^n,X^{n-1};R)$ of the pair $X^n,X^{n-1}$, with the boundary maps being the compositions of the connecting homomorphism $H_n^S(X^n,X^{n-1};R) \rightarrow H_{n-1}^S(X^{n-1};R)$ and the map induced by a projection $H_{n-1}^S(X^{n-1};R)\rightarrow H_{n-1}^S(X^{n-1},X^{n-2};R)$.

    In the introduction to Chapter 2.6 of \cite{geoghegan2008}, there is a proof that this complex is isomorphic to the cellular complex previously defined here. Then, Theorem 2.35 of \cite{hatcher2002} shows that this third complex is such that its homology groups are isomorphic to the singular homology groups of $X$ with coefficients in $R$.
\end{demonstracao}

The cellular complex for a CW-complex is often preffered for being considerably simpler than the singular one, allowing for easier calculations of homology groups.

\begin{exemplo}
    The Klein bottle $K$ has a CW-complex structure with one 0-cell, two 1-cells and a single 2-cell. To understand the attaching maps used to calculate the incidence numbers, it is useful to view it as a quotient space. In Figure \ref{fig:bottle} below, the four vertices of the square are glued together and thus represent the single point of the 0-skeleton. The 1-skeleton will consist of the edges $a,b$ once glued as described in the figure.

    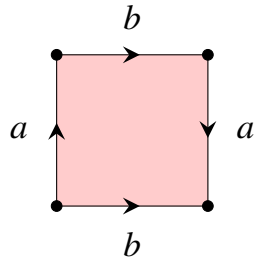
\begin{figure}[H]
        \centering
        \caption{CW-complex structure for Klein bottle $K$ via quotient topology}
        \begin{tikzpicture}[scale=0.5]
            \draw[red, ultra thin, fill, opacity=0.2] (-2,-2) -- (2,-2) -- (2,2) -- (-2,2) -- cycle;
            \draw[black] (-2,-2)--(2,-2) [arrow inside={end=stealth,opt={black,scale=2}}{0.55}];
            \node[scale=1.1, black] at (0,-3) {$b$};
            \draw[black] (2,2)--(2,-2) [arrow inside={end=stealth,opt={black,scale=2}}{0.55}];
            \node[black, scale=1.1] at (3,0) {$a$};
            \draw[black] (-2,2)--(2,2) [arrow inside={end=stealth,opt={black,scale=2}}{0.55}];
            \node[black,scale=1.1] at (0,3) {$b$};
            \draw[black] (-2,-2)--(-2,2) [arrow inside={end=stealth,opt={black,scale=2}}{0.55}];
            \node[black,scale=1.1] at (-3,0) {$a$};
            \draw [fill=black] (-2,-2) circle (4pt);
            \draw [fill=black] (-2,2) circle (4pt);
            \draw [fill=black] (2,-2) circle (4pt);
            \draw [fill=black] (2,2) circle (4pt);
        \end{tikzpicture}
        \label{fig:bottle}
        \fautor
    \end{figure}

    Note that the edges marked $a$ and $b$ are both attached to the single point as a wedge of circles, so the incidence numbers are 0.

    The 2-cell on the left can be seen as attached going clockwise around the edge (this is just a choice of orientation) which, abusing notation, corresponds to a path $abab^{-1}$, that is, is winds around $b$ once in each direction, cancelling out, and twice around $a$ in the same direction. Thus, if we call this cell $e^2$ and the cells associated to $a$ and $b$ as $e^1_a,e^1_b$, then the incidence numbers are $[e^2:e^1_a] = 2$, $[e^2:e^1_b]=0$. Thus, the cellular chain complex $C_\bullet(K,\z)$ is $$0 \rightarrow \z \rightarrow \z\oplus\z \rightarrow \z \rightarrow 0$$ where the first map is defined by $1\mapsto (2,0)$ and the second is trivial. Then, we conclude that $$H_n(K,\z) \simeq \begin{cases}
        \z,& \text{ if } n=0 \\
        \z\oplus(\z/2\z),& \text{ if }n=1 \\
        0,& \text{ otherwise}
    \end{cases}.$$
\end{exemplo}

\section{Trees}

CW-complexes of dimension 1 are particularly useful in this thesis. These complexes are in one-to-one correspondence with graphs, and the contractible dimension 1 complexes are in one-to-one correspondence with trees, that is, graphs without cycles.

In this section, we recall a few important results from the theory of trees. For further information on this, we reccomend \cite{serre1980}. This is the main reference for this section.

\begin{definicao}
    A \textit{graph} $\Gamma$ consists of a set $X = \text{vert}(\Gamma)$, a set $Y = \text{edge}(\Gamma)$ and two maps $Y\rightarrow X\times X$, $Y\rightarrow Y$ such that $y\mapsto (o(y),t(y))$, $y\mapsto \overline{y}$ satisfying the conditions that $\overline{\overline{y}}=y$ and $o(y)=t(\overline{y})$.
\end{definicao} 

In this definition, the first map associates an edge to its origin and terminus and the second map associates an edge to its inverse.

Any graph can be identified with its \textit{geometrical realization}, which is formed as the disjoint union of $X$ and $Y \times [0, 1]$ (where $X$ and $Y$ are provided with the discrete topology) with an equivalence relation for which $(y,0) \sim o(y)$, $(y, 1)\sim t(y)$ and $(y, t) \sim (\overline{y}, 1 - t)$ for $y\in Y, t\in [0,1]$. 
    
Note that this defines a CW-complex structure of dimension 1, where the 0-skeleton is the set $X$ and the 1-skeleton is constructed by adding the copies $D^1_y  \simeq y\times[0,1]$ to the 0-skeleton via the relations described above. Conversely, any 1-dimensional CW-complex defines a graph where the verices are 0-cells, the edges are 1-cells, and the gluing maps describe the origin and terminus of each vertex.

\begin{definicao}
    A \textit{path} in a graph is a sequence $y_1\dots y_n$ of edges where $t(y_i)=o(y_{i+1})$ for all $i\in\{1,\dots,n-1\}$. A \textit{loop} is a path $y_1\dots y_n$ where $t(y_n)=o(y_1)$. A path is said to \textit{backtrack} if $y_i = \overline{y_{i+1}}$ for some $i$.
\end{definicao}

\begin{definicao}
    A graph $\Gamma$ is said to be connected if given any vertices $P,Q$ of $\Gamma$, there exists a path $y_1\dots y_n$ in $\Gamma$ where $o(y_1)=P,t(y_n)=Q$.
\end{definicao}

\begin{definicao}
    A \textit{tree} is a connected graph in which there are no loops without backtracking.
\end{definicao}

\begin{proposicao}
    A graph $\Gamma$ is a tree if and only if given any two vertices $P,Q$ of $\Gamma$, there is a unique path without backtracking from $P$ to $Q$.
\end{proposicao}

\begin{demonstracao}
    For this proof, paths are assumed to not backtrack.
    
    Note that trees are path-connected. Thus there is always a path between two vertices of a tree. If there are two distinct paths $y_1\dots y_n$ and $y_1'\dots y_k'$ from $P$ to $Q$, then $y_1\dots y_n\overline{y_1'}\dots\overline{y_k'}$ is a loop at $P$, which would be a contradiction.

    On the other hand, if there is always a path between two vertices, the graph is path-connected and there can be no loops, since a path from $P$ to itself can be broken as a path from $P$ to some other vertex $Q$ and another path from $Q$ back to $P$, so that these determine two distinct paths from $P$ to $Q$.
\end{demonstracao}

This implies that there is a well-defined notion of distance between vertices of a tree, given by the number of edges in the unique path connecting two vertices.

If $P$ is a vertex of a tree $\Gamma$, for each integer $n\geq 0$, let $X_n$ be the set of vertices at distance less or equal to $n$ from $P$. Note that $X_0=\{P\}$. If $Q \in X_n$, with $n\geq 1$, there is a unique path $(y_1,\dots,y_i)$ from $P$ to $Q$ with $i\leq n$. Let $Q'$ be $o(y_i)$. Then $Q'$ is the unique vertex adjacent to $Q$ at distance strictly less then $n$ from $P$. This defines a map $f_n: X_n \rightarrow X_{n-1}$ that maps each $Q$ to the corresponding $Q'$. Then, we have an inverse system $$\dots\rightarrow X_n \overset{f_n}{\longrightarrow} X_{n-1} \rightarrow \dots \rightarrow X_1 \overset{f_1}{\longrightarrow}X_0$$ which allows us to reconstruct $\Gamma$, since the vertices of $\Gamma$ are the elements in $\bigcup_{n\geq 0} X_n$ and the edges are $\{Q,f_n(Q)\}$ with $Q\in X_n$.

\begin{figure}[H]
    \centering
    \caption{Tree as an inverse limit}
    \begin{tikzpicture}[scale=0.7]
        \draw[thick] (0,0)--(8,2);
        \draw[thick] (-4,5)--(8,2);
        \draw[thick] (0,2)--(4,1);
        \draw[thick] (0,4)--(-4,3);
        \draw [fill=white] (0,0) circle (4pt);
        \draw [fill=white] (4,1) circle (4pt);
        \draw [fill=white] (0,2) circle (4pt);
        \draw [fill=white] (0,4) circle (4pt);
        \draw [fill=white] (4,3) circle (4pt);
        \draw [fill=white] (8,2) circle (4pt);
        \draw [fill=white] (-4,3) circle (4pt);
        \draw [fill=white] (-4,5) circle (4pt);

        \node[black,scale=1.1] at (-4,-2) {$X_3$};
        \node[black,scale=1.1] at (0,-2) {$X_2$};
        \node[black,scale=1.1] at (4,-2) {$X_1$};
        \node[black,scale=1.1] at (8,-2) {$X_0$};

        \draw[thick, ->] (-3,-2)--(-1,-2);
        \draw[thick, ->] (1,-2)--(3,-2);
        \draw[thick, ->] (5,-2)--(7,-2);
    \end{tikzpicture}
    \label{fig:inv lim tree}
    \fautor
\end{figure}
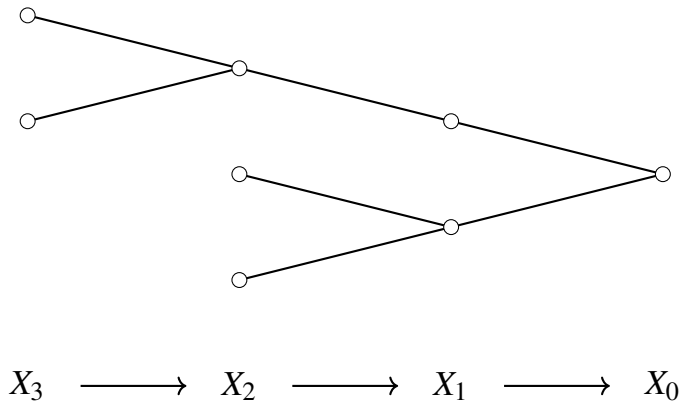

\begin{teorema}
    There is a one-to-one correspondence between trees and contractible CW-complexes of dimension $\leq 1$.
\end{teorema}

\begin{demonstracao}
    Recall that any CW-complex defines a graph where the edges are given by 1-cells, the vertices by 0-cells and the gluing map describes the origin and terminus of each cell. Conversely, any graph defines a CW-complex in the same way.

    The previous depiction of a tree as an inverse system shows that trees are contractible complexes, since each edge $\{Q,f_n(Q)\}$ with $Q\in X_n$ can be contracted into the vertex $f_n(Q)$, implying that the entire tree can be contracted into the vertex $P$ associated with $X_0$.

    On the other hand, if a CW-complexes of dimension $\leq 1$ is made into a graph, then contractibility implies no trivial loops, so the graph will be a tree.
\end{demonstracao}

\begin{definicao}
    Let $\Gamma$ be a graph, $P$ be a vertex of $\Gamma$ and let $Y_P$ be the set of edges $y$ such that $P = t(y)$. The number of elements in $Y_P$ is called the \textit{index} of $P$. If $n = 0$, we call $P$ \textit{isolated}. If $n\leq 1$, we call $P$ a \textit{terminal vertex}.
\end{definicao}

Note that if $\Gamma$ is a connected graph with more than one vertex, then any vertex $P$ is connected to some other vertex $Q$ by a path $y_1\dots y_n$ where either $t(y_n)=P$ or $o(y_1) = P$ (implying $ t(\overline{y_1})=P$). In other words, such a graph has no isolated vertices.

Denote by $\Gamma-P$ the subgraph of $\Gamma$ whose vertex set is $\text{vert}(\Gamma)\backslash\{P\}$ and edge set is $\text{edge}(\Gamma)\backslash(Y_p\cup\overline{Y_p})$.

If $P$ is a non-isolated terminal vertex, there is a single edge connecting $P$ to "the rest of the graph", that is, one can obtain $\Gamma$ from $\Gamma-P$ by adding $P$ and this single edge. 

Note also from our construction of trees as inverse systems that if a tree has more than 1 vertex, at least one vertex will be terminal (the one at greatest distance from a fixed edge $P$).

\begin{proposicao}
    Let $P$ be a non-isolated terminal vertex of a graph $\Gamma$. Then
    \begin{enumerate}
        \item $\Gamma$ is connected if and only if $\Gamma-P$ is connected.
        
        \item Every non-backtracking loop of $\Gamma$ is contained in $\Gamma-P$.
        
        \item $\Gamma$ is a tree if and only if $\Gamma-P$ is a tree. 
    \end{enumerate}
\end{proposicao}

\begin{demonstracao}
    Denote by $y_P$ the single edge that $P$ is a terminal vertex of.
    
    The first item is obvious. For the second item, if $y_1\dots y_n$ is a loop in $\Gamma$ not contained in $\Gamma-P$, then $y_i = y_P$ for some $i$, implying $o(y_{i+1})=t(y_i) = y_P$. But this implies that $t(\overline{y_{i+1}}) = y_P$ so $\overline{y_{i+1}}=y_P=y_i$, hence the loop backtracks. The third item follows from the second. 
\end{demonstracao}

\begin{corolario}\label{cor: edges}
    If $\Gamma$ is a tree with $n$ vertices, then $\Gamma$ has $n-1$ edges (counting $y$ and $\overline{y}$ as the same edge).
\end{corolario}

\begin{demonstracao}
    Note that if there are 2 vertices, $\Gamma$ being a tree forces there to be a single path between these two vertices, so there is exactly one edge which connects them.

    Suppose the statement is true for all $n< k$ and let $\Gamma$ be a tree with $k$ vertices. Then at least one vertex $P$ of $\Gamma$ is terminal and non-isolated (again, this can be constructed by fixing any vertex and picking the furthest vertex from it). From the last proposition, $\Gamma-P$ is also a tree and it has $k-1$ vertices. Thus from the induction hypothesis, $\Gamma-P$ has $k-2$ edges. Then $\Gamma$ has the $k-2$ edges of $\Gamma-P$ plus the single edge connecting $P$ to it, that is, $\Gamma$ has $k-1$ edges.
\end{demonstracao}

\begin{definicao}
    We say that $T$ is a \textit{tree in $X$} if $T$ is a tree and is a subcomplex of $X$; $T$ is a \textit{maximal tree} if there is no larger subcomplex of $X$ which is a tree.
\end{definicao}

\begin{proposicao}
    $T$ is a maximal tree of a connected non-empty graph $\Gamma$ if and only if $T$ contains all the vertices of $\Gamma$.    
\end{proposicao}

\begin{demonstracao}
    If $T$ doesn't contain all vertices of $\Gamma$, there is an edge $y$ starting at a vertex of $T$ and ending outside of $T$. Adding this vertex and the pair $y,\overline{y}$ is adding a terminal vertex, so this remains a tree and will contain $T$, contradicting maximality.

    On the other hand, if $T$ already contains all vertices of $\Gamma$, then only edges could be added to $T$, but that would create two distinct paths between two vertices, so any subgraph containing $T$ is not a tree.
\end{demonstracao}

Note that the proof of the last proposition also gives a way of constructing a maximal tree inside a graph by adding vertices.

The next proposition allows us to characterize the fundamental group of a CW-complex $X$ by identifying the 1-skeleton $X^1$, which is itself a complex of dimension 1, with its corresponding graph.

For this proposition, the notation $\langle W \mid R,\rho\rangle$ describes a presentation of a group, that is, $W$ is a set of generators, $R$ is an indexing set for relations and $\rho$ is a map from $R$ to the free group $F$ generated by the set $W$, whose elements we call "words". Denote by $N(R)$ the normal subgroup of $F$ generated by $\rho(R)$. Then $\langle W \mid R,\rho\rangle = F/N(R)$.

As an example, let $W = \{a\}$, $R = \{1\}$ and let $\rho$ be defined by $\rho(1) = a^3$. Then $F/N(R) = \langle a \rangle/\langle a^3 \rangle \simeq \z/3$, the cyclic group of order three. This corresponds to the presentation $\langle a \mid a^3=1\rangle$ of $\z/3$.

Given a 2-cell $e^2_\gamma$ of $X$, its border $\partial e^2_\gamma$ is given in $C_1(X;\z)$ as a sum of 1-cells, which can be identified with a path $e^1_1\dots e^1_n$ in $X^1$, which itself corresponds to a word in the free group generated by the 1-cells.

\begin{proposicao}
    Let $X$ be an oriented path connected CW complex, let $T$ be a maximal tree in $X^1$ and let $P$ be a vertex of $X^1$ (that is, $P$ is a 0-cell of $X$). Let $F$ be the free group generated by the set $W$ of (oriented) 1-cells of X. Let $R$ be the set of (oriented) 2-cells of $X$ and let $S$ be the set of (oriented) 1-cells of $T$. For each $e^2_\gamma \in R$, let $\tau(\gamma)$ be a path in $X^1$ representing $\partial(e^2_\gamma)$. Let $\rho : R \sqcup S \rightarrow F$ take each $e^2_\gamma \in R$ to the word in $F$ that corresponds to $\tau(\gamma)$, and take each $\sigma \in S$ to the one-letter word $\sigma\in F$. Then $\langle W \mid R \sqcup S, \rho\rangle$ is a presentation of $\pi_1(X, P)$.
\end{proposicao}

\begin{demonstracao}
    See Theorems 3.1.16 and 3.4.1 of \cite{geoghegan2008}. Theorem 3.1.16 proves that $\langle W \mid R \sqcup S, \rho\rangle$ is a representation of $\pi_1(X,P)$ in the combinatorial definition given by the author, and Theorem 3.4.1 shows that this definition coincides with the usual topological definition.
\end{demonstracao}

\begin{teorema}\label{teo: homology graph}
    Let $X$ be an oriented path-connected CW-complex of dimension 1. If $X$ has $v$ 0-cells and $e$ 1-cells, then $H_1(X,\z)\simeq \z^{e-v+1}$.
\end{teorema}

\begin{demonstracao}
    From the previous proposition, if $X$ has dimension 1, then it has no 2-cells and $\pi_1(X,P)$ is simply given as the free group generated by the cells outside of a maximal tree $T$. From Corollary \ref{cor: edges}, since $X$ has $v$ 0-cells, a maximal tree will have $v-1$ edges, so there are $e-(v-1)$ remaining edges generating $\pi_1(X,P)$. The theorem then follows from the fact that $H_1(X,\z)$ is the abelianization of $\pi_1(X,P)$ (see Theorem 3.1.19 of \cite{geoghegan2008}).
\end{demonstracao}

\section{Homology of groups}

We now turn to the question of homology of groups. The reader familiar with algebraic topology may have heard of the homology groups of a space, which can be defined either as the homology of a singular, simplicial or cellular chain complex, all of which will coincide. For the case of groups, one can also define homology via an associated chain complex. The main reference used for this section is \cite{brown1994}.

\begin{definicao}\label{def: action}
    Let $G$ be a group with neutral element $e$ and $X$ a set. An \textit{action} of $G$ on $X$ is a map $\varphi(G,X)\rightarrow X$ such that $\varphi(e,x)=x$ and $\varphi(g,\varphi(h,x)) = \varphi(gh,x)$ for all $x\in X$ and all $g,h\in G$.  
\end{definicao}

Given $x\in X$, the set $$G_x = \{g\in G \mid \varphi(g,x)=x\}\subseteq G$$ is called the \textit{stabilizer} of $x$. The set $$Gx = \{\varphi(g,x)\mid g\in G\}\subseteq X$$ is called the \textit{orbit} of $x$. One can show that $x\sim y \Leftrightarrow y\in Gx$ is an equivalence relation, that is, the orbits under an action form a partition of $X$. 

\begin{definicao}
    Let $G$ be a group. The group ring $\z G$ is defined as the ring whose underlying abelian group structure is $\displaystyle\bigoplus_{g\in G}\z$, that is, elements in $\z G$ can be understood as sums $\displaystyle\sum_{g\in G}n_gg$ with $n_g\in \z$, $n_g=0$ for all but a finite subset of $G$, where addition is coordinate-wise, as expected, $$\sum_{g\in G}n_gg+\sum_{g\in G}m_gg = \sum_{g\in G}(n_g+m_g)g$$ and multiplication is induced by the group operation, that is, $$\left(\sum_{g\in G}n_gg\right)\left(\sum_{h\in G}m_hh\right) := \sum_{g,h\in G}n_gm_h gh = \sum_{k\in G}\left(\sum_{g\in G}n_gm_{g^{-1}k}\right)k$$
\end{definicao}

We leave it to the reader to check that $\z G$ is indeed a ring. Note that $\z G$ will be commutative if and only if $G$ is abelian.

\begin{exemplo}
    Let $G$ be a finite cyclic group with $n$ elements. Then $G = \{1,t,\dots,t^{n-1}\}$ for some generator $t$, allowing us to write $\z G = \displaystyle\bigoplus_{i=0}^{n-1}t^i\z$. Then $\z G \simeq \displaystyle\frac{\z[T]}{\langle T^n-1\rangle}$ via the isomorphism $\displaystyle\sum_{i=0}^{n-1}a_it^i \mapsto \overline{\sum_{i=0}^{n-1}a_iT^i}$.
\end{exemplo}

\begin{definicao}
    Given any group $G$, the augmentation map $\varepsilon: \z G \rightarrow \z$ is the map defined by $$\sum_{g\in G}n_g g \mapsto \sum_{g\in G}n_g.$$ Since this sum is always finite, this map is well-defined.
\end{definicao}

\begin{exemplo}
    Let $G$ be a group. Define an action of $\z G$ on $\z$ by $$\left(\displaystyle\sum_{g\in G} n_gg\right) m = \varepsilon\left(\displaystyle\sum_{g\in G} n_gg\right) m = \left(\displaystyle\sum_{g\in G} n_g\right) m.$$ This action makes $\z$ into a left $\z G$-module. One can similarly define $\z$ as a right $\z G$-module. This is called the trivial action of $G$ on $\z$.
\end{exemplo}

\begin{definicao}
    Let $G$ be a group. The $n$-th integral homology of $G$ (or the $n$-th homology of $G$ with coefficients in $\z$) is defined as the abelian group $H_n(G,\z):=\text{Tor}_n^{\z G}(\z,\z)$. For simplicity, we usually denote $H_n(G,\z)$ by $H_n(G)$. For the defition of $\text{Tor}$ and its properties, see \cite{rotman2009}.
\end{definicao}

\begin{exemplo}\label{exemplo1}
    Let $G$ be a cyclic group with $n$ elements, so that $\z G\simeq \z[T]/\langle T^n-1\rangle$ and let $t$ be a generator for $G$. Note that $$(t-1)(a_0+a_1t+\dots+a_{n-1}t^{n-1}) = (a_{n-1}-a_0) + (a_0-a_1)t+ \dots + (a_{n-2}-a_{n-1})t^{n-1}.$$ In particular, $(1+t+\dots+t^{n-1})(t-1)=0$ and if $\varepsilon$ is the augmentation map $\sum a_it^i \mapsto \sum a_i$, then $\varepsilon((t-1)(a_0+a_1t+\dots+a_{n-1}t^{n-1})) = a_{n-1}-a_0+a_0-a_1+\dots+a_{n-2}-a_{n-1} = 0$.
    
    Denote by $N$ the element $1+t+\dots+t^{n-1}\in\z G$. Then the sequence \begin{center}
        \begin{tikzcd}
            F: & \dots \arrow[r, "t-1"] & \z G \arrow[r, "N"]& \z G \arrow[r, "t-1"]& \z G\arrow[r, "\varepsilon"] & \z \arrow[r] & 0
        \end{tikzcd}
    \end{center} forms a chain complex. We show that this complex is exact, thus being a projective resolution of $\z$ over $\z G$ suitable to calculate $\Tor_n^{\z G}(\z,\z)$.

    Let $\sum a_it^i\in\ker\varepsilon$, that is, $\sum a_i = 0$. In particular, $a_0+a_1+\dots+a_{n-2}=-a_{n-1}$. Let $a':= -a_0-(a_0+a_1)t-\dots-(a_0+\dots+a_{n-2})t^{n-2}$. Then \begin{align*}
        (t-1)a' &= (t-1)(-a_0-(a_0+a_1)t-\dots-(a_0+\dots+a_{n-2})t^{n-2}) \\
        &= a_0+a_1t+\dots+a_{n-2}t-(a_0+\dots+a_{n-2})t^{n-1} \\
        &= a_0+a_1t+\dots+a_{n-2}t+a_{n-1}t^{n-1}. 
    \end{align*} So that $\ker\varepsilon\subseteq \im(t-1)$.

    If $a = \sum_{i=0}^{n-1} a_it^i\in\ker(t-1)$, then $a_{n-1}-a_0 = a_0-a_1 = \dots = a_{n-2}-a_{n-1}=0$, so that $a_i=a_0$ for all $i$. This implies that $a = a_0+a_0t+\dots+a_0t^{n-1} = N\cdot a_0$, so $\ker(t-1)\subseteq \im\,N$. 

    Also, 
    
    \vspace{-20pt}\begin{align*}
        N\cdot a &= (1+t+\dots+t^{n-1})(a_0+a_1t+\dots+a_{n-1}t^{n-1}) \\
        &= a_0+a_1t+\dots+a_{n-2}t^{n-2}+a_{n-1}t^{n-1} \\
        &+a_{n-1}+a_0t+\dots+a_{n-3}t^{n-2}+a_{n-2}t^{n-1} \\
        &+ a_{n-2}+a_{n-1}t+\dots+a_{n-4}t^{n-2}+a_{n-3}t^{n-1} \\
        &\hspace{100pt}\vdots \\
        &+ a_1+a_2t+\dots+a_{n-1}t^{n-2}+a_{0}t^{n-1} \\
        &= \varepsilon(a)+\varepsilon(a)t+\dots+\varepsilon(a)t^{n-2}+\varepsilon(a)t^{n-1}.
    \end{align*} So $$a\in\ker(N) \implies \varepsilon(a)=0 \implies a \in \ker\varepsilon \implies a \in \im(t-1).$$ Thus, our chain complex is indeed exact.

    To use $F$ to calculate $\Tor_n^{\z G}(\z,\z)$, we consider the deleted resolution $F\otimes_{\z G}\z$. Let $\phi$ be the isomorphism $\z G\otimes_{\z G}\z \simeq \z$ given by $(a\otimes n)\mapsto \varepsilon(a)n$. We have that $$(\phi\circ(t-1)\otimes\id)(a\otimes m) = \phi((t-1)a\otimes m) = \varepsilon((t-1)a)\cdot m = 0$$ $$(\phi\circ N\otimes\id)(a\otimes m) = \phi(\varepsilon(a)N\otimes m) = n\cdot \varepsilon(a)\cdot m = n\cdot \phi(a\otimes m).$$ This implies, then, that we have the following commutative diagram: \begin{center}
        \begin{tikzcd}[column sep=8ex]
            \dots \arrow[r, "(t-1)\otimes\id"] & \z G\otimes_{\z G}\z \arrow[d, "\phi"] \arrow[r, "N\otimes\id"] & \z G\otimes_{\z G}\z \arrow[d, "\phi"] \arrow[r, "(t-1)\otimes\id"] & \z G\otimes_{\z G}\z \arrow[d, "\phi"] \arrow[r] & 0 \\
            \dots \arrow[r, "0"] & \z \arrow[r, "n"]& \z \arrow[r, "0"] & \z \arrow[r] & 0
        \end{tikzcd}
    \end{center} so that $\Tor_n^{\z G}(\z,\z)$, the $n$-th homology of the upper complex, coincides with the $n$-th homology of the inferior complex. Thus, $$H_n(G) \simeq \begin{cases}
        \z, &\text{ if } n=0, \\
        \z/n\z,& \text{ if } n \text{ odd}, \\
        0, &\text{ otherwise.}
    \end{cases}$$
\end{exemplo}

\begin{teorema}\label{teo: H0 and H1}
    Let $G$ be a group. Then $H_0(G) \simeq\z$ and $H_1(G) \simeq G^{\text{ab}} = G/[G,G]$, the abelianization of $G$.
\end{teorema}

\begin{demonstracao}
    See Chapter II.3 of \cite{brown1994}.
\end{demonstracao}

\begin{proposicao}\label{prp: surjection}
    If $f: G\rightarrow K$ is a surjective homomorphism, then the induced map $f_*: H_1(G) \rightarrow H_1(K)$ is also surjective.
\end{proposicao}

\begin{demonstracao}
    Since $H_1(G) \simeq  G^{\text{ab}} = G/[G,G]$ for any group $G$, it follows that $f_*: H_1(G)\rightarrow H_1(K)$ is simply the map $$f_*: \frac{G}{[G,G]}\rightarrow\frac{K}{[K,K]}$$ induced on the abelianizations, that is, $f_*(\overline{a}) = \overline{f(a)}$ so surjectivity of $f$ gives surjectivity of $f_*$.
\end{demonstracao}

\begin{definicao}
    An abelian group $M$ is said to be a \textit{$G$-module} if there is a $G$-action on $M$.
\end{definicao}

Observe that an action of $G$ on $M$ induces an action of $\z G$ on $M$ and vice-versa, so $G$-modules and $\z G$-modules coincide. 

In particular, the augmentation map $\varepsilon: \z G \rightarrow \z$ induces the trivial $G$-action on $\z$ making $\z$ a $G$-module.

\begin{definicao}
    Let $M$ be a left $G$-module and consider $\z$ as a right $G$-module with trivial action. We define the homology of $G$ with coefficients in $M$ as the abelian group $H_n(G,M) := \text{Tor}_n^{\z G}(\z,M)$. 
\end{definicao}

Note that for homology with coefficients, the difference will be that instead of tensoring a projective resolution $F$ of $\z$ over $\z G$ with the $G$-module $\z$, we tensor $F$ with $M$. 

\begin{proposicao}\label{prp: coinvariants}
    Let $G$ be a group and $M$ a $G$-module. Define the group of co-invariants $M_G$ by $$M_G = \displaystyle\frac{M}{\langle gm-g\mid g\in G,m\in M\rangle},$$ the largest quotient of $M$ on which $G$ acts trivially. Then $\z\otimes_{\z G}M\simeq M_G$.
\end{proposicao}

\begin{demonstracao}
    Recall that $\z$ is considered a $G$-module with trivial action. Then for all $g\in G,m\in M$, we have that $1\otimes gm = 1\cdot g \otimes m = 1\otimes m$. This allows us to define a map $M_G\rightarrow \z\otimes_{\z G}M$ by $\overline{m}\mapsto 1\otimes m$ which will be well defined since the generators $gm-m$ are such that $1\otimes m = 1\otimes gm$. On the other hand, the map $\z\times M \rightarrow M_G$ given by $(a,m)\mapsto a\overline{m}$ is bilinear, so it defines a map $\z\otimes_{\z G}M \rightarrow M_G$. These two maps are inverses of each other, proving the isomorphism.
\end{demonstracao}

\begin{proposicao}
    Let $f: G \rightarrow H$ be a group homomorphism. Then $f$ induces for all $n\geq 0$ and all $H$-modules $M$ a map $f_*: H_n(G,M)\rightarrow H_n(H,M)$.
\end{proposicao}

\begin{demonstracao}
    Note that $f$ induces a map $\z G\rightarrow \z H$ where each component $n_gg$ of a finite sum is mapped to $n_{g}f(g)$. Thus, any $H$-module can be made into a $G$-module via change of rings. So if $F,F'$ are resolutions of $\z$ over $\z G$ and $\z H$, respectively, then $F'$ can also be regarded as a  complex of $G$-modules. Since $F$ is a projective resolution, this will induce a chain map $f': F \rightarrow F'$ such that $f'(gx) = f(g)f'(x)$ (see Section II.6 of \cite{brown1994}). This induces a map on the invariants  $f'':(F_p)_G\rightarrow (F_p')_{H}$ given by $[x]_G \mapsto [f'(x)]_H$, since $[gx-x]_H\mapsto [f(g)f'(x)-f'(x)]_H = [0]_H$. Then $f''$ can be considered as a map $f'':F_p\otimes_{\z G}\z \rightarrow F_p'\otimes_{\z H}\z$. Furthermore, since $F_p\otimes_{\z G} M \simeq F_p\otimes_{\z G}(\z\otimes_{\z} M) \simeq (F_p\otimes_{\z G}\z)\otimes_{\z} M$, we have that $f''$ induces a well defined map $f_{\#}: F_p\otimes_{\z G} M \rightarrow F_p'\otimes_{\z H}M$. Descending to homology gives us $f_*: H_n(G,M)\rightarrow H_n(H,M)$.
\end{demonstracao}

If $i: H\rightarrow G$ is a group inclusion, the induced map $i_*: H_n(H) \rightarrow H_n(G)$ is called a \textit{corestriction map} and is often denoted by $i_* = \text{cor}^G_H$. 

\begin{proposicao}\label{prp: transfer}
    If $[G:H]<\infty$, there exists a \textit{transfer map} $\text{tr}^G_H: H_n(G)\rightarrow H_n(H)$, also called a \textit{restriction map} and denoted by $\text{res}^G_H$ such that $\text{cor}^G_H\circ\text{res}^G_Hz = [G:H]z$.
\end{proposicao}

\begin{demonstracao}
    There are actually five possible ways of defining such a map, one of which will be explored in the section on Equivariant Homology. To see all these constructions, we recommend Section III.9 of \cite{brown1994}.
\end{demonstracao}

\begin{exemplo}\label{exe: induced}
    We wish to show that the map induced on homology by the inclusion $i: \z/n\z \rightarrow \z/(nm)\z$ concides with multiplication by $m$ on odd homology and with the identity map otherwise. 

    Let $G = \langle t \mid t^{mn}=1\rangle \simeq \z/nm$ and $H= \langle s \mid s^n=1\rangle\simeq \z/n$, so that the inclusion $i$ maps $s$ to $t^m$, the element of order $n$ in $G$. This induces an inclusion $\z H \rightarrow \z G$ where $\sum a_is^i \mapsto \sum a_it^{mi}$.
    
    Let $N = 1+s+\dots+s^{n-1}$, $MN = 1+t+\dots+t^{mn-1}$ as in the resolutions for $G$ and $H$ discussed in Example \ref{exemplo1}. Let $M = 1+t+\dots+t^n\in\z G$ so that the map in the diagram below denoted by $M$ is the map that sends $p\in\z H$ to $i(p)(1+t+\dots+t^n)$ in $\z G$. \begin{center}
        \begin{tikzcd}
            F_{\bullet}: \arrow[d] &\dots \arrow[r, "s-1"] & \z H \arrow[r, "N"] \arrow[d, "i"] & \z H \arrow[d, "M"] \arrow[r, "s-1"] & \z H \arrow[r, "\varepsilon"] \arrow[d, "i"] & \z \arrow[d, "\id"] \\
            F'_{\bullet}: & \dots \arrow[r, "t-1"] & \z G \arrow[r, "MN"] & \z G \arrow[r, "t-1"] & \z G \arrow[r, "\varepsilon"] & \z \\
        \end{tikzcd}
    \end{center}

    Note that the first square obviously commutes, since the inclusion $i:\z H \rightarrow \z G$ doesn't affect coefficients. Then \begin{align*}
        M(t-1) &= (1+t+\dots+t^{m-1})(t-1) \\
        &= t+t^2+\dots+t^m-1-t-\dots-t^{m-1} \\
        &= t^m-1
    \end{align*} which is the image of $s-1$ under the inclusion $i$. Also, \begin{align*}
        (M\circ N)(p) &= (1+t+\dots+t^{m-1})(1+t^m+\dots+t^{m(n-1)})i(s) \\
        &= (1+t+\dots+t^{m-1}+t^m+t^{m+1}+\dots+t^{2m-1}+\dots+t^{m(n-1)}+\dots+t^{mn-1})i(s) \\
        &= MNi(s),
    \end{align*} so the diagram commutes. Thus, the maps $M,i$ form a chain map $\tau$ such that for all $h\in H, p\in\z H$, $\tau(hp) = i(h)\tau(p)$, that is, $\tau$ is an $H$-map and so the map $$\tau': F_\bullet \otimes_{\z H}\z \rightarrow F'_\bullet \otimes_{\z G}\z$$ is a well-defined chain map. Abusing notation, consider the commutative diagram: \begin{center}
        \begin{tikzcd}[column sep=8ex]
            \dots \arrow[r, "0"] & \z \arrow[r, "n"] \arrow[d] & \z\arrow[d] \arrow[r, "0"] & \z \arrow[d] \arrow[r] & 0 \\
            \dots \arrow[r, "(s-1)\otimes\id"] & \z H\otimes_{\z H}\z \arrow[d, "i\otimes\id"] \arrow[r, "N\otimes\id"] & \z H\otimes_{\z H}\z \arrow[d, "M\otimes\id"] \arrow[r, "(s-1)\otimes\id"] & \z H\otimes_{\z H}\z \arrow[d, "i\otimes\id"] \arrow[r] & 0 \\
            \dots \arrow[r, "(t-1)\otimes\id"] & \z G\otimes_{\z G}\z \arrow[d, ] \arrow[r, "MN\otimes\id"] & \z G\otimes_{\z G}\z \arrow[d, ] \arrow[r, "(t-1)\otimes\id"] & \z G\otimes_{\z G}\z \arrow[d, ] \arrow[r] & 0 \\
            \dots \arrow[r, "0"] & \z \arrow[r, "mn"]& \z \arrow[r, "0"] & \z \arrow[r] & 0
        \end{tikzcd}
    \end{center} so the chain map we are looking for is simply the composition of the vertical maps described. For any integer $a$, $$a \longmapsto 1\otimes_H a \longmapsto i(1)\otimes_G a = 1\otimes_G a \longmapsto \varepsilon(1)\cdot a = a,$$ $$a \longmapsto 1\otimes_H a \longmapsto M(1)\otimes_G a = (1+t+\dots+t^{m-1})\otimes_G a \longmapsto \varepsilon(1+t+\dots+t^{m-1})a = ma.$$

    Thus this chain map is multiplication by 1 or $M$, depending on the index, so that for a homology class, $[z]=[z]$ or $[z] = m[z]$. For the non-zero even homologies, the chain map is $\id$ but the homologies are trivial (see Example \ref{exemplo1}). For dimension 0, $i_*: H_0(H)\rightarrow H_0(G)$ will be the identity $\id:\z\rightarrow\z$. For the odd homologies, the chain map is multiplication by $m$, so the induced map $H_n(H)\rightarrow H_n(G)$ is also multiplication by $m$, as expected.

    As an addendum, note that

    \vspace{-20pt}\begin{align*}
        \sum_{i=0}^{mn-1}a_it^i &= \sum_{k=0}^{n-1}a_{km}t^{km}+t\left(\sum_{k=0}^{n-1}a_{km+1}t^{km}\right)+\dots+t^{m-1}\left(\sum_{k=0}^{n-1}a_{km+(m-1)}t^{km}\right) \\
        &= \sum_{k=0}^{n-1}a_{km}s^{k}+t\left(\sum_{k=0}^{n-1}a_{km+1}s^{k}\right)+\dots+t^{m-1}\left(\sum_{k=0}^{n-1}a_{km+(m-1)}s^{k}\right)
    \end{align*} which shows that $\z G = \displaystyle\bigoplus_{i=0}^{m-1}t^i\z H$.
\end{exemplo}

\begin{observacao}\label{zg free over zh}
    The addendum from the last example can be generalized: if $H$ is a normal subgroup of $G$, then $\z G = \displaystyle\bigoplus_{g_i\in E}g_i\z H$ where $E$ is a set of coset representatives (see Exercise 1 in Section I.3 of \cite{brown1994}).
\end{observacao}

\begin{definicao}\label{def: induced}
    Let $G$ be a group, $H$ a subgroup of $G$ and $M$ an $H$-module. The induced module $\Ind^G_HM$ is defined to be the $G$-module $$\Ind^G_HM \coloneq \z G\otimes_{\z H}M.$$
\end{definicao}

\begin{lema}[Shapiro]\label{shapiro}
    If $H$ is a subgroup of $G$ and $M$ is an $H$-module, then $H_*(H, M)\simeq H_*(G, \Ind^G_HM)$.
\end{lema}

\begin{demonstracao}
    Let $F$ be a projective resolution of $\z$ over $\z G$. Since $\z G$ is a direct sum of copies of $\z H$ as in the last example, then $F$ is also a projective resolution of $\z$ over $\z H$, which implies $H_*(H,M) = H_*(F\otimes_{\z H}M)$. The result then follows from the isomorphism $$F\otimes_{\z H}M \simeq F\otimes_{\z G}(\z G\otimes_{\z H} M)\simeq F\otimes_G\Ind^G_HM.$$ 
\end{demonstracao}

\begin{teorema}\label{teo: rational homology}
    For any group $G$, $H_n(G,\q) \simeq H_n(G)\otimes_{\z} \q$. In particular, $H_n(G,\q)$ is a $\q$-vector space.
\end{teorema}

\begin{demonstracao}
    This follows from the fact that $\q$ is flat as a $\z$-module, thus preserving exactness of sequences and so $H_n(C_\bullet\otimes_{\z}\q) \simeq H_n(C_\bullet)\otimes_{\z}\q$ for any complex $C_\bullet$ of $\z$-modules (see Section III.2 of \cite{brown1994}). Then, if $F$ is a deleted projective resolution of $\z$ over $\z G$, $$H_n(G,\q)=H_n(F\otimes_{\z G}\q) \simeq H_n((F\otimes_{\z G}\z)\otimes_{\z}\q) \simeq H_n((F\otimes_{\z G}\z))\otimes_{\z}\q = H_n(G)\otimes_{\z}\q.$$ 
\end{demonstracao}

We recall an important result from the theory of modules that will be used for the next lemma and throughout the last chapter of this dissertation.

\begin{teorema}[Structure theorem]\label{teo: structure}
    Let $R$ be a PID and $M$ a finitely generated module over $R$. Then there exists an integer $k$ and elements $r_1,r_2,\dots,r_n\in R$ with $r_1\mid r_2\mid\cdots\mid r_n$ such that $$M\simeq R^k\oplus \frac{R}{(r_1)}\oplus\frac{R}{(r_2)}\oplus\cdots\oplus \frac{R}{(r_n)}$$ and the ideals $(r_1),(r_2),\dots,(r_n)$ are uniquely determined.
\end{teorema}

\begin{demonstracao}
    See Theorems 7.3 and 7.7 in Chapter III.7 of \cite{lang2002}. Theorem 7.3 shows that $M$ is decomposed as a free part and a torsion part and Theorem 7.7 characterizes the torsion part of $M$ as described above.
\end{demonstracao}

\begin{lema}
    If $G$ is a finitely generated abelian group, then $\rank\, G = \dim_\q (G\otimes_{\z}\q)$.
\end{lema}

\begin{demonstracao}
    By the structure theorem for finitely generated modules over PIDs, there exists $r\geq 0$, $n_i\in\z$ such that $$G\simeq \z^r\oplus \frac{\z}{n_1\z} \oplus \frac{\z}{n_2\z} \oplus \dots\oplus \frac{\z}{n_k\z}$$ where $r = \rank\, G$. Then \begin{align*}
        G\otimes\q &\simeq \left(\z^r\oplus \frac{\z}{n_1\z} \oplus \frac{\z}{n_2\z} \oplus \dots\oplus \frac{\z}{n_k\z}\right)\otimes_{\z} \q \\
        &\simeq (\z^r\otimes_{\z} \q)\oplus \left(\frac{\z}{n_1\z}\otimes_{\z} \q\right) \oplus \left(\frac{\z}{n_2\z}\otimes_{\z} \q\right) \oplus \dots\oplus \left(\frac{\z}{n_k\z}\otimes_{\z} \q\right) \\
        &\simeq \q^r \oplus 0 \oplus 0 \oplus\dots\oplus 0 \\
        &\simeq \q^r.
    \end{align*} 
\end{demonstracao}

\begin{teorema}
    Let $G$ be a group. If $H_n(G)$ is finitely generated, then the rank of $H_n(G)$ equals $\dim_\q(H_n(G,\q))$.
\end{teorema}

\begin{demonstracao}
    From the previous results, $\rank\, H_n(G) = \dim_\q(H_n(G)\otimes_{\z}\q) = \dim_\q(H_n(G,\q))$.
\end{demonstracao}

\chapter{Spectral sequences}
\label{chapter:sequences}
A spectral sequence is a generalization of an exact sequence that provides means of calculating homology via succesive approximations. There are many such sequences that are commonly used in algebra, like the Serre spectral sequence associated to a fibration or the Lyndon-Hochschild-Serre spectral sequence relating the homologies of a normal subgroup H and the quotient group G/H to the homology of the total group G.

\section{Spectral sequences associated to double complexes}

In this section, we focus on sequences that arise from filtrations of double complexes, which will be necessary the next section on Equivariant Homology. The main reference used for this section is \cite{rotman2009}.

\begin{definicao}
    A spectral sequence is a collection $\{E^r_{p,q},d^r_{p,q}\}_{p,q\in\z,r\geq a}$ of $R$-modules and homomorphisms  $d^r_{p,q}: E^r_{p,q} \rightarrow E^r_{p-r,q+r-1}$ such that $d^r_{p-r,q+r-1}\circ d^r_{p,q} = 0$ and $E^{r+1}_{p,q} \simeq \displaystyle\frac{\ker(d^r_{p,q})}{\im (d^r_{p-r,q+r-1})}$. Here $a$ is some integer, not necessarily positive.
\end{definicao}

Intuitively, a spectral sequence can be thought of as a collection of "pages" where the $i$-th page consists of the modules $E^i_{p,q}$ with upper index $i$ and all the maps $d^i_{p,q}$ with upper index $i$ which map to and from modules of the $i$-th page. The $i$-th page is often represented as $E^i$. In the definition, the condition $r\geq a$ indicates that the sequence starts at the page $E^a$ of index $a$.

The condition $d^r_{p-r,q+r-1}\circ d^r_{p,q} = 0$ implies that the maps $d$ form chain complexes around the modules, and the condition $E^{r+1}_{p,q} \simeq \displaystyle\frac{\ker(d^r_{p,q})}{\im (d^r_{p-r,q+r-1})}$ tells us that to find $E^{r+1}_{p,q}$, one must go back one page to the module of same $p,q$ index and take the homology around that module.

As an example, letting $r=0$ gives that the maps in the $0$-th page map $E^0_{p,q}$ to $E^0_{p,q-1}$, that is, the maps point down. In the first page, $d^1_{p,q}: E^1_{p,q}\rightarrow E^1_{p-1,q}$, so the maps point left. Figure \ref{fig: spectral} below illustrates the first pages of a first quadrant spectral sequence.

\begin{definicao}
    A \textit{first quadrant spectral sequence} is a spectral sequence such that $E^r_{p,q}=0$ if either $p$ or $q$ is negative.
\end{definicao}

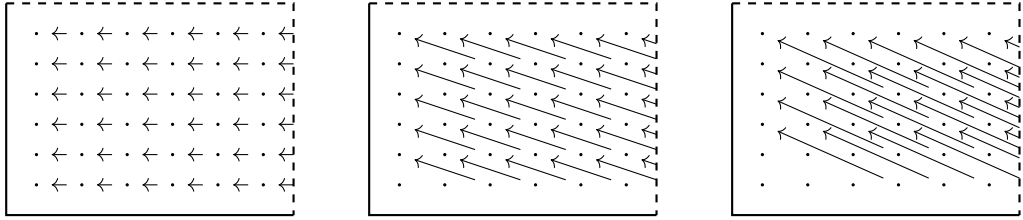
\begin{figure}[H]
    \centering
    \caption{Pages $E^1,E^2$ and $E^3$ of a first quadrant spectral sequence}
    \begin{tikzpicture}[scale=0.2]
        \draw [fill=black, shift={(0,0)}] (0,0) circle (2pt);
        \draw [fill=black, shift={(0,0)}] (3,0) circle (2pt);
        \draw [fill=black, shift={(0,0)}] (6,0) circle (2pt);
        \draw [fill=black, shift={(0,0)}] (9,0) circle (2pt);
        \draw [fill=black, shift={(0,0)}] (12,0) circle (2pt);
        \draw [fill=black, shift={(0,0)}] (15,0) circle (2pt);

        \draw[ ->, shift={(0,0)}] (2,0)--(1,0);
        \draw[ ->, shift={(0,0)}] (5,0)--(4,0);
        \draw[ ->, shift={(0,0)}] (8,0)--(7,0);
        \draw[ ->, shift={(0,0)}] (11,0)--(10,0);
        \draw[ ->, shift={(0,0)}] (14,0)--(13,0);
        \draw[ ->, shift={(0,0)}] (17,0)--(16,0);

        \draw [fill=black, shift={(0,2)}] (0,0) circle (2pt);
        \draw [fill=black, shift={(0,2)}] (3,0) circle (2pt);
        \draw [fill=black, shift={(0,2)}] (6,0) circle (2pt);
        \draw [fill=black, shift={(0,2)}] (9,0) circle (2pt);
        \draw [fill=black, shift={(0,2)}] (12,0) circle (2pt);
        \draw [fill=black, shift={(0,2)}] (15,0) circle (2pt);

        \draw[ ->, shift={(0,2)}] (2,0)--(1,0);
        \draw[ ->, shift={(0,2)}] (5,0)--(4,0);
        \draw[ ->, shift={(0,2)}] (8,0)--(7,0);
        \draw[ ->, shift={(0,2)}] (11,0)--(10,0);
        \draw[ ->, shift={(0,2)}] (14,0)--(13,0);
        \draw[ ->, shift={(0,2)}] (17,0)--(16,0);

        \draw [fill=black, shift={(0,4)}] (0,0) circle (2pt);
        \draw [fill=black, shift={(0,4)}] (3,0) circle (2pt);
        \draw [fill=black, shift={(0,4)}] (6,0) circle (2pt);
        \draw [fill=black, shift={(0,4)}] (9,0) circle (2pt);
        \draw [fill=black, shift={(0,4)}] (12,0) circle (2pt);
        \draw [fill=black, shift={(0,4)}] (15,0) circle (2pt);

        \draw[ ->, shift={(0,4)}] (2,0)--(1,0);
        \draw[ ->, shift={(0,4)}] (5,0)--(4,0);
        \draw[ ->, shift={(0,4)}] (8,0)--(7,0);
        \draw[ ->, shift={(0,4)}] (11,0)--(10,0);
        \draw[ ->, shift={(0,4)}] (14,0)--(13,0);
        \draw[ ->, shift={(0,4)}] (17,0)--(16,0);

        \draw [fill=black, shift={(0,6)}] (0,0) circle (2pt);
        \draw [fill=black, shift={(0,6)}] (3,0) circle (2pt);
        \draw [fill=black, shift={(0,6)}] (6,0) circle (2pt);
        \draw [fill=black, shift={(0,6)}] (9,0) circle (2pt);
        \draw [fill=black, shift={(0,6)}] (12,0) circle (2pt);
        \draw [fill=black, shift={(0,6)}] (15,0) circle (2pt);

        \draw[ ->, shift={(0,6)}] (2,0)--(1,0);
        \draw[ ->, shift={(0,6)}] (5,0)--(4,0);
        \draw[ ->, shift={(0,6)}] (8,0)--(7,0);
        \draw[ ->, shift={(0,6)}] (11,0)--(10,0);
        \draw[ ->, shift={(0,6)}] (14,0)--(13,0);
        \draw[ ->, shift={(0,6)}] (17,0)--(16,0);

        \draw [fill=black, shift={(0,8)}] (0,0) circle (2pt);
        \draw [fill=black, shift={(0,8)}] (3,0) circle (2pt);
        \draw [fill=black, shift={(0,8)}] (6,0) circle (2pt);
        \draw [fill=black, shift={(0,8)}] (9,0) circle (2pt);
        \draw [fill=black, shift={(0,8)}] (12,0) circle (2pt);
        \draw [fill=black, shift={(0,8)}] (15,0) circle (2pt);

        \draw[ ->, shift={(0,8)}] (2,0)--(1,0);
        \draw[ ->, shift={(0,8)}] (5,0)--(4,0);
        \draw[ ->, shift={(0,8)}] (8,0)--(7,0);
        \draw[ ->, shift={(0,8)}] (11,0)--(10,0);
        \draw[ ->, shift={(0,8)}] (14,0)--(13,0);
        \draw[ ->, shift={(0,8)}] (17,0)--(16,0);

        \draw [fill=black, shift={(0,10)}] (0,0) circle (2pt);
        \draw [fill=black, shift={(0,10)}] (3,0) circle (2pt);
        \draw [fill=black, shift={(0,10)}] (6,0) circle (2pt);
        \draw [fill=black, shift={(0,10)}] (9,0) circle (2pt);
        \draw [fill=black, shift={(0,10)}] (12,0) circle (2pt);
        \draw [fill=black, shift={(0,10)}] (15,0) circle (2pt);

        \draw[ ->, shift={(0,10)}] (2,0)--(1,0);
        \draw[ ->, shift={(0,10)}] (5,0)--(4,0);
        \draw[ ->, shift={(0,10)}] (8,0)--(7,0);
        \draw[ ->, shift={(0,10)}] (11,0)--(10,0);
        \draw[ ->, shift={(0,10)}] (14,0)--(13,0);
        \draw[ ->, shift={(0,10)}] (17,0)--(16,0);


        \draw [fill=black, shift={(24,0)}] (0,0) circle (2pt);
        \draw [fill=black, shift={(24,0)}] (3,0) circle (2pt);
        \draw [fill=black, shift={(24,0)}] (6,0) circle (2pt);
        \draw [fill=black, shift={(24,0)}] (9,0) circle (2pt);
        \draw [fill=black, shift={(24,0)}] (12,0) circle (2pt);
        \draw [fill=black, shift={(24,0)}] (15,0) circle (2pt);

        \draw[ ->, shift={(24,0)}] (5,1/3)--(1,5/3);
        \draw[ ->, shift={(24,0)}] (8,1/3)--(4,5/3);
        \draw[ ->, shift={(24,0)}] (11,1/3)--(7,5/3);
        \draw[ ->, shift={(24,0)}] (14,1/3)--(10,5/3);
        \draw[ ->, shift={(24,0)}] (17,1/3)--(13,5/3);
        \draw[ ->, shift={(24,0)}] (17,4/3)--(16,5/3);

        \draw [fill=black, shift={(24,2)}] (0,0) circle (2pt);
        \draw [fill=black, shift={(24,2)}] (3,0) circle (2pt);
        \draw [fill=black, shift={(24,2)}] (6,0) circle (2pt);
        \draw [fill=black, shift={(24,2)}] (9,0) circle (2pt);
        \draw [fill=black, shift={(24,2)}] (12,0) circle (2pt);
        \draw [fill=black, shift={(24,2)}] (15,0) circle (2pt);

        \draw[ ->, shift={(24,2)}] (5,1/3)--(1,5/3);
        \draw[ ->, shift={(24,2)}] (8,1/3)--(4,5/3);
        \draw[ ->, shift={(24,2)}] (11,1/3)--(7,5/3);
        \draw[ ->, shift={(24,2)}] (14,1/3)--(10,5/3);
        \draw[ ->, shift={(24,2)}] (17,1/3)--(13,5/3);
        \draw[ ->, shift={(24,2)}] (17,4/3)--(16,5/3);

        \draw [fill=black, shift={(24,4)}] (0,0) circle (2pt);
        \draw [fill=black, shift={(24,4)}] (3,0) circle (2pt);
        \draw [fill=black, shift={(24,4)}] (6,0) circle (2pt);
        \draw [fill=black, shift={(24,4)}] (9,0) circle (2pt);
        \draw [fill=black, shift={(24,4)}] (12,0) circle (2pt);
        \draw [fill=black, shift={(24,4)}] (15,0) circle (2pt);

        \draw[ ->, shift={(24,4)}] (5,1/3)--(1,5/3);
        \draw[ ->, shift={(24,4)}] (8,1/3)--(4,5/3);
        \draw[ ->, shift={(24,4)}] (11,1/3)--(7,5/3);
        \draw[ ->, shift={(24,4)}] (14,1/3)--(10,5/3);
        \draw[ ->, shift={(24,4)}] (17,1/3)--(13,5/3);
        \draw[ ->, shift={(24,4)}] (17,4/3)--(16,5/3);

        \draw [fill=black, shift={(24,6)}] (0,0) circle (2pt);
        \draw [fill=black, shift={(24,6)}] (3,0) circle (2pt);
        \draw [fill=black, shift={(24,6)}] (6,0) circle (2pt);
        \draw [fill=black, shift={(24,6)}] (9,0) circle (2pt);
        \draw [fill=black, shift={(24,6)}] (12,0) circle (2pt);
        \draw [fill=black, shift={(24,6)}] (15,0) circle (2pt);

        \draw[ ->, shift={(24,6)}] (5,1/3)--(1,5/3);
        \draw[ ->, shift={(24,6)}] (8,1/3)--(4,5/3);
        \draw[ ->, shift={(24,6)}] (11,1/3)--(7,5/3);
        \draw[ ->, shift={(24,6)}] (14,1/3)--(10,5/3);
        \draw[ ->, shift={(24,6)}] (17,1/3)--(13,5/3);
        \draw[ ->, shift={(24,6)}] (17,4/3)--(16,5/3);

        \draw [fill=black, shift={(24,8)}] (0,0) circle (2pt);
        \draw [fill=black, shift={(24,8)}] (3,0) circle (2pt);
        \draw [fill=black, shift={(24,8)}] (6,0) circle (2pt);
        \draw [fill=black, shift={(24,8)}] (9,0) circle (2pt);
        \draw [fill=black, shift={(24,8)}] (12,0) circle (2pt);
        \draw [fill=black, shift={(24,8)}] (15,0) circle (2pt);

        \draw[ ->, shift={(24,8)}] (5,1/3)--(1,5/3);
        \draw[ ->, shift={(24,8)}] (8,1/3)--(4,5/3);
        \draw[ ->, shift={(24,8)}] (11,1/3)--(7,5/3);
        \draw[ ->, shift={(24,8)}] (14,1/3)--(10,5/3);
        \draw[ ->, shift={(24,8)}] (17,1/3)--(13,5/3);
        \draw[ ->, shift={(24,8)}] (17,4/3)--(16,5/3);

        \draw [fill=black, shift={(24,10)}] (0,0) circle (2pt);
        \draw [fill=black, shift={(24,10)}] (3,0) circle (2pt);
        \draw [fill=black, shift={(24,10)}] (6,0) circle (2pt);
        \draw [fill=black, shift={(24,10)}] (9,0) circle (2pt);
        \draw [fill=black, shift={(24,10)}] (12,0) circle (2pt);
        \draw [fill=black, shift={(24,10)}] (15,0) circle (2pt);


        \draw [fill=black, shift={(48,0)}] (0,0) circle (2pt);
        \draw [fill=black, shift={(48,0)}] (3,0) circle (2pt);
        \draw [fill=black, shift={(48,0)}] (6,0) circle (2pt);
        \draw [fill=black, shift={(48,0)}] (9,0) circle (2pt);
        \draw [fill=black, shift={(48,0)}] (12,0) circle (2pt);
        \draw [fill=black, shift={(48,0)}] (15,0) circle (2pt);

        \draw[->,shift={(48,0)}] (8,4/9)--(1,32/9);
        \draw[->,shift={(48,0)}] (11,4/9)--(4,32/9);
        \draw[->,shift={(48,0)}] (14,4/9)--(7,32/9);
        \draw[->,shift={(48,0)}] (17,4/9)--(10,32/9);
        \draw[->,shift={(48,0)}] (17,16/9)--(13,32/9);
        \draw[->,shift={(48,0)}] (17,28/9)--(16,32/9);

        \draw [fill=black, shift={(48,2)}] (0,0) circle (2pt);
        \draw [fill=black, shift={(48,2)}] (3,0) circle (2pt);
        \draw [fill=black, shift={(48,2)}] (6,0) circle (2pt);
        \draw [fill=black, shift={(48,2)}] (9,0) circle (2pt);
        \draw [fill=black, shift={(48,2)}] (12,0) circle (2pt);
        \draw [fill=black, shift={(48,2)}] (15,0) circle (2pt);

        \draw[->,shift={(48,2)}] (8,4/9)--(1,32/9);
        \draw[->,shift={(48,2)}] (11,4/9)--(4,32/9);
        \draw[->,shift={(48,2)}] (14,4/9)--(7,32/9);
        \draw[->,shift={(48,2)}] (17,4/9)--(10,32/9);
        \draw[->,shift={(48,2)}] (17,16/9)--(13,32/9);
        \draw[->,shift={(48,2)}] (17,28/9)--(16,32/9);

        \draw [fill=black, shift={(48,4)}] (0,0) circle (2pt);
        \draw [fill=black, shift={(48,4)}] (3,0) circle (2pt);
        \draw [fill=black, shift={(48,4)}] (6,0) circle (2pt);
        \draw [fill=black, shift={(48,4)}] (9,0) circle (2pt);
        \draw [fill=black, shift={(48,4)}] (12,0) circle (2pt);
        \draw [fill=black, shift={(48,4)}] (15,0) circle (2pt);

        \draw[->,shift={(48,4)}] (8,4/9)--(1,32/9);
        \draw[->,shift={(48,4)}] (11,4/9)--(4,32/9);
        \draw[->,shift={(48,4)}] (14,4/9)--(7,32/9);
        \draw[->,shift={(48,4)}] (17,4/9)--(10,32/9);
        \draw[->,shift={(48,4)}] (17,16/9)--(13,32/9);
        \draw[->,shift={(48,4)}] (17,28/9)--(16,32/9);

        \draw [fill=black, shift={(48,6)}] (0,0) circle (2pt);
        \draw [fill=black, shift={(48,6)}] (3,0) circle (2pt);
        \draw [fill=black, shift={(48,6)}] (6,0) circle (2pt);
        \draw [fill=black, shift={(48,6)}] (9,0) circle (2pt);
        \draw [fill=black, shift={(48,6)}] (12,0) circle (2pt);
        \draw [fill=black, shift={(48,6)}] (15,0) circle (2pt);

        \draw[->,shift={(48,6)}] (8,4/9)--(1,32/9);
        \draw[->,shift={(48,6)}] (11,4/9)--(4,32/9);
        \draw[->,shift={(48,6)}] (14,4/9)--(7,32/9);
        \draw[->,shift={(48,6)}] (17,4/9)--(10,32/9);
        \draw[->,shift={(48,6)}] (17,16/9)--(13,32/9);
        \draw[->,shift={(48,6)}] (17,28/9)--(16,32/9);

        \draw [fill=black, shift={(48,8)}] (0,0) circle (2pt);
        \draw [fill=black, shift={(48,8)}] (3,0) circle (2pt);
        \draw [fill=black, shift={(48,8)}] (6,0) circle (2pt);
        \draw [fill=black, shift={(48,8)}] (9,0) circle (2pt);
        \draw [fill=black, shift={(48,8)}] (12,0) circle (2pt);
        \draw [fill=black, shift={(48,8)}] (15,0) circle (2pt);

        \draw [fill=black, shift={(48,10)}] (0,0) circle (2pt);
        \draw [fill=black, shift={(48,10)}] (3,0) circle (2pt);
        \draw [fill=black, shift={(48,10)}] (6,0) circle (2pt);
        \draw [fill=black, shift={(48,10)}] (9,0) circle (2pt);
        \draw [fill=black, shift={(48,10)}] (12,0) circle (2pt);
        \draw [fill=black, shift={(48,10)}] (15,0) circle (2pt);

        \draw [fill=white, shift={(48,12)}, opacity=0] (0,0) circle (2pt);

        \draw [thick] (-2,12)--(-2,-2)--(17,-2);
        \draw[thick, dashed] (-2,12)--(17,12)--(17,-2);

        \draw [thick, shift={(24,0)}] (-2,12)--(-2,-2)--(17,-2);
        \draw[thick, dashed, shift={(24,0)}] (-2,12)--(17,12)--(17,-2);

        \draw [thick, shift={(48,0)}] (-2,12)--(-2,-2)--(17,-2);
        \draw[thick, dashed, shift={(48,0)}] (-2,12)--(17,12)--(17,-2);
    \end{tikzpicture}
    \label{fig: spectral}
    \fautor
\end{figure}

Suppose $\{E^r_{p,q},d^r_{p,q}\}$ is a spectral sequence starting at page 0. Then the modules at the page $E^1$ are all homologies of the complexes defined by $d^0$, thus, quotients of submodules of $E^0$. 

Let $Z^1$ be the set of "cycles" of $E^0_{p,q}$, that is, elements in the modules of $E^0_{p,q}$ that are annihilated by $d^0_{p,q}$, and let $B^1$ be the corresponding set of "boundaries" of $E^0_{p,q}$, elements that are image of some element under $d^0_{p,q-1}$. Then, $E^1_{p,q}$ is the quotient $Z^1/B^1$.

Now consider $Z^2$ and $B^2$ similarly defined. By the correspondence theorem, $Z^2 = A'/B^1$ and $B^2 = B'/B^1$ where $A'$ and $B'$ contain $B^1$. For simplicity, we will abuse notation and denote $A'$ and $B'$ by their quotients $Z^2$ and $B^2$. Note that then, we have that $$B^1 \subseteq B^2 \subseteq Z^2\subseteq Z^1\subseteq E^0.$$

This generalizes nicely: identifying $Z^r$ and $B^r$ with submodules of $E^0$ whose successive quotients lead to $Z^r$ and $B^r$, we have that $$B^1\subseteq B^2 \subseteq \dots \subseteq B^r\subseteq Z^r \subseteq \dots \subseteq Z^2\subseteq Z^1\subseteq E^0.$$

\begin{definicao}
    Given a spectral sequence $\{E^r_{p,q},d^r_{p,q}\}$, denote by $Z^\infty_{p,q} = \bigcap Z^r_{p,q}$ and $B^\infty_{p,q} = \bigcup B^r_{p,q}$. From the remarks above, $B^\infty_{p,q} \subseteq Z^\infty_{p,q}$. We define the limit term $E^\infty_{p,q}$ by $$E^\infty_{p,q} \coloneq \frac{Z^\infty_{p,q}}{B^\infty_{p,q}}.$$  
\end{definicao}

\begin{proposicao}
    Let $\{E^r_{p,q},d^r_{p,q}\}$ be a spectral sequence. Denote by $Z^r_{p,q}$ and $B^r_{p,q}$ the set of cycles and boundaries of $E^{r-1}_{p,q}$, respectively, so that $E^r_{p,q} = Z^r_{p,q}/B^r_{p,q}$. Then \begin{enumerate}
        \item[i.] $E^{r+1}_{p,q} = E^r_{p,q}$ if and only if $Z^{r+1}_{p,q} = Z^r_{p,q}$ and $B^{r+1}_{p,q} = B^r_{p,q}$.
        
        \item[ii.] If $E^{r+1}_{p,q} = E^r_{p,q}$ for all $r \geq s$, then $E^s_{p,q} = E^\infty_{p,q}$.
    \end{enumerate} 
\end{proposicao}

\begin{demonstracao}
    We will omit indices $p,q$ in this proof.

    For the first item, one of the implications is immediate, since $$E^{r+1} = Z^{r+1}/B^{r+1} = Z^r/B^r = E^r.$$ For the other implication, notice that the inclusions $B^{r+1}\subseteq Z^{r+1}\subseteq E^r$ gives us that $Z^{r+1}/B^{r+1} = E^r$. This implies  $B^{r+1}=0$ in $E^r = Z^r/B^r$ (which implies $B^{r+1}=B^r$) and so $$E^{r+1} =\frac{Z^{r+1}}{B^{r+1}} = \frac{Z^{r+1}}{B^r} = E^r = \frac{Z^r}{B^r}.$$ This shows that $Z^{r+1}=Z^r$.

    For the second item, if $E^r = E^{r+1}$ for $r \geq s$, then the previous item gives $Z^s = Z^r$ for all $r \geq s$, so that $Z^s = \bigcap_{r\geq s}Z^r = Z^\infty$. Similarly, $B^s = B^r$ for all $r \geq s$, so $B^s = \bigcup_{r\geq s}B^r = B^\infty$. Thus $$E^s = \frac{Z^s}{B^s} = \frac{Z^\infty}{B^\infty} = E^\infty.$$
\end{demonstracao}

The proposition above is especially useful for first quadrant spectral sequences: Let $\{E^r_{p,q},d^r_{p,q}\}$ be such a sequence and fix an index $p,q$. If $s>p$, then $p-s<0$, implying that $d^s_{p,q}: E^s_{p,q} \rightarrow E^s_{p-s,q+s-1}$ points to a zero module, thus its kernel is the entire set $E^s_{p,q}$. Similarly, if $t>q+1$, then $q-t+1<0$ and the map $d^t_{p+t,q-t+1}: E^t_{p+t,q-t+1} \rightarrow E^t_{p,q}$ has trivial domain, thus trivial image. This implies that taking $r$ large enough, the maps pointing to and out of $E^r_{p,q}$ are trivial, thus further homologies will preserve the module. 

In other words, if $\{E^r_{p,q},d^r_{p,q}\}$ is a first quadrant spectral sequence, then $E^\infty_{p,q}$ always coincides with $E^r_{p,q}$ for some large enough $r$.

\begin{definicao}\label{def: bounded}
    A filtration of a module $M$ is a family $(F^pM)_{p\in\z}$ of submodules of $M$ such that $$\dots \subseteq F^{p-1}M\subseteq F^pM \subseteq F^{p+1}M\subseteq \dots.$$

    Let $M=(M_n)$ be a graded module. A filtration $FM = (FM_n)$ of $M$ is said to be bounded if for all $n$ there exists $s=s(n),t=t(n)$ such that $F^sM_n = 0$ and $F^tM_n=M_n$.
\end{definicao}

\begin{definicao}\label{def: convergence}
    A spectral sequence is said to \textit{converge} to a graded module $H = (H_n)$, which is denoted by $$E^r_{p,q} \Longrightarrow H_{p+q}$$ for some $r$, if there is some bounded filtration $FH$ of $H$ with $E^{\infty}_{p,n-p} \simeq \displaystyle\frac{F^pH_n}{F^{p-1}H_n}$ for all $n,p$.
\end{definicao}

The remainder of this section will be dedicated to the construction of two spectral sequences associated with a double complex.

\begin{definicao}
    A \textit{bigraded module} is a family $M=(M_{p,q})_{p,q\in\z}$ of $R$-modules doubly indexed by $\z$.
\end{definicao}

\begin{definicao}
    Let $M$ and $N$ be bigraded $R$-modules, and let $(a,b)\in\z\times\z$. A \textit{bigraded map} of bidegree $(a, b)$ is a family of homomorphisms $f = ( f_{p,q} : M_{p,q} \rightarrow N_{p+a,q+b})_{(p,q)\in\z\times\z}$. The pair $(a,b)$ is said to be the \textit{bidegree} of $f$.
\end{definicao}

\begin{definicao}
    A \textit{double complex} (or \textit{bicomplex}) is an ordered triple $(M, d', d'')$, where $M = (M_{p,q})$ is a bigraded module, $d', d'' : M \rightarrow M$ are bigraded maps of bidegrees $(-1, 0)$ and $(0, -1)$, such that $d'\circ d'=0$, $d''\circ d''=0$ and $d'_{p,q-1}d''_{p,q} + d''_{p-1,q}d'_{p,q} = 0$.
\end{definicao}

The conditions $d'\circ d'=0$ and $d''\circ d''=0$ imply that fixing one of the $p,q$ coordinates of $M$ gives a chain complex. Fixing $p$, we get a "vertical" chain complex with differential $d''$ and fixing $q$ gives a "horizontal" chain complex with differential $d'$. The third condition  $d'_{p,q-1}d''_{p,q} =- d''_{p-1,q}d'_{p,q}$ implies that squares formed by the maps are anti-commutative.

\begin{center}
    \begin{tikzcd}[column sep=10ex]
        M_{p-1,q} \arrow[dd, "d''_{p-1,q}", swap] & M_{p,q} \arrow[l, "d'_{p,q}", swap] \arrow[dd, "d''_{p,q}"] \\ & \\
        M_{p-1,q-1} & M_{p,q-1} \arrow[l, "d'_{p,q-1}"]
    \end{tikzcd}
\end{center}

Given a bigraded module $M$ with differentials $d',d''$ with the correct bidegrees whose squares comute (that is, $d'_{p,q-1}d''_{p,q} = d''_{p-1,q}d'_{p,q}$), one can create a bicomplex via a sign change: Let $\partial'_{p,q} = d'_{p,q}$ and let $\partial''_{p,q} = (-1)^pd''_{p,q}$. Then $\partial',\partial''$ will still be differentials with same bidegrees as $d',d''$ and \begin{align*}
    \partial'_{p,q-1}\partial''_{p,q}+\partial''_{p-1,q}\partial'_{p,q} &= (-1)^pd'_{p,q-1}d''_{p,q}+(-1)^{p-1}d''_{p-1,q}d'_{p,q} \\ &= (-1)^p(d'_{p,q-1}d''_{p,q}-d''_{p-1,q}d'_{p,q}) \\ &=0,
\end{align*} thus making $(M,\partial',\partial'')$ a double complex.

\begin{exemplo}
    Let $(A,d')$ and $(B,d'')$ be two chain complexes of $R$-modules. Let $M = (M_{p,q}) = (A_p\otimes_R B_q)$ and consider the maps $d'\otimes_R \id$ and $\id\otimes_R d''$. Then the squares in the lattice formed by $M$ will commute, so that adding a factor $(-1)^p$ to $\id\otimes_Rd''$ makes $M$ into a double complex.
\end{exemplo}

\begin{definicao}
    If $(M,d',d'')$ is a double complex, then its \textit{total complex}, denoted by $\Tot(M)$ is a chain complex whose $n$-term is given by $$\Tot(M)_n = \bigoplus_{p+q=n}M_{p,q}$$ and whose differentials $D_n: \Tot(M)_n\rightarrow \Tot(M)_{n-1}$ are given by $$D_n = \sum_{p+q=n}(d'_{p,q}+d''_{p,q}).$$
\end{definicao}

The notation above for $D_n$ is to indicate that its action on a fixed $M_{p,q}$ summand of $M_{n}$, $n=p+q$ is given by $d'_{p,q}+d''_{p,q}$. Note that $d'_{p,q}$ maps into $M_{p-1,q}$ and $d''_{p,q}$ maps into  $M_{p,q-1}$, both with indices that sum to $n-1 = p+q-1$, thus both modules that lie in $\Tot(M)_{n-1}$. 

To see that $D_n$ is indeed a differential, let $x$ be any element in some $M_{p,q}$. Then $D_n$ will map it to $d'_{p,q}(x)+d''_{p,q}(x)$, where the first summand lives in $M_{p-1,q}$ and the second lives in $M_{p,q-1}$. Then \begin{align*}
    D_{n-1}D_n(x) &= D_{n-1}(d'_{p,q}(x)+d''_{p,q}(x)) \\
    &= (d'_{p-1,q}+d''_{p-1,q})(d'_{p,q}(x))+(d'_{p,q-1}+d''_{p,q-1})(d''_{p,q}(x)) \\
    &= d'_{p-1,q}d'_{p,q}(x)+(d''_{p-1,q}d'_{p,q}+d'_{p,q-1}d''_{p,q})(x)+(d''_{p,q-1}d''_{p,q})(x) \\
    &=0.
\end{align*} This follows from the fact that $d',d''$ are horizontal and vertical differentials implies that $d'_{p-1,q}d'_{p,q}=0$, $d''_{p,q-1}d''_{p,q}=0$ and the middle map $d''_{p-1,q}d'_{p,q}+d'_{p,q-1}d''_{p,q}$ is zero because $(M,d',d'')$ is a bicomplex.  

This can be thought of in a more geometric way: the sum $d'+d''$ maps left and downwards. Applying it to some $x$ leads to $d'(x)$ which points left and to $d''(x)$ which points down. Applying it again to $d'(x)$ leads to a point which is twice left from $x$ (thus zero, since horizontal lines form a complex) and to a point which is left, then down. Similarly, applying it again to $d''(x)$ leds to a point which is twice down from $x$, thus zero, and one which is down, then left. Anti-commutativity of the squares implies these two remaining points sum to zero.

\begin{definicao}
    Given a double complex $(M,d',d'')$, the homology groups $H_n(\Tot(M)_\bullet)$ of its total complex will be referred to as the homology groups of the double complex.
\end{definicao}

\begin{definicao}
    Let $(M,d',d'')$ be a double complex. The \textit{transpose complex} of $(M,d',d'')$ is the double complex $(M^t,d'',d')$ where $M^t_{p,q} = M_{q,p}$.
\end{definicao}

\begin{observacao}\label{rem: total complex}
    Note that $M$ and $M^t$ have the same total complexes, since $\Tot(M)_n = \displaystyle\bigoplus_{n=p+q}M_{p,q} = \bigoplus_{n=q+p}M_{q,p} = \Tot(M^t)_n$ and $d'+d'' = d''+d'$, so their differentials also agree. This implies $M$ and $M^t$ have the same homology groups.
\end{observacao}

One way in which spectral sequences arise is from a filtration of given module. Here, we intend to give two filtrations to a double complex.

\begin{definicao}
    Let $(M,d',d'')$ be a double complex. The first filtration of $\Tot(M)$ is given by $$\leftindex^IF^p\Tot(M)_n = \bigoplus_{i\leq p}M_{i,n-i}.$$
\end{definicao}

\begin{definicao}
    Let $(M,d',d'')$ be a double complex. The second filtration of $\Tot(M)$ is given by $$\leftindex^{II}F^p\Tot(M)_n = \bigoplus_{j\leq p}M_{n-j,j}.$$
\end{definicao}

These filtrations can be thought of as follows: For the first one, start at $M_{0,n}$ and progressively add $M_{1,n-1},M_{2,n-2},\dots$ going downwards in the diagonal $n=p+q$. For the second one, start at $M_{n,0}$ and progressively add $M_{n-1,1},M_{n-2,2},\dots$ going upwards in the diagonal $n=p+q$.

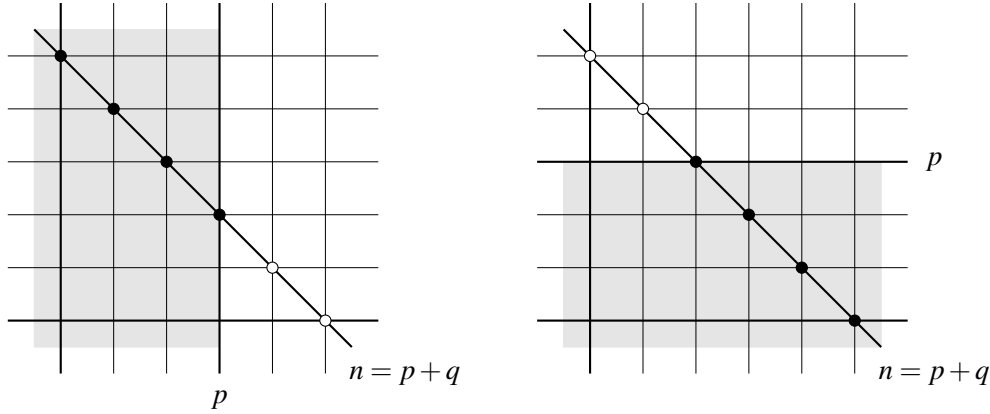
\begin{figure}[H]
    \centering
    \caption{First and second filtrations for a total complex}
    \begin{tikzpicture}[scale=0.7]
        \draw[gray, ultra thin, fill, opacity=0.2, shift={(0,0)}] (-0.5,-0.5)--(3,-0.5)--(3,5.5)--(-0.5,5.5) -- cycle;
        \draw[thick, shift={(0,0)}] (-1,0)--(6,0);
        \draw[ultra thin, shift={(0,0)}] (-1,1)--(6,1);
        \draw[ultra thin, shift={(0,0)}] (-1,2)--(6,2);
        \draw[ultra thin, shift={(0,0)}] (-1,3)--(6,3);
        \draw[ultra thin, shift={(0,0)}] (-1,4)--(6,4);
        \draw[ultra thin, shift={(0,0)}] (-1,5)--(6,5);
        \draw[thick, shift={(0,0)}] (0,-1)--(0,6);
        \draw[ultra thin, shift={(0,0)}] (1,-1)--(1,6);
        \draw[ultra thin, shift={(0,0)}] (2,-1)--(2,6);
        \draw[thick, shift={(0,0)}] (3,-1)--(3,6);
        \draw[ultra thin, shift={(0,0)}] (4,-1)--(4,6);
        \draw[ultra thin, shift={(0,0)}] (5,-1)--(5,6);
        \draw[thick, shift={(0,0)}] (-0.5,5.5)--(5.5,-0.5);

        \draw [fill=black, shift={(0,0)}] (0,5) circle (3pt);
        \draw [fill=black, shift={(0,0)}] (1,4) circle (3pt);
        \draw [fill=black, shift={(0,0)}] (2,3) circle (3pt);
        \draw [fill=black, shift={(0,0)}] (3,2) circle (3pt);
        \draw [fill=white, shift={(0,0)}] (4,1) circle (3pt);
        \draw [fill=white, shift={(0,0)}] (5,0) circle (3pt);

        \node[black,scale=0.9] at (3,-1.5) {$p$};
        \node[black,scale=0.9] at (6.5,-1) {$n=p+q$};


        \draw[gray, ultra thin, fill, opacity=0.2, shift={(10,0)}] (-0.5,-0.5)--(5.5,-0.5)--(5.5,3)--(-0.5,3) -- cycle;
        \draw[thick, shift={(10,0)}] (-1,0)--(6,0);
        \draw[ultra thin, shift={(10,0)}] (-1,1)--(6,1);
        \draw[ultra thin, shift={(10,0)}] (-1,2)--(6,2);
        \draw[thick, shift={(10,0)}] (-1,3)--(6,3);
        \draw[ultra thin, shift={(10,0)}] (-1,4)--(6,4);
        \draw[ultra thin, shift={(10,0)}] (-1,5)--(6,5);
        \draw[thick, shift={(10,0)}] (0,-1)--(0,6);
        \draw[ultra thin, shift={(10,0)}] (1,-1)--(1,6);
        \draw[ultra thin, shift={(10,0)}] (2,-1)--(2,6);
        \draw[ultra thin, shift={(10,0)}] (3,-1)--(3,6);
        \draw[ultra thin, shift={(10,0)}] (4,-1)--(4,6);
        \draw[ultra thin, shift={(10,0)}] (5,-1)--(5,6);
        \draw[thick, shift={(10,0)}] (-0.5,5.5)--(5.5,-0.5);

        \draw [fill=white, shift={(10,0)}] (0,5) circle (3pt);
        \draw [fill=white, shift={(10,0)}] (1,4) circle (3pt);
        \draw [fill=black, shift={(10,0)}] (2,3) circle (3pt);
        \draw [fill=black, shift={(10,0)}] (3,2) circle (3pt);
        \draw [fill=black, shift={(10,0)}] (4,1) circle (3pt);
        \draw [fill=black, shift={(10,0)}] (5,0) circle (3pt);

        \node[black,scale=0.9] at (16.5,3) {$p$};
        \node[black,scale=0.9] at (16.5,-1) {$n=p+q$};
    \end{tikzpicture}
    \label{fig:filtrations}
    \fautor
\end{figure}

\begin{definicao}
    An \textit{exact couple} is a 5-tuple $(D, E, \alpha, \beta, \gamma)$, where $D$ and $E$ are bigraded modules, $\alpha,\beta,\gamma$ are bigraded maps making the diagram below commute, and there is exactness at each vertex: $\ker\alpha=\im\gamma, \ker\beta=\im\alpha$ and $\ker\gamma=\im\beta$.

    \begin{center}
        \begin{tikzcd}
            D \arrow[rr, "\alpha"] & & D \arrow[dl, "\beta"] \\
            & E \arrow[ul, "\gamma"] &
        \end{tikzcd}
    \end{center}
\end{definicao}

\begin{proposicao}
    Every filtration $(F^pC)_{p\in\z}$ of a complex $C$ determines an exact couple $(D, E, \alpha, \beta, \gamma)$ where the bidegrees of $\alpha, \beta, \gamma$ are $(1,-1),(0,0)$ and $(-1,0)$, respectively.
\end{proposicao}

\begin{demonstracao}
    Abbreviate $F^pC$ as $F^p$ and recall that $F^{p-1}\subset F^p$, from the definition of a filtration. So we have a short exact sequence $$0\rightarrow F^{p-1}\rightarrow F^p \rightarrow F^p/F^{p-1}\rightarrow 0$$ where the maps are an inclusion and a projection, leading to a long exact sequence of homology (see Theorem 6.10 of \cite{rotman2009}) $$\cdots \rightarrow H_n(F^{p-1})\rightarrow H_n(F^p) \rightarrow H_n(F^p/F^{p-1})\rightarrow H_{n-1}(F^{p-1}) \rightarrow \cdots$$ Write $n=p+q$, so that the sequence becomes $$\cdots \rightarrow H_{p+1}(F^{p-1})\rightarrow H_{p+1}(F^p) \rightarrow H_{p+1}(F^p/F^{p-1})\rightarrow H_{p+q-1}(F^{p-1}) \rightarrow \cdots$$ Let $D_{p,q} = H_{p+q}(F^p)$ and $E_{p,q} = H_{p+q}(F^p/F^{p-1})$. Then the sequence is written as $$\cdots \rightarrow D_{p-1,q-1} \rightarrow D_{p,q}\rightarrow E_{p,q} \rightarrow D_{p-1,q} \rightarrow \cdots$$ Note that the maps involved have bidegrees $(1,-1),(0,0)$ and $(-1,0)$ and they are, respectively, the maps induced by the inclusion $F^{p-1}\hookrightarrow F^p$, induced by the projection $F^p\rightarrow F^p/F^{p-1}$, and the associated connecting map.
\end{demonstracao}

\begin{definicao}
    If $M$ is a bigraded module and $d: M \rightarrow M$ is a bigraded map of bidegree $(a,b)$ such that $dd=0$ (that is, $d$ is a differential), we call the pair $(M,d)$ a \textit{differential bigraded module}, and its homology $H(M,d)$ is the bigraded module whose $p,q$ term is $$H(M,d)_{p,q} \coloneq \frac{\ker(d_{p,q})}{\im(d_{p-a,q-b})}.$$
\end{definicao}

\begin{proposicao}\label{bidegree}
    If $(D, E, \alpha, \beta, \gamma)$ is an exact couple, then $d^1 \coloneq \beta\gamma$ is a differential $d^1 : E \rightarrow E$, and there is an exact couple $(D^2, E^2, \alpha^2, \beta^2, \gamma^2)$, called the \textit{derived couple}, with $E^2 = H(E, d^1)$. If $\alpha,\beta,\gamma$ have bidegrees $(a,a'),(b,b'),(c,c')$, respectively, then $\alpha^2,\beta^2,\gamma^2$ have bidegrees $(a,a'),(b-a,b'-a'),(c,c')$, respectively.
\end{proposicao}

\begin{demonstracao}
    Let $d^1 := \beta\gamma$. This composition makes sense, since $\gamma: E \rightarrow D$ and $\beta: D\rightarrow E$, making $d^1$ a map from $E$ to $E$. Furthermore, exactness of the original couple gives $\gamma\beta=0$, so that $d^1d^1 = \beta(\gamma\beta)\gamma = 0$, making it a differential.

    Let $E^2 = H(E,d^1)$ and let $D^2 = \im\alpha \subseteq D$.

    Define $\alpha^2:D^2\rightarrow D^2$ as the restriction of $\alpha$ to $D^2$. If $x\in D^2 = \im\alpha$, then $x = \alpha y$ for some $y\in D$ and so $$\alpha^2x = \alpha\rvert_{D^2}x = \alpha(x) = \alpha(\alpha y).$$

    To define $\beta^2: D^2\rightarrow E^2$, note again that $x\in D^2$ is $x=\alpha y$ for some $y\in D$. Note that $d^1\beta y = \beta\gamma\beta y=0$, so $\beta y$ is a cycle for the differential $d^1$, thus, it makes sense to talk about $[\beta y] \in E^2 = H(E,d^1)$. Furthermore, if $x=\alpha y=\alpha y'$, then $y-y'\in\ker\alpha=\im\gamma$ so that $y-y' = \gamma z$ and $\beta(y)-\beta(y') = \beta\gamma z = d^1z$, and so, $[\beta(y)] = [\beta(y')]$. Abusing notation, $$\beta^2(x) = \beta^2(\alpha y) = [\beta \alpha^{-1}\alpha y]  = [\beta(y)]$$ is well defined. 

    Finally, we define $\gamma^2: E^2 \rightarrow D^2$. If $[z]\in E^2$, then $z\in E$ and $0=d^1z = \beta\gamma z$, so that $\gamma z\in\ker\beta = \im\alpha = D^2$. If $[z]=[z']$, then $z-z' =d^1w$ for some $w$ and $\gamma(z)-\gamma(z') = \gamma d^1w = \gamma\beta\gamma w = 0$, thus $\gamma(z)=\gamma(z')$, so that $$\gamma^2([z]) = \gamma z$$ is well defined.

    Now, we show that $(D^2,E^2,\alpha^2,\beta^2,\gamma^2)$ is an exact couple.

    The inclusions of the type $\im\subset\ker$ follow from the fact that $\beta\alpha,\gamma\beta,\alpha\gamma$ are all zero, since $$\beta^2\alpha^2(x) = \beta^2(\alpha \alpha y) = [\beta\alpha y ] = 0,$$ $$\gamma^2\beta^2(y) = \gamma^2[\beta y] = \gamma\beta y = 0,$$ $$\alpha^2\gamma^2([z]) = \alpha^2\gamma z = \alpha\gamma z = 0,$$ since $\gamma z\in \ker\beta=\im\alpha = D^2$.

    Let $x\in\ker\alpha^2$. Then $x\in D^2 = \im\alpha=\ker\beta$ and in particular, $\alpha^2x = \alpha x$ since $\alpha^2$ is just the restriction of $\alpha$ to $D^2$. Thus, $x\in\ker\alpha^2 $ which implies that $ x\in\ker\alpha = \im\gamma$, and $x = \gamma z$ for some $z$. Since $x\in\ker\beta$, we have that $0=\beta x = \beta\gamma z = d^1z$, so $z$ is a cycle and it makes sense to consider $[z]$. Then $x = \gamma z = \gamma^2[z]$ and $x\in\im\gamma^2$. This proves that $$\ker\alpha^2\subseteq\im\gamma^2.$$

    Let $x\in\ker\beta^2 \subset D^2 = \im\alpha$.  Then $x=\alpha y$ for some y and $0 = \beta^2x = [\beta y]$ implying $\beta y$ is a cycle, that is, $\beta y = d^1z = \beta\gamma z$ for some $z$. Then $y-\gamma z \in \ker\beta = \im\alpha = D^2$, so that $\alpha^2(y-\gamma z) = \alpha(y-\gamma z) = \alpha y -\alpha\gamma z = \alpha y = x$, thus $x\in \im \alpha^2$. This proves that $$\ker\beta^2\subseteq \im\alpha^2.$$ 

    Let $[z]\in\ker\gamma^2$. Then $\gamma[z] = \gamma z=0$, so $z\in\ker\gamma=\im\beta$, that is, $z = \beta w$ for some $w$. Hence $\beta^2(\alpha w) = [\beta\alpha^{-1}\alpha w] = [\beta w] = [z]$ and so, $z\in\im\beta^2$. This proves $$\ker\gamma^2\subseteq\im\beta^2.$$

    As a final remark, if $\alpha,\beta,\gamma$ have bidegrees $(a,a'),(b,b'),(c,c')$, since $\alpha^2$ is a restriction of $\alpha$, it also has bidegree $(a,a')$. Similarly, $\gamma^2$ takes a class $[z]$ in $E^2_{p,q}$ where $z\in E^1_{p,q}$ and maps it to $\gamma z \in D^1_{p+c,q+c'}$ such that $\gamma z$ is understood to lie in $D^2_{p+c,q+c'}$, thus, $\gamma^2$ also has bidegree $(c,c')$. Finally, $\beta^2$ will take $x\in D^2_{p,q} \subset D^1_{p,q}$, find $y\in D^1_{p-a,q-a'}$ such that $x = \alpha y$, and this $y$ is then mapped to $\beta y \in E^1_{p-a+b,q-a'+b'}$, so that $\beta^2x=[\beta y] \in E^2_{p-a+b,q-a'+b'}$. Thus, $\beta^2$ has bidegree $(b-a,b'-a')$.

    In short, if $\alpha,\beta,\gamma$ have bidegrees $(a,a'),(b,b'),(c,c')$ then $\alpha^2,\beta^2,\gamma^2$ have bidegrees $(a,a'),(b-a,b'-a'),(c,c')$.
\end{demonstracao}

Since $(D^2, E^2, \alpha^2, \beta^2, \gamma^2)$ is also an exact couple, one could also consider its derived couple. More generally, it is possible to iterate this process to obtain the $r$-th derived couple. 

\begin{definicao}
    Let $(D, E, \alpha, \beta, \gamma)$ be an exact couple. Its $(r + 1)$-st derived couple $(D^{r+1}, E^{r+1},$ $ \alpha^{r+1}, \beta^{r+1}, \gamma^{r+1})$ is the derived couple of $(D^{r}, E^{r}, \alpha^{r}, \beta^{r}, \gamma^{r})$.
\end{definicao}

\begin{corolario}\label{derived}
    Let $(D, E, \alpha, \beta, \gamma)$ be  be the exact couple arising from a filtration $F^p$ of a complex $C$. Then the $r$-th derived couple has the following properties:

    \begin{enumerate}
        \item[i. ] The bigraded maps $\alpha^r,\beta^r,\gamma^r$ have bidegrees $(1,-1)$, $(1-r,r-1)$ and $(-1,0)$, respectively;

        \item[ii. ] The differential $d^r$ has bidegree $(-r,r-1)$ and it is induced by $\beta\alpha^{-r+1}\gamma$;

        \item[iii. ] $E^{r+1}_{p,q} = \ker (d^r_{p,q})/\im (d^r_{p+r,q-r+1})$;

        \item[iv. ] $D^r_{p,q} = \im(\alpha_{p-1,q+1})(\alpha_{p-2,q+2})\dots(\alpha_{p-r+1,q+r-1})$. In particular, for the exact couple arising from a filtration, $D^r_{p,q}$ is the image of the induced homomorphism $$(j^{p-1}j^{p-2}\dots j^{p-r+1})_*: H_n(F^{p-r+1}) \rightarrow H_n(F^p).$$
    \end{enumerate}
\end{corolario}

\begin{demonstracao}
    Proposition \ref{bidegree} shows that if $\alpha,\beta,\gamma$ have bidegrees $(a,a'),(b,b'),(c,c')$ then $\alpha^2,\beta^2,\gamma^2$ have bidegrees $(a,a'),(b-a,b'-a'),(c,c')$. Induction shows that $\alpha^r,\beta^r,\gamma^r$ have bidegrees $(a,a'),(b-(r-1)a,b'-(r-1)a'),(c,c')$. Since the exact couple arising from a filtration has bidegrees $(1,-1),(0,0)$ and $(-1,0)$, the maps $\alpha^r,\beta^r,\gamma^r$ will have bidegrees $(1,-1),(1-r,r-1),(-1,0)$, which shows (i).

    For (ii), note that $d^1 = \beta\gamma$, so $d^2 = \beta^2\gamma^2$. Then $$d^2[z] = \beta^2\gamma^2[z] = \beta^2\gamma z = [\beta\alpha^{-1}\gamma z]$$ and by induction, $$\beta^r\gamma^r[z] = [\beta\alpha^{-r+1}\gamma z].$$ To see the claim for the bidegree, recall that $\beta^r,\gamma^r$ have bidegrees $(1-r,r-1),(-1,0)$ so the composition has bidegree $(-r,r-1)$.

    For the third statement, $$E^{r+1}_{p,q} = H(E^r_{p,q},d^{r}_{p,q}) = \frac{\ker (d^r_{p,q})}{\im(d^r_{p+r,q-r+1})}.$$

    Finally, note that each $\alpha^r$ is a restriction of $\alpha^{r-1}$ to $D^{r-1}$ and since the bidegree of $\alpha$ is $(1,-1)$, $$D^2_{p,q} = \alpha_{p-1,q+1}D_{p-1,q+1}$$ and so $$D^3_{p,q} = \alpha_{p-1,q+1} D_{p-1,q+1} = \alpha_{p-1,q+1}\alpha_{p-2,q+2}D_{p-2,q+2}.$$ So by induction, $$D^r_{p,q} = \alpha_{p-1,q+1}\alpha_{p-2,q+2}\dots\alpha_{p-r+1,q+r-1}D_{p-r+1,q-r+1}.$$ Finally, recall that this $\alpha_{p,q}$ is induced by $j^p: H_{p+q}(F^p)\rightarrow H_{p+q}(F^{p+1})$, that is, $\alpha_{p,q} = j^p_*: D_{p,q}\rightarrow D_{p+1,q-1}$. Then \begin{align*}
        D^r_{p,q} &= \alpha_{p-1,q+1}\alpha_{p-2,q+2}\dots\alpha_{p-r+1,q+r-1}D_{p-r+1,q-r+1} \\
        &= j_*^{p-1}j_*^{p-2}\dots j_*^{p-r+1}D_{p-r+1} \\
        &= (j^{p-1}j^{p-2}\dots j^{p-r+1})_*D_{p-r+1,q-r+1}.
    \end{align*}
\end{demonstracao}

\begin{teorema}\label{sequence}
    Every filtration of a complex yields a spectral sequence.
\end{teorema}

\begin{demonstracao}
    Every filtration yields an exact couple $(D,E,\alpha,\beta,\gamma)$. Let $E^r, d^r$ be the bigraded module $E^r_{p,q}$ and the differential $d^r$ from the $r$-th derived couple. Then $d^r$ has bidegree $(-r,r+1)$ and $E^{r+1}_{p,q} = \ker(d^r_{p,q})/\im(d^r_{p+r,q-r+1})$, as in Corollary \ref{derived}.
\end{demonstracao}

\begin{definicao}
    If $(F^pC)$ is a filtration of a complex $C$ and $i^p : F^pC \rightarrow C$ are inclusions, define the induced filtration of $H_n(C)$ to be $\Phi^pH_n(C) = \im(i^p_*)$.
\end{definicao}

Note that $i^p_*: H_n(F^pC) \rightarrow H_n(C)$, so that $\im(i^p_*)$ is always a submodule of $H_n(C)$.

Let $j^p: F^pC \rightarrow F^{p+1}C$ be inclusions, so that $i^{p+1}\circ j^p = i^p$, implying that $i^{p+1}_*\circ j^p_* = i^p_*$, so $$i^{p+1}_*\circ j^p_*(H_n(F^pC)) = i^p_*(H_n(F^pC)).$$ Then $j^p_*(H_n(F^pC)) \subseteq H_n(F^{p+1}C)$ and so $$i_*^p(H_n(F^pC)) = i_*^{p+1}(j^p_*(H_n(F^pC))) \subseteq i_*^{p+1}(H_n(F^{p+1}C)).$$ In other words, $\Phi^p(H_n(C)) \subseteq \Phi^{p+1}(H_n(C))$.

Thus, $\Phi^pH_n(C)$ is indeed a filtration of $H_n(C)$.

\begin{proposicao}\label{prp: bounds}
    If $(F^pC)$ is a bounded filtration (as in Definition \ref{def: bounded}), the induced filtration $\Phi^pH_n(C)$ is also bounded with same bounds.
\end{proposicao}

\begin{demonstracao}
    If $F^pC$ is bounded, that is, $F^sC=0$ and $F^tC=C$ for some $s,t$, then $H_n(F^sC)=0$ and so $i_*^s: H_n(F^sC)\rightarrow H_n(C)$ is a map with null image, that is, $\Phi^sH_n(C) = 0$. Also, $H_n(F^tC) = H_n(C)$, so $i^t = \id$ and $\Phi^tH_n(C) = H_n(C)$. Thus the bounds for $\Phi^pH_n(C)$ are at most the bounds for $F^pC$.
\end{demonstracao}

\begin{teorema}\label{convergence}
    Let $(F^pC)$ be a bounded filtration of a complex $C$, and let $\{E^r,d^r\}$ be the associated spectral sequence (as in Theorem \ref{sequence}). Then for each $p,q$ we have $E^\infty_{p,q} = E^r_{p,q}$ for a large enough $r$ and $E^2_{p,q} \Longrightarrow H_{p+q}(C)$.
\end{teorema}

\begin{demonstracao}
    Let $n=p+q$ and $s(n),t(n)$ be the bounds for $F^pC$, which will also be bounds for $\Phi^nH_n(C)$.

    Note that if $p> t(n)$, then $F^{p-1}C=F^pC = C$, so $E_{p,q} = H_{p+q}(F^p/F^{p-1}) = H_{p+1}(0) = 0$. Further quotients will also be zero, so $E^r_{p,q} = 0$ for all $r$. Similarly, if $p<s(t)$, then $F^p=0$ and $E^r_{p,q}=0$ for all $r$.

    Now, if $r> p-s(t)$, we have that $p-r<s(t)$, so $E^r_{p-r,q+r-1} = 0$, and the differential $d^r_{p,q}: E^r_{p,q}\rightarrow E^r_{p-r,q+r-1}$ has kernel equal to $E^r_{p,q}$. Similarly, if $r> t(n)-p$, then $r+p> t(n)$ and $E^r_{p+r,q-r+1}=0$, so the differential $d^r_{p+r,q-r+1}: E^r_{p+r,q-r+1} \rightarrow E^r_{p,q}$ has trivial image. Thus, taking $r$ large enough implies that $$E^{r+1}_{p,q} \simeq \frac{\ker d^r_{p,q}}{\im d^r_{p+r,q-r+1}} \simeq E^r_{p,q}$$ so that $E^r_{p,q} = E^\infty_{p,q}$.   

    To show that $E^2_{p,q}\Longrightarrow H_n(C)$, recall that the couple $(D^r,E^r,\alpha^r,\beta^r,\gamma^r)$ forms an exact sequence $$D^r_{p+r-2,q-r+2} \overset{\alpha^r}{\longrightarrow} D^r_{p+r-1,q-r+1}\overset{\beta^r}{\longrightarrow} E^r_{p,q} \overset{\gamma^r}{\longrightarrow} D^r_{p-1,q}$$ and, as in Corollary \ref{derived}, we have that $$D^r_{p,q} = \im(j^{p-1}j^{p-2}\dots j^{p-r+1})_*: H_n(F^{p-r+1})\rightarrow H_n(F^p).$$ Hence $$D^r_{p+r-1,q-r+1} = \im(j^{p+r-2}\dots j^p)_*\subset H_n(F^{p+r-1}),$$ $$D^r_{p+r-2,q-r+2} = \im(j^{p+r-3}\dots j^{p-1})_*\subset H_n(F^{p+r-2}),$$ and for large enough $r$, we have that $F^{p+r-1}C=C$. Thus the composition $j^{p+r-2}\dots j^p$ will coincide with the inclusion $i^p: F^pC\rightarrow C$, so that $$D^r_{p+r-1,q-r+1} = \im(i^p_*) = \Phi^p(H_n(C))$$ and similarly, $D^r_{p+r-2,q+r-2} = \Phi^{p-1}(H_n(C))$.

    Furthermore, since $D^r_{p-1,q}$ is the image of $H_n(F^{p-r}C)$ under a homomorphism, taking $r$ large enough leads to $F^{p-r}C=0$, so this image will be trivial. Thus, the exact sequence from the derived couple becomes $$\Phi^{p-1}(H_n(C)) \rightarrow \Phi^p(H_n(C)) \rightarrow E^r_{p,q} \rightarrow 0$$ implying that $$E^r_{p,q}\simeq \frac{\Phi^p(H_n(C))}{\Phi^{p-1}(H_n(C))}.$$ Hence the sequence converges to $H_n(C)$.
\end{demonstracao}

\begin{teorema}
    Let $M$ be a first quadrant double complex (that is, $M_{p,q}=0$ if either $p$ or $q$ is negative) and let $\leftindex^IE^r, \leftindex^{II}E^r$ be the spectral sequences determined by the first and second filtrations of $\Tot(M)$. Then \begin{enumerate}
        \item[i. ] The first and second filtrations are bounded;

        \item[ii. ] For all $p,q$ we have $\leftindex^{I}E^\infty_{p,q} = \leftindex^{I}E^r_{p,q}$ and $\leftindex^{II}E^\infty_{p,q} = \leftindex^{II}E^r_{p,q}$ for large enough $r$;

        \item[iii. ] $ \leftindex^{I}E^2_{p,q}\Longrightarrow H_{p+q}(\Tot(M))$ and $ \leftindex^{II}E^2_{p,q}\Longrightarrow H_{p+q}(\Tot(M))$. 
    \end{enumerate}
\end{teorema}

\begin{demonstracao}
    Recall that the first filtration is given by starting at $M_{0,n}$ and progressively adding terms $M_{r,n-r}$, so that lower and upper bounds for this filtration is given by $s(n)=-1$ and $t(n)=n$, respectively. The argument for the second filtration is the same.

    Notice that the second filtration of a double complex coincides with the first filtration of the transposed complex, implying their total complexes are the same, since transposing doesn't affect the total complex (see Remark \ref{rem: total complex}).

    Since these filtrations are bounded, the second and third statements follow from the previous theorem.
\end{demonstracao}

Consider a double complex $(M,d',d'')$. As mentioned earlier, both $d'$ and $d''$ form complexes once one of the indexes $p,q$ is fixed. 

Denote by $M_{p,*}$ the chain complex where the first index is some fixed integer $p$. Taking homology of this "vertical" complex gives a bigraded module whose $p,q$ term
is $H_q(M_{p,*})$. Now, the horizontal maps $d': M_{p,*} \rightarrow M_{p-1,*}$ induce maps $d'_*: H_q(M_{p,*}) \rightarrow H_q(M_{p-1,*})$. Since $d'd'=0$, we have that $d'_*d'_*=0$, thus the maps $d'_*$ are also differentials forming "horizontal" complexes, allowing us to take homology once again. We denote this by $H_p'H_q''(M)$, emphasizing that homology was taken first with respect to $d''$, then to $d'_*$. We call this the first iterated homology of the double complex. The second iterated homology is defined similarly to be $H_p''H_q'(M)$.

\begin{teorema}\label{teo: first seq}
    If $M$ is a first quadrant double complex, then $$\leftindex^{I}E^1_{p,q} = H_q(M_{p,*})$$ and $$\leftindex^{I}E^2_{p,q} = H_p'H_q''(M) \Longrightarrow H_n(\Tot(M)).$$
\end{teorema}

\begin{demonstracao}
    The fact that $\leftindex^{I}E^2_{p,q}  \Longrightarrow H_{p+q}(\Tot(M))$ has already been proven and follows from the fact that the associated filtration is bounded.

    We will omit the superscprit I and let $n=p+q$ throughout the proof.

    Recall that for a sequence arising from a filtration, $E^1_{p,q} = H_n(F^p/F^{p-1})$. Also, $$(F^p)_n = \bigoplus_{i\leq p}M_{i,n-i} = \dots \oplus M_{p-2,q+2} \oplus M_{p-1,q+1}\oplus M_{p,q},$$ $$(F^{p-1})_n = \bigoplus_{i\leq p-1}M_{i,n-i} = \dots \oplus M_{p-2,q+2} \oplus M_{p-1,q+1},$$ so that the $n$-th term of $F^p/F^{p-1}$ is simply given by $M_{p,q}$. Then, if $a_n\in M_{p,q}$, the differential $\overline{D}$ of $F^p/F^{p-1}$ induced by the total differential $D$ can be described as $$\overline{D}_n(a_n+(F^{p-1})n) = D_na_n+(F^{p_1})_{n-1} = d'a_n+d''a_n+(F^{p_1})_{n-1}$$ but $d''_{p,q}$ maps into $M_{p-1,q}\subset (F^{p-1})_n$, so $$\overline{D}_n(a_n+(F^{p-1})n) = d''a_n+(F^{p_1})_{n-1}.$$ Thus we get $$E^1_{p,q} = H_n\left(\frac{F^p}{F^{p-1}}\right) = \frac{\ker \overline{D}_n}{\im\overline{D}_{n+1}} \simeq \frac{\ker d''_{p,q}}{\im d''_{p,q+1}} = H_q(M_{p,*})$$ which shows the first part of the statement.

    We can then see elements of $E^1_{p,q}$ as classes $[z]$ where $z\in M_{p,q}$ and $d''z=0$. To prove that $E^2_{p,q}=H'_pH''_q(M)$, we will show that $d^1_{p,q}$ maps $[z]$ to $[d'z]$. 

    Consider the short exact sequence $$0 \rightarrow \frac{F^{p-1}}{F^{p-2}} \rightarrow \frac{F^p}{F^{p-2}} \rightarrow \frac{F^{p}}{F^{p-1}}\rightarrow 0$$ which induces a long exact sequence (see Theorem 6.10 of \cite{rotman2009}) $$\cdots \rightarrow H_n\left(\frac{F^p}{F^{p-2}}\right) \rightarrow H_n\left(\frac{F^{p}}{F^{p-1}}\right) \overset{\partial}{\longrightarrow} H_{n-1}\left(\frac{F^{p-1}}{F^{p-2}}\right) \rightarrow \cdots $$ where this connecting map $\partial$ is defined by lifting a cycle in $F^p/F^{p-1}$ to one in $F^p/F^{p-2}$, applying the boundary map, then finding an element in $F^{p-1}/F^{p-2}$ whose inclusion coincides with this boundary. This is precisely what the composition $\beta\gamma$ does to a class (just moving through a different lifting), since $\gamma$ is the connecting homomorphism $H_n(F^p/F^{p-1})\rightarrow H_{n-1}(F^{p-1})$ and $\beta$ is induced by the projection $\pi: F^{p-1}\rightarrow F^{p-1}/F^{p-2}$, so the connecting homomorphism $\partial$ coincides with the differential $d^1$. In other words, $d^1$ is the map that arises from the diagram \begin{center}
        \begin{tikzcd}
            & & (F^p/F^{p-2})_n \arrow[r, "\pi"] \arrow[d, "\overline{D}"] & (F^p/F^{p-1})_n \\
            0 \arrow[r] & (F^{p-1}/F^{p-2})_{n-1} \arrow[r, "i"] & (F^p/F^{p-2})_{n-1}. &
        \end{tikzcd}
    \end{center} Let $z\in M_{p,q} = (F^p/F^{p-1})_n$ be a cycle for $\overline{D}$, that is, $d''z=0$. Choose $\pi^{-1}(z) = (0,z)$ so that $\overline{D}(0,z) = (d'_{p,q}z,0)$. Then $$d^1[z] = [i^{-1}\overline{D}\pi^{-1}z] = [d'z].$$
\end{demonstracao}

\begin{teorema}\label{teo: sec seq}
    If $M$ is a first quadrant double complex, then $$\leftindex^{II}E^1_{p,q} = H_q(M_{*,p})$$ and $$\leftindex^{II}E^2_{p,q} = H_p''H_q'(M) \Longrightarrow H_{p+q}(\Tot(M)).$$
\end{teorema}

\begin{demonstracao}
    The proof is analogous to the previous theorem.
\end{demonstracao}

To conclude this section, we summarize our main results: If $(M,d',d'')$ is a double complex, there exists two spectral sequences $\leftindex^{I}E^r_{p,q},\leftindex^{II}E^r_{p,q}$ both converging to the homology $H_{p+q}(\Tot(M))$, where the differentials $d^0,d^1$ are induced by $d'',d'$, respectively, for the first sequence and by $d',d''$, respectively, for the second. For $\leftindex^{I}E^r_{p,q}$, the first page consists of the vertical homologies $H_q(M_{p,*})$ of the original double complex and the second page consists of the iterated homology $H'_qH''_p(M)$, the horizontal homologies of the complexes arising from $H_q(M_{p,*})$, with analogous results for $\leftindex^{II}E^r_{p,q}$.

\section{Equivariant homology}

In this section, given a group $G$ acting on some CW-complex $X$, we wish to construct a spectral sequence that converges to the homology of $G$ and whose first page is given by simpler groups, namely, the homology groups of the isotropy subgroups of this action. This sequence will be the one that arises from a filtration of a double complex, as constructed in the last section. The main reference for this section is chapter VII of \cite{brown1994}.

We recall from the first section (see Definition \ref{def: action}) that an action of a group $G$ on a set $X$ is a map $\varphi: (G,X) \rightarrow X$ such that $\varphi(e,x)=x$ for all $x\in X$, where $e$ is the trivial element of $G$, and $\varphi(g,\varphi(h,x)) = \varphi(gh,x)$. Such actions are often represented by multiplication, so these properties are represented as $ex=x, g(hx)=(gh)x$. Note that this implies $gg^{-1}x = g^{-1}gx = ex=x$, so fixing a $g$ gives that $\phi(g,{-}): X \rightarrow X$ is a bijection.

Here, we are particularly interested in actions on CW-complexes, so we will demand that for each $g \in G$, the action of $g$ in $X$ gives a homeomorphism of X such that the image $g\sigma$ of any cell $\sigma$ of X is again a cell in $X$. Essentially, we are asking that the action permutes cells in $X$.

\begin{definicao}
    Let $(C, \partial), (D,\partial')$ be two chain complexes over a ring $R$. We define the \textit{tensor complex} $C\otimes_R D$ as the total complex of the bigraded module whose $p,q$ entry is $C_p\otimes_R D_q$ and whose differentials are given by $\partial\otimes 1$ and $1\otimes \partial'$.
\end{definicao}

\begin{definicao}
    Given a group $G$ and a nonnegative $G$-chain complex $C = (C_n)_{n\geq0}$ (that is, a complex of $G$-modules), we define the homology of $G$ with coefficients in $C$, denoted by $H_*(G,C)$, by $$H_*(G,C) = H_*(F\otimes_G C);$$ where $F$ is a projective resolution of $\z$ over $\z G$.
\end{definicao}

\begin{proposicao}\label{tensor-point}
    If $C$ consists of a single $G$-module $M$ concentrated in dimension 0, then $H_*(G, C)$ reduces to $H_*(G, M)$.    
\end{proposicao}

\begin{demonstracao}
    If $C_0=M$ and $C_q = 0$ otherwise, then $$(F\otimes_G C)_n = \bigoplus_{p+q=n} F_p\otimes_G C_q = F_n\otimes_G M.$$ So the total complex $F\otimes_G C$ is equivalent to the complex $F\otimes_G M$, whose homology is by definition $H_*(G,M)$.
\end{demonstracao}

\begin{teorema}\label{teo: equiv hom seq}
    Let $G$ a group, $C = (C_n)$ a non-negative $G$-chain complex. Then there exists two spectral sequences $\leftindex^I{E}, \leftindex^{II}{E}$ such that $$\leftindex^I{E}^2_{pq} = H_p(G, H_qC) \Longrightarrow H_{p+q}(G,C)$$ $$\leftindex^{II}{E}^1_{pq} = H_q(G, C_p) \Longrightarrow H_{p+q}(G,C)$$ 
\end{teorema}

\begin{demonstracao}
    These are simply the sequences previously described in Theorems \ref{teo: first seq} and \ref{teo: sec seq} that arises from the two filtrations of the total complex $F\otimes_G C$.
\end{demonstracao}

\begin{proposicao}
    If $f: C' \rightarrow C$ is a weak equivalence between complexes of right $R$-modules and $P$ is a non-negative complex of flat left $R$-modules, then $f \otimes_R \id_P: C' \otimes_R P \rightarrow C \otimes_R P$ is a weak equivalence.
\end{proposicao}

\begin{demonstracao}
    To prove this, we first recall that the suspension $(\Sigma C',\Sigma d')$ of a chain $(C',d')$ is given by $(\Sigma C)_n = C'_{n-1}$, $\Sigma d' = -d'$. The mapping cone of $f$ is then defined as the chain $C''=C\oplus\Sigma C'$ with $d''(c,c') = (dc+fc',-d'c')$. We have, then, a short exact sequence $$0 \rightarrow C \rightarrow C'' \rightarrow \Sigma C' \rightarrow 0$$ where the mappings are the inclusion of $C$ in $C''$ and the projection on the second coordinate. In the induced long exact sequence in homology, the connecting map is given by $f_*: H(C') \rightarrow H(C)$ (see Theorem 6.10 of \cite{rotman2009} for the long exact sequence induced by a short exact sequence of modules and see secion I.0 of \cite{brown1994} for this connecting map). Thus, the long exact sequence is $$\cdots \rightarrow H_{n+1}(C'') \rightarrow H_n(C') \overset{f_*}{\longrightarrow} H_n(C) \rightarrow H_n(C'') \rightarrow \cdots$$ which gives us that $f$ is a weak equivalence if and only the mapping cone $C''$ is acyclic.

    Thus, it suffices to show that the mapping cone of $f\otimes_R \id_P$ is acyclic. If $C''$ is the mapping cone of $f$, then $C''\otimes_R P$ will be the mapping cone of $f\otimes_R \id_P$.

    Let $P^{(n)}$ be the truncation $(P_i)_{i\leq n}$ of $P$. Since $P^{(n)}/P^{(n-1)}$ is a complex consisting of a single flat module $P_n$, we have that tensoring an exact complex with it preserves exactness. Assume by induction that $C'' \otimes_R P^{(n-1)}$ is acyclic (for $n=1$, it follows from flatness of $P_1$). From the exact sequence $$0 \rightarrow C'' \otimes_R P^{(n-1)} \rightarrow C'' \otimes_R P^{(n)} \rightarrow C'' \otimes_R (P^{(n)}/P^{n-1}) \rightarrow 0$$ and its associated long exact sequence we get that $C''\otimes_R P^{(n)}$ is acyclic. Since $C''\otimes_RP$ is the increasing union of the acyclic complexes $C''\otimes_R P^{(n)}$, it is acyclic.
\end{demonstracao}

\begin{proposicao}\label{weak-equiv}
    If $\tau: C \rightarrow C'$ is a weak equivalence of $G$-chain complexes, then $\tau$ induces an isomorphism $H_*(G, C) \simeq H_*(G, C')$. 
\end{proposicao}

\begin{demonstracao}
    By the previous proposition, if $\tau$ is a weak equivalence then $\id_F\otimes \tau: F\otimes_G C \rightarrow F\otimes_G C'$ is also a weak equivalence, that is, $$(\id_F\otimes\tau)_*: H_*(G,C) \rightarrow H_*(G,C')$$ is an isomorphism.
\end{demonstracao}

\begin{definicao}
    If $C(X)$ is the cellular chain complex of a $G$-complex $X$, the resulting homology groups $H_*(G, C(X))$ are denoted $H^G_*(X)$ and called the \textit{equivariant homology groups} of $(G, X)$.
\end{definicao}

The previous definition can be extended to consider coefficients in a $G$-module $M$.

\begin{definicao}
    Let $M$ be a $G$-module and consider the diagonal $G$-action on $C(X, M) := C(X) \otimes_{\z} M$. The homology groups $$H^G_*(X, M) = H_*(G, C(X, M))$$ are called the equivariant homology groups of the pair $(G,X)$ with coefficients in $M$.
\end{definicao}

\begin{proposicao}
    $H^G_*(\text{pt}., M) = H_*(G, M)$. 
\end{proposicao}

\begin{demonstracao}
    This follows from the fact that the cellular chain complex of a point consists of a single module $C_0$ in dimension 0 (see Definition \ref{def: cell complex}). So $C(\text{pt}.,M)$ consists of the single module $M$ in dimension 0. Applying Proposition \ref{tensor-point} gives the result.
\end{demonstracao}

\begin{teorema}
    If $f: X \rightarrow Y$ is a cellular map of CW-complexes on which $G$ acts (so that $C(X),C(Y)$ are $G$-complexes) such that $f_*: H_*X \rightarrow H_* Y$ is an isomorphism, then $f$ induces an isomorphism $H^G_*(X, M)\simeq H^G_*(Y, M)$ for any $G$-module $M$. 
\end{teorema}

\begin{demonstracao}
    In these conditions, $f$ is a weak equivalence between $C(X)$ and $C(Y)$. Applying Proposition \ref{weak-equiv} gives the result. The induced isomorphism is simply the map induced by $\id_F\otimes f: F\otimes_G C(X) \rightarrow F\otimes_G C(Y)$.
\end{demonstracao}

\begin{corolario}\label{equiv-contr}
    For any CW-complex $X$ on which $G$ acts and any $G$-module $M$, there is a canonical map $H^G_*(X, M) \rightarrow H_*(G, M)$. In particular, if $X$ is contractible then this map is an isomorphism $H^G_*(X, M) \simeq H_*(G, M)$.
\end{corolario}

\begin{demonstracao}
    Let $Y = \{y\}$ be a CW-complex consisting of a single point and let $G$ act trivially by $gy=y$. Then there exists a map $f:X\rightarrow Y$ given by $x\mapsto y$ for all $x\in X$, which will be cellular since each cell in $X$ gets mapped to the 0-cell $y$. This induces a map $f_*: H_*^G(X,M) \rightarrow H_*^G(Y,M) \simeq H_*(G,M)$. If $X$ is contractible, then $H_*X,H_*Y$ are both trivial and $f_*$ will be an isomorphism by the last theorem.
\end{demonstracao}

Given a contractible space $X$ on which $G$ acts, we then have that there is a spectral sequence converging to $H_*^G(X,M)\simeq H_*(G,M)$. This gives us, then, a way to approximate the homology of $G$ using its action on $X$. We now turn to analysing the terms $E^1_{pq} = H_q(G,C_p)$ of the spectral sequence of Theorem \ref{teo: equiv hom seq}.

\begin{lema}
    Let $N$ be a $G$-module whose underlying abelian group is of the form $\displaystyle\bigoplus_{i\in I}M_i$· Assume that the $G$-action permutes the summands according to some action of $G$ on $I$. Let $G_i$ be the isotropy group of $i$ and let $E$ be a set of representatives for $I$ mod $G$. Then $M_i$ is a $G_i$-module and there is a $G$-isomorphism $N\simeq \displaystyle\bigoplus_{i\in E}\text{Ind }^G_{G_i}M_i$.
\end{lema}

\begin{demonstracao}
    To say that the $G$-action permutes the summands is to say that $gM_i = M_{gi}$. Thus $M_i$ is not necessarily a $G$-module, but since $gM_i = M_{gi}=M_i$ for any $g\in G_i$, we get that $M_i$ is a $G_i$ module.
    
    Notice that, for any $H$-module $M$, we have (see Definition \ref{def: induced} and Remark \ref{zg free over zh}) $$\Ind^G_H M \simeq \z G\otimes_{\z H}M \simeq \left(\displaystyle\bigoplus_{g\in G/H}g\z H\right)\otimes_{\z H}M \simeq \displaystyle\bigoplus_{g\in G/H}g(\z H\otimes_{\z H}M)$$ which gives us that $\Ind^G_H M \simeq \displaystyle\bigoplus_{g\in G/H}gM$. Since we can choose to represent the coset $1H$ with $1$, this implies that $M$ is an $H$-submodule of $\Ind^G_HM$.  

    Now, since $I$ is the disjoint union of the orbits $Gi$, we get that $N = \displaystyle\bigoplus_{i\in I}M_i \simeq \displaystyle\bigoplus_{i\in E}\bigoplus_{j\in Gi} M_j$. On the other hand, the map $G/G_i \rightarrow G_i$ that takes a coset $gG_i$ and maps it to $gi$ is well-defined, since $gG_i = g'G_i$ implies that $g^{-1}g\in G_i$ so $gi=g'i$. This map is also a surjection, since any element in $Gi$ is of the form $gi$, and it is an injection, since $gi=g'i$ would imply $g^{-1}g\in G_i$.

    Then, for each $j$ in the orbit $Gi$, $M_j = M_{gi} = gM_i$ for some $g_j$ (unique for each $j$) so $\displaystyle\bigoplus_{j\in Gi} M_j$ can be written as $\displaystyle\bigoplus_{g\in G/G_i} gM_i$, which is isomorphic to $\Ind^G_{G_i}M_i$. Then $$N = \displaystyle\bigoplus_{i\in I}M_i \simeq \displaystyle\bigoplus_{i\in E}\left(\bigoplus_{j\in Gi} M_j\right) \simeq \displaystyle\bigoplus_{i\in E}\Ind^G_{G_i}M_i.$$
\end{demonstracao}

\begin{proposicao} Let $G_\sigma$ be the stabilizers of a cell $\sigma$, $\Sigma_p$ be a set of representatives for the $G$-orbits of $p$-cells, and let $\z_\sigma$ be the infinite cyclic group whose two generators correspond to the two orientations of $\sigma$ (so $g\in G_\sigma$ acts on $\z_\sigma$ as $+1$ if $g$ preserves the orientation of $\sigma$ and $-1$ otherwise). Then $C_p(X) \simeq \displaystyle\bigoplus_{\sigma\in \Sigma_p}\text{Ind }^G_{G_\sigma}\z_\sigma$.
\end{proposicao}

\begin{demonstracao}
    This follows directly from the last lemma, since $C_p(X)$ has an abelian group structure of a sum of one copy of $\z$ for each $p$-cell. Given a $p$-cell $\sigma$, denoting by $\z_\sigma$ the copy of $\z$ associated with $\sigma$ allows us to make explicit what the action of $g\in G_\sigma$ is, since $g\sigma$ must be either $\sigma$ or $-\sigma$. 
\end{demonstracao}

\begin{teorema}\label{spec-seq}
    There is a spectral sequence $E^1_{pq} = \displaystyle\bigoplus_{\sigma\in\Sigma_p}H_q(G_\sigma,M_\sigma) \Longrightarrow H^G_{p+q}(X,M)$, where $M_\sigma = \z_\sigma\otimes M$.
\end{teorema}

\begin{demonstracao}
    This spectral sequence is the one associated with the second filtration of $F\otimes_G C(X,M)$, as previously discussed.

    We have a G-module decomposition $$C_p(X,M) = C_p(X)\otimes M \simeq \left(\displaystyle\bigoplus_{\sigma\in \Sigma_p}\text{Ind }^G_{G_\sigma}\z_\sigma\right)\otimes M \simeq \displaystyle\bigoplus_{\sigma\in\Sigma_p}\text{Ind}^G_{G_\sigma}M_\sigma.$$ Shapiro's lemma (see Lemma \ref{shapiro}) now yields $$H_q(G,C_p(X,M))\simeq\displaystyle\bigoplus_{\sigma\in\Sigma_p}H_q(G_\sigma,M_\sigma)$$ so that the sequence takes the form $$E^1_{pq} = H_q(G,C_q(X,M)) = \displaystyle\bigoplus_{\sigma\in\Sigma_p}H_q(G_\sigma,M_\sigma) \Longrightarrow H^G_{p+q}(X,M).$$
\end{demonstracao}

\begin{corolario}\label{seq}
    If $X$ is a contractible $CW$-complex on which $G$ acts, there is a spectral sequence $E^1_{p,q} = \displaystyle\bigoplus_{\sigma\in\Sigma_p}H_q(G_\sigma,M_\sigma) \Longrightarrow H_{p+q}(G,M)$.
\end{corolario}

\begin{demonstracao}
    By Corollary \ref{equiv-contr}, $H^G_{*}(X,M)\simeq H_*(G,M)$, so the spectral sequence of Proposition \ref{spec-seq} converges to $H_*(G,M)$.
\end{demonstracao}

Since this spectral sequence arises from the double complex $F\otimes_G C(X,M)$ where $F$ is a projective resolution of $\z$ over $\z G$ and $C(X,M) = C(X)\otimes_{\z} M$, we know that the differential $d^1_{p,q}$ is induced by the boundary map of the CW-complex $X$. We wish to describe this differential in terms of the decomposition given for $E^1_{p,q}$. 

From here on, denote by $\partial: C_p(X,M) \rightarrow C_{p-1}(X,M)$ the boundary homomorphism of $X$ and for each pair $\sigma,\tau$ of a $p$-cell and a $(p-1)$-cell, let $\partial_{\sigma\tau}$ be the $(\sigma,\tau)$-component of $\partial$.

\begin{proposicao}
    Let $\mathcal{F}_\sigma = \{\tau: \partial_{\sigma\tau}\neq 0\}$. This is a finite set of $(p -1)$-cells and it is $G_\sigma$-invariant.
\end{proposicao}

\begin{demonstracao}
    This set is finite because each $p$-cell can only have a finite amount of $(p-1)$-cells in its boundary. Also, note that if $g\in G_\sigma$, then $g\sigma=\pm\sigma$ and from the description of $\partial$ in Definition \ref{def: cell complex}, we obtain $$\pm\sum_{\tau \in \mathcal{F_\sigma}}[\sigma:\tau]\tau=\partial(\pm\sigma)=\partial(g\sigma) = \sum_{g\tau \in\mathcal{F}_{g\sigma}}[g\sigma:g\tau]g\tau = \sum_{g\tau\in\mathcal{F_\sigma}}\pm[\sigma:g\tau]g\tau.$$ So for each $\tau$ in the boundary of $\sigma$, $g\tau$ must be some other component in that boundary. Thus $\mathcal{F}_\sigma$ is $G_\sigma$-invariant.
\end{demonstracao}

\begin{proposicao}
    Let $G_{\sigma\tau} = G_{\sigma}\cap G_\tau$. Then $[G_{\sigma}:G_{\sigma\tau}]<\infty$ and there is a transfer map $t_{\sigma\tau}: H_*(G_\sigma,M_\sigma) \rightarrow H_*(G_{\sigma\tau},M_\sigma).$
\end{proposicao}

\begin{demonstracao}
    The map describred here is the same transfer map mentioned in Proposition \ref{prp: transfer}. In addition to the construction described here, Section III.9 of \cite{brown1994} gives four other ways to define this map.
    
    If $g^{-1}h\in G_{\sigma\tau}$, denote by $\tau_1$ the cell $h\tau$. Then $(g^{-1}h)\tau = \tau$ implies $g^{-1}\tau_1=\tau$, that is, $\tau_1=g\tau$. We have shown that elements in the same coset of $G_\sigma/G_{\sigma\tau}$ must have the same action on $\tau$. The converse statement is also true, obviously: if $g\tau = h\tau$, then $\tau=g^{-1}h\tau$. Thus there are only as many cosets in $G_\sigma/G_{\sigma\tau}$ as there are orbits of $\tau$ under the action of $G_\sigma$, which is less than the cardinality of $F_\sigma$. We conclude that $(G_\sigma : G_{\sigma\tau})<\infty$ for $\tau\in\mathcal{F}_\sigma$.

    In general, given a group $G$, a subgroup $H$ of $G$ with finite index and a $G$-module $N$, we can define a map on the invariants $tr:N_G\rightarrow N_H$ by $\overline{n}\mapsto \sum_{g\in G/H}\overline{\overline{gn}}$. This is well-defined since given $g,g'$ in the same coset of $G/H$, the fact that the action of $H$ in $N_H$ is trivial implies that $\overline{\overline{g'n}}=\overline{\overline{(g'g^{-1})gn}} = \overline{\overline{gn}}$.

    In particular, we have a map $t_{\sigma\tau}: F\otimes_{G_\sigma}M_\sigma \rightarrow F\otimes_{G_{\sigma\tau}}M_\sigma$ given by $$x\otimes m \mapsto \displaystyle\sum_{g\in G_\sigma/G_{\sigma\tau}}g^{-1}x\otimes g^{-1}m$$ and the induced map on homology is the desired one. Here we are choosing $g^{-1}$ as a representative for the coset to make further calculations easier.
\end{demonstracao}

\begin{proposicao}
    The map $\partial_{\sigma\tau}$ and the inclusion $G_{\sigma\tau}\hookrightarrow G_\tau$ induce a map $u_{\sigma\tau}:H_*(G_{\sigma\tau},M_\sigma) \rightarrow H_*(G_\tau,M_\tau)$.
\end{proposicao}

\begin{demonstracao}
    This map is simply the induced map on homology given by the composition of $1\otimes\partial_{\sigma\tau}: F\otimes_{G_{\sigma\tau}}M_\sigma \rightarrow F\otimes_{G_{\sigma\tau}}M_\tau$ with the coset enlargment map $F\otimes_{G_{\sigma\tau}}M_\tau \rightarrow F\otimes_{G_{\tau}}M_\tau$. At a chain level, the action of this map is given by $$x\otimes m \mapsto x\otimes\partial_{\sigma\tau}m.$$
\end{demonstracao}

\begin{proposicao}
    Let $\tau_0\in\Sigma_{p-1}$ be the representative of the $G$-orbit of $\tau$, and choose $g(\tau) \in G$ such that $g(\tau)\tau = \tau_0$. Then there is an induced isomorphism $\nu_\tau: H_*(G_\tau, M_\tau)\rightarrow H_*(G_{\tau_0},M_{\tau_0}).$     
\end{proposicao}

\begin{demonstracao}
    In general, given two groups $G,G'$, a $G$-module $N$, a $G'$-module $N'$, a group isomorphism $\alpha: G\rightarrow G'$ and a map $f: N\rightarrow N'$ compatible with $\alpha$ in the sense that $\alpha(g)f(n) = f(gn)$, we can induce a map $N_G \rightarrow N'_{G'}$ by $\overline{n} \mapsto \overline{\overline{f(n)}}$. This is well defined in the quotient $G'$ since $g'\overline{\overline{f(n)}} = \overline{\overline{g'f(n)}} = \overline{\overline{\alpha(g)f(n)}}=\overline{\overline{f(gn)}} = \overline{\overline{f(n)}}$.

    For our case in particular, consider the map from $M_\tau = \z_\tau\otimes_{\z} M$ to $M_{\tau_0} = \z_{\tau_0}\otimes_{\z} M$ given by multiplication by $g(\tau)$ and consider the group isomorphism $G_\tau\rightarrow G_{\tau_0}$ given by $g\mapsto g(\tau)gg(\tau)^{-1}$. These maps are compatible, so we get a map $F\otimes_{G_\tau}M_\tau \rightarrow F\otimes_{G_{\tau_0}}M_{\tau_0}$ given by $$x\otimes m \mapsto g(\tau)x\otimes g(\tau)m$$ and the desired map $\nu_\tau$ is the map induced on homology.
\end{demonstracao}

\begin{proposicao}\label{prp: varphi}
    Let $$\varphi: \bigoplus_{\sigma\in\Sigma_p}H_*(G_\sigma, M_\sigma) \rightarrow \bigoplus_{\tau\in\Sigma_{p-1}}H_*(G_\tau, M_\tau)$$ be the map given by $$\varphi\rvert_{H_*(G_\sigma, M_\sigma)} = \sum_{\tau\in\mathcal{F}'_{\sigma}} \nu_\tau u_{\sigma\tau}t_{\sigma\tau}$$ where $\mathcal{F}'_\sigma$ is a set of representatives for $\mathcal{F}_\sigma/G_\sigma$. Then, up to sign, the map $\varphi$ is the differential $d^1: E^1_{p,*}\rightarrow E^1_{p-1,*}$.
\end{proposicao}

\begin{demonstracao}
    To prove this, we will show that the following diagram commutes:

\begin{center}
    \begin{tikzcd}
        \displaystyle\bigoplus_{\sigma\in\Sigma_p} F\otimes_{G_\sigma} M_\sigma \arrow[r, "\varphi"] \arrow[d, "\theta"] & \displaystyle\bigoplus_{\tau\in\Sigma_{p-1}} F\otimes_{G_\tau} M_\tau \arrow[d, "\psi"] \\
        F\otimes_G C_p(X,M) \arrow[r, "1\otimes\partial"] & F\otimes_G C_{p-1}(X,M)
    \end{tikzcd}
\end{center} so that the corresponding diagram on the homology groups commutes. Here $\theta$ and $\psi$ are the chain isomorphisms given in Shapiro's lemma. One can check from the proof of that lemma that these maps are simply inclusions. We have \begin{align*}
    \sum_{\tau\in\mathcal{F}'_{\sigma}} \nu_\tau u_{\sigma\tau}t_{\sigma\tau}(x\otimes m) &= \sum_{\tau\in\mathcal{F}'_{\sigma}} \nu_\tau u_{\sigma\tau}\left(\sum_{g\in G_\sigma/G_{\sigma\tau}}g^{-1}x\otimes g^{-1}m\right) \\
    &= \sum_{\tau\in\mathcal{F}'_{\sigma}} \nu_\tau \left(\sum_{g\in G_\sigma/G_{\sigma\tau}}g^{-1}x\otimes \partial_{\sigma\tau}g^{-1}m\right) \\
    &= \sum_{\tau\in\mathcal{F}'_{\sigma}} \nu_\tau \left(\sum_{g\in G_\sigma/G_{\sigma\tau}}g^{-1}x\otimes g^{-1}\partial_{\sigma,g\tau}m\right) \\
    &= \sum_{\tau\in\mathcal{F}'_{\sigma}} \sum_{g\in G_\sigma/G_{\sigma\tau}}g(\tau)g^{-1}x\otimes g(\tau)g^{-1}\partial_{\sigma,g\tau}m
\end{align*}. 

    Once we include the sum in $F\otimes_G C_{p-1}(X,M)$, we will be in a set where the $G$-action is trivial, so multiplying by $g(\tau)g^{-1}$ will have no effect. Also, since the first sum ranges over the orbits of $F_\sigma$ (which are in correspondence with the cosets of $G_\sigma/G_{\sigma\tau}$, as seen earlier) and the second sum ranges over the cosets, we get $$\psi\varphi(x\otimes m) = \sum_{\tau\in \mathcal{F}_\sigma}x\otimes\partial_{\sigma,\tau}m = x\otimes\partial m,$$ as desired.
\end{demonstracao}

\begin{teorema}\label{teo: maps}
    Let $X$ be a finite-dimensional CW-complex with a $G$ action such that for every cell $\sigma$, $G_\sigma$ fixes $\sigma$ pointwise. Let $(\Sigma_p)_{p\in\z}$ be a collection of representatives for the $G$-orbits with the property that for all $\sigma\in \Sigma_p$ the cells in $\partial\sigma$ are represented in $\Sigma_{p-1}$. Then the differential $d^1: E^1_{p,*}\rightarrow E^1_{p-1,*}$ is given by $$d^1\rvert_{H_*(G_\sigma,M_\sigma)}= \sum_{\tau\in\Sigma_{p-1}}[\sigma:\tau](i_{\sigma\tau})_*.$$ 
\end{teorema}

\begin{demonstracao}
    Note that the fact that $X$ is finite-dimensional implies that such a collection of representatives can always be constructed starting from the largest skeleton and choosing representatives using the boundary map.
    
    The property that every $G_\sigma$ fixes $\sigma$ pointwise implies that it is possible to orient each cell of $X$ in such a way that the $G$-action preserves orientations, thus making the action of $G_\sigma$ on $\z_\sigma$ trivial. This determines for each $\sigma$ an isomorphism $M_\sigma \simeq M$ of $G_\sigma$-modules, thus we can simply use coefficients in $M$ for all the groups $H_*(G_\sigma,M)$.

    We also have $G_{\sigma\tau} = G_\sigma$ for $\tau\in\mathcal{F}_\sigma$, since such a $\tau$ is necessarily in the closure of $\sigma$ and hence is fixed by $G_\sigma$. Because of this, the transfer map $t_{\sigma\tau}$ will just be the identity map. This also implies that $\mathcal{F}_\sigma'=\mathcal{F}_\sigma$ since the action of $G_\sigma$ on $\mathcal{F}_\sigma$ is trivial. 

    Since in this case the inclusion $G_{\sigma\tau}\hookrightarrow G_\tau$ is also the identity, we have that the map $u_{\sigma\tau}$ is simply induced by $\partial_{\sigma\tau}$.
    
    Furthermore, the choice of representatives allows the maps $\nu_\tau$ to also be identity maps, since there is no need to map $\tau \in \partial\sigma$ to a different representative.

    From Definition \ref{def: cell complex}, we have that the $\tau$-component of $\partial\sigma$ is given by the incidence number $[\sigma:\tau]$, so the map $\partial_{\sigma\tau}: M \rightarrow M$ takes the generator $\sigma$ of $M_\sigma \simeq M$ to $[\sigma:\tau]\tau \in M_\tau\simeq M$. Thus, if $i_{\sigma\tau}: G_{\sigma}\hookrightarrow G_\tau$ is the inclusion map, the map $\varphi$ in Proposition \ref{prp: varphi} is such that $$\varphi\rvert_{H_*(G_\sigma, M_\sigma)} = \sum_{\tau\in\mathcal{F}'_{\sigma}} \nu_\tau u_{\sigma\tau}t_{\sigma\tau} = \sum_{\tau\in\mathcal{F}_{\sigma}} u_{\sigma\tau} = \sum_{\tau\in\mathcal{F}_{\sigma}}(\partial_{\sigma\tau})_* = \sum_{\tau\in\mathcal{F}_{\sigma}}[\sigma:\tau](i_{\sigma\tau})_*.$$ The choice of $\Sigma_{p-1}$ also allows us to replacle $\mathcal{F}_{\sigma}$ by $\Sigma_{p-1}$ in the sum.
\end{demonstracao} 

\chapter{The homology groups of \texorpdfstring{$\Gamma_0(n)$}{Lg}}
\label{chapter:gamma0n}
In this chapter, we calculate the homology groups of the congruence subgroups $\Gamma_0(n)$ of $\SL_2(\z)$ which appear in the spectral sequence converging to $H_k(\SL_2(\z[1/m]))$. This is also done via the spectral sequence obtained by its action on a contractible 1-dimensional CW-complex, as constructed in the last chapter.

\section{The tree of \texorpdfstring{$\Gamma_0(n)$}{Lg}}

Here, the tree used for the calculations of the homology groups $H_k(\Gamma_0(n))$ is due to \cite{ww1998}. In this section, we elucidate the construction of this tree.

For a natural number $n$, $\Gamma_0(n)$ is defined as the subgroup of $\SL_2(\z)$ given by
\[
\Gamma_0(n) :=\bigg\{ {\mtxx a b c d} \in \SL_2(\z): n\mid c \bigg\}.
\]
In particular, $\Gamma_0(1) = \SL_2(\z)$.

Let $X'$ be the subset of the 
integer lattice $\z^2$ given by 
\[
X' := \left\{\binom{w}{y} \in \z^2: \gcd(w,y)=1\right\}.
\]
The group $\SL_2(\z)$ acts on $X'$ by multiplication on the left, that is, given $\begin{pmatrix}
    a & b \\ c & d
\end{pmatrix}\in \Gamma_0(n),\begin{pmatrix}
    w \\ y
\end{pmatrix}\in X'$, the action is given by $$\begin{pmatrix}
    a & b \\ c & d
\end{pmatrix}\cdot  \begin{pmatrix}
    w \\ y
\end{pmatrix} = \begin{pmatrix}
    aw+by \\ cw+dy
\end{pmatrix}.$$ Define $X$ to be the quotient of $X'$ by the 
action of $\pm I_2$:
\[
X:=X'/\{\pm I_2\}.
\]
We identify, by abuse of notation, $X$ with the subset of $X'$ consisting of points 
$\displaystyle\binom{w}{y}$ with $w \geq 1$, together with the point $\displaystyle\pm\binom{0}{1}$.

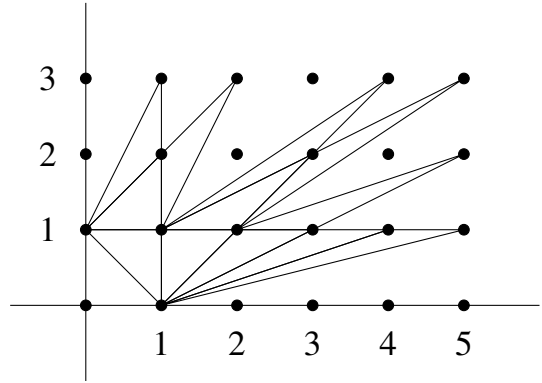
\begin{figure}[H]
    \centering
    \caption{Triangulation formed by $\Delta$}
    \begin{tikzpicture}
        \draw (0,-1)--(0,4);
        \draw (-1,0)--(6,0);
        \draw [fill=black] (0,0) circle (2pt);
        \draw [fill=black] (1,0) circle (2pt);
        \draw [fill=black] (2,0) circle (2pt);
        \draw [fill=black] (3,0) circle (2pt);
        \draw [fill=black] (4,0) circle (2pt);
        \draw [fill=black] (5,0) circle (2pt);

        \draw [fill=black] (0,1) circle (2pt);
        \draw [fill=black] (1,1) circle (2pt);
        \draw [fill=black] (2,1) circle (2pt);
        \draw [fill=black] (3,1) circle (2pt);
        \draw [fill=black] (4,1) circle (2pt);
        \draw [fill=black] (5,1) circle (2pt);

        \draw [fill=black] (0,2) circle (2pt);
        \draw [fill=black] (1,2) circle (2pt);
        \draw [fill=black] (2,2) circle (2pt);
        \draw [fill=black] (3,2) circle (2pt);
        \draw [fill=black] (4,2) circle (2pt);
        \draw [fill=black] (5,2) circle (2pt);

        \draw [fill=black] (0,3) circle (2pt);
        \draw [fill=black] (1,3) circle (2pt);
        \draw [fill=black] (2,3) circle (2pt);
        \draw [fill=black] (3,3) circle (2pt);
        \draw [fill=black] (4,3) circle (2pt);
        \draw [fill=black] (5,3) circle (2pt);

        \draw (0,1)--(1,1)--(1,0)-- cycle;
        \draw (2,1)--(1,1)--(1,0)-- cycle;
        \draw (2,1)--(3,1)--(1,0)-- cycle;
        \draw (3,1)--(4,1)--(1,0)-- cycle;
        \draw (4,1)--(5,1)--(1,0)-- cycle;
        \draw (0,1)--(1,1)--(1,2)-- cycle;
        \draw (2,1)--(1,1)--(3,2)-- cycle;
        \draw (2,1)--(3,1)--(5,2)-- cycle;
        \draw (0,1)--(1,2)--(1,3)-- cycle;
        \draw (1,1)--(1,2)--(2,3)-- cycle;
        \draw (1,1)--(3,2)--(4,3)-- cycle;
        \draw (2,1)--(3,2)--(5,3)-- cycle;

        \node[black,scale=1.1] at (1,-0.5) {$1$};
        \node[black,scale=1.1] at (2,-0.5) {$2$};
        \node[black,scale=1.1] at (3,-0.5) {$3$};
        \node[black,scale=1.1] at (4,-0.5) {$4$};
        \node[black,scale=1.1] at (5,-0.5) {$5$};
        \node[black,scale=1.1] at (-0.5,1) {$1$};
        \node[black,scale=1.1] at (-0.5,2) {$2$};
        \node[black,scale=1.1] at (-0.5,3) {$3$};
    \end{tikzpicture}
    \label{fig:triang}
    \fautor
\end{figure}

Let $\Delta$ be the set of \textit{triangles} whose vertices are vectors
\[
\binom{w}{y},\binom{x}{z} \text{ and }\binom{t}{u} \in X
\]
for which 
\[
t = w + x, \ \ u = y + z,\ \text{and}\ wz - xy = \pm 1.
\]
Note that, since $zt-xu=wz - xy = \pm 1$, we have $\displaystyle\binom{t}{u}\in X$. 
We shall represent such a triangle in matrix form, 
\[
T = \begin{pmatrix}w & x & t \\ y & z & u\end{pmatrix}.
\]
Without loss of generality, in $T$, we may require that 
\[
\det\begin{pmatrix}
w & x\\
y &z 
\end{pmatrix}=1.
\]
Triangles of this form will turn out to provide a triangulation of a subset of the right half-plane, as seen in Figure \ref{fig:triang}.

The action of $\SL_2(\z)$ (and thus of $\Gamma_0(n)$) on $\Delta$ is given by multiplication 
on the left.

\begin{proposicao}\label{LC}
Linear combinations of vertices in $X$ with coprime coefficients will also be vertices in $X$.
\end{proposicao}
\begin{demonstracao}Let $a,b$ be integers and suppose $r$ and $s$ are both linear combinations of $a,b$ 
with integer coefficients. If $d \mid a$ and $d\mid b$, then $d\mid r$ and $d\mid s$. More 
importantly, if $r$ and $s$ are coprime, then $a$ and $b$ will also be coprime. 

Now, consider a triangle with vertices $\displaystyle\binom{w}{y},\binom{x}{z},\binom{w+x}{y+z}$ 
and such that $wz-yx=1$. Let $r,s$ be coprime integers. Note that 
\[
r\begin{pmatrix}
    w \\ y
\end{pmatrix}+s\begin{pmatrix}
    x \\ z
\end{pmatrix}=\begin{pmatrix}
    rw+sx \\ ry+sz
\end{pmatrix} = \begin{pmatrix}
    w & x \\ y & z
\end{pmatrix}\begin{pmatrix}
    r\\ s
\end{pmatrix} \hspace{20pt} \implies \hspace{20pt} \begin{pmatrix}
    r\\ s
\end{pmatrix} = \begin{pmatrix}
    w & x \\ y & z
\end{pmatrix}^{-1}\begin{pmatrix}
    rw+sx \\ ry+sz
\end{pmatrix},
\]
that is, $r$ and $s$ are linear combinations of $rw+sx$ and $ry+sz$ with integer coefficients. 
Since $r,s$ are coprime, then $rw+sx$ and $ry+sz$ are also coprime.
\end{demonstracao}

Now, we turn to the question of equivalence classes under the action of $\Gamma_0(n)$.
Given a triangle 
$T = \begin{pmatrix}
w & x & t \\ y & z & u
\end{pmatrix}$, 
define three associated elements of $\SL_2(\z)$ by 
\[
T_1 = \begin{pmatrix}
w & x \\ y & z
\end{pmatrix}, \
T_2 = \begin{pmatrix}
t& -w \\ u & -y
\end{pmatrix}, \ 
T_3 = \begin{pmatrix}
x & -t \\ z & -u
\end{pmatrix}.
\]

Define 
\[
\epsilon:\SL_2(\mathbb{Z})\rightarrow \Delta, \ \ \ \ 
\begin{pmatrix}
a & b \\ c & d
\end{pmatrix}
\mapsto 
\begin{pmatrix}
a & b & a+b \\ c & d & c+d
\end{pmatrix}.
\]
This induces a projection 
\[
\pi: \SL_2(\z)/\Gamma_0(n) \rightarrow \Delta/\Gamma_0(n).
\]

\begin{proposicao}\label{II}
The inverse image of a class of an element $T$ by $\pi$ consists of either one or three elements.
\end{proposicao}

\begin{demonstracao}
If a triangle has vertices $\displaystyle\binom{w}{y},\binom{x}{z},\binom{t}{u}$, there are 
$6\cdot 4 \cdot 2 = 48$ matrices who could represent it. If we insist on representing it with 
the third column being the sum of the first two, we narrow down to 12 representations, and if 
we insist that the determinant of the first minor be 1, we further reduce to 6 representations, 
which contains 3 pairs of matrices equivalent under the action of $-I_2$, which are precisely 
$\pm \epsilon(T_1),\pm \epsilon(T_2)$ and $\pm \epsilon(T_3)$.

To exemplify, consider we choose $\displaystyle\binom{x}{z}$ as the third column. Then the first 
two must be $\displaystyle\binom{t}{u}$ and $\displaystyle-\binom{w}{y}$, in this order. Thus, our 
representation is $\epsilon(T_2)$. Choosing $\displaystyle-\binom{x}{z}$ would lead to $-\epsilon(T_3)$.
This implies, then, that $\epsilon^{-1}(T) = \{\pm T_1,\pm T_2,\pm T_3\}$.

Note that the projection $\pi$ takes a class $[X]$ and maps it to $\epsilon(X)$, which is well 
defined in $\Delta/\Gamma_0(n)$ since matrices in the same coset will give $\Gamma_0(n)$-equivalent triangles. 

If $S$ and $T$ are $\Gamma_0(n)$-equivalent, then $T_i=\Gamma_0(n)S_j$ for some indexes $i,j$, so that each column 
in $\epsilon(S_j)$ corresponds to a column in $\epsilon(T_i)$, thus each coset $[S_1],[S_2],[S_3]$ 
corresponds to a coset $[T_1],[T_2],[T_3]$. For example, if $T_1=\Gamma_0(n)S_2$, then $[T_2]=[S_3]$ and 
$[T_3]=[S_1]$.

Furthermore, we have 
\[
T_2T^{-1}_1 = 
\begin{pmatrix}
(w+x)z+wy & -(w+x)x-w^2 \\ z^2+yz+y^2 & -(y+z)x-wy
\end{pmatrix},
\]
\[
T_3T^{-1}_1 =
\begin{pmatrix}
xz+(w+x)y & -x^2-w(w+x) \\ z^2+yz+y^2 & -xz-(y+z)w
\end{pmatrix}.
\]
So the three cosets $[T_1],[T_2],[T_3]$ will be pairwise distinct or all the same if 
$y^2+yz+z^2 \equiv 0 \pmod{n}$.
\end{demonstracao}

\begin{definicao}
A triangle $T$ will be called of Type 2 if $[T_1]=[T_2]=[T_3]$ and of Type 1 otherwise.
\end{definicao}

\begin{lema}\label{T-S}
Given triangles $T$, $S$, if there is $A\in \Gamma_0(n)$ such that the action of $A$ takes 
$T$ to $S$, then $A$ must be of the form $\pm S_i(T_1)^{-1}$ for some $i=1,2,3$. 
\end{lema}

\begin{demonstracao}
Since the action is independent of representation, we can choose the representative
$\epsilon(T_1)$ to represent $T$, so that by $AT=S$, we mean $A\epsilon(T_1)$ is a 
representative of the triangle $S$. Since the third column of $\epsilon(T_1)$ is 
the sum of the first two, so it will be for $A\epsilon(T_1)$. Note also that 
$(A\epsilon(T_1))_1 = AT_1$, so $\det(A)=\det(T_1)=1$ implies that 
$\det((A\epsilon(T_1))_1)=1$. 

Essentially, we have that $A\epsilon(T_1)$ must be a representative of $S$ whose 
third column is the sum of the first two and whose first minor has determinant 1, 
that is, it must be of the form $\pm\epsilon(S_i)$. 
Thus, $AT_1=\pm S_i$ and we obtain $A = \pm S_iT_1^{-1}$ for some $i$, as desired.
\end{demonstracao}

\begin{definicao}
For every $k\in\mathbb{Z}$, let $R^{(k)}$ denote the triangle 
\[
R^{(k)}:=\begin{pmatrix}
1 & 0 & 1 \\ k & 1 & 1+k
\end{pmatrix}.
\]
\end{definicao}

\begin{observacao}
    Note that if $k=qn+r$, then \[
\begin{pmatrix}
1 & 0 \\ qn & 1
\end{pmatrix}\cdot \begin{pmatrix}
1 & 0 & 1 \\ r & 1 & 1+r
\end{pmatrix} = \begin{pmatrix}
1 & 0 & 1 \\ k & 1 & 1+k
\end{pmatrix}
\] implying that a matrix $R^{(k)}$ is $\Gamma_0(n)$-equivalent to $R^{(r)}$ where $r$ is its remainder mod $n$. 
\end{observacao}

\begin{proposicao}\label{R^k}
If $R^{(k)}$ and $R^{(r)}$ are $\Gamma_0(n)$-equivalent, there is only one matrix 
(up to sign) in $\Gamma_0(n)$ that exchange them.
\end{proposicao}
\begin{demonstracao}
Suppose $0\leq k,r<n$ are such that $R^{(k)}$ and $R^{(r)}$ are equivalent. So there is 
$A\in \Gamma_0(n)$ such that $(AR^{(k)})_1$ is either $\pm(R^{(r)})_1$, $\pm(R^{(r)})_2$ or 
$\pm(R^{(r)})_3$ so that $A$ is one of the following matrices: 
\[
\pm\begin{pmatrix}
1 & 0 \\ r-k & 1
\end{pmatrix}, \hspace{20pt} \pm\begin{pmatrix}
1+k & -1 \\ 1+r+rk & -r
\end{pmatrix}, \hspace{20pt} \pm\begin{pmatrix}
k & -1 \\ 1+k+rk & -1-r
\end{pmatrix}.
\]
Note that the first matrix is never in $\Gamma_0(n)$, and the other two can't be simultaneously in 
$\Gamma_0(n)$, since $1+r+rk \equiv 0$ and $1+k+rk\equiv 0$ would imply $r\equiv k \pmod{n}$, which 
is impossible.
\end{demonstracao}

\begin{proposicao}
Given a triangle $T$, $T$ is $\Gamma_0(n)$-equivalent to $R^{(k)}$ if and only if an 
integer $k$ exists such that at least one of $u+ky, z+ku$ or $y-kz$ is a multiple of $n$.
\end{proposicao}

\begin{demonstracao}
Note that 
\[
\begin{pmatrix}
a & b \\ nc & d
\end{pmatrix}\cdot \begin{pmatrix}
1 & 0 \\ k & 1
\end{pmatrix} = \begin{pmatrix}
a+bk & b \\ nc+kd & d
\end{pmatrix}
\]

In other words, if there exists $A\in\Gamma_0(n)$ such that $A(R^{(k)})_1$ is $T_1, T_2$ or $T_3$, then either $y,z$ or $u$ must be given by $nc+kd$ where $d$ is the lower right entry of $T_i$.

If $A(R^{(k)})_1 = T_1$, this condition becomes $y=nc+kz$, so $n\mid (y-kz)$. If $A(R^{(k)})_1 = T_2$, the condition is $u=nc-ky$ and we get $n\mid (u+ky)$ and finally, if $A(R^{(k)})_1 = T_3$, the condition is $z=(nc-ku)$ which gives $n\mid z+ku$.

For the converse statement, if $n\mid (y-kz)$, then $y=nc+kz$ for some $c$ and \[
\begin{pmatrix}
w-xk & x \\ nc & z
\end{pmatrix}\cdot \begin{pmatrix}
1 & 0 \\ k & 1
\end{pmatrix} = \begin{pmatrix}
w & x \\ nc+kd & z
\end{pmatrix}= \begin{pmatrix}
w & x \\ y & z
\end{pmatrix}.
\] So $T$ is $\Gamma_0(n)$-equivalent to $R^{(k)}$. The cases $n\mid (u+ky)$ and $n\mid (z+ku)$ will give the same result analogously.
\end{demonstracao}

\begin{corolario}\label{yzu coprime to n}
If at least one of $y,z,u$ is coprime to $n$, then $T$ is $\Gamma_0(n)$-equivalent to some 
$R^{(k)}$.
\end{corolario}

\begin{demonstracao}
If $\gcd(y,n)=1$, then $ay+bn=1$ for some $a,b$ and thus, $uay+ubn=u$, so that taking $c=ub$ 
and $k=-ua$, we obtain $cn = u+ky$, as desired. The cases for $z$ and $u$ are analogous.
\end{demonstracao}

\begin{proposicao}\label{type2}
Every Type $2$ triangle is $\Gamma_0(n)$-equivalent to a unique triangle $R^{(i)}$.
\end{proposicao}

\begin{demonstracao}
Let $T = \begin{pmatrix}
w & x & t \\ y & z & u
\end{pmatrix}$ 
be a Type 2 triangle. Then $y^2+yz+z^2 \equiv 0\pmod{n}$ and we can show that 
both $y$ and $z$ are coprime to $n$.

Indeed, let $p$ be a prime that divides $n$. Since $wz-xy=1$, we have that $\gcd(y,z)=1$ 
and thus we can't have both $y\equiv z \equiv 0 \pmod{p}$. Suppose, then, that $z\not\equiv 0$ 
(the other case is symmetric) and consider $r = y/z$ in $\mathbb{F}_p$. Then 
\[
y^2+yz+z^2 = 0 \implies z^2(r^2+r+1)=0 \implies r^2+r+1=0
\]
in $\mathbb{F}_p$, so that $r\neq 0$ and thus $y\not\equiv 0 \pmod{p}$. We conclude, then, that 
none of the primes that divide $n$ will divide $y$ and $z$, so that $\gcd(y,n)=\gcd(z,n)=1$.

Furthermore, if $r$ and $k$ are distinct roots of $X^2+X+1$ mod $n$, then $R^{(r)}$ and $R^{(k)}$ 
are not exchanged. Recall that for these matrices to be exchanged, it is necessary for either 
$1+r+rk \equiv 0\pmod{n}$ or $1+k+rk \equiv 0\pmod{n}$, which are mutually exclusive conditions if $r\neq k$ (see Proposition \ref{R^k}).

Indeed, suppose $1+k+rk\equiv 0 \pmod{n}$, so that $k(1+r)\equiv -1$ and $k(1+r)^2\equiv -(1+r)$. 
But $(1+r)^2 = 1+2r+r^2 \equiv 1+2r-(1+r) = r$, and substituting we get $kr \equiv -(1+r)$, so 
that $1+r+kr \equiv 0$. By the previous remarks, this is absurd. In a similar way, if we suppose 
$1+r+rk\equiv 0 \pmod{n}$, we will conclude $1+k+rk\equiv 0 \pmod{n}$.
\end{demonstracao}

\begin{definicao}
Let $b(n)$ be the number of roots of the polynomial $X^2 + X +1$ in $\z/n$.
\end{definicao}

\begin{corolario}\label{cor: type 2}
There are exactly $b(n)$ equivalence classes of Type $2$ triangles.
\end{corolario}

\begin{demonstracao}
If $T$ is a Type 2 triangle, then $T$ is $\Gamma_0(n)$-equivalent to a unique $R^{(k)}$. The condition 
for a triangle $R^{(k)}$ to be Type 2 is for $k$ to be a root of $X^2+X+1$ mod $n$, and none 
of the different roots give equivalent triangles.

Thus, the number of equivalence classes of Type 2 triangles is equal to the number of roots 
of $X^2+X+1$ mod $n$, which we defined as $b(n)$. 
\end{demonstracao}

Let $Y$ be the dual graph of this triangulation. More precisely, the vertices of $Y$ are the 
triangles of $\Delta$ and a pair of vertices shares an edge if the corresponding triangles 
are adjacent. 

\begin{proposicao}\label{cubic-tree}
The graph $Y$ is the cubic tree.
\end{proposicao}
\begin{demonstracao}
Note that, given an edge with vertices $\displaystyle u=\binom{w}{y}$ and $\displaystyle v = \binom{x}{y}$, 
if $u,v,r$ are the vertices of a triangle, then either $u+v=r$, $u+r=v$ or $r+v=u$ (up to antipodes). This 
implies, then, that $r$ must be either $u+v$ or $u-v$ (again, up to antipodes). Thus, an edge $u,v$ 
exclusively determines two triangles, namely, $u,v,u+v$ and $u,v,u-v$. Furthermore, since a triangle 
has three edges, it must have exactly three neighboring triangles.

Let us construct a path in our graph. Suppose we start with a vertex corresponding to a triangle $T$ 
with vertices $u,v,u+v$.

\begin{figure}[H]
    \centering
    \caption{Paths of lenght 1 starting on triangle $u,v,u+v$}
    \begin{tikzpicture}[scale=0.8]
        \draw [fill=black] (0,0) circle (2pt);
        \node[scale=1] at (-0.3,0.3) {$u$};
        \draw [fill=black] (3,0) circle (2pt);
        \node[scale=1] at (3.8,0.3) {$u+v$};
        \draw [fill=black] (1.5,2.598) circle (2pt);
        \node[scale=1] at (1.5,3) {$2u+v$};
        \draw [fill=black] (1.5,-2.598) circle (2pt);
        \node[scale=1] at (1.3,-3) {$v$};
        \draw [fill=black] (4.5,-2.598) circle (2pt);
        \node[scale=1] at (5.2,-3) {$u+2v$};
        \draw [fill=black] (-1.5,-2.598) circle (2pt);
        \node[scale=1] at (-1.8,-3) {$u-v$};

        \draw (0,0)--(3,0)--(1.5,2.598)--(0,0);
        \draw (0,0)--(3,0)--(1.5,-2.598)--(0,0);
        \draw (0,0)--(1.5,-2.598)--(-1.5,-2.598)--(0,0);
        \draw (3,0)--(1.5,-2.598)--(4.5,-2.598)--(3,0);

        \draw [fill=black] (1.5,0.866) circle (2pt);
        \draw [fill=black] (0,-1.732) circle (2pt);
        \draw [fill=black] (1.5,-0.866) circle (2pt);
        \draw [fill=black] (3,-1.732) circle (2pt);

        \draw[thick] (1.5,-0.866)--(1.5,0.866);
        \draw[thick] (1.5,-0.866)--(0,-1.732);
        \draw[thick] (1.5,-0.866)--(3,-1.732);
    \end{tikzpicture}
    \label{fig:1path}
    \fautor
\end{figure}
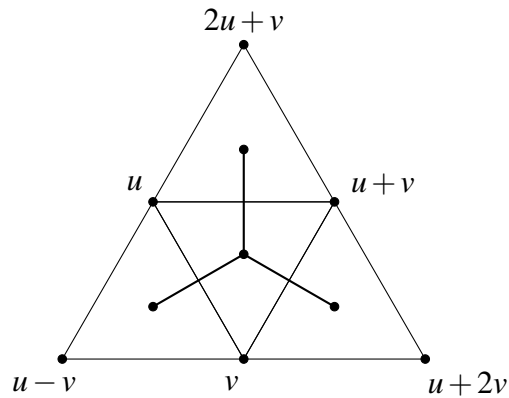

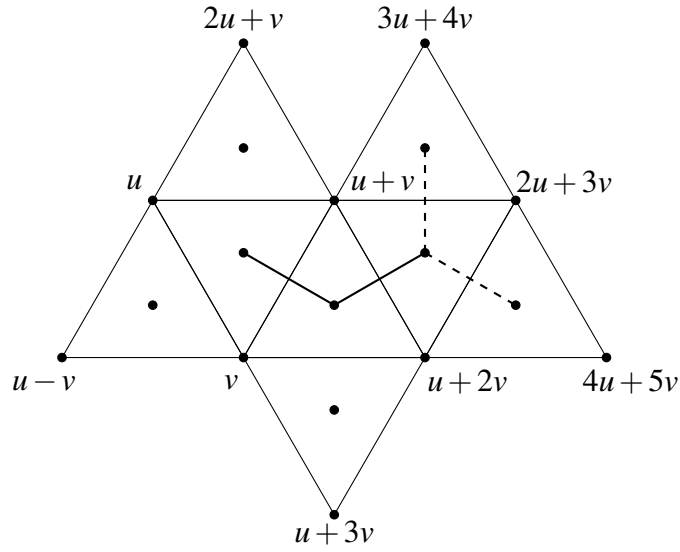
\begin{figure}[H]
    \centering
    \caption{Paths of lenght 2 starting from the triangle $u,v,u+v$ with second vertex in $u+v,v,u+2v$}
    \begin{tikzpicture}[scale=0.8]
        \draw [fill=black] (0,0) circle (2pt);
        \node[scale=1] at (-0.3,0.3) {$u$};
        \draw [fill=black] (3,0) circle (2pt);
        \node[scale=1] at (3.8,0.3) {$u+v$};
        \draw [fill=black] (6,0) circle (2pt);
        \node[scale=1] at (6.8,0.3) {$2u+3v$};
        \draw [fill=black] (1.5,2.598) circle (2pt);
        \node[scale=1] at (1.5,3) {$2u+v$};
        \draw [fill=black] (1.5,-2.598) circle (2pt);
        \node[scale=1] at (1.3,-3) {$v$};
        \draw [fill=black] (4.5,-2.598) circle (2pt);
        \node[scale=1] at (5.2,-3) {$u+2v$};
        \draw [fill=black] (3,-5.196) circle (2pt);
        \node[scale=1] at (3,-5.5) {$u+3v$};
        \draw [fill=black] (-1.5,-2.598) circle (2pt);
        \node[scale=1] at (-1.8,-3) {$u-v$};
        \draw [fill=black] (4.5,2.598) circle (2pt);
        \node[scale=1] at (4.5,3) {$3u+4v$};
        \draw [fill=black] (7.5,-2.598) circle (2pt);
        \node[scale=1] at (7.9,-3) {$4u+5v$};

        \draw (0,0)--(3,0)--(1.5,2.598)--(0,0);
        \draw (0,0)--(3,0)--(1.5,-2.598)--(0,0);
        \draw (0,0)--(1.5,-2.598)--(-1.5,-2.598)--(0,0);
        \draw (3,0)--(1.5,-2.598)--(4.5,-2.598)--(3,0);
        \draw (1.5,-2.598)--(4.5,-2.598)--(3,-5.196)--(1.5,-2.598);
        \draw (3,0)--(4.5,-2.598)--(6,0)--(3,0);
        \draw (3,0)--(6,0)--(4.5,2.598)--(3,0);
        \draw (6,0)--(4.5,-2.598)--(7.5,-2.598)--(6,0);

        \draw [fill=black] (1.5,0.866) circle (2pt);
        \draw [fill=black] (0,-1.732) circle (2pt);
        \draw [fill=black] (1.5,-0.866) circle (2pt);
        \draw [fill=black] (3,-1.732) circle (2pt);
        \draw [fill=black] (3,-3.464) circle (2pt);
        \draw [fill=black] (4.5,-0.866) circle (2pt);
        \draw [fill=black] (6,-1.732) circle (2pt);
        \draw [fill=black] (4.5,0.866) circle (2pt);

        \draw[thick] (1.5,-0.866)--(3,-1.732);
        \draw[thick] (3,-1.732)--(4.5,-0.866);
        \draw[thick, dashed] (4.5,-0.866) -- (4.5,0.866);
        \draw[thick, dashed] (4.5,-0.866) -- (6,-1.732);
    \end{tikzpicture}
    \fautor
    \label{fig:2paths}
\end{figure}

If $S$ is a neighboring triangle, then it must share an edge with $T$ (thus, a pair of vertices). If 
the pair of vertices is $u,v$, the two triangles with those vertices are $u,v,u+v = T$  and $u,v,u-v$. 
If the pair is $u,u+v$, then the triangles are $u,u+v,u-(u+v)= u,u+v,v = T$ and $u,u+v,2u+v$. If the 
pair is $v,u+v$, the triangles are $u,v,u+v$ and $v,u+v,u+2v$. Thus, we have only three possible 
1-paths as seen in Figure \ref{fig:1path}.

Suppose we chose $v,u+v,u+2v$ for our 1-path. We find the three neighboring triangles to be 
$v,u+v,u = T$, $v,u+2v,u+3v$ and $u+v,u+2v,2u+3v$. Choosing the next triangle as $u+v,u+2v,2u+3v$,
we find the neighboring triangles (that won't backtrack) to be $u+v,2u+3v,3u+4v$ and $u+2v,2u+3v,4u+5v$, as in Figure \ref{fig:2paths}.

Following this logic, we find that geodesics starting at $u,v,u+v$ will end at triangles 
$au+bv, cu+dv, (a+c)u+(b+d)v$. Recall that for these to be triangles, we must have 
\[
\det\begin{pmatrix}
aw+bx & cw+dz \\ ay+bz & cx+dz
\end{pmatrix} = \pm1 \Leftrightarrow (ad-bc)(wz-xy)=\pm1
\]
that is, $ad-bc=\pm 1$.

This implies, then, that the graph is connected, since we can take $\displaystyle u=\binom{1}{0}, 
v=\binom{0}{1}$ and find that any triangle $\displaystyle\binom{a}{b},\binom{c}{d},\binom{a+b}{c+d}$ 
is connected to $\displaystyle \binom{1}{0},\binom{0}{1},\binom{1}{1}$ by a path.

Furthermore, for there to be a cycle, we would need for $au+bv$ to be either $\pm u, \pm v$ or 
$\pm (u+v)$. Suppose, for example, that $au+bv=u$. Then 
\[
\begin{cases}
bx=(1-a)w \\ 
bz=(1-a)y
\end{cases} \implies 1=wz-xy =\frac{1-a}{b}wy-\frac{1-a}{b}wy = 0
\]
which is absurd. Similarly, the other cases will lead to contradictions of the type $1=0$. 

Thus, we have that our graph is a tree. Since each vertex is connected to three other vertices, 
we have that it is the cubic tree.
\end{demonstracao}

We now turn to the question of which edges of the cubic tree Y will get inverted by our action, so we can alter 
our tree by breaking each of these edges in the middle, adding an extra vertex. 

Note that if an edge $y$ with $o(y)=T$, $t(y)=S$ is inverted by a matrix $A\in \Gamma_0(n)$ (that is, 
$Ay=\overline{y}$) then $$o(Ay)=A(o(y)) \implies AT = o(\overline{y})=S,$$ and we can also conclude 
that $AS=T$. Thus, $T$ and $S$ are two neighboring triangles that get swapped by $A$.

\begin{proposicao}\label{swap}
A triangle $T$ is swapped with a neighboring triangle $S$ if and only if $y^2+z^2\equiv 0 \pmod{n}$.
\end{proposicao}

\begin{demonstracao}
If the common columns of $T$ and $S$ are given by vectors $u,v$,  then $u+v$ and 
$u-v$ are the only possible remaining columns in $S$ and $T$. Let 
\[
T = \begin{pmatrix}
w & x & w+x \\ y& z & y+z
\end{pmatrix}, \hspace{20pt} S=\begin{pmatrix}
-w & x & x-w \\ -y & z & z-y
\end{pmatrix}.
\]
Lemma \ref{T-S}, $A$ must be given by $\pm S_iT_1^{-1}$ where $i=1,2,3$. We can check that 
$S_2T_1^{-1}$ maps $T$ to $S$, but it won't map $S$ to $T$. Similarly, $S_3T_1^{-1}$ 
maps $T$ to $S$ but not $S$ to $T$. We get, then, that the only two matrices that swap 
$S$ and $T$ are given by 
\[
\pm S_1T_1^{-1} = \pm\begin{pmatrix}
wy+xz & -(w^2+x^2) \\ y^2+z^2 & -(wy+xz)
\end{pmatrix}.
\]
\end{demonstracao}

\begin{proposicao}
If $T$ and $S$ are neighboring triangles that gets swapped, then $T$ and $S$ are 
$\Gamma_0(n)$-equivalent to a pair $R^{(k)}, R^{(k+1)}$.
\end{proposicao}

\begin{demonstracao}
The proof is very similar to the case of Type 2 triangles (see Proposition \ref{type2}).

Let $T = \begin{pmatrix}
w & x & t \\ y & z & u
\end{pmatrix}$ 
be a Type 3 triangle and $S= \begin{pmatrix}
-w & x & x-w \\ -y & z & z-y
\end{pmatrix}$ its neighboring triangle as described above. Then $y^2+z^2 \equiv 0\pmod{n}$ and we  show that 
both $y$ and $z$ are coprime to $n$.

Let $p$ be a prime that divides $n$. Then $wz-xy=1$ implies $\gcd(y,z)=1$ 
and thus we can't have both $y\equiv z \equiv 0 \pmod{o}$. Suppose, then, that $z\not\equiv 0$ and consider $r = y/z$ in $\mathbb{F}_p$. Then 
\[
y^2+z^2 = 0 \implies z^2(r^2+1)=0 \implies r^2+1=0
\]
in $\mathbb{F}_p$, so that $r\neq 0$ and thus $y\not\equiv 0 \pmod{p}$. We conclude, then, that 
none of the primes that divide $n$ will divide $y$ and $z$, so that $\gcd(y,n)=\gcd(z,n)=1$.

Then, Corollary \ref{yzu coprime to n} implies that $T$ is equivalent to some $R^{(k)}$ with $k^2+1\equiv 0 \pmod{n}$, so $S$ is equivalent to $R^{(k-1)}$. 

If $r$ and $k$ are distinct roots of $X^2+1$ mod $n$, it may be the case that $R^{(r)}$ and $R^{(k)}$ 
are exchanged. For example, if $n=5$, then 2 and 3 are both roots of $X^2+1$ mod 5, but $R^{(2)}$ and $R^{(3)}$ are swapped because $2=3-1$. 
\end{demonstracao}

The above proposition shows that $\Gamma_0(n)$ acts by inversion on some edges of the cubic tree
$Y$. For this reason we shall change our tree by adding a new vertex at the midpoint of each edge 
that is inverted by an element of $\Gamma_0(n)$. If an edge of the tree is inverted by $\Gamma_0(n)$, 
the two triangles corresponding to the ends of this edge have two of their vertices in common which 
are exchanged by an element of $\Gamma_0(n)$.

Let $e,e'$ be edges in $Y$ with $e$ having vertices $T,S$ and $e'$ having vertices $T',S'$. If there are $A,B\in\Gamma_0(n)$ such that $Ae = \overline{e}$ and $Be = e'$, then $BAB^{-1}e'= \overline{e'}$. 

In other words, if $e'$ is equivalent to an edge $e$ inverted by $\Gamma_0(n)$, then $e'$ is also inverted by $\Gamma_0(n)$. 

So if $E$ is an extra vertex added to $Y$ to avoid inversion of an edge $e$, the action of a matrix $A\in \Gamma_0(n)$ will only map $e$ to other invertible edges, thus it makes sense to describe the action of $A\in \Gamma_0(n)$ as mapping $E$ to the vertex added in $Ae$. 

\begin{definicao}
We call the added new vertices to Y, vertices of Type 3. The vertices of $Y$ corresponding to triangles 
of Types 1 and 2, are called vertices of Types 1 and 2, respectively.
\end{definicao}

\begin{definicao}
Let $c(n)$ be the number of roots of the polynomial $X^2 +1$ in $\z/n$.
\end{definicao}

\begin{corolario}\label{type3}
There are exactly $c(n)$ equivalence classes of Type $3$ vertices.
\end{corolario}

\begin{demonstracao}
Each Type 3 vertex corresponds to an added vertex between a pair of neighbors that are 
swapped, which themselves correspond to a pair $R^{(k)}, R^{(k-1)}$ where $k$ is a root of $X^2+1$, so there are at most $c(n)$ equivalence classes of Type 3 vertices.

Suppose $r,k$ are two distinct roots of $X^2+1$ with $R^{(r)}$ and $R^{(k)}$ equivalent and denote by $E_r,E_k$ the corresponding Type 3 vertices added between $R^{(r)}, R^{(r-1)}$ and $R^{(k)}, R^{(k-1)}$, respectively. 

By the description of the action of $\Gamma_0(n)$ on Type 3 vertices, we have that $E_r$ and $E_k$ are equivalent iff there exists a matrix $A\in\Gamma_0(n)$ such that $AR^{(k)}=R^{(r)}$ and $AR^{(k-1)} = R^{(r-1)}$. From Proposition \ref{R^k}, $A$ must be given by one of the two following matrices: $$\pm\begin{pmatrix}
1+k & -1 \\ 1+r+rk & -r
\end{pmatrix}, \hspace{20pt} \pm\begin{pmatrix}
k & -1 \\ 1+k+rk & -1-r
\end{pmatrix}.$$

It follows that neither of these matrices will map $R^{(k-1)}$ to $R^{(r-1)}$. Thus, the added vertices $E_r,E_k$ are not equivalent, so each root of $X^2+1$ gives a distinct equivalence class.
\end{demonstracao}

\begin{definicao}
    Let $a(n)$ be the index of $\Gamma_0(n)$ in $\SL_2(\z)$, that is, $a(n)\coloneq [\SL_2(\z):\Gamma_0(n)]$.
\end{definicao}

\begin{corolario}\label{v(n)}
The number of  $\Gamma_0(n)$-equivalence classes of vertices of $Y$ is given by 
\[
v(n) := \frac{a(n)+2b(n)}{3}+c(n).
\]
\end{corolario}

\begin{demonstracao}
Let $d(n)$ be the number of equivalence classes 
of Type 1 triangles. We know that each of these classes corresponds to 3 cosets in $\SL_2(\z)/\Gamma_0(n)$, and 
each class of a Type 2 triangle corresponds to a single coset in $\SL_2(\z)/\Gamma_0(n)$, so that $3d(n)+b(n)=a(n)$. 
Thus $d(n) = (a(n)-b(n))/3$, so that the total number of equivalence classes of vertices is 
\[
d(n)+b(n)+c(n) = \frac{a(n)-b(n)}{3} +b(n)+c(n) = \frac{a(n)+2b(n)}{3}+c(n) = v(n).
\]
\end{demonstracao}

\begin{proposicao}\label{e(n)}
The number of $\Gamma_0(n)$-equivalence classes of edges of $Y$ is given by 
\[
e(n) := \frac{a(n)+c(n)}{2}.
\]
\end{proposicao}

\begin{demonstracao}
Recall that in the dual graph, an edge $e$ connects two vertices if their corresponding 
triangles are adjacent, that is, if they have an edge in common. Thus, if triangles $T,S$ 
have common vertices $\displaystyle\binom{w}{y},\binom{x}{z}$, we can represent the 
corresponding edge by a matrix in $\Gamma_0(n)$.

If an edge $e$ has vertices $\displaystyle\binom{w}{y},\binom{x}{z}$ there are four 
representatives of it that lie in $\Gamma_0(n)$, specifically, 
\[
\pm\begin{pmatrix}
w & x \\
y & z
\end{pmatrix} 
\hspace{20pt}\text{and}\hspace{20pt} \pm
\begin{pmatrix}
x & -w \\
z & -y
\end{pmatrix}.
\]
Note that 
\[
\begin{pmatrix}
x & -w \\
z & -y
\end{pmatrix}\begin{pmatrix}
w & x \\
y & z
\end{pmatrix}^{-1} = \begin{pmatrix}
x & -w \\
z & -y
\end{pmatrix}\begin{pmatrix}
z & -x \\
-y & w
\end{pmatrix} = \begin{pmatrix}
wy+xz & -(x^2+w^2) \\
y^2+z^2 & -(wy+xz)
\end{pmatrix}.
\]

Consider, then, the quotient map $\SL_2(\z)/\Gamma_0(n) \rightarrow \Delta'/\Gamma_0(n)$, 
where $\Delta'$ is the set of edges of the triangles of $\Delta$, defined similarly as in the 
case for triangles. It is obviously surjective and each edge lifts to the equivalence classes 
of the two matrices listed above, which correspond to the same coset if the edge in question 
is inverted, that is, corresponds to a vertex of Type 3.

Let $d'(n)$ be the number of equivalences classes of non-inverted edges. Each of those 
corresponds to two classes in $\SL_2(\z)/\Gamma_0(n)$, so that $2d'(n)+c(n) = a(n)$. Thus, 
the number of equivalence classes of inverted edges in total is given by 
\[
d'(n)+c(n) = \frac{a(n)-c(n)}{2}+c(n) = \frac{a(n)+c(n)}{2}=e(n).
\]
\end{demonstracao}

\begin{proposicao}\label{stab}
The stabilizers of Type $1$, Type $2$ and Type $3$ vertices of $Y$ are isomorphic to $\z/2$, 
$\z/6$ and $\z/4$, respectively. The stabilizers of edges are isomorphic to $\z/2$.
\end{proposicao}

\begin{demonstracao}
Given a triangle with vertices $\displaystyle\binom{w}{y},\binom{x}{z},\binom{w+x}{y+z}$, 
we are looking for a matrix $A\in \Gamma_0(n)$ such that $AT$ is a matrix representing $T$. 
As in the Lemma \ref{T-S}, 
\[
A \in \{\pm I_2, \pm T_2T_1^{-1},\pm T_3T_1^{-1}\}.
\]
Note that $\pm I_2$ will stabilize any triangle and $\pm T_2T_1^{-1},\pm T_3T_1^{-1}$ will
stabilize triangles of Type 2 (since those being
in $\Gamma_0(n)$ would imply $T$ being of Type 2).

Recall that the characteristic polynomial of a $2\times 2$ matrix $M$ is given by 
$\lambda^2-\text{tr}(M)\lambda+\det(M)$ and that any matrix is annihilated by their 
characteristic polynomials (see Theorem 4 in Section 6.3 of \cite{hoffmankunze1971}), giving the following expressions: 
\[
M^2 = \text{tr}(M)M-\det(M)I_2,
\]
\[
M^3 = \text{tr}(M)M^2-\det(M)M = ((\text{tr}(M))^2-\det(M))M-\text{tr}(M)\det(M)I_2.
\]
Since $\det(T_1)=\det(T_2)=\det(T_3)=1$, we have $\det(\pm T_2T_1^{-1})=\det(\pm T_3T_1^{-1})=1$, $\text{tr}(T_2T_1^{-1}) = wz-xy=1$, $\text{tr}(T_3T_1^{-1}) = xy-wz=-1$ so that 
$(T_2T_1^{-1})^3=-I_2$ and $(T_3T_1^{-1})^3=I_2$. 

One can check that $T_2T_1^{-1}-T_3T_1^{-1} = I_2$, so that by letting $S = T_2T_1^{-1}$, we get 
\[
S^2 = S-I_2 = T_3T_1^{-1},
\]
\[
S^3 =  -I_2,
\]
\[
S^4 = -S = -T_2T_1^{-1},
\]
\[
S^5 = I_2-S = -T_3T_1^{-1},
\]
\[
S^6 = I_2,
\]
and thus the 6 matrices that stabilize a Type 2 triangle $T$ are the ones in 
$\langle T_2T_1^{-1}\rangle \simeq \z/6$.

Furthermore, the matrices $A$ that stabilize type 3 vertices are those who interchange 
adjacent triangles. From Proposition \ref{e(n)}, these matrices are of type $\begin{pmatrix}
    wy+xz & -(x^2+w^2) \\ y^2+z^2 & -(wy+xz) 
\end{pmatrix}$, so their trace is null and their determinant is 1, giving that $A^2=-I_2$ and their order is 4.

Recall that if an edge $e'$ with vertices $\displaystyle\binom{w'}{y'},\binom{x'}{z'}$ is mapped to 
an edge $e$ with vertices $\displaystyle\binom{w}{y},\binom{x}{z}$ by $A\in \Gamma_0(n)$, our only 
possibilities are that 
\[
A\begin{pmatrix}
    w' & x' \\
    y' & z'
\end{pmatrix} = \pm\begin{pmatrix}
    w & x \\
    y & z
\end{pmatrix} \hspace{20pt}\text{or}\hspace{20pt}A\begin{pmatrix}
    w' & x' \\
    y' & z'
\end{pmatrix} = \pm\begin{pmatrix}
    x & -w \\
    z & -y
\end{pmatrix}
\]
since those are the representations of the edge $e$ that correspond to matrices in 
$\Gamma_0(n)$ (here we are choosing a representative for $e'$ that is in $\Gamma_0(n)$).

Thus, if $A$ stabilizes $\displaystyle\binom{w}{y},\binom{x}{z}$, it must be given by 
\[
\pm\begin{pmatrix}
    w & x \\
    y & z
\end{pmatrix}\begin{pmatrix}
    w & x \\
    y & z
\end{pmatrix}^{-1} = \pm\begin{pmatrix}
     1 & 0 \\
     0 & 1
\end{pmatrix}
=\pm I_2
\]
or 
\[
\pm\begin{pmatrix}
    x & -w \\
    z & -y
\end{pmatrix}\begin{pmatrix}
    w & x \\
    y & z
\end{pmatrix}^{-1} = \pm\begin{pmatrix}
    x & -w \\
    z & -y
\end{pmatrix}\begin{pmatrix}
    z & -x \\
    -y & w
\end{pmatrix} = \pm\begin{pmatrix}
     wy+xz & -(x^2+w^2) \\
     y^2+z^2 & -(wy+xz)
\end{pmatrix}
\]
but matrices of the second kind, if in $\Gamma_0(n)$, would invert the corresponding edge, which cannot 
happen as Type 3 vertices were added to avoid that. Thus only $\pm I_2\simeq \z/2$ will stabilize edges.
\end{demonstracao}

\section{Homology calculations}

In the previous section, we constructed for all positive integers $n$ a tree on which $\Gamma_0(n)$ acts without inversion. This tree has $e(n)$ edges and $v(n)$ total vertices, with $v(n)-b(n)-c(n)$ vertices of Type 1, $b(n)$ vertices of Type 2 and $c(n)$ vertices of Type 3. The stabilizers of edges are isomorphic to $\z/2$ and the stabilizers of vertices are isomorphic to $\z/2,\z/6$ and $\z/4$, respectively. This action allows us to calculate the homology of $\Gamma_0(n)$ via the equivariant homology spectral sequence obtained in the last chapter.

\begin{teorema}\label{exact-seq}
Let a group $\Gamma$ act without inversion on a tree $Y$. Let $\Sigma_0$ (resp. $\Sigma_1$) denote 
a system of representatives of the vertices (resp. the edges) of $Y$, and for each vertex $x$ (resp. 
edge $y$) let $\Gamma_x$ (resp. $\Gamma_y$) be its stabilizer in $\Gamma$. For each $\Gamma$-module 
$M$, there exists an exact homology sequence
\[
\cdots \rightarrow H_{i+1}(\Gamma, M) \rightarrow \bigoplus_{y\in\Sigma_1}H_i(\Gamma_y,M)
\rightarrow \bigoplus_{x\in\Sigma_0}H_i(\Gamma_x,M) \rightarrow H_i(\Gamma,M) \arr \cdots
\]
\end{teorema}

\begin{demonstracao}
Since trees are contractible CW-complexes, we know from Corollary \ref{seq} that there is a spectral sequence $E^1_{p,q} = \displaystyle\bigoplus_{\sigma\in\Sigma_p}H_q(\Gamma_\sigma,M_\sigma) \Longrightarrow H_{p+q}(\Gamma,M)$. Since $\Gamma$ acts without inversion, \ref{teo: maps} shows that $M_\sigma\simeq M$ for all cells $\sigma$. Since there are only 0-cells and 1-cells, our first page will look like 

\begin{center}
        \begin{tikzcd}
            \vdots & \vdots & \vdots & \vdots & \\
            \displaystyle\bigoplus_{x\in\Sigma_0}H_2(\Gamma_x,M) & \arrow[l, "d^1_{1,2}"] \displaystyle\bigoplus_{x\in\Sigma_1}H_2(\Gamma_y,M) & \arrow[l] 0 & \arrow[l] 0 & \arrow[l]\hdots \\
            \displaystyle\bigoplus_{x\in\Sigma_0}H_1(\Gamma_x,M) & \arrow[l, "d^1_{1,1}"] \displaystyle\bigoplus_{x\in\Sigma_1}H_1(\Gamma_y,M) & \arrow[l] 0 & \arrow[l] 0 & \arrow[l]\hdots \\
            \displaystyle\bigoplus_{x\in\Sigma_0}H_0(\Gamma_x,M) & \arrow[l, "d^1_{1,0}"] \displaystyle\bigoplus_{x\in\Sigma_1}H_0(\Gamma_y,M) & \arrow[l] 0  & \arrow[l] 0 & \arrow[l]\hdots \\
        \end{tikzcd}
    \end{center}

    so that the spectral sequence collapses at the second page (that is, $E^\infty_{p,q}\simeq E^2_{p,q}$ for all $p,q$), which is given by 

    \begin{center}
        \begin{tikzcd}
            \vdots & \vdots & \vdots & \\
            \displaystyle\frac{\displaystyle\bigoplus_{x\in\Sigma_0}H_0(\Gamma_x,M)}{\text{im } d^1_{1,2}} &  \ker d^1_{1,2} & 0 & \hdots \\
            \displaystyle\frac{\displaystyle\bigoplus_{x\in\Sigma_0}H_0(\Gamma_x,M)}{\text{im } d^1_{1,1}} & \ker d^1_{1,1} & \arrow[llu] 0 & \hdots \\
            \displaystyle\frac{\displaystyle\bigoplus_{x\in\Sigma_0}H_0(\Gamma_x,M)}{\text{im } d^1_{1,0}} &   \ker d^1_{1,0} & \ \arrow[llu]0 &\hdots \\
        \end{tikzcd}
    \end{center}

    From the definition of convergence of spectral sequences (see Definition \ref{def: convergence}), there exists a filtration $\Phi H_*(\Gamma_0(n),M)$ such that $$E^\infty_{p,q}\simeq \frac{\Phi^pH_{p+q}(\Gamma_0(n),M)}{\Phi^{p-1}H_{p+q}(\Gamma_0(n),M)}.$$ From Proposition \ref{prp: bounds}, we also have that $\Phi$ has bounds $m$ and $-1$, that is, $\Phi^mH_m(\Gamma_0(n),M)=H_m(\Gamma_0(n,M))$ and $\Phi^{-1}H_m(\Gamma_0(n),M)=0$.

    Since $E^\infty_{p,q}=0$ for $q\geq 2$, we have that for all $m$, one can show by induction that $\Phi^{1}H_m(\Gamma_0(n),M) = H_m(\Gamma_0(n),M)$ which gives $$\frac{\bigoplus H_m(\Gamma_x,M)}{\im\,d^1_{1,m}} = E^2_{0,m}\simeq E^\infty_{0,m} \simeq \frac{\Phi^0H_m(\Gamma_0(n),M)}{\Phi^{-1}H_m(\Gamma_0(n),M)} = \Phi^0H_m(\Gamma_0(n),M)$$ and $$\ker d^1_{1,m-1} = E^2_{1,m-1}\simeq E^\infty_{1,m-1} \simeq \frac{\Phi^1H_m(\Gamma_0(n),M)}{\Phi^0H_m(\Gamma_0(n),M)} \simeq \frac{H_m(\Gamma_0(n),M)}{\Phi^0H_m(\Gamma_0(n),M)}.$$

    This gives short exact sequences $$0 \rightarrow \frac{\bigoplus H_m(\Gamma_x,M)}{\im\,d^1_{1,m}} \rightarrow H_m(\Gamma,M) \rightarrow \ker d^1_{1,m-1}\rightarrow 0.$$

    Let $\alpha',\beta'$ be the maps in the above sequence and let $\pi: \bigoplus H_m(\Gamma_x,M) \rightarrow \frac{\bigoplus H_m(\Gamma_x,M)}{\im\,d^1_{1,m}}$ be the natural projection, $i: \ker d^1_{1,m-1} \rightarrow H_{m-1}(\Gamma_y,M)$ the natural inclusion. Then $\beta = \pi\circ \beta'$ and $\alpha = \alpha'\circ i$ are such that $\im(\beta) = \im(\beta') = \ker(\alpha')=\ker(\alpha)$ giving exactness of the sequence $$ \bigoplus_{x\in\Sigma_0} H_m(\Gamma_x,M) \overset{\beta}{\longrightarrow} H_m(\Gamma,M) \overset{\alpha}{\longrightarrow} \bigoplus_{y\in\Sigma_1} H_{m-1}(\Gamma_y,M).$$ Furthermore, $\ker(\beta) = \ker(\pi) = \im(d^1_{1,m})$ and $\im(\alpha) = \im(i) = \ker(d^1_{1,m-1})$, so that the sequence $$\cdots \rightarrow H_{i+1}(\Gamma, M) \overset{\alpha}{\longrightarrow} \bigoplus_{y\in\Sigma_1}H_i(\Gamma_y,M) \overset{d^1}{\longrightarrow} \bigoplus_{x\in\Sigma_0}H_i(\Gamma_x,M) \overset{\beta}{\longrightarrow} H_i(\Gamma,M) \arr \cdots$$ is exact. Note that the central homomorphism is given precisely by the differential $d^1$.

    By Theorem \ref{teo: maps}, if $y\in\Sigma_1$ and $x_1,x_0 \in\Sigma_0$ are such that $\partial y = x_1-x_0$ and $[e] \in H_i(\Gamma_y,M)$, then $d^1[e] = i_{1,*}[e]-i_{0,*}[e]$ where $i_{1,*},i_{0,*}$ are the maps induced by the inclusions of $\Gamma_y$ in $\Gamma_{x_1}$ and $\Gamma_{x_0}$.
\end{demonstracao}

\begin{definicao}\label{def: r(n)}
    Given a positive integer $n$, let $r(n)\coloneq \displaystyle\frac{a(n)}{6}-\frac{2b(n)}{3}-\frac{c(n)}{2}+1$.
\end{definicao}

\begin{teorema}\label{H_n}
For any positive integer $n$, we have 
\[
H_k(\Gamma_0(n))\simeq \begin{cases}
\z & \text{if $k=0$}\\
\z^{r(n)}\oplus G & \text{if $k=1$}\\
(\z/2)^{r(n)} & \text{if $k > 1$ is even,}\\
G & \text{if $k>1$ is odd}
\end{cases}
\]
where
\[
G = \begin{cases}
(\z/2)^{c(n)-1}\oplus (\z/3)^{b(n)}\oplus \z/4 & \text{if $c(n)>0$} \\
\z/2\oplus (\z/3)^{b(n)} & \text{if $c(n)=0$}.
\end{cases}
\]
\end{teorema}
\begin{demonstracao}

Let $Y$ be the tree constructed on the previous section. We have seen that $\Gamma_0(n)$ acts without inversion on this tree. 

Denote by $d(n)$ the number of $\Gamma_0(n)$-equivalence classes of Type 1 vertices, that is, $d(n) \coloneq v(n)-b(n)-c(n)$.

Applying our results on stabilizers on the previous theorem (see Corollaries \ref{e(n)}, \ref{cor: type 2}, \ref{type3} and \ref{v(n)} for the cardinality of the equivalence classes of edges and each type of vertex and Proposition \ref{stab} for the respective stabilizers) and using the description given in Example \ref{exemplo1} of the homologies of finite cyclic groups, we obtain exact sequences 
\begin{center}
\begin{tikzcd}[column sep=2ex]
0 \arrow[r] & H_2(\Gamma_0(n)) \arrow[r] & (\z/2)^{e(n)} \arrow[r, "i"] \arrow[d, phantom, ""{coordinate, name=Z}] & 
(\z/2)^{d(n)}\oplus(\z/6)^{b(n)}\oplus(\z/4)^{c(n)} 
\arrow[dlll, rounded corners, to path={ -- ([xshift=2ex]\tikztostart.east)|- 
(Z) [near end]\tikztonodes-| ([xshift=-2ex]\tikztotarget.west)-- (\tikztotarget)}] \\ 
H_1(\Gamma_0(n)) \arrow[r] & \z^{e(n)} \arrow[r, "i_0"] & \z^{v(n)} \arrow[r] & \z \arrow[r] & 0 &
\end{tikzcd}
\end{center}
and for $k\geq 2$,
\begin{center}
\begin{tikzcd}
0 \arrow[r] & H_{2k}(\Gamma_0(n)) \arrow[r] & (\z/2)^{e(n)} \arrow[r, "i"] 
\arrow[d, phantom, ""{coordinate, name=Z}] & (\z/2)^{d(n)}\oplus(\z/6)^{b(n)}\oplus(\z/4)^{c(n)} 
\arrow[dl, rounded corners, to path={ -- ([xshift=2ex]\tikztostart.east)|- 
(Z) [near end]\tikztonodes-| ([xshift=-2ex]\tikztotarget.west)-- (\tikztotarget)}]  \\ & &   
H_{2k-1}(\Gamma_0(n)) \arrow[r] & 0
\end{tikzcd}
\end{center}

From the first exact sequence, we obtain 
\[
0 \rightarrow \coker(i) \rightarrow H_1(\Gamma_0(n)) \rightarrow \ker(i_0) \rightarrow 0
\]
\[
H_2(\Gamma_0(n)) = \ker(i)
\]
From the second exact sequence, we obtain for $k\geq 2$, 
\[
H_{2k-1}(\Gamma_0(n)) = \coker(i)
\]
\[
H_{2k}(\Gamma_0(n)) = \ker(i)
\]


Theorem \ref{exact-seq} gives the following description of the map $i$: if $y\in\Sigma_1$ and $x_1,x_0\in\Sigma_0$ are such that $\partial y =x_1-x_0$ and $[e]\in H_i(\Gamma_y)$, then $i[e] = i_{1,*}[e]-i_{0,*}[e]$ where $i_{1,*},i_{0,*}$ are induced by the inclusions of $\Gamma_y$ in $\Gamma_{x_1}$ and $\Gamma_{x_0}$. 

Since each of the stabilizer subgroups $\Gamma_y,\Gamma_x$ is isomorphic to a finite cyclic group, denote the generator of the copy of $\Gamma_\sigma$ inside $\displaystyle\bigoplus_{s\in\Sigma_i}\Gamma_s$ by its cell $\sigma$.

For dimension 0, the induced maps $H_0(\Gamma_y)\rightarrow H_0(\Gamma_x)$ are all $\id_{\z}$ (see Example \ref{exe: induced}), so the map $i: \z^{e(n)} \rightarrow \z^{v(n)}$ is such that the generator $y$ is mapped to $x_1-x_0$, that is, $i$ coincides with the boundary map of the quotient graph $Y/\Gamma_0(n)$ and $\ker(i_0)\simeq \z^{e(n)-(v(n)-1)} = \z^{r(n)}$ by Theorem \ref{teo: homology graph}. In particular, $\ker(i_0)$ is a free $\z$-module and the short exact sequence for $H_1(\Gamma_0(n))$ splits giving  $H_1(\Gamma_0(n)) = \z^{r(n)}\oplus \,\coker(i)$ (see Proposition 2.28 of \cite{rotman2009}). It remains to show that $\ker(i) = (\z/2)^{r(n)}$ 
and $\coker(i) = G$.

For even non-zero dimensions, the inclusions $i_{1,*},i_{0,*}$ are induced by the inclusions of $\z/2$ in $\z/2,\z/6$ or $\z/4$, depending on the type of the vertices $x_1,x_0$, so the map $i$ becomes $y\mapsto a_1x_1-a_0x_0$ where $a_1,a_0$ are either $1,3$ or $2$, the generators of the subgroup $\z/2$ contained in $\z/2,\z/6$ and $\z/4$. Let $j$ be the inclusion of $(\z/2)^{v(n)}$ in $(\z/2)^{d(n)}\oplus (\z/6)^{b(n)}\oplus (\z/4)^{c(n)}$ and $d$ be the map defined by $y\mapsto x_1+x_0$. Then the diagram below commutes: \begin{center}
    \begin{tikzcd}
        (\z/2)^{e(n)} \arrow[rr, "d"] \arrow[dr, "i"] & & (\z/2)^{v(n)} \arrow[ld, "j"]\\
        & (\z/2)^{d(n)}\oplus (\z/6)^{b(n)}\oplus (\z/4)^{c(n)}  &
    \end{tikzcd}.
\end{center}

Note that $d$ is the analogue of the boundary map of $Y/\Gamma_0(n)$ for coefficients in $\z/2$. This implies that the dimension of the kernel of $i$ as a $(\z/2)$-vector space corresponds to the number of non-homologous cycles in the quotient graph $Y\backslash\Gamma_0(n)$, that is, the rank of $H_1(Y\backslash\Gamma_0(n),\z)$.  Thus, $$\ker(i) \simeq (\z/2)^{e(n)-(v(n)-1)}\simeq(\z/2)^{r(n)}.$$

Now, from the fact that $(\z/2)^{e(n)}/\ker(i)\simeq \im(i)$, we obtain that $$\im(i)\simeq \frac{(\z/2)^{e(n)}}{(\z/2)^{e(n)-(v(n)-1)}} \simeq (\z/2)^{v(n)-1}.$$

If $c(n)=0$, then $d(n)=v(n)-b(n)$ and we obtain the desired result that \begin{align*}
    \coker(i) &\simeq \frac{(\z/2)^{d(n)}\oplus (\z/6)^{b(n)}}{(\z/2)^{v(n)-1}} \\ &\simeq \frac{(\z/2)^{v(n)}\oplus (\z/3)^{b(n)}}{(\z/2)^{v(n)-1}}\\ &\simeq \z/2\oplus(\z/3)^{b(n)}.
\end{align*}

For the case $c(n)>0$, we note that $$(\z/2)^{d(n)}\oplus(\z/6)^{b(n)}\oplus(\z/4)^{c(n)} \simeq (\z/2)^{v(n)-c(n)}\oplus(\z/3)^{b(n)}\oplus(\z/4)^{c(n)}$$ and we identify $i$ with the corresponding map under this isomorphism.

Proposition \ref{prp: numbers} at the end of this chapter shows that 1 and 2 are the only integers $n$ such that $c(n)=1$. We calculate the cokernel of $i$ directly for these cases in Examples \ref{graph 1} and \ref{graph 2} below. Suppose, then, that $c(n)>1$ so there are at least two Type 3 vertices.

Let $A = (\z/2)^{v(n)-c(n)}\oplus(0)\oplus(\z/4)^{c(n)}$ be a submodule of the codomain of $i$ and note that $A$ contains $\im(i)$ since 3 in $\z/6$ corresponds to $(1,0)\in(\z/2)\oplus(\z/3)$. Let $B=\im(i)\cap (2A) \simeq \im(i)\cap (\z/2)^{c(n)}$ and let $s(n) = \dim_{\mathbb{F}_2}B$. In other words, $s(n)$ is the amount of copies of $\z/2$ that sit inside the $(\z/4)^{c(n)}$-component of $\im(i)$. Then, $$\im(i) \simeq (\z/2)^{v(n)-1} \simeq (\z/2)^{v(n)-1-s(n)}\oplus(\z/2)^{s(n)}$$ so that \begin{align*}
    \coker(i) &\simeq \frac{(\z/2)^{v(n)-c(n)}\oplus (\z/3)^{b(n)}\oplus(\z/4)^{c(n)}}{(\z/2)^{v(n)-1-s(n)}\oplus(\z/2)^{s(n)}} \\ &\simeq (\z/2)^{2s(n)-c(n)+1}\oplus(\z/3)^{b(n)}\oplus(\z/4)^{c(n)-s(n)}.
\end{align*}

As vector spaces over $\mathbb{F}_2$, we have $$s(n) = \dim_{\mathbb{F}_2}(\im(i)\cap(2A)) = \dim_{\mathbb{F}_2}(\im(i))+\dim_{\mathbb{F}_2}(2A)-\dim_{\mathbb{F}_2}(\im(i)+2A).$$ We claim that $\dim_{\mathbb{F}_2}(\im(i)+2A) = v(n)$, so that $s(n)=c(n)-1$.

Indeed, let $x$ be a vertex in the quotient graph $Y\backslash\Gamma_0(n)$ and let $x'$ be any vertex of Type 3 different from $x$. Since our graph is connected, there is a sequence $e_1,e_2,\dots,e_k$ of consecutive edges such that $e_1$ has origin $x$ and $e_k$ has terminus $x'$. Then $i(e_1+e_2+\dots+e_k) = ax+2x'$ with $a=1$ or $a=2$ depending on the type of the vertex $x$. 

So $ax+2x'$ is in $\im(i)$ and $2x'\in 2A$ since $x'$ represents the generator of the copy of $\z/4$ indexed by $x'$. This implies that $ax=(ax+2x')+2x' \in \im(i)+2A$. Thus there are at least $v(n)$ linearly independent vectors in $\im(i)+2A$. This implies that $\dim_{\mathbb{F}_2}(\im(i)+2A)\geq v(n)$.

On the other hand, neither $\im(i)$ nor $2A$ has vectors with entries 1 or 3 in the coordinates corresponding to $(\z/4)$, implying that the sum $\im(i)+2A$ is contained in the set $(\z/2)^{v(n)-c(n)}\oplus 2(\z/4)^{c(n)} \simeq (\z/2)^{c(n)}$. Thus, $\dim_{\mathbb{F}_2}(\im(i)+2A) \leq v$.

Plugging in the obtained value for $s(n)$, we get that for $c(n)>0$, \begin{align*}
    \coker(i)&\simeq (\z/2)^{2s(n)-c(n)-1}\oplus(\z/3)^{b(n)}\oplus(\z/4)^{c(n)-s(n)}\\ & = (\z/2)^{c(n)+1}\oplus(\z/3)^{b(n)}\oplus(\z/4).
\end{align*}

\end{demonstracao}

Given a set $R$ of coset representatives for $\SL_2(\z)/\Gamma_0(n)$, Proposition \ref{e(n)} and Lemma \ref{T-S} encode the means to construct the quotient graph $Y/\Gamma_0(n)$, since each equivalence class of an edge can be represented by a single coset (if this is an added edge to prevent inversion) or by a pair of cosets.

For the case where $n$ is a prime integer, the matrices $(R^{(k)})_1 = \begin{pmatrix}
    1 & 0 \\ k & 1
\end{pmatrix}$ with $0\leq k < n$, along with the matrix $S = \begin{pmatrix}
    0 & -1 \\ 1 & 0
\end{pmatrix}$ form a set of coset representatives for $\SL_2(\z)/\Gamma_0(n)$. Proposition \ref{e(n)} gives that $S$ and $I = (R^{(0)})_1$ represent equivalent edges and $(R^{(k)})_1$ represents the same edge as $(R^{(r)})_1$ where $1+kr\equiv 0 \pmod{n}$.

The edge represented by $(R^{(k)})_1$ is the edge given by $R^{(k)}$ and $R^{(k-1)}$ (or the extra edge created to avoid inversion of this edge if $k^2+1\equiv 0 \pmod{n}$) and Lemma \ref{T-S} gives the conditions to decide when $R^{(k)}$ and $R^{(r)}$ are equivalent, allowing us to reduce the vertices to a set of representatives and reconstruct the quotient graph.

Appendix \ref{chapter:codigo} contains a python code that computes these equivalences and produces a list of pairs $[i,j]$ describing the distinct edges from $R^{(i)}$ to $R^{(j)}$ or $[i,E_j]$ to describe an added edge at a Type 3 vertex by a root $j$ of $X^2+1\pmod{n}$. This code can also be accessed by the link \url{https://github.com/isavpicinini/codes/blob/main/graph.py}.

Using this code for $n=19$ and $n=5$ outputs, respectively, the lists $$[[0, 0], [0, 1], [1, 2], [2, 3], [3, 4], [4, 3], [2, 7], [1, 10], [10, 11], [10, 4]],$$ $$[[0, 0], [0, 1], [1, E_2], [1, E_3]].$$ Examples \ref{graph19} and \ref{graph5} below arranges these in graphs.

\begin{exemplo}\label{graph 1}
    Since $\Gamma_0(1)=\SL_2(\z)$, we get that the action of $\Gamma_0(1)$ on the original tree is transitive. Thus, every triangle is equivalent and so in particular every triangle is of Type 2 and every edge is inverted. 
    
    Indeed, $a(1)=b(1)=c(1)=1$ giving $v(1)=2$ and $e(1)=1$. So $Y/\Gamma_0(1)$ has one vertex of Type 2, one vertex of Type 3 and a single edge connecting the vertices.

    The map $i: \z/2 \rightarrow (\z/6)\oplus(\z/4)$ is then given by $1\mapsto (3,2)$. This map is injective and has cokernel isomorphic to $\z/12$, proving the last theorem for the case $n=1$.

    The reader familiar with amalgamated products will recognize that this implies $\SL_2(\z) = \Gamma_0(1)$ is isomorphic to the amalgamated product $\z/4*_{\z/2}\z/6$. For more information on amalgamated products, we reccomend Section II.7 from \cite{brown1994} and the appendix following it.
\end{exemplo}

\begin{exemplo}\label{graph 2}
    For $n=2$, we get $v(2)=e(2)=2$. The graph $Y/\Gamma_0(2)$ then has one vertex of Type 2 denoted by $R^{(0)}$, one added vertex of Type 3 denoted by $E^{(1)}$ (here 1 is the unique root of $X^2+1\pmod{2}$) and one terminal edge connecting both vertices. The remaining edge must then be a loop at the Type 2 vertex. Figure \ref{fig:n2} below illustrates this graph.

    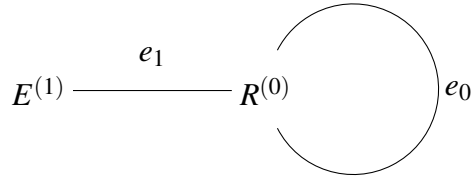
\begin{figure}[H]
        \centering
        \caption{Quotient graph of $Y/\Gamma_0(2)$}
        \begin{tikzpicture}[scale=1.5]
        \node[scale=1] at (-1,0) {$E^{(1)}$};
        \node[scale=1] at (1,0) {$R^{(0)}$};
        
        \draw (-0.7,0)--(0.7,0);
        \draw (1.102,-0.332) arc [x radius=0.752, y radius=0.752, start angle=207.2, end angle=512.8] ;

        \node[scale=1] at (0,0.3) {$e_1$};
        \node[scale=1] at (2.7,0) {$e_0$};
    \end{tikzpicture}
        \label{fig:n2}
        \fautor
    \end{figure}

    The map $i: (\z/2)^2\rightarrow (\z/6)\oplus(\z/4)$ is then given by $e_0\mapsto 0$ and $e_1\mapsto (3,2)$. This map has kernel generated by $e_0$, the only cycle in the graph, and its cokernel is isomorphic to $\z/12$, proving the last theorem for the case $n=2$.
\end{exemplo}

\begin{exemplo}\label{graph5}
    Consider the case $n=5$. We have 4 equivalence classes of vertices, 2 of which are of Type 1 and 2 of which are of Type 3, and 4 equivalence classes of edges. Choosing an $R^{(k)}$ representative for the Type 1 vertices and denoting the added Type 3 vertices by $E^{(k)}$ for the integers $k$ such that $k^2+1\equiv 0 \pmod{5}$, the quotient graph is as shown in Figure \ref{fig:n5}.

    \begin{figure}[H]
        \centering
        \caption{Quotient graph of $Y/\Gamma_0(5)$}
        \begin{tikzpicture}[scale=1.5]
        \node[scale=1] at (-1,0) {$R^{(1)}$};
        \node[scale=1] at (-3,-1) {$E^{(3)}$};
        \node[scale=1] at (-3,1) {$E^{(2)}$};
        \node[scale=1] at (1,0) {$R^{(0)}$};
        
        \draw (-0.7,0)--(0.7,0);
        \draw (-1.3,0)--(-2.8,0.8);
        \draw (-1.3,0)--(-2.8,-0.8);
        \draw (1.102,-0.332) arc [x radius=0.752, y radius=0.752, start angle=207.2, end angle=512.8] ;

        \node[scale=1] at (0,0.3) {$e_1$};
        \node[scale=1] at (-2,0.7) {$e_2$};
        \node[scale=1] at (-2,-0.7) {$e_3$};
        \node[scale=1] at (2.7,0) {$e_0$};
    \end{tikzpicture}
        \label{fig:n5}
        \fautor
    \end{figure}
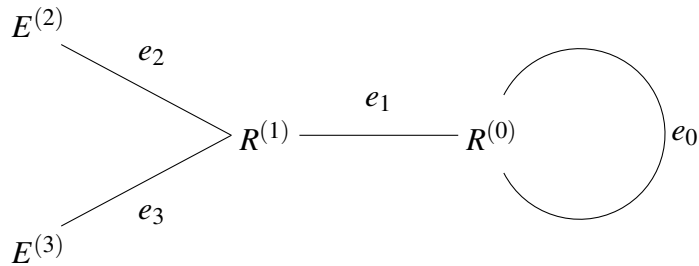

    Our map $i: (\z/2)^4 \rightarrow (\z/2)^2\oplus(\z/4)^2$ is given by \begin{align*}
        e_0 &\longmapsto 0 \\
        e_1 &\longmapsto (1,1,0,0) \\
        e_2 &\longmapsto (0,1,2,0) \\
        e_3 &\longmapsto (0,1,0,2)
    \end{align*} 

    The only possible choice of generators for $\ker(i)$ is given by $e_0$, representing the only cycle in the graph.

    Letting $A=(\z/2)^2\oplus(\z/4)^2$ as in the proof of the last theorem, we see that $2A$ is generated by $(0,0,2,0)$ and $(0,0,0,2)$. Fix the Type 3 vertex $E^{(3)}$. Then each vertex in the tree can be reached from $E^{(3)}$: $R^{(0)}$ can be reached by the path $e_1+e_3$, $R^{(1)}$ can be reached by the path $e_3$, and $E^{(2)}$ can be reached by $e_2+e_3$. Thus, \begin{align*}
        (1,0,0,0) &= (1,1,0,0)+(0,1,0,2)+(0,0,0,2) = i(e_1+e_2)+2(0,0,0,1), \\
        (0,1,0,0) &= (0,1,0,2)+(0,0,0,2) = i(e_2)+2(0,0,0,1), \\
        (0,0,2,0) &= (0,1,2,0)+(0,1,0,2)+(0,0,0,2) = i(e_2+e_3)+2(0,0,0,1),\\
        (0,0,0,2) &= (0,1,2,0)+(0,1,0,2)+(0,0,2,0) = i(e_2+e_3)+2(0,0,1,0),
    \end{align*}  and we conclude that $\im(i)+2A$ has four linearly independent vectors, one for each vertex as in the proof. Note that the element $(0,0,0,2)$ corresponding to the vertex $E^{(3)}$ we fixed was obtained by connecting $E^{(3)}$ to the other Type 3 vertex $E^{(2)}$ by the path $e_2+e_3$. This is where the condition $c(n)>1$ becomes important.
\end{exemplo}

\begin{exemplo}\label{graph19}
    Consider the case $n=19$. We have 8 equivalence classes of vertices and 10 equivalence classes of edges. Choosing an $R^{(k)}$ representative for each vertex, we get that the quotient graph is as shown in Figure \ref{fig:n19}.

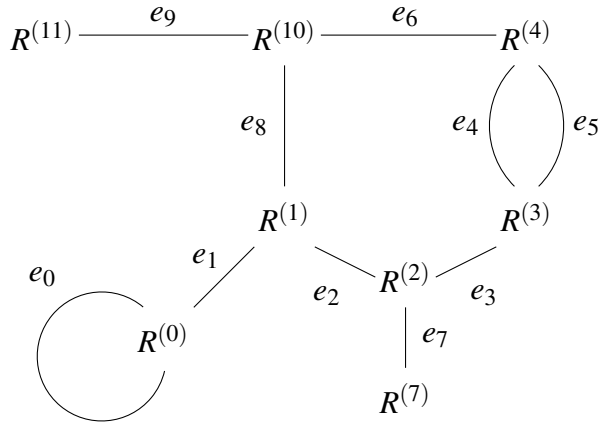
\begin{figure}[H]
    \centering
    \caption{Quotient graph $Y/\Gamma_0(19)$}
    \begin{tikzpicture}[scale=0.8]
        \node[scale=1] at (0,0) {$R^{(1)}$};
        \node[scale=1] at (0,3) {$R^{(10)}$};
        \node[scale=1] at (4,0) {$R^{(3)}$};
        \node[scale=1] at (4,3) {$R^{(4)}$};
        \node[scale=1] at (2,-1) {$R^{(2)}$};
        \node[scale=1] at (-4,3) {$R^{(11)}$};
        \node[scale=1] at (-2,-2) {$R^{(0)}$};
        \node[scale=1] at (2,-3) {$R^{(7)}$};
        
        \draw (0,0.5)--(0,2.5);
        \draw (0.6,3)--(3.5,3);
        \draw (-0.6,3)--(-3.4,3);
        \draw (-0.5,-0.5)--(-1.5,-1.5);
        \draw (0.5,-0.5)--(1.5,-1);
        \draw (3.5,-0.5)--(2.5,-1);
        \draw (2,-1.5)--(2,-2.5);
        \draw (3.8,0.5) edge [bend left=45] (3.8,2.5);
        \draw (4.2,0.5) edge [bend right=45] (4.2,2.5);
        \draw (-2.337,-1.495) arc [x radius=1.066, y radius=1.066, start angle=50, end angle=347] ;

        \node[scale=1] at (-1.3,-0.7) {$e_1$};
        \node[scale=1] at (0.7,-1.3) {$e_2$};
        \node[scale=1] at (3.3,-1.3) {$e_3$};
        \node[scale=1] at (3,1.5) {$e_4$};
        \node[scale=1] at (5,1.5) {$e_5$};
        \node[scale=1] at (-0.5,1.5) {$e_8$};
        \node[scale=1] at (2,3.3) {$e_6$};
        \node[scale=1] at (-2,3.3) {$e_9$};
        \node[scale=1] at (2.5,-2) {$e_7$};
        \node[scale=1] at (-4,-1) {$e_0$};
    \end{tikzpicture}
    \label{fig:n19}
    \fautor
\end{figure}

Here, the vertices $R^{(7)}$ and $R^{(11)}$ are of Type 2 and the rest of the vertices are of Type 1.

The edge $e_7$, for example, has the vertices $R^{(2)}$ and $R^{(7)}$ as its boundaries, so the induced map on homologies takes the vector $(0,0,0,0,0,0,1,0,0,0)$, which we identify with $e_7$, to $(0,0,1,0,0,0,3,0)$, since this map is induced by the inclusions of $\z/2$ in $\z/2$ and $\z/6$, which maps 1 to 1 and 3, respectively. Here we are letting the first 6 coordinates of $(\z/2)^6\oplus(\z/6)^2$ correspond to the first 6 Type-1 vertices in order and the last two coordinates correspond to the vertices $R^{(7)}$ and $R^{(11)}$.

Thus, our map $i: (\z/2)^{10} \rightarrow (\z/2)^{6}\oplus(\z/6)^{2}$ is given by \begin{align*}
    e_0 &\longmapsto 0 \\
    e_1 &\longmapsto (1,1,0,0,0,0,0,0) \\
    e_2 &\longmapsto (0,1,1,0,0,0,0,0) \\
    e_3 &\longmapsto (0,0,1,1,0,0,0,0) \\
    e_4 &\longmapsto (0,0,0,1,1,0,0,0) \\
    e_5 &\longmapsto (0,0,0,1,1,0,0,0) \\
    e_6 &\longmapsto (0,0,0,0,1,1,0,0) \\
    e_7 &\longmapsto (0,0,1,0,0,0,3,0) \\
    e_8 &\longmapsto (0,1,0,0,0,1,0,0) \\
    e_9 &\longmapsto (0,0,0,0,0,1,0,3) \\
\end{align*}

Notice that a choice of generators for the kernel of $i$ is given by $e_0$, $e_4+e_5$ and $e_2+e_3+e_4+e_6+e_8$, which are the three non-homologous cycles in our quotient graph.
\end{exemplo}


\begin{teorema}\label{H_nP}
Let $\PGamma_0(n)$ be the subgroup of $\PSL_2(\z) = \SL_2(\z)/\{\pm I\}$ of matrices whose lower left element is congruent to 0 mod $n$, that is, $\PGamma_0(n) = \PGamma_0(n)/\{\pm I\}$. For any positive integer $n$, we have 
\[
H_k(\PGamma_0(n))\simeq \begin{cases}
\z & \text{if $k=0$}\\
\z^{r(n)}\oplus (\z/3)^{b(n)}\oplus (\z/2)^{c(n)} & \text{if $k=1$}\\
(\z/3)^{b(n)}\oplus (\z/2)^{c(n)} & \text{if $k > 1$ is odd.}\\
0 & \text{otherwise}
\end{cases}
\]
\end{teorema}

\begin{demonstracao}
    Denote by $\overline{A}$ the class in $\PSL_2(\z)$ of a matrix $A\in\SL_2(\z)$.
    
    Note that since vertices of triangles are identified with their antipodes, $T$ and $-T$ will coincide as triangles in $\Delta$. This implies that the action of $\SL_2(\z)$ on the triangles $\Delta$ induces a well-defined action of $\PSL_2(\z)$ on $\Delta$ by $\overline{A}T = AT$.

    Note that two triangles $T$ and $S$ are $\Gamma_0(n)$-equivalent if and only if they are $\PGamma_0(n)$-equivalent, and a matrix $A\in\Gamma_0(n)$ stabilizes a triangle $T$ if and only if its class $\overline{A}\in \PGamma_0(n)$ stabilizes $T$. Thus, the action of $\PGamma_0(n)$ on $\Delta$ has the same set of orbits as the action of $\Gamma_0(n)$ and the stabilizers for $\PGamma_0(n)$ correspond to the classes of the stabilizers under $\Gamma_0(n).$ 

    For the stabilizers, recall that under the action of $\Gamma_0(n)$, the stabilizers of the vertices of Type 1, Type 2 and Type 3 are isomorphic to $\z/2,\z/6$ and $\z/4$, respectively, and only $\pm I$ will stabilize an edge. The stabilizers under the action of $\PGamma_0(n)$ will then be $\{0\},\z/3$ and $\z/2$ for triangles of Types 1, 2 and 3 respectively and trivial for edges.

    Applying Theorem \ref{exact-seq}, we obtain the exact sequences 
    \begin{center}
    \begin{tikzcd}[column sep=2ex]
    0 \arrow[r] & H_2(\PGamma_0(n)) \arrow[r] & 0 \arrow[r, "i"] \arrow[d, phantom, ""{coordinate, name=Z}] & (\z/3)^{b(n)}\oplus(\z/2)^{c(n)} 
    \arrow[dlll, rounded corners, to path={ -- ([xshift=2ex]\tikztostart.east)|- 
    (Z) [near end]\tikztonodes-| ([xshift=-2ex]\tikztotarget.west)-- (\tikztotarget)}] \\ 
    H_1(\PGamma_0(n)) \arrow[r] & \z^{e(n)} \arrow[r, "i_0"] & \z^{v(n)} \arrow[r] & \z \arrow[r] & 0 &
    \end{tikzcd}
    \end{center}
    and for $k\geq 2$, $$0 \longrightarrow H_{2k}(\PGamma_0(n)) \longrightarrow 0 \overset{i}{\longrightarrow}     (\z/3)^{b(n)}\oplus(\z/2)^{c(n)} \overset{\partial}{\longrightarrow}
    H_{2k-1}(\PGamma_0(n)) \longrightarrow 0$$

    We immediately obtain that for $k\geq 1$, the connecting homomorphism $\partial$ is an isomorphism, so $H_{2k+1}(\PGamma_0(n))\simeq (\z/3)^{b(n)}\oplus(\z/2)^{c(n)}$ and the even homologies $H_{2k}(\PGamma_0(n))$ are trivial for all $k$. It remains to check the case for $H_1(\PGamma_0(n))$. From the first sequence, we obtain $$0\rightarrow \coker(i) \rightarrow H_1(P\Gamma_0(n))\rightarrow \ker(i_0) \rightarrow 0$$ where $\coker(i) = (\z/3)^{b(n)}\oplus(\z/2)^{c(n)}$ since the domain of $i$ is 0. The calculation of $\ker(i_0)$ is exactly the same as in the case for $\Gamma_0(n)$, giving $\ker(i_0)\simeq \z^{r(n)}$, which is a free $\z$-module, so the exact sequence splits and $H_1(\PGamma_0(n))= \z^{r(n)}\oplus(\z/3)^{b(n)}\oplus(\z/2)^{c(n)}$, as expected.
\end{demonstracao}

For the remainder of this section, we wish to characterize the integers $a(n),b(n),c(n)$.

\begin{definicao}
    Given $n$ an integer, let $$\Gamma(n) = \left\{\begin{pmatrix}
        a & b \\ c & d
    \end{pmatrix}\in\SL_2(\z) \mid a\equiv d\equiv 1, b\equiv c \equiv 0 \pmod{n}\right\},$$
    $$\Gamma_1(n) = \left\{\begin{pmatrix}
        a & b \\ c & d
    \end{pmatrix}\in\SL_2(\z) \mid a\equiv d\equiv 1, c \equiv 0 \pmod{n}\right\}.$$
\end{definicao}

Note that $\Gamma(n)\subset \Gamma_1(n)\subset \Gamma_0(n) \subset \SL_2(\z)$.

\begin{lema}
    The index $[\Gamma_0(n):\Gamma(n)]$ is equal to $n\cdot\varphi(n)$ where $\varphi$ is the Euler totient function given by $\varphi(n) = n\displaystyle\prod_{p\mid n}\left(1-\frac{1}{p}\right)$.
\end{lema}

\begin{demonstracao}
    Note that the group homomorphism $\Gamma_1(n)\rightarrow \z/n$ given by $\begin{pmatrix}
        a & b \\ c & d
    \end{pmatrix} \mapsto \overline{b}$ is surjective and has kernel $\Gamma(n)$, so that $\Gamma_1(n)/\Gamma(n)\simeq \z/n$ and $[\Gamma_1(n):\Gamma(n)] = n$.

    Now, consider the group homomorphism $\Gamma_0(n)\rightarrow (\z/n)^\times$ given by $\begin{pmatrix}
        a & b \\ c & d
    \end{pmatrix}\mapsto \overline{d}$. Since $ad-bc=1$ and $c\equiv 0 \pmod{n}$, we have that $ad\equiv 1 \pmod{n}$, so $d$ is invertible and the map makes sense. This map is also surjective and has kernel $\Gamma_1(n)$, so that $[\Gamma_0(n):\Gamma_1(n)] = |(\z/n)^\times| = \varphi(n)$. Then, $$[\Gamma_0(n):\Gamma(n)] = [\Gamma_0(n):\Gamma_1(n)][\Gamma_1(n):\Gamma(n)]=n\cdot\varphi(n).$$
\end{demonstracao}

\begin{lema}
    If $A$ and $B$ are commutative rings, then $\SL_2(A\times B)\simeq \SL_2(A)\times\SL_2(B)$.
\end{lema}

\begin{demonstracao}
    The map $\SL_2(A)\times\SL_2(B)\rightarrow\SL_2(A\times B)$ given by $$\left(\begin{pmatrix}
        a_{11} & a_{12} \\ a_{21} & a_{22}
    \end{pmatrix},\begin{pmatrix}
        b_{11} & b_{12} \\ b_{21} & b_{22}
    \end{pmatrix}\right) \mapsto \begin{pmatrix}
        (a_{11},b_{11}) & (a_{12},b_{12}) \\ (a_{21},b_{21}) & (a_{22},b_{22})
    \end{pmatrix}$$ is an isomorphism, as is easily checked.
\end{demonstracao}

\begin{lema}
    The order of the group $|\SL_2(\z/n)|$ is equal to $n^3\displaystyle\prod_{p\mid n}\left(1-\frac{1}{p^2}\right)$.
\end{lema}

\begin{demonstracao}
    We first consider the case where $n$ has a single prime factor, that is, $n=p^k$ for some $k$.

    Note that if a matrix $\begin{pmatrix}
        a & b \\ c & d
    \end{pmatrix} \in \SL_2(\z/p^k)$ is such that $a\in(\z/p^k)^\times$, then $(1+bc)a^{-1}$ uniquely determines $d$ for any choices of $b,c\in\z/p^k$. Thus, each $a\in(\z/p^k)^\times$ contributes with $p^{2k}$ matrices.

    On the other hand, if $a\notin(\z/p^k)^\times$, then $c\in(\z/p^k)^\times$, since if both $a$ and $c$ are non-invertible, then $p\mid a$ and $p\mid c$, so $p\mid ad-bc=1$, which is absurd. Then, $b$ is uniquely determined by $(ad-1)c^{-1}$ for any choice of $d$, thus each $a\notin(\z/p^k)^\times$ contributes with $p^k\cdot\varphi(p^k)$ based on the choices of $c$ and $d$. 

    We have, then, that each of the $\varphi(p^k)$ choices of invertible $a$ contributes with $p^{2k}$ matrices and each of the $p^k-\varphi(p^k)$ choices of non-invertible $a$ contributes with $p^k\cdot\varphi(p^k)$ matrices, so $$|\SL_2(\z/p^k)| = \varphi(p^k)\cdot p^{2k}+(p^k-\varphi(p^k))\cdot p^k\cdot\varphi(p^k) = p^{3k}\left(1-\frac{1}{p^2}\right).$$ 

    For the general case, if $n=p_1^{k_1}\dots p_r^{k_r}$ is the prime factorization for $n$, then $\z/n \simeq \displaystyle\prod_{i=1}^r\z/p_i^{k_i}\z$, so that $\SL_2(\z/n) \simeq \displaystyle\prod_{i=1}^r\SL_2(\z/p_1^{k_i})$ and $$|\SL_2(\z/n)| = \prod_{i=1}^k p_i^{3k_i}\left(1-\frac{1}{p_i^2}\right) = n^3\prod_{p\mid n}\left(1-\frac{1}{p^2}\right).$$
\end{demonstracao}

\begin{teorema}
    \label{teo: GE2}
    If $R$ is a semi-local ring or an euclidian domain, then $\SL_2(R)$ is generated by elementary matrices.
\end{teorema}

\begin{demonstracao}
    See page 245 of \cite{Silvester1982} for the case of semi-local rings and Section 2 of \cite{cohn1966} for the case of euclidian domains. 
\end{demonstracao}

This result is particularly useful for the case of projections $\z \rightarrow \z/n$ since $\z/n$ is semi-local and an elementary matrix in $\SL_2(\z/n)$ is either $\begin{pmatrix}
    \overline{1} & \overline{a} \\ \overline{0} & \overline{1}
\end{pmatrix}$ or $\begin{pmatrix}
    \overline{1} & \overline{0} \\ \overline{b} & \overline{1}
\end{pmatrix}$. For either case, the matrices $\begin{pmatrix}
    1 & a \\ 0 & 1
\end{pmatrix},\begin{pmatrix}
    1 & b \\ 0 & 1
\end{pmatrix}$ are always in $\SL_2(\z)$, implying that the induced map $\SL_2(\z)\rightarrow \SL_2(\z/n)$ is surjective.

\begin{lema}\label{numb a(n)}
    The index $[\SL_2(\z):\Gamma_0(n)]$ is equal to $n\displaystyle\prod_{p\mid n}\left(1+\frac{1}{p}\right)$.
\end{lema}

\begin{demonstracao}
    Note that the map $\SL_2(\z)\rightarrow \SL_2(\z/n)$ given by projecting each coordinate to $\z/n$ is a surjective map with kernel equal to $\Gamma(n)$, so that $$[\SL_2(\z):\Gamma(n)] = |\SL_2(\z/n)| = n^3\prod_{p\mid n}\left(1-\frac{1}{p^2}\right).$$ Then $[\SL_2(\z):\Gamma(n)] = [\SL_2(\z):\Gamma_0(n)][\Gamma_0(n):\Gamma(n)]$, from which we conclude that $$[\SL_2(\z):\Gamma_0(n)] = \frac{[\SL_2(\z):\Gamma(n)]}{[\Gamma_0(n):\Gamma(n)]} = \frac{n^3\prod_{p\mid n}\left(1-\frac{1}{p^2}\right)}{n^2\prod_{p\mid n}\left(1-\frac{1}{p}\right)} = n\prod_{p\mid n}\left(1+\frac{1}{p}\right)$$
\end{demonstracao}

\begin{teorema}[Hensel lifting]
    Let $f(x)\in\z[x]$ and $k$ a positive integer. If $r$ is an integer such that $f(r)\equiv 0 \pmod{p^k}$ and $f'(r)\neq 0 \pmod{p}$ then for every $m>0$ there exists an integer $s$, unique modulo $p^{k+m}$ such that $f(s)\equiv 0 \pmod{p^{k+m}}$ and $r\equiv s \pmod{p^{k}}$.
\end{teorema}

\begin{demonstracao}
    See Theorem 4.5.2 of \cite{gouvea2020}, which is a more general form of this statement using the $p$-adic integers $\z_p$. The first section of the next chapter reviews some important facts on $\z_p$.
    
    To obtain this version, note that $\z\subset\z_p$, so if $f(x)\in\z[x] $ then $ f(x)\in\z_p[x]$ and $f(r)\equiv 0 \pmod{p^k}. $ Hence $ f(r)\equiv 0 \pmod{p}$ and so the theorem guarantees a unique $p$-adic integer $\alpha$ which is an actual root of $f$, that is, $f(\alpha)=0$.

    Write $\alpha=\sum_{i\geq 0}a_ip^i$ with $a_i\in\{0,\dots,p-1\}$ which converges considering the $p$-adic topology. Any $p$-adic integer has a unique such form.

    Given $m>0$, let $s$ be the truncation $s:=\sum_{i=0}^{k+m-1}a_ip^i$. Then $\alpha-s\in p^{k+m}\z_p$, that is $\alpha \equiv s \pmod{p^{k+m}}$ and so $0 = f(\alpha)\equiv f(s) \pmod{p^{k+m}}$. Uniqueness of $\alpha$ gives uniqueness of $s$. Note that $s$ is simply the $\alpha_{k+m}$ derived in the proof, giving $r\equiv s \pmod{p^k}$.
\end{demonstracao}

Hensel's lifting is particularly interesting for us in the case $k=0$. Since $f(x)\equiv 0 \pmod{p^m}$ implies $f(x)\equiv 0 \pmod{p}$, the Theorem allows us to determine wether some polynomials have roots mod $p^m$ simply by judging whether it has roots mod $p$.

Since $\z/n\z \simeq \prod \z/p_i^{l_i}\z$, one can determine whether a polynomial has roots mod $n$ by judging wheter is has roots mod $p^k$, and for our polynomials of interest, by looking at roots mod $p$. In particular, if $a \in \z/n\z$ is equivalent to $(a_1,\dots,a_k)$ under this isomorphism, then $a$ is a root of a polynomial if and only if each $a_i$ is a root.

\begin{lema}\label{numb b(n)}
    Let $n=3^{l_0}p_1^{l_1}\cdots p_k^{l_k}$ where $3,p_1,\dots,p_k$ are distinct primes. Then the polynomial $f(x)=x^2+x+1 \in (\z/n)[x]$ has $2^k$ roots if $l_0\leq 1$ and all $p_i\equiv 1 \pmod{6}$ and 0 roots otherwise.
\end{lema}

\begin{demonstracao}
    Note that $f'(x)=2x+1$, so if $f(r)\equiv 0 \pmod{p}$ and $f'(r)\equiv 0\pmod{p}$, we have that $r^2+r+1 = r^2+2r+1-r\equiv 0$ and $2r+1\equiv 0$, so $r^2-r\equiv 0$ implies $r\equiv 1$ (since $r\equiv 0$ would contradict $f(x)\equiv 0$). But $f(1) = 3$, so $1$ is only a root mod $n$ if $n=3$. Thus, if $l_0>1$, this root can't be lifted. We have, then, that 3 contributes with a single root and any other prime would contribute with 2 roots.

    Recall that the discriminant of a degree 2 polynomial $ax^2+bx+c$ is given by $b^2-4ac$. Since the discriminant of $f(x)$ is $-3$, we have that a prime $p\neq 3$ that divides $n$ will contribute with two roots if and only if $-3$ is a square mod $p$, that is, if $p\equiv 1\pmod{6}$.

    Thus, if $l_0\leq 1$ and each of the other $k$ primes that divide $n$ are congruent to 1 mod 6, there are $2^k$ possible combinations of $p_i^{l_i}$-roots to create a root mod $n$. Any other case would lead to one of the $a_i$ terms in $(a_0,\dots,a_l)$ never being 0 for $f(a_0,\dots,a_l)$, so $f$ would not have a root. 
\end{demonstracao}

\begin{lema}\label{numb c(n)}
    Let $n=2^{l_0}p_1^{l_1}\dots p_k^{l_k}$ where $2,p_1,\dots,p_k$ are distinct primes. Then the polynomial $g(x)=x^2+1 \in (\z/n)[x]$ has $2^k$ roots if $l_0\leq 1$ and all $p_i\equiv 1 \pmod{4}$ and 0 roots otherwise.
\end{lema}

\begin{demonstracao}
    The proof is analogous to the previous case.

    Since $g'(x)=2x$, if $g(r)\equiv 0 \pmod{p}$ and $g'(r)\equiv 0\pmod{p}$, we have that $(r+1)^2-2r = r^2+1 \equiv 0$ and $2r\equiv 0$, so $(r-1)^2\equiv 0$ implies $r\equiv 1$ (since $r\equiv 0$ would contradict $g(x)\equiv 0$). But $g(1) = 2$, so $1$ is only a root mod $n$ if $n=2$. Thus, if $l_0>1$, this root can't be lifted  any other prime would contribute with 2 roots.

    The discriminant of $g(x)$ is $-4$, so we have that a prime $p\neq 2$ that divides $n$ will contribute with two roots if and only if $-4$ is a square mod $p$, that is, if $p\equiv 1\pmod{4}$. The rest of the argument is analogous to the previous case.
\end{demonstracao}

Finally, we summarize our results:

\begin{proposicao}\label{prp: numbers}
    Let $n=p_0^{l_0}p_1^{l_1}\cdots p_k^{l_k}$ be the prime decomposition of a positive 
    integer $n$. then
    \begin{align*}
    a(n)&=[\SL_2(\z):\Gamma_0(n)]=n\prod_{p\mid n}\bigg(\frac{p+1}{p}\bigg),\\
    b(n)&=\begin{cases}
    2^k & \text{if $p_0=3$, $l_0\leq 1$ and all $p_i\equiv 1\!\!\!\!\! \pmod{6}$,}\\
    0 & \text{otherwise}
    \end{cases}\\
    c(n)&=\begin{cases}
    2^k & \text{if $p_0=2$, $l_0\leq 1$ and all $p_i\equiv 1\!\!\!\!\! \pmod{4}$.}\\
    0 & \text{otherwise}
    \end{cases}
    \end{align*}
\end{proposicao}

\begin{demonstracao}
    This follows from Lemmas \ref{numb a(n)}, \ref{numb b(n)} and \ref{numb c(n)}, since $a(n)$ is the index of $\Gamma_0(n)$ in $\SL_2(\z)$, $b(n)$ is the number of roots of $x^2+x+1$ mod $n$ and $c(n)$ is the number of roots of $x^2+1$ mod $n$.    
\end{demonstracao}

\chapter{The homology groups of \texorpdfstring{$\text{SL}_2(\mathbb{Z}[1/n])$}{Lg}}
\label{chapter:gammapq}
This chapter introduces for each integer $n$ a spectral sequence converging to the homology groups of $\SL_2(\z[1/n])$. If $n = p_1^{e_1}p_2^{e_2}\cdots p^{e_k}_k$ is the prime factorization of $n$, then $\SL_2(\z[1/n])=\SL_2(\z[1/p_1p_2\cdots p_k])$ and the desired spectral sequence is obtained via an action of this group on a contractible CW-complex which will be a product of trees. 

This spectral sequence allows us to obtain some results about the groups $H_s(\SL_2(\z[1/n]))$ and $H_s(\PSL_2(\z[1/n]))$, in particular, that they are finitely generated. We also obtain an upper bound for the ranks of this groups and a description of $H_1(\SL_2(\z[1/n]))$. 

\section{Bruhat-Tits Building associated with a prime}

In this section, we will construct for each $p$ prime a tree $B_p$ on which $\SL_2(\z[1/p])$ acts. This will allow us to consider the action of $\SL_2(\z[1/p_1p_2\cdots p_k])$ on the product of trees $B_{p_1}\times B_{p_2}\times\cdots\times B_{p_k}$. This is a particular case of the tree of $\SL_2$ over a local field as constructed in \cite{serre1980}.

First, we need some results about the topology of $\GL_2(\q)$ induced by the $p$-adic topology of $\q$. For further information on $p$-adics, we recommend \cite{gouvea2020}.

Let $p$ be a prime. For all $x\in\q$ non-zero, there is a unique way of writing $x$ as $p^k\displaystyle\frac{a}{b}$ where $k\in\z$ and $a,b$ are coprime integers and $p\nmid a,p\nmid b$. This defines a function $\nu_p: \q\rightarrow \z$ where $x\mapsto k$. This function is called the $p$-adic valuation and can be extended to all of $\q$ by letting $\nu_p(0) = \infty$.

Note that integers always have non-negative valuations.

\begin{proposicao}
    For all $x,y\in \q$, the following statements hold: \begin{itemize}
        \item $\nu_p(xy) = \nu_p(x)+\nu_p(y)$, 
        \item $\nu_p(x+y) \geq \min\{\nu_p(x),\nu_p(y)\}$.
    \end{itemize}
\end{proposicao}

\begin{demonstracao}
    Note that if either $x$ or $y$ is 0 the statements follow trivially. Suppose then that neither are 0 and write $x=p^{\nu_p(x)}\displaystyle\frac{a}{b}$ and $y=p^{\nu_p(y)}\displaystyle\frac{a'}{b'}$ where $p$ doesn't divide any of the $a,b,a',b'$. Then $p\nmid aa'$, $p\nmid bb'$ and so $$\nu_p(xy) = \nu_p\left(p^{\nu_p(x)}\frac{a}{b}p^{\nu_p(y)}\frac{a'}{b'}\right) = \nu_p\left(p^{\nu_p(x)+\nu_p(y)}\frac{aa'}{bb'}\right) = \nu_p(x)+\nu_p(y).$$

    For the second statement, suppose $\nu_p(x)\leq \nu_p(y)$, so that $0\leq \nu_p(y)-\nu_p(x)$. Then $ab'+p^{\nu_p(y)-\nu_p(x)}a'b$ is an integer, thus has non-negative valuation.  Hence$$\nu_p(x+y) = \nu_p\left(p^{\nu_p(x)}\frac{a}{b}+p^{\nu_p(y)}\frac{a'}{b'}\right) = \nu_p\left(p^{\nu_p(x)}\frac{ab'+p^{\nu_p(y)-\nu_p(x)}a'b}{bb'}\right)\geq \nu_p(x).$$ 
    
    Thus $\nu_p(x+y)\geq \min\{\nu_p(x),\nu_p(y)\}$. 
\end{demonstracao}

\begin{definicao}
    A non-archimedean norm on a field $\mathbb{K}$ is a function $|\phantom{x}|:\mathbb{K}\rightarrow\mathbb{R}_+$ that satisfies: \begin{enumerate}
        \item $|x|=0 \Leftrightarrow x=0$;
        \item $|xy|=|x||y|$ for all $x,y\in\mathbb{K}$;
        \item $|x+y|\leq \max\{|x|,|y|\}$ for all $x,y\in\mathbb{K}$.
    \end{enumerate}
\end{definicao}

Note that the third condition implies that this function is a norm, since $\max\{|x|,|y|\} \leq |x|+|y|$. 

Recall that a norm always defines a distance function $d: \mathbb{K}\times\mathbb{K}\rightarrow \mathbb{R}_+$ where $d(x,y) = |x-y|$, thus allowing us to induce a topology on $\mathbb{K}$ making it a metric space.

\begin{teorema}
    For all primes $p$, the function $|\phantom{x}|_p: \q \rightarrow \mathbb{R}_+$ given by $$|x|_p \coloneq p^{-\nu_p(x)}$$ is a non-archimedean norm.
\end{teorema}

\begin{demonstracao}
    Note here that $|0|_p = p^{-\infty}$ is understood to be 0, as expected. No other powers of $p$ can be zero, so the first condition is satisfied.

    For the second condition, recall that $\nu_p(xy) = \nu_p(x)+\nu_p(y)$, so that $$|xy|_p = p^{-\nu_p(xy)} = p^{-\nu_p(x)-\nu_p(y)} = p^{-\nu_p(x)}p^{-\nu_p(y)} = |x|_p|y|_p.$$

    For the last condition, let $x,y\in\q$. If either is zero, the condition is satisfied trivially, so let $x,y\neq 0$. Suppose $\max\{|x|_p,|y|_p\} = |x|_p$, so that $|y|_p\leq |x|_p$ implies $p^{-\nu_p(y)}\leq p^{-\nu_p(x)}$ and thus, $-\nu_p(y)\leq -\nu_p(x)$, leading us to conclude that $\nu_p(y)\geq \nu_p(x)$. Since $\nu_p(x+y)\geq \min\{\nu_p(x),\nu_p(y)\} = \nu_p(x)$, we have $-\nu_p(x+y)\leq -\nu_p(x)$ and thus $p^{-\nu_p(x+y)}\leq p^{-\nu_p(x)}$. Then, $$|x+y|_p = p^{-\nu_p(x+y)} \leq p^{-\nu_p(x)} = |x|_p = \max\{|x|_p,|y|_p\}.$$
\end{demonstracao}

The previous theorem shows that $\q$ is a metric space with the topology given by the $p$-adic metric $d_p(x,y) = |x-y|_p$.

The completion of $\q$ with respect to the $p$-adic metric is denoted by $\q_p$ and is called the set of $p$-adic rational.

\begin{observacao}\label{inf sum}
    Any $p$-adic rational (in particular, any rational number in the copy of $\q$ inside $\q_p$) can be written uniquely as a sum $\displaystyle\sum_{k\geq n_0} a_kp^k$ with each $a_k$ being an integer between 0 and $p-1$ and $n_0$ not necessarily positive. Since $\nu_p(p^k)=k$, we get that $|a_kp^k|_p = p^{-k}$, so the terms in the sequence get progressively smaller in the $p$-adic metric, which is why the sequence converges. See Corollary 4.3.4 of \cite{gouvea2020} for the proof of this. In particular, if a number has a representation where $n_0\geq 0$, we call it a $p$-adic integer. This condition is equivalent to having a non-negative $p$-valuation.
\end{observacao}

Let $$\z_{(p)} \coloneq \left\{\frac{a}{b}\in\q \hspace{5pt}\middle|\hspace{5pt} \gcd(a,b)=1, p\nmid b \right\}$$ be the localization of $\z$ on the complement of the ideal $(p)$. Observe that $$\z_{(p)} = \{x \in \q \mid \nu_p(x)\geq 0\},$$ that is, $\z_{(p)}$ coincides with the set of rational numbers with positive $p$-valuation (which is embedded on $\q_p$ as a subset of the $p$-adic integers). Note that this set is a PID and a local ring, where the ideals are given by powers of $p$, and $\z_{(p)}/p^n\z_{(p)} \simeq \z/p^n\z$.

We give $\q^4$ the product topology considering the $p$-adic topology on each component, and we identify $\q^4$ with the set $M_2(\q)$ of $2\times 2$ matrices with coefficients in $\q$. We will refer to this as the $p$-adic topology on $M_2(\q)$.

\begin{proposicao}
    $\GL_2(\q)$ is open in the $p$-adic topology.
\end{proposicao}

\begin{demonstracao}
    Recall that addition and multiplication are continuous functions in metric spaces, so that the determinant function is continuous in the $p$-adic topology. Furthermore, metric spaces are Hausdorff (if $x\neq y$ and $d=\frac{1}{2}|x-y|$, then the open balls $B(x,d),B(y,d)$ are disjoint), implying that $\q\backslash\{0\}$ is open and thus its pre-image under the determinant function, $\GL_2(\q)$, is an open set.
\end{demonstracao} 

As we remarked on the previous proposition, multiplication being a continuous function on metric spaces implies that if $U$ is open, then the sets $gU$ and $Ug$ will also be open for all fixed $g\in \GL_2(\q)$, since these translations are homeomorphisms.

\begin{teorema}\label{teo: GL2Zp open}
    $\GL_2(\z_{(p)})$ is open in the $p$-adic topology.
\end{teorema}

\begin{demonstracao}
    Consider the subset $\z_{(p)}^4 \subset \q^4$. Under the identification $\q^4\simeq M_2(\q)$, this set corresponds to $M_2(\z_{(p)})$ and the intersection of this set with $\GL_2(\q)$ (which is open) is precisely the set $\GL_2(\z_{(p)})$. Thus, we prove that $\z_{(p)}^4$ is open with the $p$-adic topology.

    Note that \begin{align*}
        \z_{(p)}^4 &= \{(x_1,x_2,x_3,x_4)\in\q^4 \mid \nu_p(x_i)\geq 0 \text{ for all } x_i\} \\
        &= \{(x_1,x_2,x_3,x_4)\in\q^4 \mid -\nu_p(x_i)\leq 0 \text{ for all } x_i\} \\
        &= \{(x_1,x_2,x_3,x_4)\in\q^4 \mid p^{-\nu_p(x_i)}\leq 1 \text{ for all } x_i\} \\
        &= \{(x_1,x_2,x_3,x_4)\in\q^4 \mid |x|_p\leq 1 \text{ for all } x_i\} 
    \end{align*} in other words, $\z_{(p)}^4$ is the product of the closed balls $\overline{B}(0,1)$ of $\q$. Since we are considering the product topology on $\q^4$, it suffices to show that $\overline{B}(0,1)$ is open.

    Let $y\in \overline{B}(0,1)$ and let $s<1$. We claim that the open ball $B(y,s)$ is contained in $\overline{B}(0,1)$, thus showing that every $y$ is an interior point and thus that $\overline{B}(0,1)$ is open. Let $x\in B(y,s)$. Then $|x-y|_p<s$ and $$|x|_p  = |(x-y)+(y-0)|_p \leq \max\{|x-y|_p,|y|_p\} \leq \max\{s,1\} = 1$$ that is, $|x|_p\leq 1$ and so $x\in \overline{B}(0,1)$.
\end{demonstracao}

\begin{observacao}
    The previous result is counterintuitive, since it implies closed balls in $p$-adic topology are also open. In general, the topology of non-archimedean fields also has other weird features, such as every point in a ball being the center of the ball, or every triangle being an isosceles triangle. See Corollary 2.3.5 and Proposition 2.3.7 of \cite{gouvea2020} for proofs of these statements.
\end{observacao}

\begin{teorema}\label{teo: gamma0 open}
    Let $\Gamma_0(p,\z_{(p)})$ be the subgroup of $\GL_2(\z_{(p)})$ consisting of matrices whose lower left-entry is congruent to 0 mod $p$, that is, $$\Gamma_0(p,\z_{(p)}) \coloneq \left\{\begin{pmatrix}
        a & b \\ c & d
    \end{pmatrix}\in \GL_2(\z_{(p)})\mid c\in p\z_{(p)}\right\}.$$ Then $\Gamma_0(p,\z_{(p)})$ is an open set in the $p$-adic topology.
\end{teorema}

\begin{demonstracao}
    Let $\pi_p = \begin{pmatrix}
        1 & 0 \\ 0 & p
    \end{pmatrix} \in \GL_2(\q)$. Then the subgroup we are interested in is the intersection $(\GL_2(\z_{(p)}))\cap(\pi_p\GL_2(\z_{(p)})\pi_p^{-1})$ of two open sets, and thus, it is open. 
\end{demonstracao}

\begin{lema}
    Let $G$ be a subgroup of $\GL_2(\q)$ generated by a subset $R$. If the closure of a subgroup $X$ of $\GL_2(\q)$ contains $R$, then it contains $G$.
\end{lema}

\begin{demonstracao}
    For this proof, $\|\phantom{A}\|_p$ will indicate the norm induced on $M_2(\q)$ by the $p$-adic topology.
    
    Let $A$ be some matrix in $G$. Then $A = R_1R_2\cdots R_n$ for some $R_i\in R$. This lemma will be proven by induction on the integer $n$. If $n=1$, then $A\in R \subset \overline{X}$, so the base case is trivial.

    Suppose the statement is true for all $n<k$ for some integer $k$ and let $R_1\cdots R_k \in G$, $\varepsilon>0$. Then there exists $S,S'\in X$ such that $\|S-R_1\cdots R_{k-1}\|p<\varepsilon$, $\|S'-R_k\|_p<\varepsilon$, so that \begin{align*}
        \|SS' - R_1\cdots R_k\|_p &\leq \|SS'-SR_k\|_p+\|SR_k-R_1\cdots R_k\|_p \\
        &= \|S\|_p\|S'-R_k\|_p+\|R_k\|_p\|SR_k-R_1\cdots R_{k-1}\|_p \\
        &< (\|S\|_p+\|R_k\|_p)\varepsilon.
    \end{align*}
\end{demonstracao}

\begin{definicao}
    Given $n$ an integer, let $\mathbb{Z}[1/n] := \{a/n^r: a \in \mathbb{Z}, r \in \mathbb{Z}^{\geq 0}\}.$
\end{definicao}

\begin{teorema}
    The closure of $\SL_2(\z[1/p])$ contains $\SL_2(\q)$ in the $p$-adic topology.
\end{teorema}

\begin{demonstracao}
    Recall from Remark \ref{inf sum} that every rational number $q$ can be  written in $\q_p$ as a sum $\displaystyle\sum_{k\geq n_0}a_kp^k$ with each $a_k$ an integer between 0 and $p-1$. Given $\varepsilon>0$, let $n$ be an integer such that $p^n> \displaystyle\frac{1}{\varepsilon}$. Let $q_{n-1} = \displaystyle\sum_{k=n_0}^{n-1}a_kp^k = p^{n_0}\displaystyle\sum_{k=n_0}^{n-1}a_kp^{k-n_0}$. Since $k-n_0\geq 0$ for all $k=n_0,\dots,n-1$, we have that $q_{n-1}p^{-n_0}\in\z$, so $q_{n-1}\in\z[1/p]$. Then $$|q-q_{n-1}|_p = \left|\sum_{k=n_0}^{\infty}a_kp^k-\sum_{k=n_0}^{n-1}a_kp^k\right|_p = \left|\sum_{k=n}^{\infty}a_kp^k\right|_p=\left|p^n\sum_{k=n}^{\infty}a_kp^{k-n}\right|_p\leq p^{-n}<\varepsilon.$$ Thus, $\z[1/p]$ is dense in $\q$. 

    Since $\q$ is a field (in particular, a semi-local ring) Theorem \ref{teo: GE2} shows that the matrices $$\begin{pmatrix}
        1 & a \\ 0 & 1
    \end{pmatrix}, \begin{pmatrix}
        1 & 0 \\ a & 1
    \end{pmatrix}$$ where $a\in\q$ generate $\SL_2(\q)$.

    Since each $a\in \q$ can be approximated by a sequence $(a_i)\subset \z[1/p]$, then each of the above matrices can be approximated by a matrix in $\SL_2(\z[1/p])$, so that the closure of $\SL_2(\z[1/p])$ contains all elementary matrices and thus, all of $\SL_2(\q)$.
\end{demonstracao}

\begin{lema}\label{lema: qz[1/p] dense}
    Let $q$ be a prime different than $p$ and let $q\z[1/p]$ be the ideal of $\z[1/p]$ generated by $q$, that is, $$q\z[1/p] = \left\{qx\hspace{5pt}\middle|\hspace{5pt} x \in \z[1/p]\right\} = \{x \in \z[1/p] \mid \nu_q(x)\geq 1\}.$$ Then $q\z[1/p]$ is dense in $\z[1/p]$ (and thus, in $\q$) in the $p$-adic topology.
\end{lema}

\begin{demonstracao}
    We remark that if $y\in\z[1/p]$, then $\nu_q(y)\geq 0$ since $y$ only admits powers of $p$ as denominators, and thus $\nu_q(qy)  =\nu_q(q)+\nu_q(y) = 1+\nu_q(y) \geq 1$, justifying our choice of naming this set.

    Let $x\in\z[1/p]$. We show that there is a sequence in $q\z[1/p]$ converging to $x$ in the $p$-adic topology. 

    Since $x\in\z[1/p]$, there exists $k\geq 0$ such that $x':=p^kx \in \z_{(p)}$ (this $k$ can be taken as $\nu_p(x)$).

    Let $n\geq 0$ and recall that $\z_{(p)}/p^n\z_{(p)} \simeq \z/p^n\z$. Since $\gcd(q,p^n)=1$, we have that $q$ is invertible mod $p^n$, so the equation $qy_n'=x'$ has solution $\overline{y_n} = \overline{q^{-1}x}$ in $\z/p^n\z$. Then $qy_n' \equiv x' \mod p^n$ implies that $p^n$ divides $qy_n'-x'$, that is, $qy_n'-x' = p^na$ for some $a\in\z_{(p)}$, implying that $q(p^ky_n')-x = p^n(p^ka)$. Take $y_n := p^ky_n'$. Then $qy_n-x \in p^n\z_{(p)}$, that is, $\nu_p(qy_n-x) \geq n$, implying that $|qy_n-x|_p \leq p^{-n}$, and so, $qy_n\in B(x,p^{-n})$.

    Thus, the sequence $qy_n$ is contained in $q\z[1/p]$ and converges to $x$ in the $p$-adic topology.
\end{demonstracao}

\begin{teorema}\label{teo: closure gammapq}
    The closure of $\Gamma_0(q,p) \coloneq \left\{\displaystyle\begin{pmatrix}
        a & b \\ c & d
    \end{pmatrix} \in\SL_2(\z[1/p]) \mid c\equiv 0 \mod q\right\}$ contains $\SL_2(\z[1/p])$.
\end{teorema}

\begin{demonstracao}
    From Theorem \ref{teo: GE2}, $\SL_2(\z[1/p])$ is generated by elementary matrices. Since $q\z[1/p]$ is dense in $\z[1/p]$, the only elementary matrices in $\SL_2(\z[1/p])$ not in $\Gamma_0(q,p)$ are of the type $\begin{pmatrix}
        1 & 0 \\ a & 1
    \end{pmatrix}$ with $a\in\z[1/p]$. Since $a$ can be approximated by a sequence $(a_i)\subset q\z[1/p]$, those matrices can be approximated by the sequence $\begin{pmatrix}
        1 & 0 \\ a_i & 1
    \end{pmatrix}$ in $\Gamma_0(q,p)$.
\end{demonstracao}

We now turn to the construction of our tree.

Given linearly independent vectors $v_1,v_2\in \mathbb{Q}^2$, the lattice defined by these vectors is the set $$L_{v_1v_2} = \{a_1v_1+a_2v_2 \mid a_1,a_2\in\z_{(p)}\} = \z_{(p)} v_1\oplus\z_{(p)} v_2.$$ Here, we use the notation $L_{v_1v_2} = \langle v_1,v_2\rangle$ to indicate that $L$ is generated by integer sums of $v_1$ and $v_2$. Note that the vectors $e_1 = (1,0)$ and $e_2=(0,1)$ generate the lattice $\z_{(p)}^2$. 

We now consider a relation $\sim$ on the set of lattices where two lattices $L$ and $L'$ are related if they are homothetic, that is $$L\sim L' \Longleftrightarrow L = \alpha L' \text{ for some }\alpha\in\mathbb{Q}.$$ This is obviously an equivalence relation. 

Note that since lattices are composed by two linearly independent vectors, every lattice $L$ can be identified as the set $M\z_{(p)}^2$ where $M\in \GL_2(\q)$. Furtermore, since $U\z_{(p)}^2 = \z_{(p)}^2$ for $U\in\GL_2(\z_{(p)})$, we have that $M$ and $MU$ represent the same lattice.

Furthermore, $\GL_2(\q)$ acts on lattices via the standard representation, that is, if $g\in\GL_2(\q)$ and $L$ is a lattice, $$gL = \{gv \mid v \in L\}$$ where the multiplication $gv$ is understood to be considering the coordinates of $v$ in the canonical basis $\{e_1,e_2\}$ of $\q^2$.

\begin{teorema} \label{invariance}
    Let $L,L'$ be any two lattices. Then there exists a basis $\{u_1,u_2\}$ of $L$ such that $\{p^au_1,p^bu_2\}$ is a basis of $L'$ for some integers $a,b$.
\end{teorema}

\begin{demonstracao}
    Let $L = \langle v_1,v_2\rangle$ and $L'=\langle w_1,w_2\rangle$. Since both $\{v_1,v_2\}$ and $\{w_1,w_2\}$ generate $\q^2$, there exists integers $a_{ij}\in \q$ satisfying $$v_1 = a_{11}w_1+a_{21}w_2,$$ $$v_2 = a_{12}w_1+a_{22}w_2.$$ Let $n$ be large enough that $p^na_{ij}\in\z_{(p)}$ for all $i,j$ (one can take $n = -\min\{\nu_p(a_{ij})\}$). Then $p^nL\subset L'$. A similar argument shows that $p^mL'\subset L$ for some $m\geq 0$.

    Note that $p^{n+s}L\subset p^n L \subset L'$, so any number $l$ greater than $n$ satisfies the condition $p^lL\subset L'$. If there is no minimal $n$ satisfying this condition, then in particular, $p^{-m}L\subset L'$ implies $L\subset p^mL'\subset L$ and thus, $L = p^{m}L'$ and any basis $\{u_1,u_2\}$ of $L$ will suffice, since $\{p^{-m}u_1,p^{-m}u_2\}$ will be a basis of $L'$. 

    For the sake of defining a distance funtion later on, if there is a minimal $n$ satisfying this condition, choose $n$ to be minimal.

    Consider now the quotient $L'/nL$. Since $p^{n+m}L' \subset p^nL$, we have that this quotient has no free elements, thus the structure theorem for finitely generated modules over PIDs (see Theorem \ref{teo: structure}) implies that $$\frac{L'}{p^nL}\simeq \frac{\z_{(p)}}{p^k\z_{(p)}}\oplus \frac{\z_{(p)}}{p^r\z_{(p)}}.$$ Notice that under this isomorphism, the class of $p^kw_1$ is send to 0, implying that $p^kw_1 \in p^nL$, and so, that $w_1 \in p^{n-k}L$. If $w_1 = p^{n-k}(\lambda_1v_1+\lambda_2v_2)$, let $u_1 = \lambda_1v_1+\lambda_2v_2$. Then $w_1 = p^{n-k}u_1$. Similarly, one finds $u_2$ such that $w_2 = p^{n-r}u_2$.

    To finish, take $a=n-k$ and $b=n-r$, so that $L = \langle u_1,u_2\rangle$ and $L' = \langle p^au_1,p^bu_2\rangle$.
\end{demonstracao}

\begin{definicao}
    Let $\Lambda,\Lambda'$ be two equivalence classes of vertices. Choose representatives $L$ of $\Lambda$ and $L'$ of $\Lambda'$. Using the previous theorem, let $\{u,v\}$ be a basis of $L$ such that $\{p^au,p^bv\}$ is a basis of $L'$. Define the distance $d(\Lambda,\Lambda')$ of the two classes $\Lambda,\Lambda'$ to be $|a-b|$. 
\end{definicao}

\begin{teorema}\label{dist}
    The distance $d(\Lambda,\Lambda')$ of two classes is well-defined.
\end{teorema}

\begin{demonstracao}
    Note that on the proof of the last theorem, once $n$ is fixed, the integers $k,r$ are uniquely determined by the structure theorem and so $a,b$ depend only on $n$. 

    If there is no $n$ satisfying a minimal condition on $p^nL\subset L'$, then for any basis $\{u_1,u_2\}$ for $L$, we have $\{p^{-m}u_1,p^{-m}u_2\}$ is a basis for $L'$. Any choice of $m$ leads to the conclusion that $d(\Lambda,\Lambda') = 0$, reflecting the fact that $p^mL=L'$ implies the classes are the same. 

    Suppose, then, that there exists a minimal $n$. To show that this distance also does not depend on the choice of representatives, replace $L$ and $L'$ by  $\alpha L$ and $\beta L'$ with $\alpha,\beta\in\q^\times$.

    Let $c = \nu_p(\alpha/\beta) \in\q^\times$. We claim that $n-c$ satisfies minimality for the pair $\alpha L,\beta L'$.

    Indeed, suppose $k<n-c$ and suppose $p^k\alpha L\subset \beta L'$. Then $p^{k+c}L = (\alpha/\beta)p^kL\subset L'$ implies $k+c$ satisfies the condition, and thus $k+c \geq n$, which is absurd.

    Thus, $n-c$ is the minimal element for the pair $\alpha L,\beta L'$. Then $$\frac{\beta L'}{p^{n-c}\alpha L} = \frac{p^{\nu_p(\beta)}L'}{p^{n-c}p^{\nu_p(\alpha)}L} \simeq \frac{L'}{p^nL}\simeq \frac{\z_{(p)}}{p^k\z_{(p)}}\oplus \frac{\z_{(p)}}{p^r\z_{(p)}}$$ so that the pair $(a',b')$ for $\alpha L,\beta L'$ is given by $a' = n-c-k, b' = n-c-r$, that is, $(a',b') = (a-c,b-c)$, which doesn't affect the distance.
\end{demonstracao}

\begin{definicao}\label{def: Bp}
    Let $B_p$ be the graph whose vertices are be equivalence classes $\Lambda$ of lattices. We say that two vertices are adjacent if their distance is 1.
\end{definicao}

\begin{observacao}\label{rem: iso lattice}
    If $L'\subset L$, the pair $(a,b)$ for $L,L'$ will be determined by $$\frac{L}{L'}\simeq \frac{\z_{(p)}}{p^a\z_{(p)}}\oplus\frac{\z_{(p)}}{p^b\z_{(p)}}$$ with $d(\Lambda,\Lambda') = |b-a|$, so that if one of the factors $(a,b)$ is 0, the other must be $d(\Lambda,\Lambda')$. Note also that these values $a,b$ must be non-negative since they are describing ideals of $\z_{(p)}$.
\end{observacao}

\begin{observacao}\label{rem: only p affects L}
    If $\nu_p(a)=0$, then $\nu_p(a^{-1})=0$ and $aL = L$ for any lattice $L$. This was used previously to argue that the distance function is well defined. This also implies that $\alpha L = p^{\nu_p(\alpha)}L$, so that only homothety by a power of $p$ can affect a lattice. In other words, if $L'\in\Lambda'$, then any other lattice in $\Lambda'$ is of the type $p^kL'$.
\end{observacao}

\begin{proposicao}\label{comp series}
    Fix a lattice $L$ in a class $\Lambda \in B_p$. Then each class $\Lambda'$ has exactly one representative $L'$ satisfying the following equivalent conditions:
    \begin{enumerate}
        \item[(i)] $L' \subset L$ and $L'$ is maximal (in $\Lambda'$) with this property;
        \item[(ii)] $L' \subset L$ and $L' \not\subset pL$;
        \item[(iii)] $L' \subset L$ and $L/L' \simeq \z_{(p)}/p^n\z_{(p)}$ with $n = d(\Lambda,\Lambda')$.
    \end{enumerate}
\end{proposicao}

\begin{demonstracao}
    Note that if $\Lambda\neq\Lambda'$, then for any $L''\in\Lambda'$ there is a minimal $n$ satisfying $p^nL''\subset L$ (see the proof of Theorem \ref{dist}). This allows us to choose $L' = p^nL''$ and this will be maximal in that class, since any other representative of $\Lambda'$ will be of type $p^k L''$ for some $k\in\z$ (see Remark \ref{rem: only p affects L}) and $p^nL'' \subset p^kL''$ implies $L''\subset p^{k-n}L''$ so $k-n<0$, that is, $k<n$ and so $p^kL''\not\subset L$ from minimality of $n$. 

    It remains to show that the conditions are equivalent. 

    We will first prove that (i) and (ii) are equivalent then show they are both equivalent to (iii).

    (i) $\implies$ (ii) Suppose $L'\subset pL$. Then $p^{-1}L'\subset L$ and $L'\subset p^{-1}L'$, so the maximal element of $\Lambda'$ is $p^{-1}L'$, contradicting (i).

    (ii) $\implies$ (i) Suppose $L'$ is not maximal in $\Lambda'$, that is, there exists $L''\in\Lambda'$ with $L'\subset L''\subset L$. Then $L'' = p^{-k}L'$ with $k\geq 0$ implying that $L' = p^kL''\subset p^kL \subset pL$, which contradicts (ii).

    (iii) $\implies$ (ii) Let $\phi$ be the isomorphism described above and let $w_1,w_2$ be a basis of $L'$, so that $\phi(\overline{w_1}),\phi(\overline{w_2})$ generate the invariant factors. If $L'\subset pL$, then $\phi(\overline{pw_1}) = p\phi(\overline{w_1}) = (0,0)$, that is, $p\phi(\overline{w_1}) \in p^a\z_{(p)}$, implying $\phi(\overline{w_1}) \in p^{a-1}\z_{(p)}$, forcing $a$ to be at least 1. Analogously, we conclude $b\geq 1$. This contradicts (iii). 

    (ii) $\implies$ (iii) Since (ii) and (i) are equivalent, we will use both for this proof. Note that $L'$ being maximal for the property $L'\subset L$ implies that for any other $p^kL'$ in $\Lambda'$, $$p^kL'\subset L \implies p^kL'\subset L' \implies k\geq 0,$$ that is, 0 is the minimal integer for the pair $L,L'$. Thus the pair $(a,b)$ in the distance function is determined by $$\frac{L}{L'}\simeq \frac{\z_{(p)}}{p^a\z_{(p)}}\oplus\frac{\z_{(p)}}{p^b\z_{(p)}}.$$ Then, there is a basis $\{u_1,u_2\}$ of $L$ such that $\{p^au_1,p^bu_2\}$ is a basis of $L'$. Since $L'\not\subset pL$, either $p^au_1$ or $p^bu_2$ is not in $pL$. Suppose $p^au_1\notin pL$. Then $p^{a-1}u_1\notin L$, which implies $a-1<0$. But $a\geq 0$, which leads to $a=0$, forcing $b=d(\Lambda,\Lambda')$. 
\end{demonstracao}

\begin{teorema}\label{neighbors}
    If $\Lambda$ is a class and $L\in\Lambda$ has a basis $\{u_1,u_2\}$, then the neighbors of $\Lambda$ are represented by the classes $\langle u_1,pu_2\rangle,\langle pu_1,au_1+u_2\rangle$ with $0\leq a \leq p-1$.
\end{teorema}

\begin{demonstracao}
    If $\Lambda$ and $\Lambda'$ form an edge, then $d(\Lambda,\Lambda')=1$. Choose $L\in\Lambda, L''\in\Lambda'$. Then there exists a basis $\{u_1,u_2\}$ of $L$ such that $\{p^au_1,p^bu_2\}$ is a basis of $L''$. Since $|b-a|=1$, suppose $b-a=1$ (invert the vectors if necessary). Then $\{p^au_1,p^{a+1}u_2\}$ is the given basis. Let $L':=p^{-a}L'' = \langle u_1,pu_2\rangle$, so that $pL\subset L' \subset L$. This implies that $L'/pL$ is a non-trivial submodule of $L/pL \simeq \mathbb{F}_p^2$, that is, it is generated by a line in $\mathbb{F}_p^2$. In other words, neighborhs of $L$ are in correspondence with 1-dimensional subspaces of $\mathbb{F}_p^2$, which can be taken as $\langle(1,0)\rangle$  $\langle(a,1)\rangle$ with $a\in\mathbb{F}_p$. Thus, the pre-images of these lines under the isomorphism $L/pL$ will give the other neighbors of $L$.

    Taking $\{u_1,u_2\}$ as a basis of $L$, the pre-image of $\langle(a,1)\rangle$ is given by vectors $x\overline{u_1}+y\overline{u_2}$ such that $(\overline{x},\overline{y})$ is congruent to $k(a,1)$ for some $k$, that is, $x\equiv ka, y\equiv k$, which simplifies to $x\equiv ya$. Then $x-ya = pz$ for some $z\in\z_{(p)}$, so that $$x\overline{u_1}+y\overline{u_2} = (ya+pz)\overline{u_1}+y\overline{u_2} = z \overline{pu_1}+y\overline{au_1+u_2}.$$ In other words, $\langle(a,1)\rangle$ corresponds to the neighbor $\langle pu_1,au_1+u_2\rangle$. A similar analysis leads to conclude that $\langle(1,0)\rangle$ corresponds to the neighbor $\langle u_1,pu_2\rangle$. In particular, each class $\Lambda$ has exactly $p+1$ neighbors.
\end{demonstracao}

\begin{definicao}
    Let $M$ be a module over a commutative ring $R$. A \textit{composition series} for $M$ is a sequence of submodules $M_n$ such that $$0=M_n \subset M_{n-1}\subset \dots \subset M_1 \subset M_0 = M$$ and each quotient $M_i/M_{i-1}$ is simple, that is, has no non-trivial submodules. The integer $n$ is said to be the \textit{lenght} of the module $M$ and it is denoted by $n=l(M)$.
\end{definicao}

Note that simple modules have lenght 1. For more information on composition series for modules, we recommend \cite{atiyah1969}.

\begin{observacao}\label{rem: lenght}
    In particular, since $\z_{(p)}$ is a local ring where the only ideals are generated by powers of $p$, we have that for any positive integer $k$, $$0 \subset \frac{p^{k-1}\z_{(p)}}{p^{k}\z_{(p)}} \subset \frac{p^{k-2}\z_{(p)}}{p^{k}\z_{(p)}} \subset \dots \subset \frac{p\z_{(p)}}{p^{k}\z_{(p)}} \subset \frac{\z_{(p)}}{p^{k}\z_{(p)}}$$ is a composition series for $\displaystyle\frac{\z_{(p)}}{p^{k}\z_{(p)}}$ since $\displaystyle\frac{p^{i-1}\z_{(p)}/p^{k}\z_{(p)}}{p^{i}\z_{(p)}/p^{k}\z_{(p)}} \simeq \frac{p^{i-1}\z_{(p)}}{p^{i}\z_{(p)}} \simeq \mathbb{F}_p$ is simple. Note that this composition series has lenght $k$.
\end{observacao}

\begin{teorema}
    $B_p$ is a tree.
\end{teorema}

\begin{demonstracao}
    If $\Lambda$ and $\Lambda'$ are two vertices of $B_p$, let $L$ and $L'$ be lattices representing $\Lambda$ and $\Lambda'$ with $L' \subset L$. A composition series for $L/L'$ gives a sequence of lattices $$L' = L_n \subset L_{n-1}  \subset \dots \subset L_0 = L$$ such that the lenght of each module $L_i/L_{i-1}$ is 1, that is, $L_i/L_{i-1}\simeq \F_p$. By item (iii) of Proposition \ref{comp series}, we have that $d(L_i,L_{i-1})=1$, so the classes $\Lambda_0,\dots,\Lambda_n$ form a path from $\Lambda'$ to $\Lambda$.

    Conversely, Proposition \ref{comp series} also gives that a path between $\Lambda'$ and $\Lambda$ always gives a composition series. Thus, it suffices to prove that for a series as above without backtracking, $\Lambda_0\neq\Lambda_n$. To prove this, we will show by induction on $n$ that if there is no backtracking, then $d(\Lambda_0,\Lambda_n)=n$. The case $n=1$ is obvious.

    For $\Lambda_0$ to equal $\Lambda_n$, the lattices $L_0$ and $L_n$ must be homothetic, and since $L_n\subset L_0$, this could only happen if $L_n = p^kL_0$ for some $k\geq 1$ ($k=0$ is impossible since $l(L_n/L_0)=n\neq 0$). If we prove that $L_n\not\subset pL_0$, this will imply that $L_n\not\subset p^kL_0$ for any $k\geq 0$, and thus, will imply that $L_n$ and $L_0$ cannot be homothetic.

    Note that $pL_{n-2}$ corresponds to the same class as $L_{n-2}$, so that the lattices $L_n$ and $pL_{n-2}$ are neighbors of $L_{n-1}$, and thus are the inverse images of two lines in $L_{n-1}/pL_{n-1}$. These lines have to be distinct, since otherwise $\Lambda_{n-2}\Lambda_{n-1}\Lambda_{n}$ would correspond to a backtracking in the given path. 
    
    Then these distinct lines add up to all of $L_{n-1}/pL_{n-1}$, that is, $$L_{n-1} = L_n+pL_{n-2}.$$ From the composition series, $L_{n-2}\subset L_0 $ implies that $ pL_{n-2}\subset pL_0$. So $$L_n\equiv L_{n-1}\pmod{pL_0}.$$ By the induction hypothesis, $L_{n-1}\not\subset pL_0$, so $L_n\not\subset pL_0$.
\end{demonstracao}

\begin{definicao}
    Let $L_1,L_2$ be two lattices. Define $\chi(L_1,L_2)$ to be the integer $$\chi(L_1,L_2) \coloneq l(L_1/L_3)-l(L_2/L_3)$$ for some lattice $L_3\subset L_1\cap L_2$.
\end{definicao}

\begin{proposicao}
    The definition of $\chi$ does not depend on the choice of lattice.
\end{proposicao}

\begin{demonstracao}
    The lenght of modules is additive (see Proposition 6.9 of \cite{atiyah1969}), that is, if $0\rightarrow M' \rightarrow M \rightarrow M'' \rightarrow 0$ is an exact sequence of $R$-modules, then $l(M) = l(M')+l(M'')$. As a consequence, given $A\subset B \subset C$ a sequence of submodules, the exact sequence $0 \rightarrow B/A \rightarrow C/A \rightarrow C/B \rightarrow 0$ gives $l(C/A) = l(B/A)+l(C/B)$.

    Let $L_3,L_3'$ be two lattices contained in $L_1\cap L_2$ and let $L_4 = L_3\cap L_3'$. Then $$l(L_1/L_3')+l(L_3'/L_4) = l(L_1/L_4)=l(L_1/L_3)+l(L_3/L_4)$$ giving $$l(L_1/L_3)=l(L_1/L_4)-l(L_3/L_4)$$ $$l(L_1/L_3')=l(L_1/L_4)-l(L_3'/L_4)$$

    Similarly, for $L_2$, $$l(L_2/L_3)=l(L_2/L_4)-l(L_3/L_4)$$ $$l(L_2/L_3')=l(L_2/L_4)-l(L_3'/L_4)$$

    Then, \begin{align*}
        l(L_1/L_3)-l(L_2/L_3) &=(l(L_1/L_4)-l(L_3/L_4))-(l(L_2/L_4)-l(L_3/L_4)) \\ &=l(L_1/L_4)-l(L_2/L_4) \\
        &=(l(L_1/L_4)-l(L_3'/L_4))-(l(L_2/L_4)-l(L_3'/L_4)) \\
        &= l(L_1/L_3')-l(L_2/L_3')
    \end{align*}
\end{demonstracao}

\begin{proposicao}
    Let $L$ be a lattice and let $g\in\GL_2(\q)$. Then $\chi(L,gL) = \nu_p(\det(g))$.
\end{proposicao}

\begin{demonstracao}
    Choose a basis $\{u_1,u_2\}$ for $L$ such that $\{p^au_1,p^bu_2\}$ is a basis for $gL$. Then $g$ coincides with the product $ABA^{-1}$ where $A$ is a change of basis matrix and $B$ is a diagonal matrix with entries $p^a,p^b$. Hence $\det(g) = \det(ABA^{-1}) = p^{a+b}$.

    Let $a' = \max\{a,0\},b'= \max\{b,0\}$ and let $L_3 = \langle p^{a'}u_1,p^{b'}u_2\rangle$ so that $L_3\subset L\cap gL$.

    From additivity of lenght and from Remark \ref{rem: lenght}, $$l\left(\frac{L}{L_3}\right) = l\left(\frac{\z_{(p)}}{p^{a'}\z_{(p)}}\oplus \frac{\z_{(p)}}{p^{b'}\z_{(p)}}\right) = l\left(\frac{\z_{(p)}}{p^{a'}\z_{(p)}}\right)+l\left(\frac{\z_{(p)}}{p^{b'}\z_{(p)}}\right) = a'+b'$$ $$l\left(\frac{gL}{L_3}\right) = l\left(\frac{\z_{(p)}}{p^{a'-a}\z_{(p)}}\oplus \frac{\z_{(p)}}{p^{b'-b}\z_{(p)}}\right) = (a'+b')-(a+b)$$ so that $\chi(L,gL) = l(L/L_3)-l(gL/L_3) = a+b = \nu_p(\det(g))$.
\end{demonstracao}

\begin{corolario}
    If $\Lambda$ is a vertex and $g\in\GL_2(\q)$, then $d(\Lambda,g\Lambda) \equiv \nu_p(\det(g)) \pmod{2}$.
\end{corolario}

\begin{demonstracao}
    This follows from the fact that if $d(\Lambda,g\Lambda) = b-a$, then $\nu_p(\det(g)) = a+b$. These integers are equivalent mod 2.
\end{demonstracao}

\begin{teorema}
    If $G$ is a subgroup of $\GL_2(\q)$ such that the determinant of every matrix has even $p$-valuation, then $G$ acts on $B_p$ without inversion. 
\end{teorema}

\begin{demonstracao}
    Since $\nu_p(\det(g)) \equiv d(\Lambda,g\Lambda)\pmod{2}$ for any $g\in\GL_2(\q),\Lambda\in B_p$, we have that if the $p$-valuation of the determinant is even, then $g$ can only move a vertex to ones at even distance from it, thus it cannot move $\Lambda$ to a neighboring vertex $\Lambda'$, and so, it cannot invert an edge.
\end{demonstracao}

\begin{corolario}\label{cor: no inversion}
    $\SL_2(\z[1/p])$ acts on $B_p$ without inversion.
\end{corolario}

\begin{demonstracao}
    Any matrix in $\SL_2(\z[1/p])$ has determinant 1. Since $v_p(1)=0$, the claim follows from the previous theorem.
\end{demonstracao}

\begin{teorema}\label{teo: stabv=stabl}
    If $G$ is a subgroup of $\GL_2(\q)$ such that $\nu_p(\det(g))=0$ for all $g\in G$, then under the action of $G$, the stabilizer of a vertex coincides with the stabilizer of any of its lattices.
\end{teorema}

\begin{demonstracao}
    Let $L$ be a lattice representative of $\Lambda$. If $G_L$ is the stabilizer of a lattice $L$ and $g\in G_L$, then obviously $g$ will also stabilize the vertex $\Lambda$, since its sending a representative of $\Lambda$ to istelf.

    On the other hand, if $g\in G$ and $gL = \alpha L$ for some $\alpha\in\q$, then $gL = p^{\nu_p(\alpha)}L$ and picking a basis $\{u_1,u_2\}$ for $L$, we have that $\{p^{\nu_p(\alpha)}u_1,p^{\nu_p(\alpha)}u_2\}$ is a basis for $gL$. So  $$0=\nu_p(\det(g)) = \chi(L,gL) = 2\nu_p(\alpha),$$ implying that $\alpha$ is invertible in $\z_{(p)}$. Thus $\alpha L=L$ and $g$ also stabilizes $L$.
\end{demonstracao}

\begin{teorema}
    The group $\GL_2(\z_{(p)})$ stabilizes the lattice $\z_{(p)}^2$.
\end{teorema}

\begin{demonstracao}
    Note that the action of $g\in\GL_2(\q)$ on $\z^2_{(p)}$ takes the canonical basis to its column vectors, so $g\z^2_{(p)}= \z^2_{(p)}$ requires each entry in $g$ to be in $\z_{(p)}$.
\end{demonstracao}

\begin{teorema}\label{teo: stab zp}
    Under the action of  $\SL_2(\z[1/p])$, the group $\SL_2(\z)$ is the stabilizer of the vertex associated with $\z_{(p)}^2$.
\end{teorema}

\begin{demonstracao}
    This follows from the fact that the matrices in $\SL_2(\z_{(p)})$ have determinant 1, thus $\nu_p(\det(g))=0$ for all $g\in \SL_2(\z_{(p)})$, and the stabilizer of a vertex coincides with the stabilizer of a lattice (see Theorem \ref{teo: stabv=stabl}), so $$\text{stab}_{\SL_2(\z[1/p])}[\z_{(p)}^2] = \SL_2(\z[1/p])\cap \GL_2(\z_{(p)} = \SL_2(\z)$$
\end{demonstracao}

\begin{corolario}\label{cor: stab pzp}
    Under the action of  $\SL_2(\z[1/p])$, the group $\pi_p\SL_2(\z)\pi_p^{-1}$ is the stabilizer of the vertex associated with $\pi_p\z_{(p)}^2$.
\end{corolario}

\begin{teorema}\label{teo: stab ep}
    Under the action of $\SL_2(\z[1/p])$, the subgroup $\Gamma_0(p)$ of $\SL_2(\z)$ is the stabilizer of the edge formed by $\z_{(p)}^2$ and $\pi_p\z_{(p)}^2$.
\end{teorema}

\begin{demonstracao}
    Since $B_p$ is a tree, if a matrix fixes both vertices of an edge, it must map this edge either to itself or its inverse. Since $\SL_2(\z[1/p])$ acts without inversion (see Corollary \ref{cor: no inversion}), fixing both vertices is equivalent to fixing an edge.
    
    The stabilizers of the vertices $\z_{(p)}^2$ and $\pi_p\z_{(p)}^2$ are $\SL_2(\z)$ and $\pi_p\SL_2(\z)\pi_p^{-1}$, respectively. Their intersection is $\Gamma_0(p)$.
\end{demonstracao}

\begin{teorema}
    Let $\Lambda$ and $\Lambda'$ be two classes of lattices such that $d(\Lambda,\Lambda')$ is even. Then there is a matrix $T$ in $\SL_2(\z[1/p])$ such that $T\Lambda = \Lambda'$.
\end{teorema}

\begin{demonstracao}
    Let $L,L'$ be representatives of $\Lambda$ and $\Lambda'$, respectively. By Theorem \ref{invariance}, there is a basis $\{v_1,v_2\}$ of $L$ such that $\{p^av_1,p^bv_2\}$ is a basis for $L'$. Let $2k = b-a = d(\Lambda,\Lambda')$. Then $\langle p^kv_1,p^{-k}v_2\rangle$ is also a lattice representative of $\Lambda'$.

    Let $B$ be a diagonal matrix with entries $p^k,p^{-k}$ and let $A$ be the change of basis matrix that sends the basis $\{v_1,v_2\}$ to the chanonical basis $\{e_1,e_2\}$ of $\q^2$. Then $A^{-1}BA$ sends $\{v_1,v_2\}$ to $\{p^kf_1,p^{-k}f_2\}$ via standard representation, thus the action of $T=A^{-1}BA$ maps $v_1$ to $p^kv_1$ and $v_2$ to $p^{-k}v_2$.

    Notice that the matrix $T = A^{-1}BA$ is in $\SL_2(\q)$. We wish to find some matrix in $\SL_2(\z[1/p])$ whose action is the same as $T$.

    Since the action of $\GL_2(\q)$ is transitive on the lattices, there is some $g\in \GL_2(\q)$ such that $L=g\z^2_{(p)}$. This implies that the set $g\GL_2(\z_{(p)})g^{-1}$ stabilizes $L$. 

    Furthermore, since $\GL_2(\z_{(p)})$ is open in $\GL_2(\q)$ (see Theorem \ref{teo: GL2Zp open}) and translations are homeomorphisms, we have that $Tg\GL_2(\z_{(p)})g^{-1}$ is an open neighborhood of $T$ in $\GL_2(\q)$. 
    
    Now, since the closure of $\SL_2(\z[1/p])$ contains $\SL_2(\q)$ and $T\in\SL_2(\q)$, any open neighborhood around $T$ will intersect $\SL_2(\z[1/p])$, thus, there must be some element $S $ in $ Tg\GL_2(\z_{(p)})g^{-1} \cap \SL_2(\z[1/p])$. In other words, $S=TU$ for some $U\in g\GL_2(\z_{(p)})g^{-1}$ and $S\in\SL_2(\z[1/p])$.

    Then $SL = (TU)L = T(UL) = TL = L'$.
\end{demonstracao}

\begin{teorema}\label{teo: numb vert 1}
    There are two $\SL_2(\z[1/p])$ equivalence classes of vertices.
\end{teorema}

\begin{demonstracao}
    Let $X_0,X_1$ be the set of vertices at an even and odd distance from $\z^2_p$, respectively, so that the vertices in each set are all at even distance from each other. The previous theorem guarantees that the action of $\SL_2(\z[1/p])$ is transitive on each set. Furthermore, we have that for all $g\in\GL_2(\q)$, $\nu_p(\det(g)) \equiv d(\Lambda,g\Lambda) \pmod{2}$, so if $g$ has determinant 1, it can only move a vertex to one at even distance. Hence these sets are also invariant under the action of $\SL_2(\z[1/p])$.
\end{demonstracao}

\begin{teorema}\label{teo: num edges 1}
    The action of $\SL_2(\z[1/p])$ on the edges is transitive.
\end{teorema}

\begin{demonstracao}
    It suffices to prove that the action is transitive on the edges where one of the vertex is $\z^2_{(p)}$.

    Indeed, suppose $e,e'$ are any two edges in $B_p$. Then both $e,e'$ have a vertex at an even distance from $\z^2_{(p)}$, so by the previous theorem, there exists matrices $S,T\in\SL_2(\z[1/p])$ mapping these vertices to $\z^2_{(p)}$. In other words, $Te$ and $Se'$ are edges with one vertex being $\z^2_{(p)}$. If our initial statement is true, there is $R\in \SL_2(\z[1/p])$ such that $RTe = Se'$. Thus, $S^{-1}RTe = e'$ and we have proven the theorem.

    Recall that there is a basis $\{v_1,v_2\}$ of $\z^2_{(p)}$ such that its neighbors are given by $\langle v_1,pv_2\rangle$ and $\langle pv_1,av_1+v_2\rangle$ with $0\leq a < p$ (see Theorem \ref{neighbors}). 

    Let $T_a = \begin{pmatrix}
        1 & a \\ 0 & 1
    \end{pmatrix}, S = \begin{pmatrix}
        0 & -1 \\ 1 & 0
    \end{pmatrix}$ and let $A$ be the change of basis from $\{v_1,v_2\}$ to the canonical basis of $\q^2$. Then $A^{-1}T_aA$ maps $\langle pv_1,v_2\rangle$ to $\langle pv_1,av_1+v_2\rangle$ and $A^{-1}SA$ maps $\langle pv_1,v_2\rangle$ to $\langle v_1,pv_2\rangle$.

    Notice also that both of these matrices send $\langle v_1,v_2\rangle$ to $\langle v_1,v_2\rangle$, thus stabilizing $\z^2_{(p)}$.

    Now consider the subgroup $\Gamma$ of matrices in $\GL_2(\z_{(p)})$ who are congruent to the identity mod $p$. These matrices stabilize both the vertices $\langle e_1,e_2\rangle$ and $\langle e_1,pe_2\rangle$, since this group is the intersection of $\GL_2(\z_{(p)})$ and $\pi_p\GL_2(\z_{(p)})\pi_p^{-1}$, both of which are open, so this is an open subgroup of $\GL_2(\z_{(p)})$. Since there is $g\in\GL_2(\q)$ that maps the edge formed by $\langle v_1,v_2\rangle$ and $\langle pv_1,v_2\rangle$ to the edge formed by $\langle e_1,e_2\rangle$ and $\langle e_1,pe_2\rangle$, we have that  $g\Gamma g^{-1}$ will fix both vertices of the edge we are interested in and is also open.
    
    Then there exists matrices $U_a, U \in g\Gamma g^{-1}$ such that $(A^{-1}T_aA)U_a$ and $(A^{-1}SA)U$ are in $\SL_2(\z[1/p])$ and have the desired action on the edges.
\end{demonstracao}

\begin{corolario}
    The edge formed by the lattices associated with $\z^2_{(p)}$ and $\pi_p\z^2_{(p)}$ is a fundamental domain for the action of $\SL_2(\z[1/p])$, that is, the quotient graph $B_p/\SL_2(\z[1/p])$ is isomorphic to this edge.
\end{corolario}

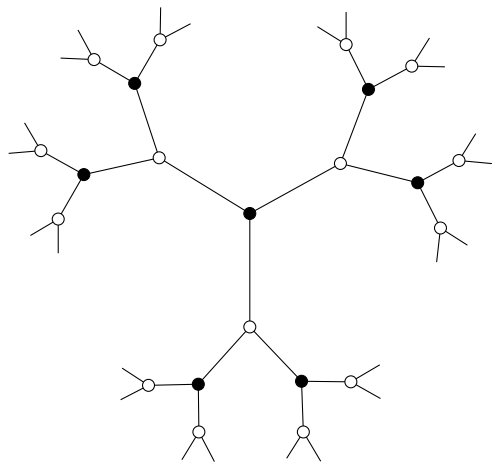
\begin{figure}[H]
    \centering
    \caption{Bruhat-Tits Building for $p=2$}
    \begin{tikzpicture}[scale=1.5]
        \draw (0,1)--(0,0);
        \draw (0,1)--(-0.799,1.491);
        \draw (0,1)--(0.799,1.441);
        \draw (-0.799,1.491)--(-1.015,2.147);
        \draw (-0.799,1.491)--(-1.463,1.343);
        \draw (0.799,1.441)--(1.046,2.094);
        \draw (0.799,1.441)--(1.48,1.271);
        \draw (0,0)--(-0.454,-0.505);
        \draw (0,0)--(0.456,-0.479);
        \draw (-0.454,-0.505)--(-0.449,-0.926);
        \draw (-0.454,-0.505)--(-0.893,-0.515);
        \draw (0.456,-0.479)--(0.472,-0.926);
        \draw (0.456,-0.479)--(0.891,-0.484);
        \draw (1.046,2.094)--(0.849,2.489);
        \draw (1.046,2.094)--(1.428,2.299);
        \draw (1.48,1.271)--(1.681,0.872);
        \draw (1.48,1.271)--(1.844,1.467);
        \draw (-1.015,2.147)--(-1.373,2.357);
        \draw (-1.015,2.147)--(-0.791,2.531);
        \draw (-1.463,1.343)--(-1.842,1.554);
        \draw (-1.463,1.343)--(-1.689,0.951);
        \draw (0.849,2.489)--(0.596,2.647);
        \draw (0.849,2.489)--(0.849,2.776);
        \draw (1.428,2.299)--(1.584,2.547);
        \draw (1.428,2.299)--(1.72,2.293);
        \draw (1.681,0.872)--(1.647,0.571);
        \draw (1.681,0.872)--(1.918,0.721);
        \draw (1.844,1.467)--(2.134,1.440);
        \draw (1.844,1.467)--(1.992,1.704);
        \draw (0.472,-0.926)--(0.325,-1.158);
        \draw (0.472,-0.926)--(0.625,-1.184);
        \draw (0.891,-0.484)--(1.152,-0.631);
        \draw (0.891,-0.484)--(1.146,-0.334);
        \draw (-0.449,-0.926)--(-0.599,-1.171);
        \draw (-0.449,-0.926)--(-0.312,-1.187);
        \draw (-0.893,-0.515)--(-1.139,-0.642);
        \draw (-0.893,-0.515)--(-1.126,-0.363);
        \draw (-1.689,0.951)--(-1.934,0.806);
        \draw (-1.689,0.951)--(-1.689,0.650);
        \draw (-1.842,1.554)--(-2.126,1.530);
        \draw (-1.842,1.554)--(-1.986,1.799);
        \draw (-1.373,2.357)--(-1.668,2.370);
        \draw (-1.373,2.357)--(-1.515,2.618);
        \draw (-0.791,2.531)--(-0.786,2.817);
        \draw (-0.791,2.531)--(-0.522,2.673);

        \draw [fill=white] (0,0) circle (1.5pt);
        \draw [fill=black] (0,1) circle (1.5pt);
        \draw [fill=white] (-0.799,1.491) circle (1.5pt);
        \draw [fill=white] (0.799,1.441) circle (1.5pt);
        \draw [fill=black] (-1.015,2.147) circle (1.5pt);
        \draw [fill=black] (-1.463,1.343) circle (1.5pt);
        \draw [fill=black] (-0.454,-0.505) circle (1.5pt);
        \draw [fill=black] (0.456,-0.479) circle (1.5pt);
        \draw [fill=black] (1.046,2.094) circle (1.5pt);
        \draw [fill=white] (0.849,2.489) circle (1.5pt);
        \draw [fill=white] (1.428,2.299) circle (1.5pt);
        \draw [fill=black] (1.48,1.271) circle (1.5pt);
        \draw [fill=white] (1.681,0.872) circle (1.5pt);
        \draw [fill=white] (1.844,1.467) circle (1.5pt);
        \draw [fill=white] (0.472,-0.926) circle (1.5pt);
        \draw [fill=white] (0.891,-0.484) circle (1.5pt);
        \draw [fill=white] (-0.449,-0.926) circle (1.5pt);
        \draw [fill=white] (-0.893,-0.515) circle (1.5pt);
        \draw [fill=white] (-1.689,0.951) circle (1.5pt);
        \draw [fill=white] (-1.842,1.554) circle (1.5pt);
        \draw [fill=white] (-1.373,2.357) circle (1.5pt);
        \draw [fill=white] (-0.791,2.531) circle (1.5pt);
    \end{tikzpicture}
    \label{fig:bp}
    \fautor
\end{figure}

Recall from Theorems \ref{dist} and \ref{neighbors} that $B_p$ is a tree (there's a well defined distance between vertices) such that every vertex has $p+1$ adjacent vertices. Figure \ref{fig:bp} illustrates what this tree looks like for the case $p=2$ and illustrates the action of $\SL_2(\z[1/p])$, where vertices of the same color (ones at even distance from each other) are equivalent under this action. 

Note that as in the case for the tree associated with $\Gamma_0$ and $\PGamma_0$, since vertices of $B_p$ are classes of lattices under homothety, the lattices $L$ and $-1\cdot L$ coincide, inducing a well defined action of $\PSL_2(\z[1/p]) = \SL_2(\z[1/p])$ on $B_p$. In the exact same way, we have that vertices will be $\PSL_2(\z[1/p])$-equivalent if and only if they are $\SL_2(\z[1/p])$-equivalent, and stabilizers under $\PSL_2(\z[1/p])$ correspond to the classes of stabilizers under $\SL_2(\z[1/p])$ mod $\pm I$. 

\begin{corolario}
    The edge formed by the lattices associated with $\z^2_{(p)}$ and $\pi_p\z^2_{(p)}$ is a fundamental domain for the action of $\PSL_2(\z[1/p])$, that is, the quotient graph $B_p/\PSL_2(\z[1/p])$ is isomorphic to this edge. The stabilizer of this edge is $\PGamma_0(p)$ and the stabilizer of the two vertices are $\PSL_2(\z),\pi_p\PSL_2(\z)\pi_p^{-1}$ respectively.
\end{corolario}

\section{Bruhat-Tits Building for a product of primes}

If $n\in\z$ with the prime decomposition of $n$ given by $n=p_1^{e_1}p_2^{e_2}\dots p_k^{e_k}$, then $\SL_2(\z[1/n]) $ coincides with $ \SL_2(\z[1/p_1p_2\dots p_k])$, so that it suffices to study $H_k(\SL_2(\z[1/n]))$ when $n$ is square free. We now consider the diagonal action of $\SL_2(\z[1/p_1\dots p_k])$ on the product $B = B_{p_1}\times\dots\times B_{p_k}$ for arbitrary $k$. Note that $B$ is contractible.

From Proposition \ref{prp: product cells}, $B_{p_1}\times\dots\times B_{p_k}$ has a CW-complex structure where the cells are products $(\sigma_1,\dots,\sigma_k)\coloneq \sigma_1\times\cdots\times\sigma_k$ with each $\sigma_i$ being a cell in $B_{p_i}$. To say that the action of $\SL_2(\z[1/p_1\dots p_k])$ on $B$ is diagonal is to say that a matrix $T \in \SL_2(\z[1/p_1\dots p_k])$ will act on a cell $(\sigma_1,\dots,\sigma_k)$ by $$T(\sigma_1,\dots,\sigma_k) = (T\sigma_1,\dots,T\sigma_k),$$ where $T\sigma_i$ represents the action of $T$ in the cell $\sigma_i$ of $B_{p_i}$. 

We denote by $\z^2_{(p_i)}$, $\pi_{p_i}\z^2_{(p_i)}$ and $e^0_{p_i}$ the 0-cells and 1-cell, respectively, of $B_{p_i}$ that compose the quotient graph $B_{p_i}\backslash \SL_2(\z[1/p_i])$.

\begin{observacao}\label{rem: stab k}
    Note that the stabilizers for $\z_{(p_k)}^2$, $\pi_p\z_{(p_k)}^2$ and $e_{p_k}^0$ under this action are $$\text{stab}(\z_{(p_k)}^2) = \SL_2(\z[1/p_1\cdots p_k]) \cap \SL_2(\z_{(p_k)}) = \SL_2(\z[1/p_1\cdots p_{k-1}]),$$ $$\text{stab}(\pi_{p_k}\z_{(p_k)}^2) = \pi_{p_k}\text{stab}(\z_{(p_k)}^2)\pi_{p_k}^{-1} = \pi_{p_k}\SL_2(\z[1/p_1\cdots p_{k-1}])\pi_{p_k}^{-1},$$ $$\text{stab}(e_{p_k}^0) = \Gamma_0(p_k,p_1\cdots p_{k-1}) = \left\{\begin{pmatrix}
    a & b \\ c & d 
\end{pmatrix}\in \SL_2(\z[1/p_1\cdots p_{k-1}])\mid p_k\mid c\right\}.$$
\end{observacao}

\begin{definicao}
    Let $\sigma$ be a cell in $B$ and let $\sigma_i$ be the $B_{p_i}$ component of $\sigma$. If for every $i\in\{1,\dots,k\}$ the component $\sigma_i$ is either $\z^2_{(p_i)}$, $\pi_{p_i}\z^2_{(p_i)}$ or $e^0_{p_i}$, the cell $\sigma$ will be called standard.
\end{definicao}

For the next proof, the set $\Gamma_0(p_{i_1}\dots p_{i_s})$ will be understood to be simply $\SL_2(\z)$ if $s=0$.

\begin{teorema} \label{stab gen}
    The stabilizer of a standard $s$-cell is given by $\pi\Gamma_0(p_{i_1}\dots p_{i_s})\pi^{-1}$ where $\pi$ is the product of matrices $\pi_q$ with $q \in \{p_1,\dots,p_k\}\backslash\{p_{i_1},\dots,p_{i_s}\}$. 
\end{teorema}

\begin{demonstracao}
    This will again be proven by induction. The case for $k=1$ has already been proven in Theorems \ref{teo: stab ep} and \ref{teo: stab zp}, since the stabilizers for $0$-cells are $\SL_2(\z)$ and $\pi_p\SL_2(\z)\pi_p^{-1}$ and the stabilizer for the single 1-cell is given by $\Gamma_0(p)$.
    
    Consider again the projection $\pi: B_{p_1}\times\dots\times B_{p_k} \rightarrow B_{p_1}\times\dots\times B_{p_{k-1}}$ which maps cells to cells and let $\sigma':=\pi(\sigma)$.

    Suppose $\sigma = (\sigma',\Lambda_{p_k})$ for some vertex $\Lambda_{p_k} \in B_{p_k}$. Then $\sigma'$ is an $s$-cell and by induction, its stabilizer is given by a conjugation $\pi\Gamma_0(p_{i_1}\dots p_{i_s})\pi^{-1}$ and the stabilizer of $\Lambda_{p_k}$ is given by $\SL_2(\z)$ or $\pi_{p_k}\SL_2(\z)\pi_{p_k}^{-1}$, depending on the class of $\Lambda_{p_k}$. Then the stabilizer of $\sigma$ is given by the two possible intersections, which will be either $\pi\Gamma_0(p_{i_1}\dots p_{i_s})\pi^{-1}$ or $(\pi_{p_k}\pi)\Gamma_0(p_{i_1}\dots p_{i_s})(\pi_{p_k}\pi)^{-1}$, as expected.

    On the other hand, if $\sigma = (\sigma',e_{p_k})$, then $\sigma'$ is an $s-1$ cell and its stabilizer is given by $\pi\Gamma_0(p_{i_1}\dots p_{i_{s-1}})\pi^{-1}$ (where $p_{i_j}\neq p_k$) with the stabilizer of $e_{p_k}$ given by $\Gamma_0(p_{k})$. Then the stabilizer of $\sigma$ is the intersection $\pi\Gamma_0(p_{i_1}\dots p_{i_{s-1}}p_k)\pi^{-1}$, as expected.
\end{demonstracao}

\begin{corolario}
    The stabilizer of a standard $s$-cell is isomorphic to $\Gamma_0(p_{i_1}\dots p_{i_s})$ where each $p_{i_j}\in\{p_1,\dots,p_k\}$.
\end{corolario}

\begin{lema}
    Let $\{i_1,\dots,i_s\}\subset\{p_1,\dots,p_{k-1}\}$. Then $p_{i_1}\dots p_{i_s}\z[1/p_k]$ is dense in $\z[1/p_k]$ considering the $p$-adic topology of $\GL_2(\q)$, and so, $\Gamma_0(p_{i_1}\dots p_{i_s}, p_k)$ is dense in $\SL_2(\z[1/p_k])$.
\end{lema}

\begin{demonstracao}
    This proof is almost exactly like the one used to conclude that $\Gamma_0(q,p)$ is dense in $\SL_2(\z[1/p])$ in the $p$-adic topology if $p,q$ are distinct primes (see Theorem \ref{teo: closure gammapq}). 
    
    Let $x\in\z[1/p_k]$ and $n\geq 0$. There exists $m\geq 0$ such that $x':= p_k^mx \in \z$ ($m$ can be taken as $\nu_{p_k}(x)$). Since $(p_{i_1}\cdots p_{i_s},p_k)=1$, we have that the equation $p_{i_1}\cdots p_{i_s}y_n'\equiv x' \mod p_k^{n+m}$ has a solution $y_n'\in\z$, that is, $p_{i_1}\cdots p_{i_s}y_n' - x' = ap_k^{n+m}$ with $a\in \z$. Let $y_n:= p_k^{-m}y \in \z[1/p_k]$. Then $p_{i_1}\cdots p_{i_s}y_n-x = ap_k^n$, that is, $\nu_{p_k}(p_{i_1}\cdots p_{i_s}y_n-x) \geq n$ and so, $|p_{i_1}\cdots p_{i_s}y-x|_{p_k}\leq p^{-n}$. In other words, given any $x\in\z[1/p_k]$ and any open ball around $x$, there exists $p_{i_1}\cdots p_{i_s}y_n \in p_{i_1}\cdots p_{i_s}y_n\z[1/p_k]$ whithin said open ball. Thus, $p_{i_1}\cdots p_{i_s}\z[1/p_k]$ is dense in $\z[1/p_k]$.  

    From Theorem \ref{teo: GE2}, $\SL_2(\z[1/p_k])$ is generated by elementary matrices $$\begin{pmatrix}
        1 & a \\ 0 & 1
    \end{pmatrix}, \begin{pmatrix}
        1 & 0 \\ a & 1
    \end{pmatrix}$$ where $a\in\z[1/p_k]$. Matrices of the first type are already in $\Gamma_0(p_{i_1}\cdots p_{i_s},p_k)$. Since any $a\in\z[1/p_k]$ can be approximated by $p_{i_1}\cdots p_{i_s}\z[1/p_k]$, we have that matrices of the second type can be approximated by matrices in $\Gamma_0(p_{i_1}\cdots p_{i_s},p_k)$. 
\end{demonstracao}

\begin{teorema}\label{teo: standard cells}
    Let $(\sigma,\sigma_k), (\sigma',\sigma'_k)$ be a pair of cells with $\sigma,\sigma'\in B_{p_1}\times\dots\times B_{p_{k-1}}$ and $\sigma_k,\sigma_k'\in B_{p_k}$ such that $\sigma$ and $\sigma'$ are $\SL_2(\z[1/p_1\cdots p_{k-1}])$-equivalent and $\sigma_k,\sigma_k'$ are $\SL_2(\z[1/p_k])$-equivalent. Suppose also that $(\sigma',\sigma'_k)$ is a standard cell. Then there is $T\in \SL_2(\z[1/p_1\cdots p_k])$ that maps $(\sigma,\sigma_k)$ to $(\sigma',\sigma'_k)$.
\end{teorema}

\begin{demonstracao}
     We prove this by induction on $k$. If $k=1$, the statement is trivial.

     Suppose $\sigma_k$ is a 0-cell, so that there exists $R\in\SL_2(\z[1/p_k])$ mapping $\sigma_k$ to $\sigma_k'$.

     If $\sigma_k' = \z^2_{(p_k)}$, then $\text{stab}(\sigma_k') = \SL_2(\z[1/p_1\cdots p_k])$. Since $\sigma,\sigma'$ are equivalent, we have that $R\sigma$ and $\sigma'$ are also equivalent, so there is $S \in \SL_2(\z[1/p_1\cdots p_k]$ which maps $R\sigma$ to $\sigma'$ and stabilizes $\sigma_k'$. Thus $$SR(\sigma,\sigma_k) = S(R\sigma,\sigma_k') = (\sigma,\sigma_k').$$

     If $\sigma_k' = \pi_{p_k}\z^2_{(p_k)}$, then $\text{stab}(\sigma_k') = \pi_{p_k}\SL_2(\z[1/p_1\cdots p_k])\pi_{p_k}^{-1}$. Again, since $\sigma,\sigma'$ are equivalent, we have that $\pi^{-1}_{p_k}R\sigma$ and $\pi^{-1}_{p_k}\sigma'$ are also equivalent, so there is $S\in \SL_2(\z[1/p_1\dots p_k])$ mapping $\pi^{-1}_{p_k}R\sigma$ to $\pi^{-1}_{p_k}\sigma'$, and so, we have that $\pi_{p_k}S\pi^{-1}_{p_k}R\sigma = \sigma'$. Thus $$\pi_{p_k}S\pi^{-1}_{p_k}R(\sigma,\sigma_k) = \pi_{p_k}S\pi^{-1}_{p_k}(R\sigma,\sigma_k') = (\sigma',\sigma_k').$$

     Finally, we consider the case where $\sigma_k'=e^0_{p_k}$.

     As mentioned in Remark \ref{rem: stab k}, the stabilizer of $\sigma'$ under the action of the group $\SL_2(\z[1/p_1\cdots p_{k-1}])$ is of the type $\pi\Gamma_0(p_{i_1}\cdots p_{i_s})\pi^{-1}$.

     Since $\sigma\sim\sigma'$, there is $R\in\SL_2(\z[1/p_1\cdots p_{k-1}])$ such that $R\sigma = \sigma'$. From Theorem \ref{teo: gamma0 open}, $\Gamma_0(p_k,\z_{(p_k)})$ is open in in the $p_k$-adic topology of $\GL_2(\q)$ and it stabilizes the vertices of the standard edge $\z_{(p_k)}$. From transitivity of the action of $\GL_2(\q)$,  for any edge in $B_{p_k}$, there is a subgroup $\Gamma$ that is a conjugation of the subgroup of $\Gamma_0(p_k,\z_{(p_k)})$ and fixes both vertices of the edge of interest. This subgroup is therefore an open set in the $p_k$-adic topology of $\GL_2(\q)$.

     Since edges of $B_{p_k}$ are all $\SL_2(\z[1/p_k])$-equivalent, let $S\in\SL_2(\z[1/p_k])$ such that $S(\pi^{-1}R\sigma_k) = \pi^{-1}\sigma_k'$ and let $\Gamma$ be the open set fixing both vertices of $\pi^{-1}R\sigma_k$. Then $S\cdot\Gamma$ is an open neighborhood of $S$ in the $p_k$-adic topology.

     Since $\Gamma_0(p_{i_1}\cdots p_{i_s},p_k)$ is dense in $\SL_2(\z[1/p_k])$, we have that $$S\cdot\Gamma \cap \Gamma_0(p_{i_1}\cdots p_{i_s},p_k) \neq \emptyset,$$ that is, there exists $U\in \Gamma$ such that $SU\in \Gamma_0(p_{i_1}\cdots p_{i_s},p_k)$. Note that this implies $\det(U)=1$ so $U$ not only stabilizes the vertices of $\pi^{-1}R\sigma_k$ but stabilizes the edge itself.

     Then $SU\pi^{-1}R\sigma_k = S\pi^{-1}R\sigma_k = \pi^{-1}\sigma_k'$ implies that $\pi SU\pi^{-1}(R\sigma_k) = \sigma_k'$ and since $SU\in \Gamma_0(p_{i_1}\cdots p_{i_s},p_k)$, we have that $\pi SU\pi^{-1} \in \Gamma$. Thus, $$\pi SU\pi^{-1} R(\sigma,\sigma_k) = \pi SU\pi^{-1}(\sigma',R\sigma_k) = (\sigma',\sigma_k').$$
\end{demonstracao}

\begin{corolario}
    A pair $(\sigma,\sigma_k),(\sigma',\sigma_k')$ with $\sigma,\sigma'\in B_{p_1}\times\dots\times B_{p_{k-1}}$ and  $\sigma_k,\sigma'_k\in B_{p_k}$ is $\SL_2(\z[1/p_1\cdots p_k)$-equivalent if and only if $\sigma$ and $\sigma'$ are $\SL_2(\z[1/p_1\cdots p_{k-1}])$-equivalent and $\sigma_k,\sigma_k'$ are $\SL_2(\z[1/p_k])$ equivalent.
\end{corolario}

\begin{demonstracao}
    If $\sigma$ and $\sigma'$ are $\SL_2(\z[1/p_1\cdots p_{k-1}])$-equivalent and $\sigma_k,\sigma_k'$ are $\SL_2(\z[1/p_k])$-equivalent, there are standard cells $\sigma'',\sigma_k''$ such that $\sigma''$ is equivalent to both $\sigma$ and $\sigma'$ and $\sigma_k''$ is equivalent to both $\sigma_k$ and $\sigma_k$. By the previous theorem, we have that $(\sigma'',\sigma_k'')$ is $\SL_2(\z[1/p_1\cdots p_{k}])$-equivalent to both $(\sigma,\sigma_k)$ and $(\sigma',\sigma_k')$, so they are all equivalent amongst themselves.

    On the other hand, if $(\sigma,\sigma_k)$ is $\SL_2(\z[1/p_1\cdots p_{k}])$-equivalent to $(\sigma',\sigma_k')$, there exists $T\in\SL_2(\z[1/p_1\cdots p_{k}])$ such that $T(\sigma,\sigma_k)=(\sigma',\sigma_k')$.
    
    Write $\sigma = \displaystyle\prod_{i=1}^{k-1} \sigma_i$, $\sigma' = \displaystyle\prod_{i=1}^{k-1} \sigma_i'$.

    For every $i\in\{1,\dots,k\}$, we have that $T\sigma_i=\sigma_i'$. Since $\text{stab}_{\GL_2(\q)}(\sigma_i)$ is open in the $p_i$-adic topology and $\SL_2(\z[1/p_i])$ is dense in $\SL_2(\q)$, there exists $U\in \text{stab}_{\GL_2(\q)}(\sigma_i)$ such that $TU\in \SL_2(\z[1/p_i])$. Then $TU\sigma_i = T\sigma_i = \sigma_i'$. In particular, $(TU)\sigma_k = \sigma_k'$.

    From the ``if'' part of the statement, since the cells $\sigma_i$ with $1\leq 1 < k$ are all $\SL_2(\z[1/p_i])$-equivalent, it follows that $\sigma$ and $\sigma'$ are $\SL_2(\z[1/p_1\dots p_{k-1}])$-equivalent.
\end{demonstracao}

\begin{teorema}\label{classes}
    There are $2^{k-s}$ equivalence classes of $s$-cells of $B$ under the action of the group $\SL_2(\z[1/p_1\cdots p_k])$.
\end{teorema}

\begin{demonstracao}
    We prove this by induction on $k$. If $k=1,2$, the statement holds as shown in Theorems \ref{teo: numb vert 1} and \ref{teo: num edges 1}.

    If $\sigma$ is an $s$-cell, then $\sigma$ is the product of $s$ 1-cells $e_{p_i} \in B_{p_i}$ and $k-s$ 0-cells $\Lambda_{p_i}$ with $p_i\in\{p_1,\dots,_k\}$.

    Consider the projection $\pi: B_{p_1}\times\dots\times B_{p_k} \rightarrow B_{p_1}\times\dots\times B_{p_{k-1}}$ which maps cells to cells and let $\sigma':=\pi(\sigma)$. Then $\sigma = (\sigma',e_{p_k})$ or $\sigma = (\sigma',\Lambda_{p_k})$ with $e_{p_k},\Lambda_{p_k}$ an edge and a vertex in $B_{p_k}$, depending on wether $\sigma'$ is a $s-1$ or $s$ cell.

    If $\sigma = (\sigma',e_{p_k})$, then $\sigma'$ is an $s-1$ cell in $B_{p_1}\times\dots\times B_{p_{k-1}}$ and there are $2^{k-1-(s-1)} = 2^{k-s}$ equivalence classes $\sigma'$ could belong to. Since $e_{p_k}$ can be mapped to  $e_{p_k}^0$ while fixing $\sigma'$, we conclude that there are $2^{k-s}$ classes for $\sigma$, depending exclusively on the class for $\sigma'$.

    On the other hand, if $\sigma = (\sigma',\Lambda_{p_k})$, then $\sigma'$ is an $s$-cell in and there are $2^{k-1-s}$ classes it could belong to. Since $\Lambda_{p_k}$ can be mapped to either $\z_{(p_k)}^2$ or $\pi_{p_k}\z_{(p_k)}^2$ while fixing $\sigma'$, there are $2\cdot 2^{k-1-s} = 2^{k-s}$ classes $\sigma$ can belong to.
\end{demonstracao}

\begin{teorema}
    The quotient graph $Y:=B\backslash\SL_2(\z[1/p_1\dots p_k])$ is a $k$-dimensional hypercube.
\end{teorema}

\begin{demonstracao}
    We prove this by induction. The case $k=1$ has already been proven in the last section.

    Let $B'=B_{p_1}\times\dots\times B_{p_{k-1}}$, so that $B = B'\times B_{p_k}$.

    As in the previous theorem, any cell in $B$ representing an equivalence class in $Y$ can be written as $(\sigma,\sigma_k)$ where $\sigma,\sigma_k$ are cells in $B'$ and $B_{p_k}$, respectively. Theorem \ref{teo: standard cells} then shows that $\sigma,\sigma_k$ can be taken as standard cells, so $\sigma_k$ is either $e^0_{p_k}$ or its boundary 0-cells, $\z^2_{(p_k)}, \pi_p\z^2_{(p_k)}$.

    Thus, for every $s$-cell $\sigma\in B'$, we have two copies of $\sigma$ in the $s$-cells $(\sigma,\z^2_{(p_k)})$ and $(\sigma,\pi_{p_k}\z^2_{(p_k)})$, joined by the $s+1$-cell $(\sigma,e^0_{p_k})$. Since $B'$ is a $(k-1)$-hypercube, this implies that $B$ is a $k$-hypercube.
\end{demonstracao}

\begin{corolario}
    The quotient graph $Y$ is contractible.
\end{corolario}

\section{The main spectral sequence}

In this section, we use the action of $\SL_2(\z[1/p_1\cdots p_k])$ on the contractible CW-complex $B_{p_1}\times\cdots\times B_{p_k}$ described in the last section to construct a spectral sequence converging to $H_k(\SL_2(\z[1/p_1\cdots p_k]))$.

\begin{teorema}\label{teo: spec seq sl}
    Let $n$ be an integer and let $p_1,\dots,p_k$ be the distinct prime factors of $n$. There is a spectral sequence converging to $H_*(\SL_2(\z[1/n]))$ whose $E^1$-term is given by $$E_{s,t}^1 = \bigoplus H_t(\Gamma_0(p_{i_1}\cdots p_{i_s}))^{2^{k-s}},$$ where the summation is over all $s$-element subsets of $\{p_1,\dots,p_k\}$. 
\end{teorema}

\begin{demonstracao}
    This is the spectral sequence described in Theorem \ref{seq} associated to the action of $\SL_2(\z[1/n])$, which coincides with  
    $ \SL_2(\z[1/p_1\cdots p_k])$, on the contractible CW-complex $B = B_{p_1}\times\dots\times B_{p_k}$. Since $B$ is finite dimensional and $\SL_2(\z[1/n])$ fixes cells pointwise, the differential $d^1$ is as described in Theorem \ref{teo: maps}.

    Recall that for such a spectral sequence, the term $E^1_{s,t}$ is given by $\displaystyle\bigoplus_{\sigma \in \Sigma_s}H^t(\Gamma_\sigma)$ where $\Sigma_s$ is a system of representatives for the orbits of $s$-cells under the action and $\Gamma_\sigma$ is the stabilizer of the cell $\sigma$.

    Theorem \ref{stab gen} shows that the stabilizer of an $s$-cell is given by $\pi\Gamma_0(p_{i_1}\cdots p_{i_s})\pi^{-1}$ where $\pi$ is a product of matrices $\pi_q$ where $q\in\{p_1,\dots,p_k\}\backslash\{p_{i_1},\dots,p_{i_s}\}$, and thus, isomorphic to $\Gamma_0(p_{i_1}\cdots p_{i_s})$. Then, Theorem \ref{classes} gives that there are $2^{k-s}$ equivalence classes of $s$-cells.

    Example \ref{exe: two primes} below illustrates this sequence for the case where $n$ has two prime factors.
\end{demonstracao}

\begin{exemplo}\label{exe: two primes}
    Let $p,q$ be two distinct primes. Then the quotient $(B_p\times B_q)/\SL_2(\z[1/pq])$ is a 2-hypercube, which is a square. 

    Denote by $\Lambda_p$ and $\Lambda_p'$ the 0-cells associated to $\z^2_{(p)}$ and $\pi_p\z^2_{(p)}$. Figure \ref{fig:quotient} below shows the quotient graph on the left, and the stabilizers of the corresponding cells on the right.

    \begin{figure}[H]
        \centering
        \caption{Quotient graph of $B_p\times B_q/\SL_2(\z[1/pq])$}
        \begin{center}
    \begin{minipage}{0.4\textwidth} \begin{center}
            \begin{tikzpicture}[scale=0.8]
            \draw[blue, ultra thin, fill, opacity=0.1] (-2,-2) -- (2,-2) -- (2,2) -- (-2,2) -- cycle;
            \node[scale=1.1, blue] at (0,0) {$e^0_p\times e^0_q$};
            
            \draw[purple] (-2,-2)--(2,-2) [arrow inside={end=stealth,opt={purple,scale=3}}{0.55}];
            \node[scale=1.1, purple] at (0,-2.5) {$\Lambda_p'\times e^0_q$};
            \draw[purple] (2,-2)--(2,2) [arrow inside={end=stealth,opt={purple,scale=3}}{0.55}];
            \node[purple, scale=1.1] at (3.5,0) {$e^0_p\times \Lambda_q'$};
            \draw[purple] (2,2)--(-2,2) [arrow inside={end=stealth,opt={purple,scale=3}}{0.55}];
            \node[purple,scale=1.1] at (0,2.5) {$\Lambda_p\times e^0_q$};
            \draw[purple] (-2,2)--(-2,-2) [arrow inside={end=stealth,opt={purple,scale=3}}{0.55}];
            \node[purple,scale=1.1] at (-3.5,0) {$e^0_p\times\Lambda_q$};

            \draw [fill=black] (-2,-2) circle (2pt);
            \node[scale=1.1] at (-2.5,-2.5) {$\Lambda_p'\times\Lambda_q$};
            \draw [fill=black] (-2,2) circle (2pt);
            \node[scale=1.1] at (-2.5,2.5) {$\Lambda_p\times\Lambda_q$};
            \draw [fill=black] (2,-2) circle (2pt);
            \node[scale=1.1] at (2.5,-2.5) {$\Lambda_p'\times\Lambda_q'$};
            \draw [fill=black] (2,2) circle (2pt);
            \node[scale=1.1] at (2.5,2.5) {$\Lambda_p\times\Lambda_q'$};
            \end{tikzpicture}
        \end{center} \end{minipage} \hspace{40pt}\begin{minipage}{0.4\textwidth} \begin{center}\vspace{5pt}
            \begin{tikzpicture}[scale=0.8]
            \draw[blue, ultra thin, fill, opacity=0.1] (-2,-2) -- (2,-2) -- (2,2) -- (-2,2) -- cycle;
            \node[scale=1.1, blue] at (0,0) {$\Gamma_0(pq)$};
            
            \draw[purple] (-2,-2)--(2,-2) [arrow inside={end=stealth,opt={purple,scale=3}}{0.55}];
            \node[scale=1.1, purple] at (0,-2.8) {$\Gamma_0(q)^p$};
            \draw[purple] (2,-2)--(2,2) [arrow inside={end=stealth,opt={purple,scale=3}}{0.55}];
            \node[purple, scale=1.1] at (3.2,0) {$\Gamma_0(p)^q$};
            \draw[purple] (2,2)--(-2,2) [arrow inside={end=stealth,opt={purple,scale=3}}{0.55}];
            \node[purple,scale=1.1] at (0,2.5) {$\Gamma_0(q)$};
            \draw[purple] (-2,2)--(-2,-2) [arrow inside={end=stealth,opt={purple,scale=3}}{0.55}];
            \node[purple,scale=1.1] at (-3.2,0) {$\Gamma_0(p)$};

            \draw [fill=black] (-2,-2) circle (2pt);
            \node[scale=1.1] at (-2.5,-2.5) {$\Gamma_0^p$};
            \draw [fill=black] (-2,2) circle (2pt);
            \node[scale=1.1] at (-2.5,2.5) {$\Gamma_0$};
            \draw [fill=black] (2,-2) circle (2pt);
            \node[scale=1.1] at (2.5,-2.5) {$\Gamma_0^{pq}$};
            \draw [fill=black] (2,2) circle (2pt);
            \node[scale=1.1] at (2.5,2.5) {$\Gamma_0^q$};
            \end{tikzpicture}
        \end{center} \end{minipage}
\end{center}
        \label{fig:quotient}
        \fautor
    \end{figure}
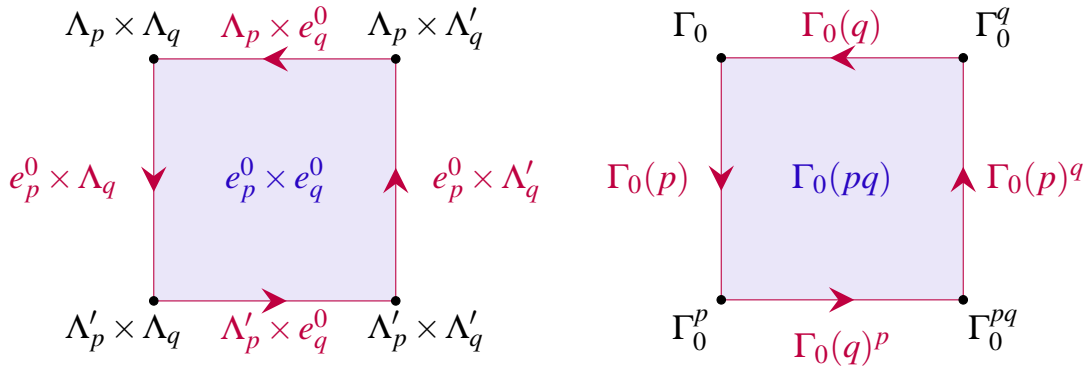

Then the first page of the associated spectral sequence is given by 

\begin{center}
    \begin{tikzcd}
        \vdots & \vdots & \vdots & \vdots  \\
        H_2(\Gamma_0)^{4} & H_2(\Gamma_0(p))^{2}\oplus H_2(\Gamma_0(q))^{2} \arrow[l]& H^2(\Gamma_0(pq)) \arrow[l]& 0  \\
        H_1(\Gamma_0)^{4} & H_1(\Gamma_0(p))^{2}\oplus H_1(\Gamma_0(q))^{2} \arrow[l]& H^1(\Gamma_0(pq)) \arrow[l]& 0  \\
        H_0(\Gamma_0)^{4} & H_0(\Gamma_0(p))^{2}\oplus H_0(\Gamma_0(q))^{2} \arrow[l]& H^0(\Gamma_0(pq))\arrow[l]& 0  
    \end{tikzcd}
\end{center} \vspace{10pt}

Note that $H_k(\Gamma_0)\!=\!H_k(\Gamma_0(1))$ is torsion for all $k\geq 1$ and the homology groups $H_k(\Gamma_0(p))$, $H_k(\Gamma_0(q))$, $H_k(\Gamma_0(pq))$ are all torsion for $k\geq 2$ (see Theorem \ref{H_n}). From the Proposition \ref{prp: line 0}, $E^2_{0,0}\simeq \z$ and $E^2_{s,0}=0$. 

This implies that only $d^1_{2,1}: E^1_{2,1} \rightarrow E^1_{1,1}$ survives tensoring with $\q$. This map $d^1_{2,1}$ is given by $$ H_k(\Gamma_0(pq)) \rightarrow H_k(\Gamma_0(p))^2\oplus H_k(\Gamma_0(q))^2$$ $$[\alpha] \longmapsto (i_{1,*}[\alpha],i_{2,*}[\alpha],j_{1,*}[\alpha],j_{2,*}[\alpha])$$ where the maps are induced by the inclusions $i_1: \Gamma_0(pq)\hookrightarrow \Gamma_0(p), i_2: \Gamma_0(pq)\hookrightarrow \Gamma_0(p)^q, j_1: \Gamma_0(pq)\hookrightarrow \Gamma_0(q), j_2: \Gamma_0(pq)\hookrightarrow \Gamma_0(q)^p$.
\end{exemplo}

\begin{proposicao}\label{prp: line 0}
    Let $E^r_{s,t}$ be the spectral sequence from Theorem \ref{teo: spec seq sl}. Then $E^2_{0,0} \simeq \z$ and $E^2_{s,0} = 0$ otherwise.
\end{proposicao}

\begin{demonstracao}
    Note that $H_0(\Gamma,\z)=\z$ for any group $\Gamma$, so that $E^1_{s,0} = \z^{2^{k-s}}$ and the differential $d^1$ is induced by the boundary map on the quotient graph, implying that the complex $E^1_{*,0}$ is simply the cellular complex of the quotient graph $Y$, so that $E^2_{s,0}$ coincides with $H_s(Y,\z)$. Since the quotient graph $Y$ is a hypercube, thus contractible, we have that $H_0(Y,\z) = \z$ and $H_s(Y,\z)=0$ otherwise, giving the desired result.
\end{demonstracao}

\begin{teorema}
    Let $n$ be an integer and let $p_1,\dots,p_k$ be the distinct prime factors of $n$. Then there is a spectral sequence converging to $H_*(\PSL_2(\z[1/n]))$ whose $E^1$-term is given by $$E_{s,t}^1 = \bigoplus H_t(\PGamma_0(p_{i_1}\cdots p_{i_s}))^{2^{k-s}},$$ where the summation is over all $s$-element subsets of $\{p_1,\dots,p_k\}$. 
\end{teorema}

\begin{demonstracao}
    This follows again from the fact that two $s$-cells $B_{p_1}\times\dots\times B_{p_k}$ are $\SL_2(\z[1/n])$-equivalent if and only if they are $\PSL_2(\z[1/n])$-equivalent, so both actions have the same orbits and the stabilizers under $\PSL_2(\z[1/n])$ are given by $\Gamma_0(p_{i_1}\cdots p_{i_s})/\{\pm I\} = \PGamma_0(p_{i_1}\cdots p_{i_s})$. 
\end{demonstracao}

\begin{teorema}\label{teo: fint gen}
    The homology groups $H_k(\SL_2(\z[1/n]))$ and $H_k(\PSL_2(\z[1/n]))$ are finitely generated for all integers $n$.
\end{teorema}

\begin{demonstracao}
    From Theorems \ref{H_n} and \ref{H_nP}, the homology groups $H_t(\Gamma_0(m))$ and $H_t(\PGamma_0(m))$ are finitely generated for all integers $t,m$. This implies that the first page of the spectral sequences converging to $H_k(\SL_2(\z[1/n]))$ and $H_k(\PSL_2(\z[1/n]))$ consist only of finitely generated abelian groups. Furthermore, since subgroups and quotients of finitely generated groups are also finitely generated, we have that every page of our two spectral sequences will consist of finitely generated abelian groups. Therefore, the stable terms $E^\infty_{p,q}$ are finitely generated.

    Let $G$ be either $\SL_2(\z[1/n])$ or $\PSL_2(\z[1/n])$ and let $E^1_{p,q}$ be the respective spectral sequence converging to the homology of $G$. From Theorem \ref{convergence}, we have that for all $k\geq 0$, there exists a filtration $\Phi$ of $H_k(G)$ such that $E^\infty_{p,q}\simeq \Phi^p H_{p+q}(G) / \Phi^{p-1}H_{p+q}(G)$, giving the group extension $$0 \rightarrow \Phi^{p-1} H_{p+q}(G) \rightarrow \Phi^{p} H_{p+q}(G) \rightarrow E^{\infty}_{p,q}\rightarrow 0.$$ Since each term $E^\infty_{p,q}$ is finitely generated and $\Phi^{-1} H_{p+q}(G)=0$ is also finitely generated, induction gives that all filtrations $\Phi^{s}H_{p+q}(G)$ are finitely generated. In particular, $\Phi^{p+q}H_{p+q}(G) = H_{p+q}(G)$ is finitely generated.
\end{demonstracao}

\begin{teorema}\label{teo: sequence sl}
    Let $E^r_{s,t}$ be the spectral sequence converging to $\SL_2(\z[1/n])$. Then the rational homology $H_*(\SL_2(\z[1/n]),\q)$ is given by the term $E^2_{s-1,1}\otimes_{\z} \q$.
\end{teorema}

\begin{demonstracao}
    From Theorem \ref{teo: rational homology}, $H_k(G,\q) \simeq H_k(G,\z)\otimes\q$ for all groups $G$. Then, the spectral sequence defined by $E^r_{s,t}\otimes_{\z}\q$ converges to $H_*(\SL_2(\z[1/n]),\z)\otimes_{\z}\q \simeq H_*(\SL_2(\z[1/n]),\q)$.

    Furthermore, the terms $E^1_{s,t}$ are all given by sums of terms of the type $H_t(\Gamma_0(m))$ for some $m$, which as shown in Theorem \ref{H_n}, are torsion if $t>1$. Thus, $E^1_{s,t}\neq 0$ only when $t=0,1$. Furthermore, the Proposition \ref{prp: line 0} gives that $E^2_{0,0}=\z$ and $E^2_{s,0}=0$ otherwise, thus, the only non-zero terms $E^2_{s,t}$ happen when $t=1$ or $(s,t)=(0,0)$. Theorem \ref{H_n} also shows that $H_t(\SL_2(\z))$ is torsion if $t>0$ (recall that $\SL_2(\z) = \Gamma_0(1)$), implying that $E^2_{0,t} = 0$ for $t>0$. 

    Thus, we have that $E^2_{s,t} \otimes_{\z} \q$ is given by 

    \begin{center}
        \begin{tikzcd}
            \vdots & \vdots & \vdots & \vdots & \vdots  \\
            0  & 0  & 0 & 0  & \dots \\
            0   & \displaystyle\coker(d^1_{2,1}\otimes\id) & \displaystyle\frac{\ker(d^1_{2,1}\otimes\id)}{\im(d^1_{3,1}\otimes\id)} \arrow[llu]& \displaystyle\frac{\ker(d^1_{3,1}\otimes\id)}{\im(d^1_{4,1}\otimes\id)} \arrow[llu]  & \dots  \\
            \q & 0 & 0 \arrow[llu] & 0 \arrow[llu] & \dots
        \end{tikzcd}
    \end{center} and the sequence collapses on the second page. Since each diagonal $s+t=k$ only contains a single non-zero element given by $(k-1,1)$ (for $k>1$), we have that the induced filtration on $H_k(\SL_2(\z[1/n]),\q)$ has one element for each $n$ and so, $H_k(\SL_2(\z[1/n]),\q) \simeq E^2_{k-1,1}\otimes_{\z}\q$.
\end{demonstracao}

\begin{corolario}\label{cor: ranks PSL and SL}
    The ranks of $H_k(\PSL_2(\z[1/n]))$ and $H_k(\SL_2(\z[1/n]))$ coincide.
\end{corolario}

\begin{demonstracao}
    From our characterizations of $H_k(\Gamma_0(n))$ and $H_k(\PGamma_0(n))$, we have that these ranks are both 1 if $k=0$, both $r(n)$ if $k=1$ and 0 otherwise.

    Denote by $PE^r_{s,t}$ the spectral sequence converging to $H_*(\PSL_2(\z[1/n]))$. Then we have the analogous result that $H_k(\PSL_2(\z[1/n])) \simeq PE^2_{k-1,1}\otimes_{\z}\q$. Since the modules $PE^r_{s,t}$ and $E^r_{s,t}$ have the same rank, they will coincide once tensored with $\q$. Thus \begin{align*}
        \rank(H_k(\PSL_2(\z[1/n]))) &= \dim_q(PE^2_{k-1,1}\otimes_{\z}\q) \\ &= \dim_q(E^2_{k-1,1}\otimes_{\z}\q) \\ &= \rank(H_k(\SL_2(\z[1/n]))).
    \end{align*}
\end{demonstracao}

As a consequence of Theorem \ref{teo: sequence sl}, we have:

\begin{corolario}\label{cor: H1 is 0}
    $H_1(\SL_2(\z[1/n]),\q) = 0$ for all integers $n$.
\end{corolario}

This shows that for all non-zero integers $n$, the group $H_1(\SL_2(\z[1/n]))$ is torsion. The structure of these groups will be calculated in the next section.

\begin{corolario}\label{cor: Hk for prime}
    If $n=p$ is a prime, then $$\rank(H_k(\SL_2(\z[1/p]))) = \begin{cases}
        1, & \text{ if }k=0 \\
        r(p), & \text{ if }k=2 \\
        0, & \text{ otherwise}.
    \end{cases}$$
\end{corolario}

\begin{demonstracao}
    If $n$ has a single prime factor, the spectral sequence has only two non-zero columns and we obtain $E^2_{0,0}\otimes_{\z}\q\simeq\q$, $E^2_{1,1}\otimes_{\z}\q = \ker(d^1_{2,1}\otimes\id)$ and $E^2_{s,t}=0$ otherwise. Furthermore, $d^1_{2,1}$ is a map from $H_1(\Gamma_0(p))$ to two copies of $H_1(\SL_2(\z))\simeq \z/12$, so $d^1_{2,1}\otimes\id = 0$ and thus $E^2_{1,1}$ has rank $r(p)$.
\end{demonstracao}

\begin{observacao}
    From Definition \ref{def: r(n)} and Proposition \ref{prp: numbers}, we obtain the following description of $r(p)$ for $p$ prime: $$r(p) = \begin{cases}
        1, & \text{ if } p=2,3\\
        \frac{p+1}{6}, &\text{ if } p\equiv 5\pmod{12}\\
        \frac{p-1}{6}, &\text{ if }p\equiv 7\pmod{12}\\
        \frac{p+7}{6}, &\text{ if }p\equiv 11\pmod{12}\\
        \frac{p-7}{6}, &\text{ if }p\equiv 1\pmod{12}
    \end{cases}.$$ In particular, $\rank(H_2(\SL_2(\z[1/p]))) = 1$ if and only if $p$ is one of the primes $2,3,5,7,13$.
\end{observacao}

\begin{corolario}\label{cor: upper bound}
    Let $n = p_1p_2\cdots p_k$. Then $H_s(\SL_2(\z[1/n]),\q) = 0$ for all $s>k+1$, and for $2\leq s \leq k+1$, $$\rank(H_s(\SL_2(\z[1/n]))) \leq 2^{k-s+1}\sum r(p_{i_1}\cdots p_{i_{s-1}})$$ where the sum is defined over all $(s-1)$-element subsets of $\{p_1,\dots,p_k\}$ and $r(m)$ is the rank of $H_1(\Gamma_0(m))$ as in Theorem \ref{H_n}.
\end{corolario}

\begin{demonstracao}
    From Theorem \ref{teo: sequence sl}, $H_s(\SL_2(\z[1/n]),\q) \simeq E^2_{s-1,1}\otimes_{\z}\q$. Since $E^1_{p,q}=0$ for $p>k$ (which follows from the fact that there are no $t$-element subsets of $\{p_1,\dots,p_k\}$, or equivalently, from the fact that the complex giving this action has cells only in dimensions less than $k$), we have that for $s>k+1$, $H_s(\SL_2(\z[1/n]),\q) \simeq E^2_{s-1,1}\otimes_{\z}\q = 0$. 
    
    For the case $2\leq s \leq k+1$, we have that $E^2_{s-1,1}$ is a subquotient of $E^1_{s-1,1}$ so the dimension of $E^2_{s-1,1}\otimes_{\z}\q$ is at most the dimension of $E^1_{s-1,1}\otimes_{\z}\q$, that is, at most the rank of $E^1_{s-1,1}$.

    Each copy of $H_1(\Gamma_0(p_{i_1}\cdots p_{i_s}))$ contributes with $r(p_{i_1}\cdots p_{i_{s-1}})$ copies of $\z$ (see Theorem \ref{H_n}), so \begin{align*}
        \rank(E^1_{s-1,1}) &= \rank\left(\bigoplus H_1(\Gamma_0(p_{i_1}\cdots p_{i_{s-1}}))^{2^{k-s+1}}\right) \\
        &= \sum (2^{k-s+1}\rank(H_1(\Gamma_0(p_{i_1}\cdots p_{i_{s-1}})))) \\
        &= 2^{k-s+1}\sum r(p_{i_1}\cdots p_{i_{s-1}})
    \end{align*} and the result follows.
\end{demonstracao}

This upper bound can be improved for $s=2,3$ using the following results:

\begin{teorema}\label{teo: r(p)}
    Let $n$ be an integer and $p$ be the prime divisor of $n$ such that $r(p)\leq r(q)$ for every other prime divisor $q$ of $n$. Then $$1 \leq \rank (H_2(\SL_2(\z[1/n])))\leq r(p).$$ In particular, if one of the primes $2,3,5,7,13$ divides $n$, then $$\rank (H_2(\SL_2(\z[1/n]))) = 1.$$
\end{teorema}

\begin{demonstracao}
See Propositions 7.1 and 7.2 of \cite{MRV2025}.
\end{demonstracao}

\begin{corolario}\label{cor: rank H3}
    Let $n = p_1\cdots p_k$ be the prime decomposition of a square-free integer $n$ and organize the primes so that $r(p_1)\leq r(p_i)$ for all $i=2,\dots,k$. Then $$\rank(H_3(\SL_2(\z[1/n]))) \leq 2^{k-2}\left(\sum_{1\leq i \leq k, i<j\leq k} r(p_ip_j)\right)-2^{k-1}\left(\sum_{1\leq i \leq k}r(p_i) \right)+r(p_1).$$ In particular, $\rank(H_3(\SL_2(\z[1/pq])))\leq r(pq)-2r(p)-r(q)$ if $r(q)\leq r(p)$.
\end{corolario}

\begin{demonstracao}
    Consider the following exact sequence of abelian groups: $$0 \rightarrow \ker(d^1_{2,1}) \rightarrow E^1_{2,1} \overset{d^1_{2,1}}{\longrightarrow} E^1_{1,1} \rightarrow \coker(d^1_{2,1}) \rightarrow 0.$$ This gives $$\rank(\ker(d^1_{2,1})) - \rank(E^1_{2,1})+\rank(E^1_{1,1})-\rank(\coker(d^1_{2,1}))=0.$$

    Since $\rank(H_3(\SL_2(\z[1/n]))) = \rank(\ker(d^1_{2,1})/\im(d^1_{3,1}))$ (see Theorem \ref{teo: sequence sl}), we obtain from Theorem \ref{teo: r(p)} \begin{align*}
        \rank(H_3(\SL_2(\z[1/n]))) &\leq \rank(\ker(d^1_{2,1})) \\
        &= \rank(E^1_{2,1})-\rank(E^1_{1,1})+\rank(\coker(d^1_{2,1})) \\
        &= \rank(E^1_{2,1})-\rank(E^1_{1,1})+\rank(H_2(\SL_2(\z[1/n]))) \\
        &\leq \rank(E^1_{2,1})-\rank(E^1_{1,1})+r(p_1) \\
        &= 2^{k-2}\left(\sum_{1\leq i \leq k, i<j\leq k} r(p_ip_j)\right)-2^{k-1}\left(\sum_{1\leq i \leq k}r(p_i) \right)+r(p_1).
    \end{align*}
\end{demonstracao}

\begin{corolario}\label{r(pq)}
Let $p,q$ be two distinct primes. If $p$ or $q$ is one of the primes $2,3,5, 7, 13$,
 then
 $$
 \rank(H_k(\SL_2(\z[1/pq])))=
 \begin{cases}
1 & \text{if $k=0,2,$}\\
r(pq)-2r(p)-2r(q)+1 & \text{if $k=3$}, \\
0 & \text{otherwise.}
 \end{cases}
 $$
\end{corolario}
\begin{demonstracao}
    Let $q\in\{2,3,5,7,13\}$ so that $r(q)=1$. Note that $r(pq)-2r(p)-2r(q)+1 = r(pq)-2r(p)-1 = r(pq)-2r(p)-r(q)$.

    If $n$ has two prime factors, then the associated spectral sequence has two non-zero columns and in particular, $E^1_{3,1}=0$, which implies $\im(d^1_{3,1}) = 0$. Theorem \ref{teo: sequence sl} then implies that $\rank(H_k(\SL_2(\z[1/pq]))) = \rank(\ker(d^1_{3,1}))$. We also have from Theorem \ref{teo: r(p)} that $\rank (H_2(\SL_2(\z[1/pq]))) = 1 = r(q)$.

    Thus, the two inequalities in the proof of Corollary \ref{cor: rank H3} become equalities, proving the statement.
\end{demonstracao}

\section{Application for the first homology}

For the remainder of this chapter, we focus on the homology groups $H_1(\SL_2(\z[1/n]))$. These groups are known (see \cite{B-E2025}), but the spectral sequence derived in Theorem \ref{teo: spec seq sl} allows us to construct a new proof for these results.

We first analyse the cases $H_1(\SL_2(\z[1/p]))$ where $p$ is a prime. The spectral sequence comes from the 1-dimensional complex $B_p$, so it has only two non-zero columns and it collapses on the second page. As in Proposition \ref{prp: line 0}, $E^2_{1,0}=0$ and $E^2_{0,1} = \coker(d^1_{1,1})$ where $d^1_{1,1}: H_1(\Gamma_0(p)) \rightarrow H_1(\SL_2(\z))\oplus H_1(\SL_2(\z))$ is induced by $A\mapsto(A,\pi_p^{-1}A\pi_p)$ where $\pi_p = \begin{pmatrix}
    1 & 0 \\ 0 & p
\end{pmatrix}$. 

Since $H_1(G)\simeq G^{\text{ab}}$ (see Theorem \ref{teo: H0 and H1}), for the cases $p=2,3$, we calculate these cokernels directly via a set of generators.

\begin{teorema}
    The group $\SL_2(\z)$ has presentation $\langle R,S \mid R^6=1,S^4=1, R^3=S^2\rangle$ where $S = \begin{pmatrix}
    0 & -1 \\ 1 & 0
\end{pmatrix},T = \begin{pmatrix}
    1 & 1 \\ 0 & 1
\end{pmatrix} $ and $ R = ST = \begin{pmatrix}
    0 & -1 \\ 1 & 1
\end{pmatrix}$.
\end{teorema}

\begin{demonstracao}
    See Example 1.5.3 of \cite{serre1980}.
\end{demonstracao}

Note that $R$ and $S$ have order 6 and 4, respectively. Theoerem \ref{H_n} gives $H_1(\SL_2(\z))\simeq \z/12$, which can be restated as $H_1(\SL_2(\z))\simeq \z/3\oplus\z/4$. For this second description, the generators correspond to the classes of $R$ and $S$ on the abelianization.

Given a set of generators of $\SL_2(\z)$ and a set of representatives for the cosets of $\SL_2(\z)/\Gamma_0(n)$, one can apply the Reidemeister-Schreier method (see Section III.6 of \cite{baumslag1993}) to obtain generators for $\Gamma_0(n)$, some of them redundant. 

We give a brief description of the method. Given $G$ a group, $H$ a subgroup, $\mathcal{S}_G$ a set of generators for $G$ and $\mathcal{T}$ a set of representatives for $G/H$, for every pair $x\in\mathcal{S}_G$ and $t\in\mathcal{T}$, we have that $tx$ must be in some coset $t'H$ for a unique $t' \in\mathcal{T}$. Denote by $s(t,x) := tx(t')^{-1}$. Then the set $$\mathcal{S}_H = \{s(t,x) \mid t\in\mathcal{T},x\in\mathcal{S}_G\}$$ is a set of generators for $H$.

\begin{proposicao}\label{prp:cosets p}
    Let $I$ denote the $2\times 2$ identity matrix. The set of matrices $ST^k$ with $0\leq k <n$, together with $I$, represent different cosets of $\SL_2(\z)/\Gamma_0(n)$. In particular, if $n=p$ prime, these matrices form a set of coset representatives for $\SL_2(\z)/\Gamma_0(p)$.
\end{proposicao}

\begin{demonstracao}
    Note that $$ST^k = \begin{pmatrix}
        0 & -1 \\ 1 & k
    \end{pmatrix}$$ which is not in $\Gamma_0(n)$, so neither of these matrices are in the coset represented by $I$. In addition, $$ST^k(ST^r)^{-1} = ST^{k-r}S^{-1} = -ST^{k-r}S = -\begin{pmatrix}
        -1 & 0 \\ k-r & -1
    \end{pmatrix}$$ is also not in $\Gamma_0(n)$, thus each $ST^k$ matrix represents a different coset.

    In particular, if $n=p$ is prime, since $[\SL_2(\z):\Gamma_0(p)] = a(p) = p+1$, this set will have $p+1$ elements representing distinct cosets, so the result follows.
\end{demonstracao}

\begin{definicao}
    Let $p$ be a prime number and let $k$ be an integer such that $0<k<p$. Denote by $k_*$ the unique integer such that $0<k_*<p$ and $1+kk_*\equiv 0 \mod p$.
\end{definicao}

\begin{proposicao}
    Let $p$ be a prime number. Then the set $$R = \{-I,T,ST^pS, -ST^kST^{-k_*}S \mid 0<k<p\}$$  is a set of generators for $\Gamma_0(p)$.
\end{proposicao}

\begin{demonstracao}
    We apply the Reidemeister-Schreier method as described earlier using $\{S,T\}$ as generators for $\SL_2(\z)$ and the set $\{I, ST^k \mid 0\leq k < p\}$ of coset representatives from the last proposition.

    Note that $T \in \Gamma_0(n)$ for all $n\in\z$ and $I,S, ST$ are always chosen as coset representatives. Thus, $$s(I,S) = ISS^{-1} = I, \hspace{30pt} s(I,T) = ITI^{-1} = T$$ $$s(S,S) = SSI^{-1}= -I, \hspace{30pt} s(S,T) = ST(ST^{-1}) = I$$ so our set of generators always includes the matrices $-I$ and $T$. 

    For the matrices $ST^k$ with $k>0$, let $r$ be an integer between 0 and $p-1$ such that $r\equiv 1+k \mod p$. Then $$s(ST^k,S) = (ST^kS)(ST^{k_*})^{-1} = -ST^kST^{-k_*}S = -\begin{pmatrix}
        k_* & 1 \\ -(1+kk_*) & -k
    \end{pmatrix},$$ $$s(ST^k,T) = (ST^kT)(ST^{r})^{-1} = -ST^{k+1-r}S = -\begin{pmatrix}
        -1 & 0 \\ k+1-r & -1
    \end{pmatrix}.$$

    Note that for $0\leq k \leq p-2$, we have that $r$ is simply given by $k+1$, so $s(ST^k,T) = I$. For $k=p-1$, we have $r = 0$, so $$s(ST^{p-1},T) = -ST^pS = \begin{pmatrix}
        1 & 0 \\ -p & 1
    \end{pmatrix}.$$
\end{demonstracao}

\begin{corolario}\label{gamma02}
    The matrices $T = \begin{pmatrix}
        1 & 1 \\ 0 & 1
    \end{pmatrix}$ and $STST^{-1}S = \begin{pmatrix}
      1 & 1 \\ -2 & -1  
    \end{pmatrix}$ generate $\Gamma_0(2)$.
\end{corolario}

\begin{demonstracao}
    The matrices obtained in the process described in the previous proposition are $-I,T$ (from the coset representatives $I,S$), the matrix $s(ST,T) = -ST^2S$ and for the unique pair $(1,1)$ of integers $(r,k)$ satisfying $1+kr\equiv 0 \mod 2$, the matrix $-STST^{-1}S$.

    One can check that $(STST^{-1}S)^2=-I$ and that $(STST^{-1}S)^{-1}T^{-1} = ST^2S$, so the generators $-I$ and $ST^2S$ are redundant. Thus it suffices to consider the generators $T$ and $STST^{-1}S$.
\end{demonstracao}

\begin{corolario}\label{gamma03}
    The matrices $T = \begin{pmatrix}
        1 & 1 \\ 0 & 1
    \end{pmatrix}$ and $-ST^2ST^{-1}S = \begin{pmatrix}
      -1 & -1 \\ 3 & 2  
    \end{pmatrix}$ generate $\Gamma_0(3)$.
\end{corolario}

\begin{demonstracao}
    The matrices obtained in the process described in the previous proposition are $-I,T$ (from the coset representatives $I,S$), the matrix $s(ST^2,T) = -ST^3S$ and for the pair $(1,2)$ of integers $(r,k)$ satisfying $1+kr\equiv 0 \mod 3$, the matrix $-ST^2ST^{-1}S$. 

    One can check that $(-ST^2ST^{-1}S)^3=-I$ and that $(-ST^2ST^{-1}S)^{-2}T^{-1} = -ST^3S$, so the generators $-I$ and $-ST^3S$ are redundant. Thus it suffices to consider the generators $T$ and $-ST^2ST^{-1}S$.
\end{demonstracao}

\begin{proposicao}\label{z[1/2]}
    $H_1(\SL_2(\z[1/2]))\simeq \z/3$.
\end{proposicao}

\begin{demonstracao}
    From Corollary \ref{gamma02}, $\Gamma_0(2)$ is generated by $T$ and $STST^{-1}S$. One can check that $\pi_2^{-1}T\pi_2 = T^2 = (RS^3)^2$ and $\pi_2^{-1}(STST^{-1}S)\pi_2 = T^{-2}S^2TST = SR^{-1}SR^{-1}S^2R^2S^3$. In terms of the generators $R,S$ of $\SL_2(\z)$, we have that $d^1_{1,1}: H_1(\Gamma_0(2)) \rightarrow H_1(\SL_2(\z))\oplus H_1(\SL_2(\z))$ is given by \begin{align*}
        \overline{T} &\mapsto (\overline{RS^3},\overline{RS^3RS^3}) = (\overline{RS^3},\overline{R^2S^2}),\\ \overline{STST^{-1}S}&\mapsto (\overline{S^3},\overline{S^7}) = (\overline{S^3},\overline{S^3}),
    \end{align*} in other words, the corresponding map $\z\oplus\z/4 \rightarrow (\z/3\oplus\z/4)^2$ is given by \begin{align*}
        (1,\overline{0})&\mapsto (\overline{1},\overline{3},\overline{2},\overline{2}), \\ (0,\overline{1})&\mapsto (\overline{0},\overline{3},\overline{0},\overline{3}),
    \end{align*} so $$(n,\overline{m})\in\ker(d^1_{1,1}) \implies \begin{cases}
        n\equiv 0 \pmod{3} \\ 2n\equiv 0 \pmod{3} \\ 3(n+m)\equiv 0 \pmod{4} \\ 2n+3m\equiv 0 \pmod{4} 
    \end{cases} \implies \begin{cases}
        n\equiv 0 \pmod{3} \\ n+m\equiv 0 \pmod{4} \\ 2n+3m\equiv 0 \pmod{4} 
    \end{cases}$$ Summing the last two equations gives $n\equiv 0 \pmod{4}$, which implies $n\equiv 0 \pmod{12}$. Plugging in either, this gives $m\equiv 0 \pmod{4}$. Thus, $$\ker(d^1_{1,1}) = \{(n,\overline{0})\mid n\equiv 0 \pmod{12}\} = 12\z\oplus 0$$ so that $$\im(d^1_{1,1}) \simeq \frac{\z\oplus\z/4}{12\z\oplus(0)} \simeq \z/12\oplus\z/4$$ and we conclude that $$H_1(\SL_2(\z[1/2])) \simeq E^2_{1,0} = \coker(d^1_{1,1}) \simeq \frac{\z/12\oplus\z/12}{\z/12\oplus\z/4} \simeq \z/3.$$ 
\end{demonstracao}

\begin{proposicao}\label{z[1/3]}
    $H_1(\SL_2(\z[1/3]))\simeq \z/4$.
\end{proposicao}

\begin{demonstracao}
    From Corollary \ref{gamma03}, $\Gamma_0(3)$ is generated by $T$ and $-ST^2ST^{-1}S$. These correspond to the generators of infinite order and order 6 of $H_1(\Gamma_0(3)) \simeq \z\oplus\z/2\oplus\z/3 \simeq \z\oplus\z/6$. One can check that $\pi_3^{-1}T\pi_3 = T^3 = (RS^3)^3$ and $\pi_3^{-1}(-ST^2ST^{-1}S)\pi_3 = T^{-1}S^3T^2S^2 = SR^{-1}S^3RS^3RS$. In terms of the generators $R,S$ of $\SL_2(\z)$, we have that $d^1_{1,1}$ is given by \begin{align*}
        \overline{T} &\mapsto (\overline{RS^3},\overline{RS^3RS^3RS^3}) = (\overline{RS^3},\overline{S^3}),\\ \overline{-ST^2ST^{-1}S}&\mapsto (\overline{S^3RS^3RSR^{-1}S},\overline{SR^{-1}S^3RS^3RS}) = (\overline{R},\overline{R}).
    \end{align*} In other words, the corresponding map $\z\oplus\z/6 \rightarrow (\z/3\oplus\z/4)^2$ is given by \begin{align*}
        (1,\overline{0})&\mapsto (\overline{1},\overline{3},\overline{0},\overline{3}), \\ (0,\overline{1})&\mapsto (\overline{1},\overline{0},\overline{1},\overline{0}),
    \end{align*} so $$(n,\overline{m})\in\ker(d^1_{1,1}) \implies \begin{cases}
        n+m\equiv 0 \pmod{3} \\ 3n\equiv 0 \pmod{4} \\ m\equiv 0 \pmod{3} \\ 3n\equiv 0 \pmod{4} 
    \end{cases} \implies \begin{cases}
        n\equiv 0 \pmod{3} \\ n\equiv 0 \pmod{4} \\ m\equiv 0 \pmod{3} 
    \end{cases}$$ So $n\equiv 0 \pmod{12}$ and $\overline{m}$ is either $\overline{0}$ or $\overline{3}$, that is, $$\ker(d^1_{1,1}) = \{(n,\overline{m})\mid n\equiv 0 \pmod{12},\, \overline{m}\in 3(\z/6)\}\} = 12\z\oplus 3(\z/6)$$ so that $$\im(d^1_{1,1}) \simeq \frac{\z\oplus\z/6}{12\z\oplus3(\z/6)} \simeq \z/12\oplus\z/3$$ and we conclude that $$H_1(\SL_2(\z[1/3])) \simeq E^2_{1,0} = \coker(d^1_{1,1}) \simeq \frac{\z/12\oplus\z/12}{\z/12\oplus\z/3} \simeq \z/4.$$ 
\end{demonstracao}

For the general case, a few results are needed:

\begin{proposicao}\label{prp: h1sl2(zn)}
    Let $p$ be a prime and $r$ a positive integer. Then $$H_1(\SL_2(\z/p^r)) \simeq \begin{cases}
        \z/2, &\text{ if } p=2,r=1 \\
        \z/4, &\text{ if } p=2,r\geq 2 \\
        \z/3, &\text{ if } p=3 \\
        0, &\text{ otherwise. } \\
    \end{cases}$$
\end{proposicao}

\begin{demonstracao}
    More generally, if $A$ is a local ring with maximal ideal $\mathfrak{m}$, then $$H_1(\SL_2(A)) \simeq \begin{cases}
        A/\mathfrak{m}^2, \hspace{5pt}\text{ if } |A/\mathfrak{m}|=2, \\
        A/\mathfrak{m}, \hspace{10pt}\text{ if } |A/\mathfrak{m}|=3, \\
        0, \hspace{27pt} \text{ if } |A/\mathfrak{m}|\geq 4.
    \end{cases}$$ For a proof of this, see Proposition 4.1 of \cite{B-E2025}. Note that the maximal ideal of the local ring $\z/p^r$ is 0 if $r=1$ and $p\z/p^r$ otherwise, so the residue fields are $\F_p$ in either case.
\end{demonstracao}

\begin{lema}\label{surjection gen}
    Let $f: A \rightarrow B$ be a surjective map and let \begin{center}
        \begin{tikzcd}
            A \arrow[r, "\psi"] & B\oplus B \arrow[r, "\phi"] & C \arrow[r] & 0
        \end{tikzcd}
    \end{center} be an exact sequence where $\psi: B\oplus B \rightarrow C$ is given by $a\mapsto (f(a),-f(a))$. Then there exists a surjection $B\rightarrow C$.
\end{lema}

\begin{demonstracao}
    Define a map $\Phi: B \rightarrow (B\oplus B)/\im\psi$ by $b\mapsto \overline{(b,0)}$. Let $\overline{(b_1,b_2)}\in B \rightarrow (B\oplus B)/\im\psi$. Since $f$ is surjective, there exists $a\in A$ such that $f(a)=b_2$. Then \begin{align*}
        \overline{(b_1,b_2)} &= \overline{(-b_2,b_2)}+\overline{(b_1+b_2,0)} \\ &= \overline{(f(-a),-f(-a))}+\overline{(b_1+b_2,0)} \\ &= \overline{\psi(a)}+\overline{(b_1+b_2,0)} \\ &= \overline{(b_1+b_2,0)} \\ &= \Phi(b_1+b_2),
    \end{align*} so $\Phi$ is surjective.

    From the first isomorphism theorem, $$\frac{B\oplus B}{\im\psi} = \frac{B\oplus B}{\ker\phi} \simeq \im\phi = C,$$ where this isomorphism is given by $\overline{\phi}$, which is defined as $\overline{\phi}(\overline{(b,b)}) = \phi(b,b)$. Define $\varphi: B \rightarrow C$ by $\varphi = \phi\circ\Phi$. Then $\varphi$ is a composition of surjections and is therefore a surjection.
\end{demonstracao}

\begin{proposicao}\label{prp: five term}
    Let $G$ be a group and $N$ a normal subgroup of $G$. Then there exists a five-term exact sequence $$H_2(G) \rightarrow H_2(G/N) \rightarrow H_1(N)_{G/N}\rightarrow H_1(G) \rightarrow H_1(G/N) \rightarrow 0$$ where $H_1(N)_{G/N}$ denoted the set of $G/N$-coinvariants of $H_1(N)$ (see Proposition \ref{prp: coinvariants}).
\end{proposicao}

\begin{demonstracao}
    See Corollary 6.4 in Chapter VII.6 of \cite{brown1994}.
\end{demonstracao}

Given a pair $(G,N)$ of a group and a normal subgroup, the exact sequence $$1 \rightarrow N \rightarrow G \rightarrow G/N \rightarrow 1$$ is said to be a group extension and the above sequence is said to be derived from this group extension. This sequence can be obtained from the Hochschild–Serre spectral sequence associated to the above extension.  

\begin{proposicao}
    \label{surjection}
    Let $p$ be a prime and $m$ an integer not divisible by $p$. Then, there exists a surjection $H_1(\SL_2(\z[1/m])) \rightarrow H_1(\SL_2(\z[1/pm])$.
\end{proposicao}

\begin{demonstracao}
    Note that $\SL_2(\z[1/pm])$ acts on the tree $B_p$ as in Definition \ref{def: Bp}. This group contains $\SL_2(\z)$, so each edge of the tree is equivalent and vertices at same distance are also equivalent. Furthermore, since each matrix has determinant 1, this group will not identify vertices at even distances, so the orbits coincide with the orbits of $\SL_2(\z[1/p])$. 
    
    The stabilizer of the base vertex is $\SL_2(\z[1/pm])\cap\SL_2(\z_{(p)}) = \SL_2(\z[1/m])$, so the other vertex has stabilizer $\pi_p\SL_2(\z[1/m])\pi_p^{-1}$ and the edges have stabilizer 
    $$\Gamma_0(p,m) = \left\{ \begin{pmatrix}
        a & b \\ c & d
    \end{pmatrix}\in \SL_2(\z[1/m]) : p \mid c \right\},$$ 
    the subgroup of matrices of $\SL_2(\z[1/m])$ whose lower left entry is divisible by $p$. Then, Theorem \ref{exact-seq} gives an exact sequence $$\cdots \rightarrow H_1(\Gamma_0(p,m)) \rightarrow H_1(\SL_2(\z[1/m]))\oplus H_1(\SL_2(\z[1/m]))\rightarrow H_1(\SL_2(\z[1/pm])) \rightarrow \cdots$$ 

    The next map in the sequence is $\z\mapsto \z\oplus\z$ given by $1\mapsto (1,-1)$, which is injective, so that $$ H_1(\Gamma_0(p,m)) \rightarrow H_1(\SL_2(\z[1/m]))\oplus H_1(\SL_2(\z[1/m]))\rightarrow H_1(\SL_2(\z[1/pm])) \rightarrow 0$$ is an exact sequence (that is, the second map is surjective).

    Let $\Gamma(p), B(\mathbb{F}_p)$ be the following groups: $$\Gamma(p) = \left\{\begin{pmatrix}
        a & b \\ c & d
    \end{pmatrix}\in\SL_2(\z) : a-1,d-1,b,c\equiv 0\pmod{p}\right\},$$ $$B(\mathbb{F}_p) = \left\{\begin{pmatrix}
        a & b \\ 0 & a^{-1}
    \end{pmatrix}\in\SL_2(\mathbb{F}_p) : a\in (\mathbb{F}_p)^\times, b\in \mathbb{F}_p\right\}.$$ Then we have the following commutative diagram of group extensions: \begin{center}
        \begin{tikzcd}
            1 \arrow[r] & \Gamma(p,m) \arrow[r] \arrow[d]  & \Gamma_0(p,m) \arrow[r] \arrow[d] & B(\mathbb{F}_p) \arrow[r] \arrow[d] & 1 \\
            1 \arrow[r] & \Gamma(p,m) \arrow[r] & \SL_2(\z[1/m]) \arrow[r] & \SL_2(\mathbb{F}_p) \arrow[r] & 1
        \end{tikzcd}
    \end{center} where the maps $\Gamma_0(p,m)\rightarrow B(\mathbb{F}_p)$, $\SL_2(\z[1/m])\rightarrow \SL_2(\F_p)$ take each entry of a matrix to its residue mod $p$ and the other maps are inclusions. These map make sense since $p\nmid m$ implies that elements of $\z[1/m]$ are fractions $a/b$ where $p\nmid b$, so the term $\overline{a}\overline{b}^{-1}$ makes sense in $\mathbb{F}_p$. 

    Surjectivity of the two maps on the left follows from the fact that $\mathbb{F}_p$ is a field (thus a local ring) and so $\SL_2(\mathbb{F}_p)$ is generated by elementary matrices (see Theorem \ref{teo: GE2}). Exactness follows from the fact that $\Gamma(p,m)$ is the subgroup of matrices congruent to the identity mod $p$, that is, the preimage of $I$ under these surjections.
    
    Each extension induces a five-term exact sequence in homology (see Proposition \ref{prp: five term}), and commutativity of the diagram implies that the following diagram also commutes: \begin{center}
        \begin{tikzcd}
             & H_1(\Gamma(p,m))_{B(\mathbb{F}_p)} \arrow[r] \arrow[d, two heads]  & H_1(\Gamma_0(p,m)) \arrow[r] \arrow[d] & H_1(B(\mathbb{F}_p)) \arrow[r] \arrow[d, two heads] & 0 \\
            0 \arrow[r] & H_1(\Gamma(p,m))_{\SL_2(\mathbb{F}_p)} \arrow[r] & H_1(\SL_2(\z[1/m])) \arrow[r] & H_1(\SL_2(\mathbb{F}_p)) \arrow[r] & 0
        \end{tikzcd}
    \end{center} 

    By Proposition \ref{prp: h1sl2(zn)}, $H_2(\SL_2(\mathbb{F}_p))=0$. The left vertical map is surjective because it is a coset enlargment, since the set of elements $gt-t$ with $t\in \Gamma(p,m), g\in B(\mathbb{F}_p)$ is a subset of the set of elements $gt-t$ with $t\in \Gamma(p,m), g\in \SL_2(\mathbb{F}_p)$. 

    To see surjectivity of the right vertical map (which is induced by the inclusion $B(\mathbb{F}_p)\hookrightarrow\SL_2(\mathbb{F}_p)$ on the abelianizations), note that since $\mathbb{F}_p$ is a local ring, $\SL_2(\mathbb{F}_p)$ is generated by elementary matrices and thus $\SL_2(\z)$ surjects to it. So by Proposition \ref{prp: surjection} there is a surjection $H_1(\SL_2(\z))\rightarrow H_1(\SL_2(\mathbb{F}_p))$. But $H_1(\SL_2(\z))\simeq \z/12$ can be seen as generated by $\overline{T}$, so via this surjection,  $H_1(\SL_2(\mathbb{F}_p))$ is also generated by $\overline{T} = \overline{\begin{pmatrix}
        \overline{1} & \overline{1} \\ \overline{0} & \overline{1}
    \end{pmatrix}}$ which is in $H_1(B(\mathbb{F}_p))$.

    The Snake Lemma (see Corollary 6.12 of \cite{rotman2009}) then implies that the map $H_1(\Gamma_0(p,m)) \rightarrow H_1(\SL_2(\z[1/m]))$ induced by the inclusion $\Gamma_0(p,m)\rightarrow \SL_2(\z[1/m])$ is surjective. If this map is $i_*$, then the map $H_1(\Gamma_0(p,m))\rightarrow H_1(\SL_2(\z[1/m]))\oplus H_1(\SL_2(\z[1/m]))$ from the exact sequence described earlier is given by $[a]\mapsto (i_*[a],-i_*[a])$. Lemma \ref{surjection gen} then gives the surjection $H_1(\SL_2(\z[1/m]))\rightarrow H_1(\SL_2(\z[1/pm]))$, as desired.
\end{demonstracao}

\begin{corolario}
    Let $p>3$ be prime. Then $H_1(\SL_2(\z[1/p]))\simeq \z/12$.
\end{corolario}

\begin{demonstracao}
    Since $p\neq 3$, the elements of $\z[1/p]$ are fractions $a/b$ where $3\nmid b$ (since $b$ is a power of $p$) so the residue of $b$ mod 3 is non-zero. In other words, the map $\z[1/p]\rightarrow \z/3$ given by mapping a fraction $a/b$ to $\overline{a}/\overline{b}$ is a well-defined surjection. In a similar way, since $p\neq 2$, there is a surjection $\z[1/p]\rightarrow \z/4$.

    Since both $\z/3$ and $\z/4$ are local rings, we have that $\SL_2(\z/3)$ and $\SL_2(\z/4)$ are generated by elementary matrices (see Theorem \ref{teo: GE2}), so the induced maps from $\SL_2(\z[1/p])$ to $ \SL_2(\z/3)$ and $\SL_2(\z/4)$ are also surjective. Then, surjectivity of these maps gives surjections from $H_1(\SL_2(\z[1/p]))$ to $H_1(\SL_2(\z/3))=\z/3$ and $H_1(\SL_2(\z/4))=\z/4$ (see Propositions \ref{surjection} and \ref{prp: h1sl2(zn)}). Let $\varphi_1,\varphi_2$ be these surjections.

    Define $\Phi: H_1(\SL_2(\z[1/p])) \rightarrow \z/3\oplus\z/4\simeq \z/12$ by $\Phi([x])= (\varphi_1([x]),\varphi_2([x]))$. Since there exists $[x]\in H_1(\SL_2(\z[1/p]))$ with $\varphi_1([x])=\overline{1}\in\z/3$, we have that $$\Phi(4[x])= (\varphi_1(4[x]),\varphi_2(4[x])) = (4\varphi_1([x]),4\varphi_2([x])) = (\overline{1},\overline{0})$$ and similarly, one finds a class $[y]$ with $\varphi_2([y])=\overline{1}\in\z/4$ and $\Phi(-3[y]) = (\overline{0},\overline{1})$. This implies $\Phi$ is surjective.

    From Proposition \ref{surjection}, taking $m=1$, there exists a surjection from $H_1(\SL_2(\z)) $ to $ H_1(\SL_2(\z[1/p]))$. The composition of this surjection with $\Phi$ gives a surjection $H_1(\SL_2(\z))\simeq  \z/12\rightarrow\z/12$ which must also be injective from finiteness of $\z/12$. This implies $\Phi$ is an isomorphism, so $H_1(\SL_2(\z[1/p]))\simeq \z/12$. 
\end{demonstracao}

\begin{teorema}\label{teo: calculation H1}
    For any positive integer $n$, we have $$H_1(\SL_2(\z[1/n])) \simeq \begin{cases}
        0, &\text{ if } 2\mid n, 3\mid n \\
        \z/3, &\text{ if } 2\mid n, 3\nmid n \\
        \z/4, &\text{ if } 2\nmid n, 3\mid n \\
        \z/12, &\text{ if } 2\nmid n, 3\nmid n. 
    \end{cases}$$
\end{teorema}

\begin{demonstracao}
    Recall that if $n=p_1^{e_1}p_2^{e_2}\cdots p_k^{e_k}$, then $\z[1/n] = \z[1/p_1p_2\cdots p_k]$, so it suffices to consider the case where $n$ is square free.

    Let $p$ be the largest prime that divides $n$ and write $n=pm$. If $p=2$, this implies $n=2$ and the result follows from Proposition \ref{z[1/2]}. If $p=3$, then either $m=1$ (so $n=3$) and the result follows from Proposition \ref{z[1/3]} or $m=2$ and so $n=6$. In this case, $\SL_2(\z[1/6])$ is an euclidean domain, so it is generated by elementary matrices $E_{1,2}(a)$ and $E_{2,1}(b)$. But $$E_{1,2}(a) = \begin{pmatrix}
        1 & a \\ 0 & 1
    \end{pmatrix} = \left[\begin{pmatrix}
        1 & -\frac{1}{3}a \\ 0 & 1
    \end{pmatrix}, \begin{pmatrix}
        2 & 0 \\ 0 & \frac{1}{2}
    \end{pmatrix}\right]$$ $$E_{2,1}(b) = \begin{pmatrix}
        1 & 0 \\ b & 1
    \end{pmatrix} = \left[\begin{pmatrix}
        1 & 0 \\ \frac{4}{3}b & 1
    \end{pmatrix}, \begin{pmatrix}
        2 & 0 \\ 0 & \frac{1}{2}
    \end{pmatrix}\right]$$ implying that $\SL_2(\z[1/6]) = [\SL_2(\z[1/6]),\SL_2(\z[1/6])]$, that is, $H_1(\SL_2(\z[1/6])) = 0$, as expected.

    We now consider the case $p>3$. We prove the claim by induction on the number of primes that divides $n = p_1\cdots p_k$ and we assume that $k>1$.

    If $2\mid n, 3\mid n$, then $2\mid m, 3\mid m$. By induction, $H_1(\SL_2(\z[1/m]))=0$, then Proposition \ref{surjection} gives a surjection $$0  =H_1(\SL_2(\z[1/m])) \rightarrow H_1(\SL_2(\z[1/n]))$$ so $H_1(\SL_2(\z[1/n]))$ = 0.

    If $2\mid n, 3\nmid n$, then $2\mid m$ and $3\nmid m$ since $n=pm$ with $p$ the largest prime dividing $n$, $p>3$.  By induction, $H_1(\SL_2(\z[1/m])) \simeq \z/3$ and Proposition \ref{surjection} gives a surjection $$\z/3 =H_1(\SL_2(\z[1/m])) \rightarrow  H_1(\SL_2(\z[1/n])).$$ On the other hand, since $3\nmid n$, there is a surjection $\z[1/n] \rightarrow \z/3$ that induces a surjection $H_1(\SL_2(\z[1/n])) \rightarrow H_1(\SL_2(\z/3))\simeq\z/3$ (see Proposition \ref{prp: h1sl2(zn)}). These two surjections then imply $H_1(\SL_2(\z[1/n]))\simeq \z/3$.

    If $2\nmid n, 3\mid n$, then once again $2\nmid m, 3\mid m$.  By induction, $H_1(\SL_2(\z[1/m])) \simeq \z/4$ and Proposition \ref{surjection} gives a surjection $$\z/4 \simeq H_1(\SL_2(\z[1/m])) \rightarrow  H_1(\SL_2(\z[1/n])).$$ From $2\nmid n$, there is a surjection $\z[1/n] \rightarrow \z/4$ that induces a surjection $H_1(\SL_2(\z[1/n])) \rightarrow H_1(\SL_2(\z/4))\simeq\z/4$ (see Proposition \ref{prp: h1sl2(zn)}). These two surjections will then imply that $H_1(\SL_2(\z[1/n]))\simeq \z/4$.

    Finally, we consider the case $2\nmid n,3\nmid n$. Then $2\nmid m,3\nmid m$ and so $H_1(\SL_2(\z[1/m]))\simeq \z/12$ by induction. Proposition \ref{surjection} then gives a surjection $$\z/12 \simeq H_1(\SL_2(\z[1/m])) \rightarrow  H_1(\SL_2(\z[1/n])).$$ The conditions $2\nmid n,3\nmid n$ give surjections from $H_1(\SL_2(\z[1/n]))$ to both $\z/4$ and $\z/3$, thus giving a surjection to $\z/12$. This implies $H_1(\SL_2(\z[1/n]))\simeq \z/12$.
\end{demonstracao}

\bookmarksetup{startatroot}%

\postextual

\bibliography{references}




\begin{apendicesenv}

    \chapter{Code for the graph of \texorpdfstring{$Y/\Gamma_0(p)$}{Lg}}
    \label{chapter:codigo}
    \begin{codigo}[caption={Code for the construction of the quotient graph $Y/\Gamma_0(p)$},  language={python},  breaklines=true]
#chose a prime number n
n = 19

E=[]
for i in range(0,n):
    w=[]
    w.append(i)
    for j in range(0,n):
        if (1+i*j)
            w.append(j)
    E.append(w)
#this creates a list of edge representatives
#for every i, it decides which j is such that R^j is the same edge as R^i
         
redE=E.copy()
for i in range(0,len(E)):
    for j in range(0,i):
        if len(E[i])==2:
            if E[i][1]==E[j][0]:
                redE.remove(E[i])
#redE is the reduction of E
#if i,j are equivalent then [i,j] and [j,i] show up in E
#the loop tests for each pair in E if there is some other vector with coordinates swapped
         
V=[]
for i in range(0,n):
    w=[]
    w.append(i)
    for j in range(0,n):
        if (1+i+i*j)
            w.append(j)
        elif (1+j+i*j)
            w.append(j)
    V.append(w)
#this creates a list of vertex representatives
#for every i, it decides what other vertices are equivalent to it
#R^0 is always equivalent to R^(n-1) and to S 
#0 and n-1 produce a pair, other integers produce triples 
        
V.remove(V[len(V)-1])

redV=V.copy()
for i in range(0,len(V)):
    if len(V[i])==3:
        for j in range(0,i):
             if V[i][1]==V[j][0]:
                 redV.remove(V[i])
                 break
             if V[i][2]==V[j][0]:
                 redV.remove([i])
                 break 
#redV is the reduction of V
                       
edges=[[0,0]]
for i in range(1,len(redE)):
    w=[]
    if redE[i][0]==redE[i][1]:
        w.append(redE[i][0])
        w.append('E'+str(redE[i][0]))
    else:
        w.append(redE[i][0]-1)
        w.append(redE[i][0])
    edges.append(w)
#each coset representative R^k is an edge from R^k to R^(k-1)
#this runs through the reduced representatives and creates edges [k,k-1] 
#if k corresponds to added edge, the loop creates [k, E]

corE=edges.copy()
for i in range(1,len(edges)):
    for j in range(0,2):
            for k in range(0,len(redV)):
                if edges[i][j]==redV[k][1]:
                    corE[i][j] = redV[k][0]
                elif len(redV[k])==3:
                    if edges[i][j]==redV[k][2]:
                        corE[i][j] = redV[k][0]
        
#corE corrects the list of edges by swapping [k,k-1] for [r,s] if k~r and k-1~s in V
                        
print(corE)
                        
                
\end{codigo}
    

\end{apendicesenv}





\end{document}